\input amstex
\documentstyle{amsppt}
\TagsOnRight
\magnification=\magstep1
\pageheight{24 truecm}
\pagewidth{16.8 truecm}
%
\loadbold
\loadeusm
\def\scr#1{{\fam\eusmfam\relax#1}}
\message{<Paul Taylor's commutative diagrams, version 3.83, 18 May 1995>}%
\ifx\diagram\undefined\else\message{WARNING: the \string\diagram\space command
is already defined and will not be loaded again}\expandafter\endinput\fi

\edef\cdrestoreat{
\noexpand\catcode\lq\noexpand\@=\the\catcode\lq\@
\noexpand\catcode\lq\noexpand\#=\the\catcode\lq\#
\noexpand\catcode\lq\noexpand\$=\the\catcode\lq\$
\noexpand\catcode\lq\noexpand\<=\the\catcode\lq\<
\noexpand\catcode\lq\noexpand\>=\the\catcode\lq\>
\noexpand\catcode\lq\noexpand\+=\the\catcode\rq53%
}\catcode\lq\@=11 \catcode\lq\#=6 \catcode\lq\<=12 \catcode\lq\>=12 \catcode
\rq53=12

\ifx\diagram@help@messages\undefined\let\diagram@help@messages y\fi

\def\cdps@Rokicki#1{\special{ps:#1}}\let\cdps@dvips\cdps@Rokicki\let
\cdps@RadicalEye\cdps@Rokicki\let\Cd@KJ\cdps@Rokicki\let\Cd@CB\cdps@Rokicki
\def\cdps@Bechtolsheim#1{\special{dvitps: Literal "#1"}}%
\let\cdps@dvitps\cdps@Bechtolsheim\let\cdps@IntegratedComputerSystems
\cdps@Bechtolsheim
\def\cdps@Clark#1{\special{dvitops: inline #1}}
\let\cdps@dvitops\cdps@Clark
\let\cdps@OzTeX\empty\let\cdps@oztex\empty\let\cdps@Trevorrow\empty
\def\cdps@Coombes#1{\special{ps-string #1}}

\count@=\year\multiply\count@12 \advance\count@\month
\ifnum\count@>2396000 
\message{***********************************************************}
\message{! THIS IS AN EXPERIMENTAL VERSION OF COMMUTATIVE DIAGRAMS *}
\message{! it expired in August 1996 and is time-bombed for January *}
\message{! You may obtain an up to date version of this package by *}
\message{! "anonymous FTP" from theory.doc.ic.ac.uk (146.169.2.27) *}
\message{***********************************************************}
\ifnum\count@>2396300 
\errhelp{You may press RETURN and carry on for the time being.}\message{! It
is embarrassing to see papers in conference proceedings}\message{! and
journals containing bugs which I had fixed years before.}\message{! It is easy
to obtain and install a new version, which will}\errmessage{! remain
compatible with your files. Please get it NOW.}\fi\fi

\def\Cd@oD{\global\let}\def\Cd@mG{\outer\def}

{\escapechar\m@ne\xdef\Cd@k{\string\{}\xdef\Cd@mC{\string\}}
\catcode\lq\&=4 \Cd@oD\Cd@N=&\xdef\Cd@P{\string\&}
\catcode\lq\$=3 \Cd@oD\Cd@h=$\Cd@oD\Cd@zC=$
\xdef\Cd@dC{\string\$}\gdef\Cd@mF{$$}
\catcode\lq\_=8 \Cd@oD\Cd@oI=_
\obeylines\catcode\lq\^=7 \Cd@oD\@super=^
\ifnum\newlinechar=10 \gdef\Cd@QG{^^J}
\else\ifnum\newlinechar=13 \gdef\Cd@QG{^^M}
\else\Cd@oD\Cd@QG\space\expandafter\message{! input error: \noexpand
\newlinechar\space is ASCII \the\newlinechar, not LF=10 or CR=13.}
\fi\fi}

\mathchardef\lessthan=\rq30474 \mathchardef\greaterthan=\rq30476

\ifx\tenln\undefined
\font\tenln=line10\relax
\fi\ifx\tenlnw\undefined\ifx\tenln\nullfont\let\tenlnw\nullfont\else
\font\tenlnw=linew10\relax
\fi\fi

\ifx\inputlineno\undefined\csname newcount\endcsname\inputlineno\inputlineno
\m@ne\message{***************************************************}\message{!
Obsolete TeX (version 2). You should upgrade to *}\message{! version 3, which
has been available since 1990. *}\message{***********************************%
****************}\fi

\newif\if@ignore

\def\cd@shouldnt#1{\Cd@FB{* THIS (#1) SHOULD NEVER HAPPEN! *}}

\def\get@round@pair#1(#2,#3){#1{#2}{#3}}
\def\get@square@arg#1[#2]{#1{#2}}
\def\Cd@lD#1{\Cd@QJ\let\Cd@tJ\Cd@kD\Cd@kD#1,],}
\def\Cd@j{[}\def\Cd@BD{]}\def\commdiag#1{{\let\enddiagram\relax\diagram[]#1%
\enddiagram}}

\def\Cd@iE{{\ifx\Cd@aG[\aftergroup\get@square@arg\aftergroup\Cd@tG\else
\aftergroup\Cd@eG\fi}}
\def\Cd@hE#1#2{\def\Cd@tG{#1}\def\Cd@eG{#2}\futurelet\Cd@aG\Cd@iE}

\def\Cd@MJ{|}

\def\Cd@KB{
\tokcase\Cd@qC:\Cd@u\break@args;\catcase\@super:\upper@label;\catcase\Cd@oI:%
\lower@label;\tokcase{~}:\middle@label;
\tokcase<:\left@label;
\tokcase>:\right@label;
\tokcase(:\Cd@wB;
\tokcase[:\optional@;
\tokcase.:\Cd@RI;
\catcase\space:\eat@space;\catcase\bgroup:\positional@;\default:\Cd@w
\break@args;\endswitch}

\def\switch@arg{
\catcase\@super:\upper@label;\catcase\Cd@oI:\lower@label;\tokcase[:\optional@
;
\tokcase.:\Cd@RI:
\catcase\space:\eat@space;\catcase\bgroup:\positional@;\tokcase{~}:%
\middle@label;
\default:\Cd@u\break@args;\endswitch}


\def\Cd@VA#1#2{\def#1{\Cd@bB{#2\Cd@XD}\Cd@oD#1\relax}}\let\Cd@wI\relax\ifx
\protect\undefined\let\protect\relax\fi\def\Cd@VF#1\repeat{\def\Cd@l{#1}%
\Cd@tE}\def\Cd@tE{\Cd@l\relax\expandafter\Cd@tE\fi}\def\Cd@TF#1\repeat{\def
\Cd@m{#1}\Cd@uE}\def\Cd@uE{\Cd@m\relax\expandafter\Cd@uE\fi}\def\Cd@UF#1%
\repeat{\def\Cd@n{#1}\Cd@vE}\def\Cd@vE{\Cd@n\relax\expandafter\Cd@vE\fi}\def
\Cd@PG#1#2#3{\def#2{\let#1\iftrue}\def#3{\let#1\iffalse}#3}\if y%
\diagram@help@messages\def\Cd@NG#1#2{\csname newtoks\endcsname#1#1=%
\expandafter{\csname#2\endcsname}}\else\csname newtoks\endcsname\no@cd@help
\no@cd@help{See the manual}\def\Cd@NG#1#2{\let#1\no@cd@help}\fi\chardef\Cd@OF
=1 \chardef\Cd@xH=2 \chardef\Cd@hG=5 \chardef\Cd@IH=6 \chardef\Cd@HH=7
\chardef\Cd@GC=9 \dimendef\Cd@uH=2 \dimendef\Cd@LF=3 \dimendef\Cd@PF=4
\dimendef\Cd@yH=5 \dimendef\Cd@zI=6 \dimendef\Cd@EI=8 \dimendef\Cd@DI=9
\skipdef\Cd@oB=1 \skipdef\Cd@sE=2 \skipdef\Cd@nB=3 \skipdef\Cd@IE=4 \skipdef
\Cd@LJ=5 \skipdef\Cd@wH=6 \skipdef\Cd@NF=7 \skipdef\Cd@BI=8 \skipdef\Cd@AI=9
\countdef\Cd@BC=9 \countdef\Cd@QD=8 \countdef\Cd@A=7 \def\sdef#1#2{\def#1{#2}%
}\def\Cd@I#1{\expandafter\aftergroup\csname#1\endcsname}\def\Cd@IC#1{%
\expandafter\def\csname#1\endcsname}\def\Cd@dD#1{\expandafter\gdef\csname#1%
\endcsname}\def\Cd@jC#1{\expandafter\edef\csname#1\endcsname}\def\Cd@QF#1#2{%
\expandafter\let\csname#1\expandafter\endcsname\csname#2\endcsname}\def\Cd@pD
#1#2{\expandafter\Cd@oD\csname#1\expandafter\endcsname\csname#2\endcsname}%
\def\Cd@CJ#1{\csname#1\endcsname}\def\Cd@bI#1{\expandafter\show\csname#1%
\endcsname}\def\Cd@dI#1{\expandafter\showthe\csname#1\endcsname}\def\Cd@aI#1{%
\expandafter\showbox\csname#1\endcsname}\def\Cd@oA{Commutative Diagram}\edef
\Cd@CH{\string\par}\edef\Cd@UC{\string\diagram}\edef\Cd@uC{\string\enddiagram
}\edef\Cd@xB{\string\\}\def\Cd@IF{LaTeX}\def\Cd@d{{\ifnum0=\lq}\fi}\def\Cd@kC
{\ifnum0=\lq{\fi}}\def\catcase#1:{\ifcat\noexpand\Cd@aG#1\Cd@wI\expandafter
\Cd@aC\else\expandafter\Cd@hI\fi}\def\tokcase#1:{\ifx\Cd@aG#1\Cd@wI
\expandafter\Cd@aC\else\expandafter\Cd@hI\fi}\def\Cd@aC#1;#2\endswitch{#1}%
\def\Cd@hI#1;{}\let\endswitch\relax\def\default:#1;#2\endswitch{#1}\ifx\at@
\undefined\def\at@{@}\fi\edef\Cd@M{\Cd@k pt\Cd@mC}\Cd@IC{\Cd@M>}#1>#2>{\Cd@v
\rTo\sp{#1}\sb{#2}\Cd@v}\Cd@IC{\Cd@M<}#1<#2<{\Cd@v\lTo\sp{#1}\sb{#2}\Cd@v}%
\Cd@IC{\Cd@M)}#1)#2){\Cd@v\rTo\sp{#1}\sb{#2}\Cd@v}
\Cd@IC{\Cd@M(}#1(#2({\Cd@v\lTo\sp{#1}\sb{#2}\Cd@v}
\def\Cd@L{\def\endCD{\enddiagram}\Cd@IC{\Cd@M A}##1A##2A{\uTo<{##1}>{##2}%
\Cd@v\Cd@v}\Cd@IC{\Cd@M V}##1V##2V{\dTo<{##1}>{##2}\Cd@v\Cd@v}\Cd@IC{\Cd@M=}{%
\Cd@v\hEq\Cd@v}\Cd@IC{\Cd@M\Cd@MJ}{\vEq\Cd@v\Cd@v}\Cd@IC{\Cd@M\string\vert}{%
\vEq\Cd@v\Cd@v}\Cd@IC{\Cd@M.}{\Cd@v\Cd@v}\let\Cd@v\Cd@N}\def\Cd@sD{\let\tmp
\Cd@tD\ifcat A\noexpand\Cd@ZG\else\ifcat=\noexpand\Cd@ZG\else\ifcat\relax
\noexpand\Cd@ZG\else\let\tmp\at@\fi\fi\fi\tmp}\def\Cd@tD#1{\Cd@QF{tmp}{\Cd@M
\string#1}\ifx\tmp\relax\def\tmp{\at@#1}\fi\tmp}\def\Cd@v{}\begingroup
\aftergroup\def\aftergroup\Cd@Q\aftergroup{\aftergroup\def\catcode\lq\@%
\active\aftergroup @\endgroup{\futurelet\Cd@ZG\Cd@sD}}\newcount\Cd@pA
\newcount\Cd@qA\newcount\Cd@rA\newcount\Cd@sA\newdimen\Cd@KA\newdimen\Cd@LA
\Cd@PG\Cd@ME\Cd@w\Cd@u\Cd@PG\Cd@NE\Cd@AA\Cd@y\newdimen\Cd@NA\newdimen\Cd@OA
\newcount\Cd@tA\newcount\Cd@uA\newdimen\Cd@MA\newbox\Cd@@A\Cd@PG\Cd@RE\Cd@YA
\Cd@XA\newcount\Cd@gG\newcount\Cd@KC\def\Cd@S#1#2{\ifdim#1<#2\relax#1=#2%
\relax\fi}\def\Cd@U#1#2{\ifdim#1>#2\relax#1=#2\relax\fi}\newdimen\Cd@sG\Cd@sG
=1sp \newdimen\Cd@nC\Cd@nC\z@\def\Cd@gI{\ifdim\Cd@nC=1em\else\Cd@qI\fi}\def
\Cd@qI{\Cd@nC1em\def\Cd@FC{\fontdimen8\textfont3 }\Cd@LI\Cd@zJ\setbox0=\vbox{%
\Cd@p\noindent\Cd@h\null\penalty-9993\null\Cd@zC\null\endgraf\setbox0=%
\lastbox\unskip\unpenalty\setbox1=\lastbox\global\setbox\Cd@jF=\hbox{\unhbox0%
\unskip\unskip\unpenalty\setbox0=\lastbox}\global\setbox\Cd@lF=\hbox{\unhbox1%
\unskip\unpenalty\setbox1=\lastbox}}}\newdimen\Cd@NH\Cd@NH=1true in \divide
\Cd@NH300 \def\Cd@MH#1{\multiply#1\tw@\advance#1\ifnum#1<\z@-\else+\fi\Cd@NH
\divide#1\tw@\divide#1\Cd@NH\multiply#1\Cd@NH}\def\MapBreadth{%
\afterassignment\Cd@tH\Cd@qE}\newdimen\Cd@qE\newdimen\Cd@@I\def\Cd@tH{\Cd@@I
\Cd@qE\Cd@S\Cd@NH{4\Cd@sG}\Cd@U\Cd@NH\p@\Cd@MH\Cd@@I\ifdim\Cd@qE>\z@\Cd@S
\Cd@@I\Cd@NH\fi\Cd@gI}\def\Cd@WI#1{\Cd@jD\count@\Cd@NH#1\ifnum\count@>\z@
\divide\Cd@NH\count@\fi\Cd@tH\Cd@zJ}\def\Cd@zJ{\dimen@\Cd@HC\count@\dimen@
\divide\count@5\divide\count@\Cd@NH\edef\Cd@mJ{\the\count@}}\def\Cd@MI{\Cd@EI
\axisheight\advance\Cd@EI-.5\Cd@@I\Cd@MH\Cd@EI\Cd@DI-\Cd@EI\advance\Cd@EI
\Cd@qE}\def\Cd@fJ{\Cd@DI\z@\Cd@EI\Cd@qE\relax}\def\horizhtdp{height\Cd@EI
depth\Cd@DI}\def\axisheight{\fontdimen22\the\textfont\tw@}\def
\script@axisheight{\fontdimen22\the\scriptfont\tw@}\def\ss@axisheight{%
\fontdimen22\the\scriptscriptfont\tw@}\def\Cd@FC{0.4pt}\def\Cd@WJ{\fontdimen3%
\textfont\z@}\def\Cd@VJ{\fontdimen3\textfont\z@}\newdimen\PileSpacing
\newdimen\Cd@hA\Cd@hA\z@\def\Cd@kA{\ifincommdiag1.3em\else2em\fi}\newdimen
\Cd@TB\def\CellSize{\afterassignment\Cd@fB\DiagramCellHeight}\newdimen
\DiagramCellHeight\DiagramCellHeight-\maxdimen\newdimen\DiagramCellWidth
\DiagramCellWidth-\maxdimen\def\Cd@fB{\DiagramCellWidth\DiagramCellHeight}%
\def\Cd@HC{3em}\newdimen\MapShortFall\def\MapsAbut{\MapShortFall\z@
\objectheight\z@\objectwidth\z@}\newdimen\Cd@cA\Cd@cA\z@\def\newarrowhead{%
\Cd@IG h\Cd@cF\Cd@hF>}\def\newarrowtail{\Cd@IG t\Cd@cF\Cd@hF>}\def
\newarrowmiddle{\Cd@IG m\Cd@cF\hbox@maths\empty}\def\newarrowfiller{\Cd@IG f%
\Cd@KE\Cd@OJ-}\def\Cd@IG#1#2#3#4#5#6#7#8#9{\Cd@IC{r#1:#5}{#2{#6}}\Cd@IC{l#1:#%
5}{#2{#7}}\Cd@IC{d#1:#5}{#3{#8}}\Cd@IC{u#1:#5}{#3{#9}}\Cd@jC{-#1:#5}{%
\expandafter\noexpand\csname-#1:#4\endcsname\noexpand\Cd@EC}\Cd@jC{+#1:#5}{%
\expandafter\noexpand\csname+#1:#4\endcsname\noexpand\Cd@EC}}\Cd@VA\Cd@EC{%
\Cd@IF\space diagonals are used unless PostScript is set}\def
\defaultarrowhead#1{\edef\Cd@vI{#1}\Cd@LI}\def\Cd@LI{\Cd@QI\Cd@vI<>ht\Cd@QI
\Cd@vI<>th}\def\Cd@QI#1#2#3#4#5{\Cd@OI{r#4}{#3}{l#5}{#2}{r#4:#1}\Cd@OI{r#5}{#%
2}{l#4}{#3}{l#4:#1}\Cd@OI{d#4}{#3}{u#5}{#2}{d#4:#1}\Cd@OI{d#5}{#2}{u#4}{#3}{u%
#4:#1}}\def\Cd@OI#1#2#3#4#5{\begingroup\aftergroup\Cd@PI\Cd@I{#1+:#2}\Cd@I{#1%
:#2}\Cd@I{#3:#4}\Cd@I{#5}\endgroup}\def\Cd@PI#1#2#3#4{\csname newbox%
\endcsname#1\def#2{\copy#1}\def#3{\copy#1}\setbox#1=\box\voidb@x}\def\Cd@vI{}%
\Cd@LI\def\Cd@PI#1#2#3#4{\setbox#1=#4}\ifx\tenln\nullfont\def\Cd@vI{vee}\else
\let\Cd@vI\Cd@IF\fi\def\Cd@XF#1#2#3{\begingroup\aftergroup\Cd@YF\Cd@I{#1#2:#3%
#3}\Cd@I{#1#2:#3}\endgroup}\def\Cd@YF#1#2{\def#1{\hbox{\rlap{#2}\kern.4\Cd@nC
#2}}}\Cd@XF rh>\Cd@XF lh>\Cd@XF rt>\Cd@XF lt>\Cd@XF rh<\Cd@XF lh<\Cd@XF rt<%
\Cd@XF lt<\def\Cd@YF#1#2{\def#1{\vbox{\vbox to\z@{#2\vss}\nointerlineskip
\kern.4\Cd@nC#2}}}\Cd@XF dh>\Cd@XF uh>\Cd@XF dt>\Cd@XF ut>\Cd@XF dh<\Cd@XF uh%
<\Cd@XF dt<\Cd@XF ut<\def\Cd@cF#1{\hbox{\mathsurround\z@\offinterlineskip
\Cd@h\mkern-1.5mu{#1}\mkern-1.5mu\Cd@zC}}\def\hbox@maths#1{\hbox{\Cd@h#1%
\Cd@zC}}\def\Cd@hF#1{\hbox to\Cd@qE{\setbox0=\hbox{\offinterlineskip
\mathsurround\z@\Cd@h{#1}\Cd@zC}\dimen0.5\wd0\advance\dimen0-.5\Cd@@I\Cd@MH{%
\dimen0}\kern-\dimen0\unhbox0\hss}}\def\Cd@dJ#1{\hbox to2\Cd@qE{\hss
\offinterlineskip\mathsurround\z@\Cd@h{#1}\Cd@zC\hss}}\def\Cd@WF#1{\hbox{%
\mathsurround\z@\Cd@h{#1}\Cd@zC}}\def\Cd@KE#1{\hbox{\kern-.15\Cd@nC\Cd@h{#1}%
\Cd@zC\kern-.15\Cd@nC}}\def\Cd@OJ#1{\vbox{\offinterlineskip\kern-.2ex\Cd@hF{#%
1}\kern-.2ex}}\def\@fillh{\xleaders\vrule\horizhtdp}\def\@fillv{\xleaders
\hrule width\Cd@qE}\Cd@QF{rf:-}{@fillh}\Cd@QF{lf:-}{@fillh}\Cd@QF{df:-}{%
@fillv}\Cd@QF{uf:-}{@fillv}\Cd@QF{rh:}{null}\Cd@QF{rm:}{null}\Cd@QF{rt:}{null%
}\Cd@QF{lh:}{null}\Cd@QF{lm:}{null}\Cd@QF{lt:}{null}\Cd@QF{dh:}{null}\Cd@QF{%
dm:}{null}\Cd@QF{dt:}{null}\Cd@QF{uh:}{null}\Cd@QF{um:}{null}\Cd@QF{ut:}{null%
}\Cd@QF{+h:}{null}\Cd@QF{+m:}{null}\Cd@QF{+t:}{null}\Cd@QF{-h:}{null}\Cd@QF{-%
m:}{null}\Cd@QF{-t:}{null}\Cd@IC{rf:}{\hbox{\kern1pt}}\Cd@QF{lf:}{rf:}\Cd@QF{%
+f:}{rf:}\Cd@IC{df:}{\vbox{\kern1pt}}\Cd@QF{uf:}{df:}\Cd@QF{-f:}{df:}\edef
\Cd@HG{\string\newarrow}\def\newarrow#1#2#3#4#5#6{\begingroup\edef\@name{#1}%
\edef\Cd@rI{#2}\edef\Cd@SD{#3}\edef\Cd@qF{#4}\edef\Cd@TD{#5}\edef\Cd@vD{#6}%
\let\Cd@rD\Cd@OG\let\Cd@HJ\Cd@YG\let\@x\Cd@XG\ifx\Cd@rI\Cd@SD\let\Cd@rI\empty
\fi\ifx\Cd@vD\Cd@TD\let\Cd@vD\empty\fi\def\Cd@ZH{r}\def\Cd@xE{l}\def\Cd@AC{d}%
\def\Cd@AJ{u}\def\Cd@zG{+}\def\@m{-}\ifx\Cd@SD\Cd@TD\ifx\Cd@qF\Cd@SD\let
\Cd@qF\empty\fi\ifx\Cd@vD\empty\ifx\Cd@SD\Cd@JE\let\@x\Cd@UG\else\let\@x
\Cd@VG\fi\fi\else\edef\Cd@X{\Cd@SD\Cd@rI}\ifx\Cd@X\empty\ifx\Cd@qF\Cd@TD\let
\Cd@qF\empty\fi\fi\fi\ifmmode\aftergroup\Cd@GG\else\Cd@w\Cd@jB rh{head\space
\space}\Cd@vD\Cd@jB rf{filler}\Cd@SD\Cd@jB rm{middle}\Cd@qF\ifx\Cd@TD\Cd@SD
\else\Cd@jB rf{filler}\Cd@TD\fi\Cd@jB rt{tail\space\space}\Cd@rI\Cd@ME\Cd@rD
\Cd@HJ\@x\Cd@JG l-2+2{lu}{nw}\NorthWest\Cd@JG r+2+2{ru}{ne}\NorthEast\Cd@JG l%
-2-2{ld}{sw}\SouthWest\Cd@JG r+2-2{rd}{se}\SouthEast\else\aftergroup\Cd@Y
\Cd@I{r\@name}\fi\fi\endgroup}\def\Cd@OG{\Cd@RG\Cd@ZH\Cd@xE rl\Horizontal@Map
}\def\Cd@YG{\Cd@RG\Cd@AC\Cd@AJ du\Vertical@Map}\def\Cd@XG{\Cd@RG\Cd@zG\@m+-%
\Vector@Map}\def\Cd@UG{\Cd@RG\Cd@zG\@m+-\Slant@Map}\def\Cd@VG{\Cd@RG\Cd@zG\@m
+-\Slope@Map}\catcode\lq\/=\active\def\Cd@RG#1#2#3#4#5{\Cd@FG#1#3#5t:\Cd@rI/f%
:\Cd@SD/m:\Cd@qF/f:\Cd@TD/h:\Cd@vD//\Cd@FG#2#4#5h:\Cd@vD/f:\Cd@TD/m:\Cd@qF/f:%
\Cd@SD/t:\Cd@rI//}\def\Cd@FG#1#2#3#4//{\edef\Cd@BG{#2}\aftergroup\sdef\Cd@I{#%
1\@name}\aftergroup{\aftergroup#3\Cd@J#4//\aftergroup}}\def\Cd@J#1/{\edef
\Cd@aG{#1}\ifx\Cd@aG\empty\else\Cd@I{\Cd@BG#1}\expandafter\Cd@J\fi}\catcode
\lq\/=12 \def\Cd@JG#1#2#3#4#5#6#7#8{\aftergroup\sdef\Cd@I{#6\@name}%
\aftergroup{\Cd@I{#2\@name}\if#2#4\aftergroup\Cd@QH\else\aftergroup\Cd@PH\fi
\Cd@I{#1\@name}
\aftergroup(\aftergroup#3\aftergroup,\aftergroup#5\aftergroup)\aftergroup}}%
\def\Cd@jB#1#2#3#4{\expandafter\ifx\csname#1#2:#4\endcsname\relax\Cd@u\Cd@bB{%
arrow#3 "#4" undefined}\fi}\Cd@NG\Cd@EE{All five components must be defined
before an arrow.}\Cd@NG\Cd@BE{\Cd@HG, unlike \string\HorizontalMap, is a
declaration.}\def\Cd@Y#1{\Cd@UA{Arrows \string#1 etc could not be defined}%
\Cd@EE}\def\Cd@GG{\Cd@UA{misplaced \Cd@HG}\Cd@BE}\def\newdiagramgrid#1#2#3{%
\Cd@IC{cdgh@#1}{#2,],}
\Cd@IC{cdgv@#1}{#3,],}}
\Cd@PG\ifincommdiag\incommdiagtrue\incommdiagfalse\Cd@PG\Cd@fE\Cd@nE\Cd@mE
\newcount\Cd@RA\Cd@RA=0 \def\Cd@LH{\Cd@RA6 }\def\Cd@JB{\Cd@RA1 \global\Cd@tA1
\Cd@oD\Cd@CF\empty}\def\Cd@CF{}\def\Cd@iB#1{\relax\Cd@xC\edef\Cd@yI{#1}%
\begingroup\Cd@XE\else\ifcase\Cd@RA\ifmmode\else\Cd@wF\Cd@E0\fi\or\Cd@LE5\or
\Cd@wF\Cd@F5\or\Cd@wF\Cd@B5\or\Cd@wF\Cd@B5\or\Cd@wF\Cd@C5\or\Cd@LE7\or\Cd@wF
\Cd@D7\fi\fi\endgroup\xdef\Cd@CF{#1}}\def\Cd@kB#1#2#3#4#5{\relax\Cd@xC\xdef
\Cd@yI{#4}\begingroup\ifnum\Cd@RA<#1 \expandafter\Cd@LE\ifcase\Cd@RA0\or#2\or
#3\else#2\fi\else\ifnum\Cd@RA<6 \Cd@wI\Cd@wF\Cd@B#2\else\Cd@wF\Cd@G#2\fi\fi
\endgroup\Cd@oD\Cd@CF\Cd@yI\ifincommdiag\let\Cd@JD#5\else\let\Cd@JD\Cd@NJ\fi}%
\def\Cd@JI{\global\Cd@tA=\ifnum\Cd@RA<5 1\else2\fi\relax}\def\Cd@cH{\Cd@RA
\Cd@tA}\def\Cd@LE#1{\aftergroup\Cd@RA\aftergroup#1\aftergroup\relax}\def
\Cd@dG{\def\Cd@iB##1{\relax}\let\Cd@kB\Cd@bG\let\Cd@LH\relax\let\Cd@JB\relax
\let\Cd@JI\relax\let\Cd@cH\relax}\def\Cd@bG#1#2#3#4#5{\ifincommdiag\let\Cd@JD
#5\else\xdef\Cd@yI{#4}\let\Cd@JD\Cd@NJ\fi}\def\Cd@wF#1{\aftergroup#1%
\aftergroup\relax\Cd@LE}\def\Cd@B{\Cd@HE\Cd@P\Cd@wD\Cd@N}\def\Cd@G{\Cd@HE{%
\Cd@mC\Cd@P}\Cd@GE\Cd@@D\Cd@N}\def\Cd@F{\Cd@HE{*\Cd@P}\Cd@AE\clubsuit\Cd@N}%
\def\Cd@C{\Cd@HE{\Cd@P*\Cd@P}\Cd@AE\Cd@N\clubsuit\Cd@N}\def\Cd@D{\Cd@HE\Cd@xB
\Cd@CE\\}\def\Cd@E{\Cd@HE\Cd@dC\Cd@@E\Cd@h}\def\Cd@NJ{\Cd@UA{\Cd@yI\space
ignored \Cd@yG}\Cd@FE}\def\Cd@qD{}\def\Cd@a{\Cd@UA{maps must never be enclosed
in braces}\Cd@yD}\def\Cd@yG{outside diagram}\def\Cd@yB{\string\HonV, \string
\VonH\space and \string\HmeetV}\Cd@NG\Cd@wD{The way that horizontal and
vertical arrows are terminated implicitly means\Cd@QG that they cannot be
mixed with each other or with \Cd@yB.}\Cd@NG\Cd@GE{\string\pile\space is for
parallel horizontal arrows; verticals can just be put together in\Cd@QG a cell%
. \Cd@yB\space are not meaningful in a \string\pile.}\Cd@NG\Cd@AE{The
horizontal maps must point to an object, not each other (I've put in\Cd@QG one
which you're unlikely to want). Use \string\pile\space if you want them
parallel.}\Cd@NG\Cd@CE{Parallel horizontal arrows must be in separate layers
of a \string\pile.}\Cd@NG\Cd@@E{Horizontal arrows may be used \Cd@yG s, but
must still be in maths.}\Cd@NG\Cd@FE{Vertical arrows, \Cd@yB\space\Cd@yG s don%
't know where\Cd@QG where to terminate.}\Cd@NG\Cd@yD{This prevents them from
stretching correctly.}\def\Cd@HE#1{\Cd@UA{"#1" inserted \ifx\Cd@CF\empty
before \Cd@yI\else between \Cd@CF\ifx\Cd@CF\Cd@yI s\else\space and \Cd@yI\fi
\fi}}\count@=\year\multiply\count@12 \advance\count@\month\ifnum\count@>2396500
\message{because this one expired in August 1996!}\expandafter\endinput\fi
\def\Horizontal@Map{\Cd@iB{horizontal map}\Cd@DC\Cd@YI\let\Cd@UD\hfdot\Cd@bD}%
\def\Cd@YI{\Cd@BB-9999 \let\Cd@JD\Cd@HD\ifincommdiag\else\Cd@gI\ifinpile\else
\skip2\z@ plus 1.5\Cd@WJ minus .5\Cd@VJ\skip4\skip2 \fi\fi\let\Cd@VD\@fillh}%
\def\Vector@Map{\Cd@JJ4}\def\Slant@Map{\Cd@JJ{\Cd@kE255\else6\fi}}\def
\Slope@Map{\Cd@JJ\Cd@mJ}\def\Cd@JJ#1#2#3#4#5#6{\Cd@DC\def\Cd@XJ{2}\def\Cd@bJ{%
2}\def\Cd@aJ{1}\def\Cd@cJ{1}\let\Horizontal@Map\Cd@zH\def\Cd@pF{#1}\def\Cd@bH
{\Cd@R#2#3#4#5#6}}\def\Cd@zH{\Cd@YI\Cd@EB\let\Cd@JD\Cd@DD\Cd@bD}\Cd@PG\Cd@VE
\Cd@mA\Cd@lA\Cd@mA\def\cds@missives{\Cd@mA}\def\Cd@DD{\Cd@bE\let\Cd@pF\Cd@mJ
\Cd@t\Cd@AF\setbox\z@\hbox{\incommdiagfalse\Cd@VH}\Cd@VE\Cd@KD\else\global
\Cd@PC\Cd@LD\fi\else\Cd@bH\Cd@aH\global\Cd@PC\Cd@ID\fi}\def\Cd@DC{\begingroup
\dimen1=\MapShortFall\dimen2=\Cd@kA\dimen5=\MapShortFall\setbox3=\box\voidb@x
\setbox6=\box\voidb@x\setbox7=\box\voidb@x\Cd@aD\mathsurround\z@\skip2\z@ plus%
1fill minus 1000pt\skip4\skip2 \Cd@OB}\Cd@PG\Cd@ZE\Cd@PB\Cd@OB\def\Cd@R#1#2#3%
#4#5{\let\Cd@rI#1\let\Cd@SD#2\let\Cd@qF#3\let\Cd@TD#4\let\Cd@vD#5\Cd@OB\ifx
\Cd@SD\Cd@TD\Cd@PB\fi}\def\Cd@bD#1#2#3#4#5{\Cd@R#1#2#3#4#5\Cd@eD}\def
\Vertical@Map{\Cd@kB433{vertical map}\Cd@MD\Cd@DC\Cd@BB-9995 \let\Cd@VD
\@fillv\let\Cd@UD\vfdot\Cd@bD}\def\break@args{\def\Cd@eD{\Cd@JD}\Cd@JD
\endgroup\aftergroup\Cd@qD}\def\Cd@TI{\setbox1=\Cd@rI\setbox5=\Cd@vD\ifvoid3
\ifx\Cd@qF\null\else\setbox3=\Cd@qF\fi\fi\Cd@aF2\Cd@SD\Cd@aF4\Cd@TD}\def
\Cd@aF#1#2{\ifx#2\Cd@VD\setbox#1=\box\voidb@x\else\setbox#1=#2\def#2{%
\xleaders\box#1}\fi}\Cd@VA\Cd@DJ{\string\HorizontalMap, \string\VerticalMap
\space and \string\DiagonalMap\Cd@QG are obsolete - use \Cd@HG\space to pre-%
define maps}\def\HorizontalMap#1#2#3#4#5{\Cd@DJ\Cd@iB{old horizontal map}%
\Cd@DC\Cd@YI\def\Cd@rI{\Cd@pG{#1}}\Cd@nG\Cd@SD{#2}\def\Cd@qF{\Cd@pG{#3}}%
\Cd@nG\Cd@TD{#4}\def\Cd@vD{\Cd@pG{#5}}\Cd@eD}\def\VerticalMap#1#2#3#4#5{%
\Cd@DJ\Cd@kB433{vertical map}\Cd@MD\Cd@DC\Cd@BB-9995 \let\Cd@VD\@fillv\def
\Cd@rI{\Cd@hF{#1}}\Cd@qG\Cd@SD{#2}\def\Cd@qF{\Cd@hF{#3}}\Cd@qG\Cd@TD{#4}\def
\Cd@vD{\Cd@hF{#5}}\Cd@eD}\def\DiagonalMap#1#2#3#4#5{\Cd@DJ\Cd@DC\def\Cd@pF{4}%
\let\Cd@VD\undefined\let\Cd@JD\Cd@ID\def\Cd@XJ{2}\def\Cd@bJ{2}\def\Cd@aJ{1}%
\def\Cd@cJ{1}\def\Cd@qF{\Cd@WF{#3}}\ifPositiveGradient\let\mv\raise\def\Cd@rI
{\Cd@WF{#5}}\def\Cd@SD{\Cd@WF{#4}}\def\Cd@TD{\Cd@WF{#2}}\def\Cd@vD{\Cd@WF{#1}%
}\else\let\mv\lower\def\Cd@rI{\Cd@WF{#1}}\def\Cd@SD{\Cd@WF{#2}}\def\Cd@TD{%
\Cd@WF{#4}}\def\Cd@vD{\Cd@WF{#5}}\fi\Cd@eD}\def\Cd@JE{-}\def\Cd@oC{\empty}%
\def\Cd@nG{\Cd@fF\Cd@KE\Cd@JE\@fillh}\def\Cd@qG{\Cd@fF\Cd@OJ\Cd@MJ\@fillv}%
\def\Cd@fF#1#2#3#4#5{\def\Cd@ZG{#5}\ifx\Cd@ZG#2\let#4#3\else\let#4\null\ifx
\Cd@ZG\empty\else\ifx\Cd@ZG\Cd@oC\else\let#4\Cd@ZG\fi\fi\fi}\def\Cd@pG#1{%
\hbox{\mathsurround\z@\offinterlineskip\def\Cd@ZG{#1}\ifx\Cd@ZG\empty\else
\ifx\Cd@ZG\Cd@oC\else\Cd@h\mkern-1.5mu{\Cd@ZG}\mkern-1.5mu\Cd@zC\fi\fi}}\def
\Cd@hD#1#2{\setbox#1=\hbox\bgroup\setbox0=\hbox{\Cd@h\labelstyle()\Cd@zC}%
\setbox1=\null\ht1\ht0\dp1\dp0\box1 \kern.1\Cd@nC\Cd@h\bgroup\labelstyle
\aftergroup\Cd@yC\Cd@iD}\def\Cd@yC{\Cd@zC\kern.1\Cd@nC\egroup\Cd@eD}\def
\Cd@iD{\futurelet\Cd@aG\Cd@pI}\def\Cd@pI{
\catcase\bgroup:\Cd@r;\catcase\egroup:\missing@label;\catcase\space:\Cd@yE;%
\tokcase[:\Cd@BF;
\default:\Cd@BJ;\endswitch}\def\Cd@r{\let\Cd@xC\Cd@Z\let\Cd@ZG}\def\Cd@BJ#1{%
\let\Cd@zE\egroup{\let\actually@braces@missing@around@macro@in@label\Cd@uG
\let\Cd@xC\Cd@lC\let\Cd@zE\Cd@@F#1%
\actually@braces@missing@around@macro@in@label}\Cd@zE}\def
\actually@braces@missing@around@macro@in@label{\let\Cd@ZG=}\def\missing@label
{\egroup\Cd@UA{missing label}\Cd@zD}\def\Cd@lC{\egroup\missing@label}\outer
\def\Cd@uG{}\def\Cd@zE{}\def\Cd@@F{\Cd@kC\Cd@zE}\def\Cd@xC{}\def\Cd@BF{\let
\Cd@K\Cd@iD\get@square@arg\Cd@lD}\Cd@NG\Cd@zD{The text which has just been
read is not allowed within map labels.}\def\Cd@Z{\egroup\Cd@UA{missing \Cd@mC
\space inserted after label}\Cd@zD}\def\upper@label{\Cd@ZD\Cd@hD6}\def
\lower@label{\def\positional@{\Cd@w\break@args}\Cd@hD7}\def\middle@label{%
\Cd@hD3}\Cd@PG\Cd@eE\Cd@aD\Cd@ZD\def\left@label{\ifPositiveGradient\Cd@wI
\expandafter\upper@label\else\expandafter\lower@label\fi}\def\right@label{%
\ifPositiveGradient\Cd@wI\expandafter\lower@label\else\expandafter
\upper@label\fi}\Cd@VA\Cd@lG{labels as positional arguments are obsolete}\def
\positional@{\Cd@lG\Cd@eE\Cd@wI\expandafter\upper@label\else\expandafter
\lower@label\fi-}\def\Cd@eD{\futurelet\Cd@aG\switch@arg}\def\eat@space{%
\afterassignment\Cd@eD\let\Cd@aG= }\def\Cd@yE{\afterassignment\Cd@iD\let
\Cd@aG= }\def\Cd@wB{\get@round@pair\Cd@fD}\def\Cd@fD#1#2{\def\Cd@XJ{#1}\def
\Cd@bJ{#2}\Cd@eD}\def\optional@{\let\Cd@K\Cd@eD\get@square@arg\Cd@lD}\def
\Cd@RI.{\Cd@oJ\Cd@eD}\def\Cd@oJ{\let\Cd@SD\Cd@UD\let\Cd@TD\Cd@UD\def\Cd@aH{%
\let\Cd@SD\dfdot\let\Cd@TD\dfdot}}\def\Cd@aH{}\def\Cd@kD#1,{\Cd@EH#1,%
\begingroup\ifx\@name\Cd@BD\Cd@lE\aftergroup\Cd@b\fi\aftergroup\Cd@nJ\else
\expandafter\def\expandafter\Cd@wE\expandafter{\csname\@name\endcsname}%
\expandafter\Cd@rJ\Cd@wE\Cd@pJ\ifx\Cd@wE\empty\aftergroup\Cd@eC\expandafter
\aftergroup\csname\Cd@AB\@name\endcsname\expandafter\aftergroup\csname\Cd@AB @%
\@name\endcsname\else\gdef\Cd@qJ{#1}\Cd@bB{\string\relax\space inserted before
`[\Cd@qJ'}\message{(I was trying to read it as an option.)}\aftergroup\Cd@lJ
\fi\fi\endgroup}\def\Cd@rJ#1#2\Cd@pJ{\def\Cd@wE{#2}}\def\Cd@nJ{\let\Cd@ZG
\Cd@K\let\Cd@K\relax\Cd@ZG}\def\Cd@lJ#1],{
\Cd@nJ\relax\def\Cd@wE{#1}\ifx\Cd@wE\empty\def\Cd@wE{[\Cd@qJ]}%
\else\def\Cd@wE{[\Cd@qJ,#1]}
\fi\Cd@wE}\def\Cd@eC#1#2{\ifx#2\undefined\ifx#1\undefined\Cd@bB{option `%
\@name' undefined}\else#1\fi\else\Cd@lE\expandafter#2\Cd@IJ\Cd@QJ\else\Cd@RJ
\fi\fi\Cd@tJ}\Cd@PG\Cd@lE\Cd@RJ\Cd@QJ\def\Cd@EH#1,{\Cd@lE\ifx\Cd@IJ\undefined
\Cd@b\else\expandafter\Cd@GH\Cd@IJ,#1,(,),(,)[]%
\fi\fi\Cd@lE\else\Cd@FH#1==,\fi}\def\Cd@b{\Cd@bB{option `\@name' needs (x,y)
value}\Cd@QJ\let\@name\empty}\def\Cd@FH#1=#2=#3,{\def\@name{#1}\def\Cd@IJ{#2}%
\def\Cd@wE{#3}\ifx\Cd@wE\empty\let\Cd@IJ\undefined\fi}%
\def\Cd@GH#1(#2,#3)#4,(#5,#6)#7[]{\def\Cd@IJ{{#2}{#3}}\def\Cd@wE{#1#4#5#6}%
\ifx\Cd@wE\empty\def\Cd@wE{#7}\ifx\Cd@wE\empty\Cd@b\fi\else\Cd@b\fi}\def
\Cd@AB{cds@}\let\Cd@K\relax\def\Cd@jD#1{\ifx\Cd@IJ\undefined\Cd@bB{option `%
\@name' needs a value}\else#1\Cd@IJ\relax\fi}\def\Cd@mD#1#2{\ifx\Cd@IJ
\undefined#1#2\relax\else#1\Cd@IJ\relax\fi}\def\cds@@showpair#1#2{\message{x=%
#1,y=#2}}\def\cds@@diagonalbase#1#2{\edef\Cd@aJ{#1}\edef\Cd@cJ{#2}}\def\Cd@RH
#1{\Cd@QF{@x}{cdps@#1}\ifx\@x\undefined\Cd@c{unknown}\else\ifx\@x\empty\Cd@c{%
cannot be used}\else\let\Cd@KJ\@x\fi\fi}\def\Cd@c#1{\Cd@bB{PostScript
translator `\Cd@IJ' #1}\Cd@UB\let\cds@PS\empty\let\cds@noPS\empty}\def\Cd@kG{%
}\def\Cd@VI{\Cd@ZA\edef\Cd@kG{\noexpand\Cd@FB{\@name\space ignored within
maths}}}\def\diagramstyle{\Cd@gI\let\Cd@K\relax\Cd@hE\Cd@lD\Cd@lD}\Cd@PG\Cd@YE\Cd@NB\Cd@MB\Cd@PG\Cd@WE\Cd@@B\Cd@zA
\Cd@PG\Cd@UE\Cd@jA\Cd@iA\Cd@PG\Cd@OE\Cd@DA\Cd@CA\Cd@DA\Cd@PG\Cd@PE\Cd@FA
\Cd@EA\Cd@PG\Cd@QE\Cd@HA\Cd@GA\Cd@PG\Cd@bE\Cd@VB\Cd@UB\Cd@PG\Cd@kE\Cd@FJ
\Cd@EJ\Cd@PG\Cd@XE\Cd@EB\Cd@DB\Cd@PG\Cd@SE\Cd@aA\Cd@ZA\Cd@PG\Cd@TE\Cd@eA
\Cd@dA\Cd@PG\Cd@gE\Cd@EG\Cd@DG\Cd@IC{cds@ }{}\Cd@IC{cds@}{}\Cd@IC{cds@1em}{%
\CellSize1\Cd@nC}\Cd@IC{cds@1.5em}{\CellSize1.5\Cd@nC}\Cd@IC{cds@2em}{%
\CellSize2\Cd@nC}\Cd@IC{cds@2.5em}{\CellSize2.5\Cd@nC}\Cd@IC{cds@3em}{%
\CellSize3\Cd@nC}\Cd@IC{cds@3.5em}{\CellSize3.5\Cd@nC}\Cd@IC{cds@4em}{%
\CellSize4\Cd@nC}\Cd@IC{cds@4.5em}{\CellSize4.5\Cd@nC}\Cd@IC{cds@5em}{%
\CellSize5\Cd@nC}\Cd@IC{cds@6em}{\CellSize6\Cd@nC}\Cd@IC{cds@7em}{\CellSize7%
\Cd@nC}\Cd@IC{cds@8em}{\CellSize8\Cd@nC}\def\cds@abut{\MapsAbut\dimen1\z@
\dimen5\z@}\def\cds@alignlabels{\Cd@EA\Cd@GA}\def\cds@amstex{\ifincommdiag
\Cd@L\else\def\CD{\diagram[amstex]}
\fi\Cd@Q\catcode\lq\@\active}\def\cds@b{\let\Cd@YB\Cd@WB}\def\cds@balance{%
\let\Cd@bA\Cd@x}\let\cds@bottom\cds@b\def\cds@center{\cds@vcentre
\cds@nobalance}\let\cds@centre\cds@center\def\cds@centerdisplay{\Cd@DA\Cd@VI
\cds@balance}\let\cds@centredisplay\cds@centerdisplay\def\cds@defaultsize{%
\Cd@mD{\let\Cd@HC}{3em}\Cd@zJ}\def\cds@displayoneliner{\Cd@zA}\let\cds@dotted
\Cd@oJ\def\cds@dpi{\Cd@WI{1truein}}\def\cds@dpm{\Cd@WI{100truecm}}\let\Cd@TA
\undefined\def\cds@eqno{\let\Cd@TA\Cd@IJ\let\Cd@xJ\empty}\def\cds@fixed{%
\Cd@lA}\def\cds@flushleft{\Cd@CA\Cd@VI\cds@nobalance\Cd@mD\Cd@hA\z@}\def
\cds@grid{\ifx\Cd@IJ\undefined\let\h@grid\relax\let\v@grid\relax\else\Cd@QF{%
h@grid}{cdgh@\Cd@IJ}\Cd@QF{v@grid}{cdgv@\Cd@IJ}\ifx\h@grid\relax\Cd@UA{%
unknown grid `\Cd@IJ'}\else\Cd@RB\fi\fi}\let\h@grid\relax\let\v@grid\relax
\def\cds@gridx{\ifx\Cd@IJ\undefined\else\cds@grid\fi\let\Cd@ZG\h@grid\let
\h@grid\v@grid\let\v@grid\Cd@ZG}\def\cds@h{\Cd@jD\DiagramCellHeight}\def
\cds@hcenter{\let\Cd@bA\Cd@WA}\let\cds@hcentre\cds@hcenter\def\cds@heads{%
\Cd@mD{\let\Cd@vI}\Cd@vI\Cd@LI\Cd@bE\else\ifx\Cd@vI\Cd@IF\else\Cd@EC\fi\fi}%
\let\cds@height\cds@h\let\cds@hmiddle\cds@balance\def\cds@htriangleheight{%
\Cd@mD\DiagramCellHeight\DiagramCellHeight\DiagramCellWidth1.73205%
\DiagramCellHeight}\def\cds@htrianglewidth{\Cd@mD\DiagramCellWidth
\DiagramCellWidth\DiagramCellHeight.57735\DiagramCellWidth}\def\cds@inline{%
\Cd@aA\let\Cd@kG\empty}\def\cds@inlineoneliner{\Cd@@B}\Cd@IC{cds@l>}{\Cd@jD{%
\let\Cd@kA}\dimen2=\Cd@kA}\def\cds@labelstyle{\Cd@jD{\let\labelstyle}}\def
\cds@landscape{\Cd@eA}\def\cds@large{\CellSize5\Cd@nC}\let\Cd@xJ\empty\def
\Cd@yJ{\refstepcounter{equation}\def\Cd@TA{\hbox{\@eqnnum}}}\def
\cds@LaTeXeqno{\let\Cd@xJ\Cd@yJ}\def\cds@lefteqno{\Cd@jA}\def
\cds@leftshortfall{\Cd@jD{\dimen1 }}\def\cds@lowershortfall{%
\ifPositiveGradient\cds@leftshortfall\else\cds@rightshortfall\fi}\def
\cds@loose{\Cd@QB}\def\cds@midhshaft{\Cd@FA}\def\cds@midshaft{\Cd@FA}\def
\cds@midvshaft{\Cd@bB{midvshaft option doesn't work}}\def\cds@moreoptions{%
\Cd@w}\let\cds@nobalance\cds@hcenter\def\cds@nohcheck{\Cd@dG}\def
\cds@nooptions{\def\Cd@RC{\Cd@FD}}\let\cds@noorigin\cds@nobalance\def
\cds@nopixel{\Cd@NH4\Cd@sG\Cd@gI}\def\cds@noPostScript{\Cd@mD\Cd@RH\empty
\Cd@UB\let\cds@PS\empty\let\cds@noPS\empty}\let\cds@noPS\Cd@UB\def
\cds@notextflow{\Cd@MB}\def\cds@noTPIC{\Cd@EJ}\def\cds@objectstyle{\Cd@jD{%
\let\objectstyle}}\def\cds@origin{\let\Cd@bA\Cd@dB}\def\cds@p{\Cd@jD
\PileSpacing}\let\cds@pilespacing\cds@p\def\cds@pixelsize{\Cd@jD\Cd@NH\Cd@tH}%
\def\cds@portrait{\Cd@dA}\def\cds@PostScript{\Cd@VB\let\cds@PS\Cd@VB\let
\cds@noPS\Cd@UB\Cd@mD\Cd@RH\empty}\let\cds@PS\Cd@VB\def
\cds@repositionpullbacks{\let\make@pbk\Cd@kJ\let\Cd@vJ\Cd@uJ}\def
\cds@righteqno{\Cd@iA}\def\cds@rightshortfall{\Cd@jD{\dimen5 }}\def
\cds@ruleaxis{\Cd@jD{\let\axisheight}}\def\cds@cmex{\let\Cd@hF\Cd@dJ\let
\Cd@MI\Cd@fJ}\def\cds@s{\cds@height\DiagramCellWidth\DiagramCellHeight}\def
\cds@scriptlabels{\let\labelstyle\scriptstyle}\def\cds@shortfall{\Cd@jD
\MapShortFall\dimen1\MapShortFall\dimen5\MapShortFall}\def\cds@showfirstpass{%
\Cd@mD{\let\Cd@YD}\z@}\def\cds@silent{\def\Cd@FB##1{}\def\Cd@bB##1{}}\let
\cds@size\cds@s\def\cds@small{\CellSize2\Cd@nC}\def\cds@t{\let\Cd@YB\Cd@aB}%
\def\cds@textflow{\Cd@NB\Cd@VI}\def\cds@thick{\let\Cd@SF\tenlnw\Cd@qE\Cd@FC
\Cd@mD\MapBreadth{2\Cd@qE}}\def\cds@thin{\let\Cd@SF\tenln\Cd@mD\MapBreadth{%
\Cd@FC}}\def\cds@tight{\Cd@RB}\let\cds@top\cds@t\def\cds@TPIC{\Cd@FJ}\def
\cds@uppershortfall{\ifPositiveGradient\cds@rightshortfall\else
\cds@leftshortfall\fi}\def\cds@vcenter{\let\Cd@YB\Cd@XB}\let\cds@vcentre
\cds@vcenter\def\cds@vtriangleheight{\Cd@mD\DiagramCellHeight
\DiagramCellHeight\DiagramCellWidth.577035\DiagramCellHeight}\def
\cds@vtrianglewidth{\Cd@mD\DiagramCellWidth\DiagramCellWidth
\DiagramCellHeight1.73205\DiagramCellWidth}\def\cds@vmiddle{\let\Cd@YB\Cd@ZB}%
\def\cds@w{\Cd@jD\DiagramCellWidth}\let\cds@width\cds@w\def\diagram{\relax
\protect\Cd@QC}\def\enddiagram{\protect\Cd@rF}\def\Cd@QC{\Cd@d\Cd@FI
\incommdiagtrue\edef\Cd@HI{\the\Cd@IB}\global\Cd@IB\z@\boxmaxdepth\maxdimen
\everycr{}\Cd@RC}\def\Cd@RC{\Cd@u\let\Cd@K\Cd@SC\Cd@hE\Cd@lD\Cd@FD}\def\Cd@SC
{\Cd@ME\expandafter\Cd@RC\else\expandafter\Cd@FD\fi}\def\Cd@FD{\let\Cd@aG
\relax\Cd@TE\Cd@bE\else\Cd@FB{landscape ignored without PostScript}\Cd@dA\fi
\fi\Cd@xJ\setbox2=\vbox\bgroup\Cd@oE\Cd@GD}\def\Cd@xG{\Cd@TE\Cd@aB\else\Cd@YB
\fi\Cd@bA\nointerlineskip\setbox0=\null\ht0-\Cd@AI\dp0\Cd@AI\wd0\Cd@wH\box0
\global\Cd@MA\Cd@NF\global\Cd@tA\Cd@SB\egroup\Cd@ZF\Cd@TE\setbox2=\hbox to\dp
2{\dp2=\Cd@MA\global\Cd@MA\ht2\ht2\wd2\global\Cd@EG\Cd@KJ{0 1 bturn}\box2%
\Cd@KJ{eturn}\hss}\Cd@zA\fi\ifnum\Cd@tA=1 \else\Cd@zA\fi\global\@ignorefalse
\Cd@SE\leavevmode\fi\ifvmode\Cd@PA\else\ifmmode\Cd@kG\Cd@UH\else\Cd@WE\Cd@aA
\fi\ifinner\Cd@aA\fi\Cd@SE\Cd@UH\else\Cd@YE\Cd@LB\else\Cd@PA\fi\fi\fi\fi
\Cd@ND}\def\Cd@ND{\global\Cd@IB\Cd@HI\relax\Cd@dE\global\Cd@vC\else
\aftergroup\Cd@cC\fi\if@ignore\aftergroup\ignorespaces\fi\Cd@kC\ignorespaces}%
\def\Cd@aB{\advance\Cd@AI\dimen1\relax}\def\Cd@ZB{\advance\Cd@AI.5\dimen1%
\relax}\def\Cd@WB{}\def\Cd@XB{\Cd@aB\advance\Cd@AI\Cd@TB\divide\Cd@AI2
\advance\Cd@AI-\axisheight\relax}\def\Cd@WA{}\def\Cd@dB{\Cd@NF\z@}\def\Cd@x{%
\ifdim\dimen2>\Cd@NF\Cd@NF\dimen2 \else\dimen2\Cd@NF\Cd@wH\dimen0 \advance
\Cd@wH\dimen2 \fi}\def\Cd@LB{\skip0\z@\relax\loop\skip1\lastskip\ifdim\skip1>%
\z@\unskip\advance\skip0\skip1 \repeat\vadjust{\prevdepth\dp\strutbox\penalty
\predisplaypenalty\vskip\abovedisplayskip\Cd@QA\penalty\postdisplaypenalty
\vskip\belowdisplayskip}\ifdim\skip0=\z@\else\hskip\skip0 \global\@ignoretrue
\fi}\def\Cd@PA{\Cd@mF\kern-\displayindent\Cd@QA\Cd@mF\global\@ignoretrue}\def
\Cd@QA{\hbox to\hsize{\setbox1=\hbox{\ifx\Cd@TA\undefined\else\Cd@h\Cd@TA
\Cd@zC\fi}\Cd@UE\Cd@OE\else\advance\Cd@MA\wd1 \fi\wd1\z@\box1 \fi\dimen0\wd2
\advance\dimen0\wd1 \advance\dimen0-\hsize\ifdim\dimen0>-\Cd@hA\Cd@DA\fi
\advance\dimen0\Cd@MA\ifdim\dimen0>\z@\Cd@FB{wider than the page by \the
\dimen0 }\Cd@DA\fi\Cd@OE\hss\else\Cd@S\Cd@MA\Cd@hA\fi\Cd@UH\hss\kern-\wd1\box
1 }}\def\Cd@UH{\Cd@gE\Cd@fE\else\Cd@JC\global\Cd@DG\fi\fi\kern\Cd@MA\box2 }%
\Cd@PG\Cd@cE\Cd@PC\Cd@OC\def\Cd@oE{\Cd@gI\ifdim\DiagramCellHeight=-\maxdimen
\DiagramCellHeight\Cd@HC\fi\ifdim\DiagramCellWidth=-\maxdimen
\DiagramCellWidth\Cd@HC\fi\global\Cd@OC\Cd@nE\let\Cd@qD\empty\let\Cd@v\Cd@N
\let\overprint\Cd@jJ\let\Cd@o\Cd@uI\let\enddiagram\Cd@rC\let\\\Cd@TC\let\par
\Cd@AH\let\Cd@xC\empty\let\switch@arg\Cd@KB\let\shift\Cd@cA\baselineskip
\DiagramCellHeight\lineskip\z@\lineskiplimit\z@\mathsurround\z@\tabskip\z@
\Cd@JB}\def\Cd@GD{\penalty-123 \begingroup\Cd@dA\aftergroup\Cd@H\halign
\bgroup\global\advance\Cd@IB1 \vadjust{\penalty1}\global\Cd@BA\z@\Cd@JB\Cd@f#%
#\Cd@qC\Cd@N\Cd@N\Cd@cH\Cd@f##\Cd@qC\cr}\def\Cd@rC{\Cd@xC\Cd@tC\crcr\egroup
\global\Cd@wC\endgroup}\def\Cd@f{\global\advance\Cd@BA1 \futurelet\Cd@aG\Cd@g
}\def\Cd@g{\ifx\Cd@aG\Cd@qC\Cd@wI\hskip1sp plus 1fil \relax\let\Cd@qC\relax
\Cd@GI\else\hfil\Cd@h\objectstyle\let\Cd@qD\Cd@a\fi}\def\Cd@qC{\Cd@xC\relax
\Cd@JI\Cd@GI\global\Cd@MA\Cd@cA\penalty-9993 \Cd@zC\hfil\null\kern-2\Cd@MA
\null}\def\Cd@TC{\cr}\def\across#1{\span\omit\mscount=#1 \global\advance
\Cd@BA\mscount\global\advance\Cd@BA\m@ne\Cd@TF\ifnum\mscount>2 \Cd@iI\repeat
\ignorespaces}\def\Cd@iI{\relax\span\omit\advance\mscount\m@ne}\def\Cd@tI{%
\ifincommdiag\ifx\Cd@SD\@fillh\ifx\Cd@TD\@fillh\ifdim\dimen3>\z@\else\ifdim
\dimen2>93\Cd@NH\ifdim\dimen2>18\p@\ifdim\Cd@qE>\z@\count@\Cd@fI\advance
\count@\m@ne\ifnum\count@<\z@\count@20\let\Cd@eI\Cd@xI\fi\xdef\Cd@fI{\the
\count@}\fi\fi\fi\fi\fi\fi\fi}\def\Cd@@G#1{\vrule\horizhtdp width#1\dimen@
\kern2\dimen@}\def\Cd@xI{\rlap{\dimen@\Cd@NH\Cd@S\dimen@{.182\p@}\Cd@MH
\dimen@\advance\Cd@EI\dimen@\Cd@@G0\Cd@@G0\Cd@@G2\Cd@@G6\Cd@@G6\Cd@@G2\Cd@@G0%
\Cd@@G0\Cd@@G2\Cd@@G6\Cd@@G0\Cd@@G0\Cd@@G2\Cd@@G2\Cd@@G6\Cd@@G0\Cd@@G0\Cd@@G2%
\Cd@@G6\Cd@@G2\Cd@@G2\Cd@@G0\Cd@@G0}}\def\Cd@fI{10}\def\Cd@eI{}\def\Cd@HD{%
\Cd@ME\Cd@OB\fi\Cd@t\Cd@AF\Cd@VH}\def\Cd@t{\Cd@MI\Cd@TI\ifvoid3 \setbox3=%
\null\ht3\Cd@EI\dp3\Cd@DI\else\Cd@S{\ht3}\Cd@EI\Cd@S{\dp3}\Cd@DI\fi\dimen3=.5%
\wd3 \ifdim\dimen3=\z@\Cd@ZE\else\dimen3-\Cd@sG\fi\else\Cd@OB\fi\Cd@S{\dimen2%
}{\wd7}\Cd@S{\dimen2}{\wd6}\Cd@tI\advance\dimen2-2\dimen3 \dimen4.5\dimen2
\dimen2\dimen4 \advance\dimen2-\wd1 \advance\dimen4-\wd5 \ifvoid2 \else\Cd@S{%
\ht3}{\ht2}\Cd@S{\dp3}{\dp2}\Cd@S{\dimen2}{\wd2}\fi\ifvoid4 \else\Cd@S{\ht3}{%
\ht4}\Cd@S{\dp3}{\dp4}\Cd@S{\dimen4}{\wd4}\fi\advance\skip2\dimen2 \advance
\skip4\dimen4 \Cd@ZE\advance\skip2\skip4 \dimen0\dimen5 \advance\dimen0\wd5
\skip3-\skip4 \advance\skip3-\dimen0 \let\Cd@TD\empty\else\skip3\z@\relax
\dimen0\z@\fi}\def\Cd@AF{\offinterlineskip\lineskip.2\Cd@nC\ifvoid6 \else
\setbox3=\vbox{\hbox to2\dimen3{\hss\box6\hss}\box3}\fi\ifvoid7 \else\setbox3%
=\vtop{\box3 \hbox to2\dimen3{\hss\box7\hss}}\fi}\def\Cd@VH{\kern\dimen1 \box
1 \Cd@eI\Cd@SD\hskip\skip2 \kern\dimen0 \ifincommdiag\Cd@PE\penalty1\fi\kern
\dimen3 \penalty\Cd@BB\hskip\skip3 \null\kern-\dimen3 \else\hskip\skip3 \fi
\box3 \Cd@TD\hskip\skip4 \box5 \kern\dimen5}\def\Cd@rE{\ifnum\Cd@gG>\Cd@KC
\Cd@S{\dimen1}\objectheight\Cd@S{\dimen5}\objectheight\else\Cd@S{\dimen1}%
\objectwidth\Cd@S{\dimen5}\objectwidth\fi}\def\Cd@V{\begingroup\ifdim\dimen7=%
\z@\kern\dimen8 \else\ifdim\dimen6=\z@\kern\dimen9 \else\dimen5\dimen6 \dimen
6\dimen9 \Cd@SI\dimen4\dimen2 \Cd@AG{\dimen4}\dimen6\dimen5 \dimen7\dimen8
\Cd@SI\Cd@ZC{\dimen2}\ifdim\dimen2<\dimen4 \kern\dimen2 \else\kern\dimen4 \fi
\fi\fi\endgroup}\def\Cd@mI{\Cd@XH\setbox\z@\hbox{\lower\axisheight\hbox to%
\dimen2{\Cd@jE\ifPositiveGradient\dimen8\ht\Cd@hG\dimen9\Cd@yH\else\dimen8\dp
3 \dimen9\dimen1 \fi\else\dimen8 \ifPositiveGradient\objectheight\else\z@\fi
\dimen9\objectwidth\fi\advance\dimen8 \ifPositiveGradient-\fi\axisheight\Cd@V
\unhbox\z@\Cd@jE\ifPositiveGradient\dimen8\dp3 \dimen9\dimen0 \else\dimen8\ht
\Cd@hG\dimen9\Cd@PF\fi\else\dimen8 \ifPositiveGradient\z@\else\objectheight
\fi\dimen9\objectwidth\fi\advance\dimen8 \ifPositiveGradient\else-\fi
\axisheight\Cd@V}}}\def\Cd@LD{\dimen6 \Cd@bJ\DiagramCellHeight\dimen7 \Cd@XJ
\DiagramCellWidth\Cd@mI\ifPositiveGradient\advance\dimen7-\Cd@aJ
\DiagramCellWidth\else\dimen7 \Cd@aJ\DiagramCellWidth\dimen6\z@\multiply
\Cd@gG\m@ne\fi\advance\dimen6-\Cd@cJ\DiagramCellHeight\setbox0=\rlap{\global
\Cd@EG\kern-\dimen7 \lower\dimen6\hbox{\Cd@OD{\the\Cd@KC\space\the\Cd@gG
\space bturn}\box0 \Cd@KJ{eturn}}}\ht0\z@\dp0\z@\raise\axisheight\box0 }\def
\Cd@pB{\advance\Cd@LF-\Cd@yH\Cd@zI\Cd@LF\advance\Cd@zI\Cd@uH\ifvoid\Cd@HH
\ifdim\Cd@zI<.1em\ifnum\Cd@QD=\@m\else\Cd@yF h\Cd@zI<.1em:objects overprint:%
\Cd@BA\Cd@QD\fi\fi\else\ifhbox\Cd@HH\Cd@TJ\else\Cd@UJ\fi\advance\Cd@zI\Cd@yH
\Cd@wG{-\Cd@yH}{\box\Cd@HH}{\Cd@zI}\z@\fi\Cd@LF-\Cd@PF\Cd@QD\Cd@BA\Cd@uH\z@}%
\def\Cd@TJ{\setbox\Cd@HH=\hbox{\unhbox\Cd@HH\unskip\unpenalty}\setbox\Cd@IH=%
\hbox{\unhbox\Cd@IH\unskip\unpenalty}\setbox\Cd@HH=\hbox to\Cd@zI{\Cd@KA\wd
\Cd@HH\unhbox\Cd@HH\Cd@LA\lastkern\unkern\ifdim\Cd@LA=\z@\Cd@PB\advance\Cd@KA
-\wd\Cd@IH\else\Cd@OB\fi\ifnum\lastpenalty=\z@\else\Cd@FA\unpenalty\fi\kern
\Cd@LA\ifdim\Cd@LF<\Cd@KA\Cd@FA\fi\ifdim\Cd@uH<\wd\Cd@IH\Cd@FA\fi\Cd@PE\Cd@uH
\Cd@zI\advance\Cd@uH-\Cd@KA\advance\Cd@uH\wd\Cd@IH\ifdim\Cd@uH<2\wd\Cd@IH
\Cd@yF h\Cd@uH<2\wd\Cd@IH:arrow too short:\Cd@BA\Cd@QD\fi\divide\Cd@uH\tw@
\Cd@LF\Cd@zI\advance\Cd@LF-\Cd@uH\fi\Cd@ZE\kern-\Cd@uH\fi\hbox to\Cd@uH{%
\unhbox\Cd@IH}\Cd@iF}}\Cd@PG\ifinpile\inpiletrue\inpilefalse\inpilefalse\def
\pile{\protect\Cd@ZI\protect\Cd@JH}\def\Cd@JH#1{\Cd@i#1\Cd@@D}\def\Cd@ZI{%
\Cd@iB{pile}\setbox0=\vtop\bgroup\aftergroup\Cd@WD\inpiletrue\let\Cd@qD\empty
\let\pile\Cd@pE\let\Cd@@D\Cd@AD\let\Cd@tC\Cd@sC\Cd@LH\baselineskip.5%
\PileSpacing\lineskip.1\Cd@nC\relax\lineskiplimit\lineskip\mathsurround\z@
\tabskip\z@\let\\\Cd@KH}\def\Cd@i{\Cd@oD\Cd@CF\empty\halign\bgroup\hfil\Cd@h
\let\Cd@qD\Cd@a##\Cd@xC\Cd@zC\hfil\Cd@N\Cd@O##\cr}\Cd@NG\Cd@xD{pile only
allows one column.}\Cd@NG\Cd@DE{you left it out!}\def\Cd@O{\Cd@@D\Cd@N\relax
\Cd@UA{missing \Cd@mC\space inserted after \string\pile}\Cd@xD}\def\Cd@AD{%
\Cd@xC\crcr\egroup\egroup}\def\Cd@tC{\Cd@xC}\def\Cd@sC{\Cd@xC\relax\Cd@@D
\Cd@UA{missing \Cd@mC\space inserted between \string\pile\space and \Cd@uC}%
\Cd@DE}\def\Cd@@D{\Cd@xC}\def\Cd@WD{\vbox{\dimen1\dp0 \unvbox0 \setbox0=%
\lastbox\advance\dimen1\dp0 \nointerlineskip\box0 \nointerlineskip\setbox0=%
\null\dp0.5\dimen1\ht0-\dp0 \box0}\ifincommdiag\Cd@wI\penalty-9998 \fi\xdef
\Cd@CF{pile}}\def\Cd@KH{\cr}\def\Cd@pE#1{#1}\def\Cd@UJ{\setbox\Cd@HH=\vbox{%
\unvbox\Cd@HH\setbox1=\lastbox\setbox0=\box\voidb@x\Cd@UF\setbox\Cd@HH=%
\lastbox\ifhbox\Cd@HH\Cd@gC\repeat\unvbox0 \global\Cd@MA\Cd@IE}\Cd@IE\Cd@MA}%
\def\Cd@gC{\Cd@PE\setbox\Cd@HH=\hbox{\unhbox\Cd@HH\unskip\setbox\Cd@HH=%
\lastbox\unskip\unhbox\Cd@HH}\ifdim\Cd@zI<\wd\Cd@HH\Cd@yF h\Cd@zI<\wd\Cd@HH:%
arrow in pile too short:\Cd@BA\Cd@QD\else\setbox\Cd@HH=\hbox to\Cd@zI{\unhbox
\Cd@HH}\fi\else\Cd@jI\fi\setbox0=\vbox{\box\Cd@HH\nointerlineskip\ifvoid0
\Cd@wI\box1 \else\vskip\skip0 \unvbox0 \fi}\skip0=\lastskip\unskip}\def\Cd@jI
{\penalty7 \noindent\unhbox\Cd@HH\unskip\setbox\Cd@HH=\lastbox\unskip\unhbox
\Cd@HH\endgraf\setbox\Cd@IH=\lastbox\unskip\setbox\Cd@IH=\hbox{\Cd@kF\unhbox
\Cd@IH\unskip\unskip\unpenalty}\ifcase\prevgraf\cd@shouldnt P\or\ifdim\Cd@zI<%
\wd\Cd@IH\Cd@yF h\Cd@zI<\wd\Cd@HH:object in pile too wide:\Cd@BA\Cd@QD\setbox
\Cd@HH=\hbox to\Cd@zI{\hss\unhbox\Cd@IH\hss}\else\setbox\Cd@HH=\hbox to\Cd@zI
{\hss\kern\Cd@LF\unhbox\Cd@IH\kern\Cd@uH\hss}\fi\or\setbox\Cd@HH=\lastbox
\unskip\Cd@TJ\else\cd@shouldnt Q\fi\unskip\unpenalty}\def\Cd@MD{\Cd@TI\ifvoid
3 \setbox3=\null\ht3\axisheight\dp3-\ht3 \dimen3.5\Cd@qE\else\dimen4\dp3
\dimen3.5\wd3 \setbox3=\Cd@hF{\box3}\dp3\dimen4 \ifdim\ht3=-\dp3 \else\Cd@OB
\fi\fi\setbox0=\null\Cd@ZE\dimen4=\ht\Cd@hG\advance\dimen4\dp5 \advance\dimen
4\dimen1 \let\Cd@TD\empty\else\dimen4\ht3 \fi\ht0\dimen4 \offinterlineskip
\setbox8=\vbox spread2ex{\kern\dimen5 \box1 \Cd@SD\vfill\box0}\ht8=\z@\setbox
9=\vtop spread2ex{\kern-\ht3 \box3 \Cd@TD\vfill\box5 \kern\dimen1}\dp9=\z@
\dimen0\dimen3 \advance\dimen0-.5\Cd@qE\hskip\z@ plus.0001fil \box6 \kern
\dimen0 \box9 \kern-\Cd@qE\box8 \Cd@ZE\penalty1 \fi\kern\PileSpacing\kern-%
\PileSpacing\kern-.5\Cd@qE\penalty\Cd@BB\null\kern\dimen3 \box7}\def\Cd@qH{%
\ifhbox\Cd@RA\Cd@FB{clashing verticals}\ht\Cd@hG.5\dp\Cd@RA\dp\Cd@hG-\ht5
\Cd@rB\ht\Cd@hG\z@\dp\Cd@hG\z@\fi\dimen1\dp\Cd@RA\Cd@sA\prevgraf\unvbox\Cd@RA
\Cd@rA\lastpenalty\unpenalty\setbox\Cd@RA=\null\setbox\Cd@xH=\hbox{\Cd@kF
\unhbox\Cd@xH\unskip\unpenalty\dimen0\lastkern\unkern\unkern\unkern\kern
\dimen0 \Cd@iF}\setbox\Cd@OF=\hbox{\unhbox\Cd@OF\dimen0\lastkern\unkern
\unkern\global\Cd@MA\lastkern\unkern\kern\dimen0 }\Cd@UF\ifnum\Cd@sA>4 \Cd@KI
\repeat\unskip\unskip\advance\Cd@PF.5\wd\Cd@RA\advance\Cd@PF\wd\Cd@OF\advance
\Cd@yH.5\wd\Cd@RA\advance\Cd@yH\wd\Cd@xH\setbox\Cd@RA=\hbox{\kern-\Cd@PF\box
\Cd@OF\unhbox\Cd@RA\box\Cd@xH\kern-\Cd@yH\penalty\Cd@rA\penalty\Cd@IB}\ht
\Cd@RA\dimen1 \dp\Cd@RA\z@\wd\Cd@RA\Cd@nB\Cd@pB}\def\Cd@KI{\ifdim\wd\Cd@OF<%
\Cd@MA\setbox\Cd@OF=\hbox to\Cd@MA{\Cd@kF\unhbox\Cd@OF}\fi\advance\Cd@sA\m@ne
\setbox\Cd@RA=\hbox{\box\Cd@OF\unhbox\Cd@RA}\unskip\setbox\Cd@OF=\lastbox
\setbox\Cd@OF=\hbox{\unhbox\Cd@OF\unskip\unpenalty\dimen0\lastkern\unkern
\unkern\global\Cd@MA\lastkern\unkern\kern\dimen0 }}\def\Cd@rB{\dimen1\dp
\Cd@RA\ifhbox\Cd@RA\Cd@sB\else\Cd@tB\fi\setbox\Cd@RA=\vbox{\penalty\Cd@IB}\dp
\Cd@RA-\dp\Cd@hG\wd\Cd@RA\Cd@nB}\def\Cd@tB{\unvbox\Cd@RA\Cd@rA\lastpenalty
\unpenalty\ifdim\dimen1<\ht\Cd@hG\Cd@yF v\dimen1<\ht\Cd@hG:rows overprint:%
\Cd@IB\Cd@rA\fi}\def\Cd@sB{\dimen0=\ht\Cd@RA\setbox\Cd@RA=\hbox\bgroup
\advance\dimen1-\ht\Cd@hG\unhbox\Cd@RA\Cd@sA\lastpenalty\unpenalty\Cd@rA
\lastpenalty\unpenalty\global\Cd@NA-\lastkern\unkern\setbox0=\lastbox\Cd@UF
\setbox\Cd@RA=\hbox{\box0\unhbox\Cd@RA}\setbox0=\lastbox\ifhbox0 \Cd@nI
\repeat\global\Cd@OA-\lastkern\unkern\global\Cd@MA\Cd@LJ\unhbox\Cd@RA\egroup
\Cd@LJ\Cd@MA\Cd@wG{\Cd@OA}{\box\Cd@RA}{\Cd@NA}{\dimen1}}\def\Cd@nI{\setbox0=%
\hbox to\wd0\bgroup\unhbox0 \unskip\unpenalty\dimen7\lastkern\unkern\ifnum
\lastpenalty=1 \unpenalty\Cd@PB\else\Cd@OB\fi\setbox0=\lastbox\dimen6%
\lastkern\unkern\setbox1=\lastbox\setbox0=\vbox{\unvbox0\Cd@ZE\kern-\dimen1%
\fi}\ifdim\dimen0<\ht0 \Cd@yF v\dimen0<\ht0:upper part of vertical too short:%
{\Cd@ZE\Cd@IB\else\Cd@rA\fi}\Cd@sA\else\setbox0=\vbox to\dimen0{\unvbox0}\fi
\setbox1=\vtop{\unvbox1}\ifdim\dimen1<\dp1 \Cd@yF v\dimen1<\dp1:lower part of
vertical too short:\Cd@IB\Cd@rA\else\setbox1=\vtop to\dimen1{\unvbox1}\fi\box
1 \kern\dimen6 \box0 \kern\dimen7 \Cd@iF\global\Cd@MA\Cd@LJ\egroup\Cd@LJ
\Cd@MA\relax}\countdef\Cd@q=14 \newcount\Cd@z\newcount\Cd@SB\newcount\Cd@IB
\let\Cd@GB\insc@unt\newcount\Cd@BA\newcount\Cd@fA\let\Cd@gA\Cd@SB\newcount
\Cd@HB\Cd@PG\Cd@jE\Cd@pH\Cd@oH\Cd@oH\def\Cd@YD{-1}\def\Cd@H{\ifnum\Cd@YD<\z@
\else\begingroup\scrollmode\showboxdepth\Cd@YD\showboxbreadth\maxdimen
\showlists\endgroup\fi\Cd@pH\Cd@ZF\Cd@z=\Cd@q\advance\Cd@z1 \Cd@SB=\Cd@z
\ifnum\Cd@IB=1 \Cd@FA\fi\advance\Cd@SB\Cd@IB\dimen1\z@\skip0\z@\count@=%
\insc@unt\advance\count@\Cd@q\divide\count@2 \ifnum\Cd@SB>\count@\Cd@FB{The
diagram has too many rows! It can't be reformatted.}\else\Cd@oF\Cd@kH\fi
\Cd@xG}\def\Cd@oF{\Cd@IB\Cd@z\Cd@VF\ifnum\Cd@IB<\Cd@SB\setbox\Cd@IB\box
\voidb@x\advance\Cd@IB1\relax\repeat\Cd@IB\Cd@z\skip\z@\z@\Cd@VF\Cd@BB
\lastpenalty\unpenalty\ifnum\Cd@BB>\z@\Cd@uD\repeat\ifnum\Cd@BB=-123 \Cd@wI
\unpenalty\else\cd@shouldnt D\fi\ifx\v@grid\relax\else\Cd@IB\Cd@SB\advance
\Cd@IB\m@ne\expandafter\Cd@hJ\v@grid\fi\Cd@HB\Cd@gA\Cd@nB\z@\Cd@vF\ifx\h@grid
\relax\else\expandafter\Cd@gJ\h@grid\fi\count@\Cd@SB\advance\count@\m@ne
\Cd@TB\ht\count@}\def\Cd@uD{\ifcase\Cd@BB\or\Cd@nF\else\Cd@pA-\lastpenalty
\unpenalty\Cd@qA\lastpenalty\unpenalty\setbox0=\lastbox\Cd@uF\fi\Cd@gD}\def
\Cd@gD{\skip1\lastskip\unskip\advance\skip0\skip1 \ifdim\skip1=\z@\else
\expandafter\Cd@gD\fi}\def\Cd@nF{\setbox0=\lastbox\Cd@AI\dp0 \advance\Cd@AI
\skip\z@\skip\z@\z@\advance\Cd@sE\Cd@AI\Cd@aE\ifnum\Cd@IB>\Cd@z\Cd@sE
\DiagramCellHeight\Cd@AI\Cd@sE\advance\Cd@AI-\Cd@BI\fi\fi\Cd@BI\ht0 \Cd@sE
\Cd@BI\setbox\Cd@IB\hbox{\unhbox\Cd@IB\unhbox0}\dp\Cd@IB\Cd@AI\ht\Cd@IB\Cd@BI
\advance\Cd@IB1 }\def\Cd@uF{\ifnum\Cd@pA<\z@\advance\Cd@pA\Cd@SB\ifnum\Cd@pA<%
\Cd@z\Cd@sJ\else\Cd@KA\dp\Cd@pA\Cd@LA\ht\Cd@pA\setbox\Cd@pA\hbox{\box\z@
\penalty\Cd@qA\penalty\Cd@BB\unhbox\Cd@pA}\dp\Cd@pA\Cd@KA\ht\Cd@pA\Cd@LA\fi
\else\Cd@sJ\fi}\def\Cd@sJ{\Cd@FB{diagonal goes outside diagram (lost)}}\def
\Cd@eJ{\advance\Cd@pA\Cd@SB\ifnum\Cd@pA<\Cd@z\Cd@sJ\else\ifnum\Cd@pA=\Cd@IB
\Cd@tF\else\ifnum\Cd@pA>\Cd@IB\cd@shouldnt M\else\Cd@KA\dp\Cd@pA\Cd@LA\ht
\Cd@pA\setbox\Cd@pA\hbox{\box\z@\penalty\Cd@qA\penalty\Cd@BB\unhbox\Cd@pA}\dp
\Cd@pA\Cd@KA\ht\Cd@pA\Cd@LA\fi\fi\fi}\def\Cd@kH{\Cd@p\Cd@MI\setbox\Cd@GC=%
\hbox{\Cd@h A\@super f\Cd@oI f\Cd@zC}\Cd@IE\z@\Cd@LJ\z@\Cd@wH\z@\Cd@NF\z@
\Cd@IB=\Cd@SB\Cd@sE\z@\Cd@oB\z@\Cd@VF\ifnum\Cd@IB>\Cd@z\advance\Cd@IB\m@ne
\Cd@BI\ht\Cd@IB\Cd@AI\dp\Cd@IB\advance\Cd@sE\Cd@BI\Cd@CI\advance\Cd@oB\Cd@sE
\Cd@CC\Cd@nH\Cd@s\ht\Cd@IB\Cd@BI\dp\Cd@IB\Cd@AI\nointerlineskip\box\Cd@IB
\Cd@sE\Cd@AI\setbox\Cd@IB\null\ht\Cd@IB\Cd@oB\repeat\Cd@qB\nointerlineskip
\box\Cd@IB\Cd@CG\Cd@IE\DiagramCellWidth{width}\Cd@CG\Cd@LJ\DiagramCellHeight{%
height}\Cd@RA\Cd@GB\advance\Cd@RA-\Cd@fA\advance\Cd@RA\m@ne\advance\Cd@RA
\Cd@gA\dimen0\wd\Cd@RA\Cd@EI\axisheight\dimen1\Cd@oB\advance\dimen1-\Cd@TB
\dimen2\Cd@wH\advance\dimen2-\dimen0 \advance\Cd@SB-\Cd@z\advance\Cd@GB-%
\Cd@fA}\count@\year\multiply\count@12 \advance\count@\month\ifnum\count@>%
23966000 \loop\iftrue\message{gone February 1997!}\repeat\fi\def\Cd@qB{\Cd@BI-%
\Cd@sE\Cd@AI\Cd@sE\setbox\Cd@hG=\null\dp\Cd@hG\Cd@sE\ht\Cd@hG-\Cd@sE\Cd@PF\z@
\Cd@yH\z@\Cd@fA\Cd@GB\advance\Cd@fA-\Cd@HB\advance\Cd@fA\Cd@gA\Cd@BA\Cd@GB
\Cd@RA\Cd@HB\Cd@TF\ifnum\Cd@BA>\Cd@fA\advance\Cd@BA\m@ne\advance\Cd@RA\m@ne
\Cd@nB\wd\Cd@RA\setbox\Cd@BA=\box\voidb@x\Cd@rB\repeat\Cd@s\ht\Cd@IB\Cd@BI\dp
\Cd@IB\Cd@AI}\def\Cd@CG#1#2#3{\ifdim#1>.01\Cd@nC\Cd@LA#2\relax\advance\Cd@LA#%
1\relax\advance\Cd@LA.99\Cd@nC\count@\Cd@LA\divide\count@\Cd@nC\Cd@FB{%
increase cell #3 to \the\count@ em}\fi}\def\Cd@CI{\Cd@BA=\Cd@GB\penalty4
\noindent\unhbox\Cd@IB\Cd@TF\unskip\setbox0=\lastbox\ifhbox0 \advance\Cd@BA
\m@ne\setbox\Cd@BA\hbox to\wd0{\null\penalty-9990\null\unhbox0}\repeat\Cd@fA
\Cd@BA\advance\Cd@BA\Cd@HB\advance\Cd@BA-\Cd@gA\ifnum\Cd@BA<\Cd@GB\count@
\Cd@BA\advance\count@\m@ne\dimen0=\wd\count@\count@\Cd@HB\advance\count@\m@ne
\Cd@nB\wd\count@\Cd@TF\ifnum\Cd@BA<\Cd@GB\Cd@NI\Cd@vF\dimen0\wd\Cd@BA\advance
\Cd@BA1 \repeat\fi\Cd@TF\Cd@BB\lastpenalty\unpenalty\ifnum\Cd@BB>\z@\Cd@qA
\lastpenalty\unpenalty\Cd@tF\repeat\endgraf\unskip\ifnum\lastpenalty=4
\unpenalty\else\cd@shouldnt S\fi}\def\Cd@tF{\advance\Cd@qA\Cd@fA\advance
\Cd@qA\m@ne\setbox0=\lastbox\ifnum\Cd@qA<\Cd@GB\setbox\Cd@qA\hbox{\box0%
\penalty\Cd@BB\unhbox\Cd@qA}\else\Cd@sJ\fi}\def\Cd@zF{}\Cd@PG\Cd@aE\Cd@RB
\Cd@QB\def\Cd@NI{\advance\dimen0\wd\Cd@BA\divide\dimen0\tw@\Cd@aE\dimen0%
\DiagramCellWidth\else\Cd@S{\dimen0}\DiagramCellWidth\Cd@sI\fi\advance\Cd@nB
\dimen0 }\def\Cd@vF{\setbox\Cd@HB=\vbox{}\dp\Cd@HB=\Cd@oB\wd\Cd@HB\Cd@nB
\advance\Cd@HB1 }\def\Cd@gJ#1,{\def\Cd@IJ{#1}\ifx\Cd@IJ\Cd@BD\else\advance
\Cd@nB\Cd@IJ\DiagramCellWidth\Cd@vF\expandafter\Cd@gJ\fi}\def\Cd@hJ#1,{\def
\Cd@IJ{#1}\ifx\Cd@IJ\Cd@BD\else\ifnum\Cd@IB>\Cd@z\Cd@sE\Cd@IJ
\DiagramCellHeight\advance\Cd@sE-\dp\Cd@IB\advance\Cd@IB\m@ne\ht\Cd@IB\Cd@sE
\fi\expandafter\Cd@hJ\fi}\def\Cd@sI{\Cd@cE\Cd@KA\dimen0 \advance\Cd@KA-%
\DiagramCellWidth\ifdim\Cd@KA>2\MapShortFall\Cd@FB{badly drawn diagonals (see
manual)}\let\Cd@sI\empty\fi\else\let\Cd@sI\empty\fi}\def\Cd@CC{\Cd@RA\Cd@gA
\Cd@TF\ifnum\Cd@RA<\Cd@HB\dimen0\dp\Cd@RA\advance\dimen0\Cd@sE\dp\Cd@RA\dimen
0 \advance\Cd@RA1 \repeat}\def\Cd@wG#1#2#3#4{\ifnum\Cd@BA<\Cd@GB\Cd@KA=#1%
\relax\setbox\Cd@BA=\hbox{\setbox0=#2\dimen7=#4\relax\dimen8=#3\relax\ifhbox
\Cd@BA\unhbox\Cd@BA\advance\Cd@KA-\lastkern\unkern\fi\ifdim\Cd@KA=\z@\else
\kern-\Cd@KA\fi\raise\dimen7\box0 \kern-\dimen8 }\ifnum\Cd@BA=\Cd@fA\Cd@S
\Cd@NF\Cd@KA\fi\else\cd@shouldnt O\fi}\def\Cd@s{\setbox\Cd@IB=\hbox{\Cd@BA
\Cd@fA\Cd@RA\Cd@gA\Cd@LA\z@\relax\Cd@TF\ifnum\Cd@BA<\Cd@GB\Cd@nB\wd\Cd@RA
\relax\Cd@sH\advance\Cd@BA1 \advance\Cd@RA1 \repeat}\Cd@S\Cd@wH{\wd\Cd@IB}\wd
\Cd@IB\z@}\def\Cd@sH{\ifhbox\Cd@BA\Cd@KA\Cd@nB\relax\advance\Cd@KA-\Cd@LA
\relax\ifdim\Cd@KA=\z@\else\kern\Cd@KA\fi\Cd@LA\Cd@nB\advance\Cd@LA\wd\Cd@BA
\relax\unhbox\Cd@BA\advance\Cd@LA-\lastkern\unkern\fi}\def\Cd@nH{\setbox
\Cd@HH=\box\voidb@x\Cd@RA=\Cd@HB\Cd@BA\Cd@GB\Cd@RA\Cd@gA\advance\Cd@RA\Cd@BA
\advance\Cd@RA-\Cd@fA\advance\Cd@RA\m@ne\Cd@nB\wd\Cd@RA\count@\Cd@GB\advance
\count@\m@ne\Cd@LF.5\wd\count@\advance\Cd@LF\Cd@nB\Cd@A\m@ne\Cd@QD\@m\Cd@TF
\ifnum\Cd@BA>\Cd@fA\advance\Cd@BA\m@ne\advance\Cd@LF-\Cd@nB\Cd@dH\wd\Cd@RA
\Cd@nB\advance\Cd@LF\Cd@nB\advance\Cd@RA\m@ne\Cd@nB\wd\Cd@RA\repeat\Cd@PF
\Cd@NF\Cd@yH-\Cd@PF\Cd@pB}\newcount\Cd@BB\def\Cd@o{}\def\Cd@p{\mathsurround
\z@\hsize\z@\rightskip\z@ plus1fil minus\maxdimen\parfillskip\z@\linepenalty
9000 \looseness0 \hfuzz\maxdimen\hbadness10000 \clubpenalty0 \widowpenalty0
\displaywidowpenalty0 \interlinepenalty0 \predisplaypenalty0
\postdisplaypenalty0 \interdisplaylinepenalty0 \interfootnotelinepenalty0
\floatingpenalty0 \brokenpenalty0 \everypar{}\leftskip\z@\parskip\z@
\parindent\z@\pretolerance10000 \tolerance10000 \hyphenpenalty10000
\exhyphenpenalty10000 \binoppenalty10000 \relpenalty10000 \adjdemerits0
\doublehyphendemerits0 \finalhyphendemerits0 \Cd@EA\prevdepth\z@}\newbox
\Cd@lF\newbox\Cd@jF\def\Cd@kF{\unhcopy\Cd@lF}\def\Cd@iF{\unhcopy\Cd@jF}\def
\Cd@lI{\hbox{}\penalty1\nointerlineskip}\def\Cd@dH{\penalty5 \noindent\setbox
\Cd@hG=\null\Cd@PF\z@\Cd@yH\z@\ifnum\Cd@BA<\Cd@GB\ht\Cd@hG\ht\Cd@BA\dp\Cd@hG
\dp\Cd@BA\unhbox\Cd@BA\skip0=\lastskip\unskip\else\Cd@PJ\skip0=\z@\fi\endgraf
\ifcase\prevgraf\cd@shouldnt Y \or\cd@shouldnt Z \or\Cd@fH\or\Cd@lH\else
\Cd@eH\fi\unskip\setbox0=\lastbox\unskip\unskip\unpenalty\noindent\unhbox0%
\setbox0\lastbox\unpenalty\unskip\unskip\unpenalty\setbox0\lastbox\Cd@UF
\Cd@BB\lastpenalty\unpenalty\ifnum\Cd@BB>\z@\setbox\z@\lastbox\Cd@gB\repeat
\endgraf\unskip\unskip\unpenalty}\def\Cd@cI{\Cd@pA\Cd@SB\advance\Cd@pA-\Cd@IB
\Cd@qA\Cd@BA\advance\Cd@qA-\Cd@fA\advance\Cd@qA1 \expandafter\message{%
prevgraf=\the\prevgraf at (\the\Cd@pA,\the\Cd@qA)}}\def\Cd@lH{\Cd@nD\setbox
\Cd@xH=\lastbox\setbox\Cd@xH=\hbox{\unhbox\Cd@xH\unskip\unpenalty}\unskip
\ifdim\ht\Cd@xH>\ht\Cd@GC\setbox\Cd@hG=\copy\Cd@xH\else\ifdim\dp\Cd@xH>\dp
\Cd@GC\setbox\Cd@hG=\copy\Cd@xH\else\Cd@gF\Cd@xH\fi\fi\advance\Cd@PF.5\wd
\Cd@xH\advance\Cd@yH.5\wd\Cd@xH\setbox\Cd@xH=\hbox{\unhbox\Cd@xH\Cd@iF}\Cd@wG
\Cd@PF{\box\Cd@xH}\Cd@yH\z@\Cd@rB\Cd@pB}\def\Cd@nD{\ifnum\Cd@A>0 \advance
\dimen0-\Cd@nB\Cd@cA-.5\dimen0 \Cd@A-\Cd@A\else\Cd@A0 \Cd@cA\z@\fi\setbox
\Cd@hG=\lastbox\setbox\Cd@hG=\hbox{\unhbox\Cd@hG\unskip\unskip\unpenalty
\setbox0=\lastbox\global\Cd@MA\lastkern\unkern}\advance\Cd@cA-.5\Cd@MA\unskip
\setbox\Cd@hG=\null\Cd@yH\Cd@cA\Cd@PF-\Cd@cA}\def\Cd@W{\ht\Cd@hG\Cd@EI\dp
\Cd@hG\Cd@DI}\def\Cd@gF#1{\setbox\Cd@hG=\hbox{\Cd@S{\ht\Cd@hG}{\ht#1}\Cd@S{%
\dp\Cd@hG}{\dp#1}\Cd@S{\wd\Cd@hG}{\wd#1}\vrule height\ht\Cd@hG depth\dp\Cd@hG
width\wd\Cd@hG}}\def\Cd@eH{\Cd@nD\Cd@W\setbox\Cd@xH=\lastbox\unskip\setbox
\Cd@OF=\lastbox\unskip\setbox\Cd@OF=\hbox{\unhbox\Cd@OF\unskip\global\Cd@tA
\lastpenalty\unpenalty}\advance\Cd@tA9999 \ifcase\Cd@tA\Cd@iH\or\Cd@mH\or
\Cd@hH\or\Cd@rH\or\Cd@qH\or\Cd@gH\else\cd@shouldnt9\fi}\def\Cd@iH{\Cd@gF
\Cd@xH\Cd@jH\setbox\Cd@HH=\box\Cd@OF\setbox\Cd@IH=\box\Cd@xH}\def\Cd@mH{%
\Cd@gF\Cd@OF\setbox\Cd@xH\hbox{\penalty8 \unhbox\Cd@xH\unskip\unpenalty\ifnum
\lastpenalty=8 \else\Cd@wJ\fi}\Cd@jH\setbox\Cd@OF=\hbox{\unhbox\Cd@OF\unskip
\unpenalty\global\setbox\Cd@@A=\lastbox}\ifdim\wd\Cd@OF=\z@\else\Cd@wJ\fi
\setbox\Cd@HH=\box\Cd@@A}\def\Cd@wJ{\Cd@FB{extra material in \string\pile
\space cell (lost)}}\def\Cd@jH{\Cd@rB\ifvoid\Cd@HH\else\Cd@FB{Clashing
horizontal arrows}\Cd@yH.5\Cd@LF\Cd@PF-\Cd@yH\Cd@pB\Cd@yH\z@\Cd@PF\z@\fi
\Cd@uH\Cd@LF\advance\Cd@uH-\Cd@yH\Cd@LF-\Cd@PF\Cd@BC\Cd@BA}\def\Cd@fH{\setbox
0\lastbox\unskip\Cd@cA\z@\Cd@W\ifdim\skip0>\z@\Cd@wI\Cd@A0 \else\ifnum\Cd@A<1
\Cd@A0 \dimen0\Cd@nB\fi\advance\Cd@A1 \fi}\def\VonH{\Cd@IA46\VonH{.5\Cd@qE}}%
\def\HonV{\Cd@IA57\HonV{.5\Cd@qE}}\def\HmeetV{\Cd@IA44\HmeetV{-\MapShortFall}%
}\def\Cd@IA#1#2#3#4{\Cd@kB34#1{\string#3}\Cd@CD\Cd@BB-999#2 \dimen0=#4\Cd@EI
\dimen0\advance\Cd@EI\axisheight\Cd@DI\dimen0\advance\Cd@DI-\axisheight\Cd@hE
\Cd@@C\Cd@JD}\def\Cd@@C#1{\setbox0=\hbox{\Cd@h#1\Cd@zC}\dimen0.5\wd0 \Cd@EI
\ht0 \Cd@DI\dp0 \Cd@JD}\def\Cd@CD{\setbox0=\null\ht0=\Cd@EI\dp0=\Cd@DI\wd0=%
\dimen0 \copy0\penalty\Cd@BB\box0 }\def\Cd@hH{\Cd@zB\Cd@rB}\def\Cd@rH{\Cd@zB
\Cd@pB}\def\Cd@gH{\Cd@zB\Cd@rB\Cd@pB}\def\Cd@zB{\setbox\Cd@xH=\hbox{\unhbox
\Cd@xH}\setbox\Cd@OF=\hbox{\unhbox\Cd@OF\global\setbox\Cd@@A=\lastbox}\ht
\Cd@hG\ht\Cd@@A\dp\Cd@hG\dp\Cd@@A\advance\Cd@PF\wd\Cd@@A\advance\Cd@yH\wd
\Cd@xH}\Cd@PG\ifPositiveGradient\Cd@QH\Cd@PH\Cd@QH\Cd@PG\ifClimbing\Cd@mB
\Cd@lB\Cd@mB\newcount\DiagonalChoice\DiagonalChoice\m@ne\ifx\tenln\nullfont
\Cd@wI\def\Cd@RF{\Cd@fG\ifPositiveGradient/\else\Cd@h\backslash\Cd@zC\fi}%
\else\def\Cd@RF{\Cd@SF\char\count@}\fi\let\Cd@SF\tenln\def\Use@line@char#1{%
\hbox{#1\Cd@SF\ifPositiveGradient\else\advance\count@64 \fi\char\count@}}\def
\Cd@GF{\Use@line@char{\count@\Cd@KC\multiply\count@8\advance\count@-9\advance
\count@\Cd@gG}}\def\Cd@DF{\Use@line@char{\ifcase\DiagonalChoice\Cd@KF\or
\Cd@JF\or\Cd@JF\else\Cd@KF\fi}}\def\Cd@KF{\ifnum\Cd@KC=\z@\count@\rq33 \else
\count@\Cd@KC\multiply\count@\sixt@@n\advance\count@-9\advance\count@\Cd@gG
\advance\count@\Cd@gG\fi}\def\Cd@JF{\count@\rq\ifcase\Cd@gG55\or\ifcase\Cd@KC
66\or22\or52\or61\or72\fi\or\ifcase\Cd@KC66\or25\or22\or63\or52\fi\or\ifcase
\Cd@KC66\or16\or36\or22\or76\fi\or\ifcase\Cd@KC66\or27\or25\or67\or22\fi\fi
\relax}\def\Cd@iC#1{\hbox{#1\setbox0=\Use@line@char{#1}\ifPositiveGradient
\else\raise.3\ht0\fi\copy0 \kern-.7\wd0 \ifPositiveGradient\raise.3\ht0\fi
\box0}}\def\Cd@MF#1{\hbox{\setbox0=#1\kern-.75\wd0 \vbox to.25\ht0{%
\ifPositiveGradient\else\vss\fi\box0 \ifPositiveGradient\vss\fi}}}\def\Cd@vH#%
1{\hbox{\setbox0=#1\dimen0=\wd0 \vbox to.25\ht0{\ifPositiveGradient\vss\fi
\box0 \ifPositiveGradient\else\vss\fi}\kern-.75\dimen0 }}\Cd@IC{+h:>}{%
\Use@line@char\Cd@JF}\Cd@IC{-h:>}{\Use@line@char\Cd@KF}\Cd@QF{+t:<}{-h:>}%
\Cd@QF{-t:<}{+h:>}\Cd@IC{+t:>}{\Cd@MF{\Use@line@char\Cd@JF}}\Cd@IC{-t:>}{%
\Cd@vH{\Use@line@char\Cd@KF}}\Cd@QF{+h:<}{-t:>}\Cd@QF{-h:<}{+t:>}\Cd@IC{+h:>>%
}{\Cd@iC\Cd@JF}\Cd@IC{-h:>>}{\Cd@iC\Cd@KF}\Cd@QF{+t:<<}{-h:>>}\Cd@QF{-t:<<}{+%
h:>>}\Cd@IC{+t:>>}{\Cd@MF{\Cd@iC\Cd@JF}}\Cd@IC{-t:>>}{\Cd@vH{\Cd@iC\Cd@KF}}%
\Cd@QF{+h:<<}{-t:>>}\Cd@QF{-h:<<}{+t:>>}\Cd@IC{+f:-}{\Cd@kE\null\else\Cd@GF
\fi}\Cd@QF{-f:-}{+f:-}\def\Cd@hC#1#2{\vbox to#1{\vss\hbox to#2{\hss.\hss}\vss
}}\def\hfdot{\Cd@hC{2\axisheight}{.7em}}
\def\vfdot{\Cd@hC{1.46ex}\z@}\def\Cd@FF{\hbox{\dimen0=.3\Cd@nC\dimen1\dimen0
\ifnum\Cd@gG>\Cd@KC\Cd@ZC{\dimen1}\else\Cd@AG{\dimen0}\fi\Cd@hC{\dimen0}{%
\dimen1}}}\newarrowfiller{.}\hfdot\hfdot\vfdot\vfdot\def\dfdot{\Cd@FF\Cd@EJ}%
\Cd@IC{+f:.}{\dfdot}\Cd@IC{-f:.}{\dfdot}\def\Cd@iJ#1{\hbox\bgroup\def\Cd@ZG{#%
1\egroup}\afterassignment\Cd@ZG
\count@=\rq}\def\lnchar{\Cd@iJ\Cd@RF}\def\lad{%
\Cd@iJ\xlad}\def\xlad{\setbox2=\hbox{\Cd@RF}\setbox0=\hbox to.3\wd2{\hss.\hss
}\dimen0=\ht0 \advance\dimen0-\dp0 \dimen1=.3\ht2 \advance\dimen1-\dimen0 \dp
0=.5\dimen1 \dimen1=.3\ht2 \advance\dimen1\dimen0 \ht0=.5\dimen1 \raise\dp0%
\box0}\def\lahh{\Cd@iJ\xlahh}\def\lat{\Cd@iJ\xlat}\def\xlat{\setbox0=\hbox{%
\Cd@RF}\dimen0=\ht0 \setbox1=\hbox to.25\wd0{\ifcase\DiagonalChoice\box0\hss
\or\hss\box0 \or\hss\box0 \or\box0\hss\fi}\vbox to.25\dimen0{\ifClimbing\box1%
\vss\else\vss\box1\fi\kern\z@}}\def\xlahh{\setbox0=\hbox{\Cd@RF}%
\ifPositiveGradient\Cd@wI\copy0 \kern-.7\wd0 \mv.3\ht0\box0 \else\ifClimbing
\Cd@wI\copy0 \kern-.7\wd0 \mv.3\ht0\box0 \else\mv-.3\ht0\copy0 \kern-.7\wd0
\box0 \fi\fi}\def\Cd@HF#1{\setbox#1=\hbox{\dimen5\dimen#1 \setbox8=\box#1
\dimen1\wd8 \count@\dimen5 \divide\count@\dimen1 \ifnum\count@=0 \box8 \ifdim
\dimen5<.95\dimen1 \Cd@bB{diagonal line too short}\fi\else\dimen3=\dimen5
\advance\dimen3-\dimen1 \divide\dimen3\count@\dimen4\dimen3 \Cd@AG{\dimen4}%
\ifPositiveGradient\multiply\dimen4\m@ne\fi\dimen6\dimen1 \advance\dimen6-%
\dimen3 \loop\raise\count@\dimen4\copy8 \ifnum\count@>0 \kern-\dimen6 \advance
\count@\m@ne\repeat\fi}}\def\Cd@dF#1{\Cd@kE\Cd@@J{#1}\else\Cd@HF{#1}\fi}%
\newdimen\objectheight\objectheight1.8ex \newdimen\objectwidth\objectwidth1em
\def\Cd@ID{\dimen6=\Cd@bJ\DiagramCellHeight\dimen7=\Cd@XJ\DiagramCellWidth
\Cd@SI\ifnum\Cd@gG>0 \ifnum\Cd@KC>0 \Cd@EF\else\aftergroup\Cd@MC\fi\else
\aftergroup\Cd@LC\fi}\def\Cd@MC{\Cd@UA{diagonal map is nearly vertical}\Cd@JA
}\def\Cd@LC{\Cd@UA{diagonal map is nearly horizontal}\Cd@JA}\Cd@NG\Cd@JA{Use
an orthogonal map instead}\def\Cd@EF{\Cd@TI\dimen3\dimen7\dimen7\dimen6\Cd@ZC
{\dimen7}\advance\dimen3-\dimen7 \Cd@rE\ifnum\Cd@gG>\Cd@KC\advance\dimen6-%
\dimen1\advance\dimen6-\dimen5 \Cd@ZC{\dimen1}\Cd@ZC{\dimen5}\else\dimen0%
\dimen1\advance\dimen0\dimen5\Cd@AG{\dimen0}\advance\dimen6-\dimen0 \fi\dimen
2.5\dimen7\advance\dimen2-\dimen1 \dimen4.5\dimen7\advance\dimen4-\dimen5
\ifPositiveGradient\dimen0\dimen5 \advance\dimen1-\Cd@XJ\DiagramCellWidth
\advance\dimen1 \Cd@aJ\DiagramCellWidth\setbox6=\llap{\unhbox6\kern.1\ht2}%
\setbox7=\rlap{\kern.1\ht2\unhbox7}\else\dimen0\dimen1 \advance\dimen1-\Cd@aJ
\DiagramCellWidth\setbox7=\llap{\unhbox7\kern.1\ht2}\setbox6=\rlap{\kern.1\ht
2\unhbox6}\fi\setbox6=\vbox{\box6\kern.1\wd2}\setbox7=\vtop{\kern.1\wd2\box7}%
\Cd@AG{\dimen0}\advance\dimen0-\axisheight\advance\dimen0-\Cd@cJ
\DiagramCellHeight\dimen5-\dimen0 \advance\dimen0\dimen6 \advance\dimen1.5%
\dimen3 \ifdim\wd3>\z@\ifdim\ht3>-\dp3\Cd@OB\fi\fi\dimen3\dimen2 \dimen7%
\dimen2\advance\dimen7\dimen4 \ifvoid3 \else\Cd@ZE\else\dimen8\ht3\advance
\dimen8-\axisheight\Cd@ZC{\dimen8}\Cd@U{\dimen8}{.5\wd3}\dimen9\dp3\advance
\dimen9\axisheight\Cd@ZC{\dimen9}\Cd@U{\dimen9}{.5\wd3}\ifPositiveGradient
\advance\dimen2-\dimen9\advance\dimen4-\dimen8 \else\advance\dimen4-\dimen9%
\advance\dimen2-\dimen8 \fi\fi\advance\dimen3-.5\wd3 \fi\dimen9=\Cd@bJ
\DiagramCellHeight\advance\dimen9-2\DiagramCellHeight\Cd@ZE\advance\dimen2%
\dimen4 \Cd@dF{2}\dimen2-\dimen0\advance\dimen2\dp2 \else\Cd@dF{2}\Cd@dF{4}%
\ifPositiveGradient\dimen2-\dimen0\advance\dimen2\dp2 \dimen4\dimen5\advance
\dimen4-\ht4 \else\dimen4-\dimen0\advance\dimen4\dp4 \dimen2\dimen5\advance
\dimen2-\ht2 \fi\fi\setbox0=\hbox to\z@{\kern\dimen1 \ifvoid1 \else
\ifPositiveGradient\advance\dimen0-\dp1 \lower\dimen0 \else\advance\dimen5-%
\ht1 \raise\dimen5 \fi\rlap{\unhbox1}\fi\raise\dimen2\rlap{\unhbox2}\ifvoid3
\else\lower.5\dimen9\rlap{\kern\dimen3\unhbox3}\fi\kern.5\dimen7 \lower.5%
\dimen9\box6 \lower.5\dimen9\box7 \kern.5\dimen7 \Cd@ZE\else\raise\dimen4%
\llap{\unhbox4}\fi\ifvoid5 \else\ifPositiveGradient\advance\dimen5-\ht5 \raise
\dimen5 \else\advance\dimen0-\dp5 \lower\dimen0 \fi\llap{\unhbox5}\fi\hss}\ht
0=\axisheight\dp0=-\ht0\box0 }\def\NorthWest{\Cd@PH\Cd@mB\DiagonalChoice0 }%
\def\NorthEast{\Cd@QH\Cd@mB\DiagonalChoice1 }\def\SouthWest{\Cd@QH\Cd@lB
\DiagonalChoice3 }\def\SouthEast{\Cd@PH\Cd@lB\DiagonalChoice2 }\def\Cd@KD{%
\vadjust{\Cd@pA\Cd@BA\advance\Cd@pA\ifPositiveGradient\else-\fi\Cd@aJ\relax
\Cd@qA\Cd@IB\advance\Cd@qA-\Cd@cJ\relax\hbox{\advance\Cd@pA
\ifPositiveGradient-\fi\Cd@XJ\advance\Cd@qA\Cd@bJ\box\z@\penalty\Cd@pA
\penalty\Cd@qA}\penalty\Cd@pA\penalty\Cd@qA\penalty104}}\def\Cd@jJ#1{\relax
\vadjust{\hbox@maths{#1}\penalty\Cd@BA\penalty\Cd@IB\penalty\tw@}}\def\Cd@gB{%
\ifcase\Cd@BB\or\or\Cd@wG{.5\wd0}{\box0}{.5\wd0}\z@\or\unhbox\z@\setbox\z@
\lastbox\Cd@wG{.5\wd0}{\box0}{.5\wd0}\z@\unpenalty\unpenalty\setbox\z@
\lastbox\or\Cd@sF\else\advance\Cd@BB-100 \ifnum\Cd@BB<\z@\cd@shouldnt B\fi
\setbox\z@\hbox{\kern\Cd@PF\copy\Cd@hG\kern\Cd@yH\Cd@pA\Cd@SB\advance\Cd@pA-%
\Cd@IB\penalty\Cd@pA\Cd@pA\Cd@BA\advance\Cd@pA-\Cd@fA\penalty\Cd@pA\unhbox\z@
\global\Cd@tA\lastpenalty\unpenalty\global\Cd@uA\lastpenalty\unpenalty}\Cd@pA
-\Cd@tA\Cd@qA\Cd@uA\Cd@eJ\fi}\def\Cd@sF{\unhbox\z@\setbox\z@\lastbox\Cd@pA
\lastpenalty\unpenalty\advance\Cd@pA\Cd@gA\Cd@qA\Cd@SB\advance\Cd@qA-%
\lastpenalty\unpenalty\dimen1\lastkern\unkern\setbox3\lastbox\dimen0\lastkern
\unkern\setbox0=\hbox to\z@{\dimen7\Cd@nB\advance\dimen7-\wd\Cd@pA\ifdim
\dimen7<\z@\Cd@QH\multiply\dimen7\m@ne\else\Cd@PH\fi\ifnum\Cd@qA>\Cd@IB\dimen
6\Cd@oB\advance\dimen6-\ht\Cd@qA\else\dimen6\z@\fi\Cd@mI\ifPositiveGradient
\dimen6\z@\else\Cd@gG-\Cd@gG\kern-\dimen7 \fi\global\Cd@EG\raise\dimen6\hbox{%
\Cd@OD{\the\Cd@KC\space\the\Cd@gG\space bturn}\box0 \Cd@KJ{eturn}}\hss}\ht0%
\z@\dp0\z@\Cd@wG{\z@}{\box\z@}{\z@}{\axisheight}}\def\Cd@OD#1{\expandafter
\Cd@KJ{#1}}\Cd@VA\Cd@GJ{output is PostScript dependent}\def\Cd@JC{\Cd@KJ{/%
bturn {gsave currentpoint currentpoint translate 4 2 roll neg exch atan rotate
neg exch neg exch translate } def /eturn {currentpoint grestore moveto} def}}%
\def\Cd@UI#1{\count@#1\relax\multiply\count@7\advance\count@16577\divide
\count@33154 }\def\Cd@PD#1{\expandafter\special{#1}} \def\Cd@@J#1{\setbox#1=%
\hbox{\dimen0\dimen#1\Cd@AG{\dimen0}\Cd@UI{\dimen0}\setbox0=\null
\ifPositiveGradient\count@-\count@\ht0\dimen0 \else\dp0\dimen0 \fi\box0 \Cd@pA
\count@\Cd@UI\Cd@qE\Cd@PD{pn \the\count@}\Cd@PD{pa 0 0}\Cd@UI{\dimen#1}\Cd@PD
{pa \the\count@\space\the\Cd@pA}\Cd@PD{fp}\kern\dimen#1}}\def\Cd@XH{\Cd@SI
\begingroup\ifdim\dimen7<\dimen6 \dimen2=\dimen6 \dimen6=\dimen7 \dimen7=%
\dimen2 \count@\Cd@gG\Cd@gG\Cd@KC\Cd@KC\count@\else\dimen2=\dimen7 \fi\ifdim
\dimen6>.01\p@\Cd@YH\global\Cd@MA\dimen0 \else\global\Cd@MA\dimen7 \fi
\endgroup\dimen2\Cd@MA}\def\Cd@YH{\Cd@kI\ifdim\dimen7>1.73\dimen6 \divide
\dimen2 4 \multiply\Cd@KC2 \else\dimen2=0.353553\dimen2 \advance\Cd@gG-\Cd@KC
\multiply\Cd@KC4 \fi\dimen0=4\dimen2 \Cd@xF4\Cd@xF{-2}\Cd@xF2\Cd@xF{-2.5}}%
\def\Cd@OH{\begingroup\count@\dimen0 \dimen2 45pt \divide\count@\dimen2 \ifdim
\dimen0<\z@\advance\count@\m@ne\fi\ifodd\count@\advance\count@1\Cd@w\else
\Cd@u\fi\advance\dimen0-\count@\dimen2 \Cd@ME\multiply\dimen0\m@ne\fi\ifnum
\count@<0 \multiply\count@-7 \fi\dimen3\dimen1 \dimen6\dimen0 \dimen7 3754936%
sp \ifdim\dimen0<6\p@\def\Cd@pF{4000}\fi\Cd@SI\dimen2\dimen3\Cd@AG{\dimen2}%
\Cd@kI\multiply\Cd@KC-6 \dimen0\dimen2 \Cd@xF1\Cd@xF{0.3}\dimen1\dimen0 \dimen
2\dimen3 \dimen0\dimen3 \Cd@xF3\Cd@xF{1.5}\Cd@xF{0.3}\divide\count@2 \Cd@ME
\multiply\dimen1\m@ne\fi\ifodd\count@\dimen2\dimen1\dimen1\dimen0\dimen0-%
\dimen2 \fi\divide\count@2 \ifodd\count@\multiply\dimen0\m@ne\multiply\dimen1%
\m@ne\fi\global\Cd@MA\dimen0\global\Cd@NA\dimen1\endgroup\dimen6\Cd@MA\dimen7%
\Cd@NA}\def\Cd@mJ{255}\let\Cd@pF\Cd@mJ\def\Cd@SI{\begingroup\ifdim\dimen7<%
\dimen6 \dimen9\dimen7\dimen7\dimen6\dimen6\dimen9\Cd@w\else\Cd@u\fi\dimen2%
\z@\dimen3\Cd@sG\dimen4\Cd@sG\dimen0\z@\dimen8=\Cd@pF\Cd@sG\Cd@bC\global
\Cd@tA\dimen\Cd@ME0\else3\fi\global\Cd@uA\dimen\Cd@ME3\else0\fi\endgroup
\Cd@gG\Cd@tA\Cd@KC\Cd@uA}\def\Cd@bC{\count@\dimen6 \divide\count@\dimen7
\advance\dimen6-\count@\dimen7 \dimen9\dimen4 \advance\dimen9\count@\dimen0
\ifdim\dimen9>\dimen8 \Cd@uB\else\Cd@vB\ifdim\dimen6>\z@\dimen9\dimen6 \dimen
6\dimen7 \dimen7\dimen9 \expandafter\expandafter\expandafter\Cd@bC\fi\fi}\def
\Cd@uB{\ifdim\dimen0=\z@\ifdim\dimen9<2\dimen8 \dimen0\dimen8 \fi\else
\advance\dimen8-\dimen4 \divide\dimen8\dimen0 \ifdim\count@\Cd@sG<2\dimen8
\count@\dimen8 \dimen9\dimen4 \advance\dimen9\count@\dimen0 \Cd@vB\fi\fi}\def
\Cd@vB{\dimen4\dimen0 \dimen0\dimen9 \advance\dimen2\count@\dimen3 \dimen9%
\dimen2 \dimen2\dimen3 \dimen3\dimen9 }\def\Cd@xF#1{\Cd@AG{\dimen2}\advance
\dimen0 #1\dimen2 }\def\Cd@AG#1{\divide#1\Cd@KC\multiply#1\Cd@gG}\def\Cd@ZC#1%
{\divide#1\Cd@gG\multiply#1\Cd@KC}\def\Cd@kI{\dimen6\Cd@gG\Cd@sG\multiply
\dimen6\Cd@gG\dimen7\Cd@KC\Cd@sG\multiply\dimen7\Cd@KC\Cd@SI}\ifx
\errorcontextlines\undefined\Cd@wI\let\Cd@cG\relax\else\def\Cd@cG{%
\errorcontextlines\m@ne}\fi\ifnum\inputlineno<0 \let\Cd@pC\empty\let\Cd@T
\empty\let\Cd@XD\relax\let\Cd@FI\relax\let\Cd@GI\relax\let\Cd@ZF\relax
\message{! Why not upgrade to TeX version 3? (available since 1990)}\else\def
\Cd@T{ at line \number\inputlineno}\def\Cd@XD{ - first occurred}\def\Cd@FI{%
\edef\Cd@e{\the\inputlineno}\global\let\Cd@eB\Cd@e}\def\Cd@e{9999}\def\Cd@GI{%
\xdef\Cd@eB{\the\inputlineno}}\def\Cd@eB{\Cd@e}\def\Cd@ZF{\ifnum\Cd@e<%
\inputlineno\edef\Cd@pC{\space at lines \Cd@e--\the\inputlineno}\else\edef
\Cd@pC{\Cd@T}\fi}\fi\let\Cd@pC\empty\def\Cd@UA#1#2{\Cd@cG\errhelp=#2%
\expandafter\errmessage{\Cd@oA: #1}}\def\Cd@FB#1{{\expandafter\message{!
\Cd@oA: #1\Cd@pC}}}\def\Cd@bB#1{{\expandafter\message{\Cd@oA\space Warning: #%
1\Cd@T}}}\def\Cd@yA#1#2{\Cd@bB{#1 \string#2 is obsolete\Cd@XD}}\def\Cd@wA#1{%
\Cd@yA{Dimension}{#1}\Cd@oD#1\Cd@xA\Cd@xA}\def\Cd@xA{\Cd@KA=}\def\Cd@vA#1{%
\Cd@yA{Count}{#1}\Cd@oD#1\Cd@jG\Cd@jG}\def\Cd@jG{\count@=}\def
\HorizontalMapLength{\Cd@wA\HorizontalMapLength}\def\VerticalMapHeight{\Cd@wA
\VerticalMapHeight}\def\VerticalMapDepth{\Cd@wA\VerticalMapDepth}\def
\VerticalMapExtraHeight{\Cd@wA\VerticalMapExtraHeight}\def
\VerticalMapExtraDepth{\Cd@wA\VerticalMapExtraDepth}\def\DiagonalLineSegments
{\Cd@vA\DiagonalLineSegments}\ifx\tenln\nullfont\Cd@VA\Cd@fG{\Cd@IF\space
diagonal line and arrow font not available}\else\let\Cd@fG\relax\fi\def\Cd@yF
#1#2<#3:#4:#5#6{\begingroup\Cd@LA#3\relax\advance\Cd@LA-#2\relax\ifdim.1em<%
\Cd@LA\Cd@pA#5\relax\Cd@qA#6\relax\ifnum\Cd@pA<\Cd@qA\count@\Cd@qA\advance
\count@-\Cd@pA\Cd@FB{#4 by \the\Cd@LA}\if#1v\let\Cd@ZG\Cd@LJ\else\advance
\count@\count@\if#1l\advance\count@-\Cd@A\else\if#1r\advance\count@\Cd@A\fi
\fi\advance\Cd@LA\Cd@LA\let\Cd@ZG\Cd@IE\fi\divide\Cd@LA\count@\ifdim\Cd@ZG<%
\Cd@LA\global\Cd@ZG\Cd@LA\fi\fi\fi\endgroup}\Cd@PG\Cd@dE\Cd@wC\Cd@vC\Cd@NG
\Cd@II{See the message above.}\Cd@NG\Cd@DH{Perhaps you've forgotten to end the
diagram before resuming the text, in\Cd@QG which case some garbage may be
added to the diagram, but we should be ok now.\Cd@QG Alternatively you've left
a blank line in the middle - TeX will now complain\Cd@QG that the remaining
\Cd@P s are misplaced - so please use comments for layout.}\Cd@NG\Cd@RD{You
have already closed too many brace pairs or environments; an \Cd@uC\Cd@QG
command was (over)due.}\Cd@NG\Cd@@H{\Cd@UC\space and \Cd@uC\space commands
must match.}\def\Cd@AH{\ifnum\inputlineno=0 \else\expandafter\Cd@BH\fi}\def
\Cd@BH{\Cd@xC\Cd@tC\crcr\Cd@UA{missing \Cd@uC\space inserted before \Cd@CH-
type "h"}\Cd@DH\enddiagram\Cd@bF\Cd@CH\par}\def\Cd@bF#1{\edef\enddiagram{%
\noexpand\Cd@cD{#1\Cd@T}}}\def\Cd@cD#1{\Cd@UA{\Cd@uC\space(anticipated by #1)
ignored}\Cd@II\let\enddiagram\Cd@rF}\def\Cd@rF{\Cd@UA{misplaced \Cd@uC\space
ignored}\Cd@@H}\def\Cd@cC{\Cd@UA{missing \Cd@uC\space inserted.}\Cd@RD\Cd@bF{%
closing group}}\ifx\DeclareOption\undefined\else\ifx\DeclareOption\@notprerr
\else\DeclareOption*{\let\Cd@K\relax\let\Cd@tJ\relax\expandafter\Cd@kD
\CurrentOption,}\fi\fi

\catcode\lq\$=3 
\def\vboxtoz{\vbox to\z@}

\def\scriptaxis#1{\@scriptaxis{$\scriptstyle#1$}}
\def\ssaxis#1{\@ssaxis{$\scriptscriptstyle#1$}}
\def\@scriptaxis#1{\dimen0\axisheight\advance\dimen0-\ss@axisheight\raise
\dimen0\hbox{#1}}\def\@ssaxis#1{\dimen0\axisheight\advance\dimen0-%
\ss@axisheight\raise\dimen0\hbox{#1}}

\ifx\boldmath\undefined
\let\boldscriptaxis\scriptaxis
\def\boldscript#1{\hbox{$\scriptstyle#1$}}
\else\def\boldscriptaxis#1{\@scriptaxis{\boldmath$\scriptstyle#1$}}
\def\boldscript#1{\hbox{\boldmath$\scriptstyle#1$}}
\fi

\def\raisehook#1#2#3{\hbox{\setbox3=\hbox{#1$\scriptscriptstyle#3$}%
\dimen0\ss@axisheight
\dimen1\axisheight\advance\dimen1-\dimen0
\dimen2\ht3\advance\dimen2-\dimen0%
\advance\dimen1 #2\dimen2
\raise\dimen1\box3}}

\def\shifthook#1#2#3{\setbox0=\hbox{#1$\scriptscriptstyle#3$}\dimen0\wd0%
\divide\dimen0 12\Cd@MH{\dimen0}
\dimen1\wd0\advance\dimen1-2\dimen0\advance\dimen1-\Cd@@I\Cd@MH{\dimen1}\kern
#2\dimen1\box0}

\def\@cmex{\mathchar"03}



\def\make@pbk#1{\setbox\tw@\hbox to\z@{#1}\ht\tw@\z@\dp\tw@\z@\box\tw@}\def
\Cd@kJ#1{\overprint{\hbox to\z@{#1}}}\def\Cd@vJ{\kern0.1em}\def\Cd@uJ{\kern0.%
25em}

\def\SEpbk{\make@pbk{\Cd@vJ\vrule depth 2.67ex height -2.55ex width 0.9em
\vrule height -0.46ex depth 2.67ex width .05em \hss}}

\def\SWpbk{\make@pbk{\hss\vrule height -0.46ex depth 2.67ex width .05em \vrule
depth 2.67ex height -2.55ex width .9em \Cd@vJ}}

\def\NEpbk{\make@pbk{\Cd@vJ\vrule depth -3.48ex height 3.67ex width 0.95em
\vrule height 3.67ex depth -1.39ex width .05em \hss}}

\def\NWpbk{\make@pbk{\hss\vrule height 3.6ex depth -1.39ex width .05em \vrule
depth -3.48ex height 3.67ex width .95em \Cd@vJ}}

\def\puncture{{\setbox0\hbox{A}\vrule height.53\ht0 depth-.47\ht0 width.35\ht
0 \kern.12\ht0 \vrule height\ht0 depth-.65\ht0 width.06\ht0 \kern-.06\ht0
\vrule height.35\ht0 depth0pt width.06\ht0 \kern.12\ht0 \vrule height.53\ht0
depth-.47\ht0 width.35\ht0 }}




\def\rhvee{\mkern-10mu\greaterthan}
\def\lhvee{\lessthan\mkern-10mu}
\def\dhvee{\vboxtoz{\vss\hbox{$\vee$}\kern0pt}}
\def\uhvee{\vboxtoz{\hbox{$\wedge$}\vss}}
\newarrowhead{vee}\rhvee\lhvee\dhvee\uhvee

\def\dhlvee{\vboxtoz{\vss\hbox{$\scriptstyle\vee$}\kern0pt}}
\def\uhlvee{\vboxtoz{\hbox{$\scriptstyle\wedge$}\vss}}
\newarrowhead{littlevee}{\mkern1mu\scriptaxis\rhvee}{\scriptaxis\lhvee}%
\dhlvee\uhlvee\ifx\boldmath\undefined
\newarrowhead{boldlittlevee}{\mkern1mu\scriptaxis\rhvee}{\scriptaxis\lhvee}%
\dhlvee\uhlvee\else
\def\dhblvee{\vboxtoz{\vss\boldscript\vee\kern0pt}}
\def\uhblvee{\vboxtoz{\boldscript\wedge\vss}}
\newarrowhead{boldlittlevee}{\mkern1mu\boldscriptaxis\rhvee}{\boldscriptaxis
\lhvee}\dhblvee\uhblvee
\fi

\def\rhcvee{\mkern-10mu\succ}
\def\lhcvee{\prec\mkern-10mu}
\def\dhcvee{\vboxtoz{\vss\hbox{$\curlyvee$}\kern0pt}}
\def\uhcvee{\vboxtoz{\hbox{$\curlywedge$}\vss}}
\newarrowhead{curlyvee}\rhcvee\lhcvee\dhcvee\uhcvee

\def\rhvvee{\mkern-13mu\gg}
\def\lhvvee{\ll\mkern-13mu}
\def\dhvvee{\vboxtoz{\vss\hbox{$\vee$}\kern-.6ex\hbox{$\vee$}\kern0pt}}
\def\uhvvee{\vboxtoz{\hbox{$\wedge$}\kern-.6ex \hbox{$\wedge$}\vss}}
\newarrowhead{doublevee}\rhvvee\lhvvee\dhvvee\uhvvee

\def\rhtriangle{\triangleright\mkern1.2mu}
\def\lhtriangle{\triangleleft\mkern1mu}
\def\uhtriangle{\vbox{\kern-.2ex \hbox{$\scriptscriptstyle\bigtriangleup$}%
\kern-.25ex}}
\def\dhtriangle{\vbox{\kern-.4ex \hbox{$\scriptscriptstyle\bigtriangledown$}%
\kern-.1ex}}
\def\dhblack{\vbox{\kern-.25ex\nointerlineskip\hbox{$\blacktriangledown$}}}%
\def\uhblack{\vbox{\kern-.25ex\nointerlineskip\hbox{$\blacktriangle$}}}%
\def\dhlblack{\vbox{\kern-.25ex\nointerlineskip\hbox{$\scriptstyle
\blacktriangledown$}}}
\def\uhlblack{\vbox{\kern-.25ex\nointerlineskip\hbox{$\scriptstyle
\blacktriangle$}}}
\newarrowhead{triangle}\rhtriangle\lhtriangle\dhtriangle\uhtriangle
\newarrowhead{blacktriangle}{\mkern-1mu\blacktriangleright\mkern.4mu}{%
\blacktriangleleft}\dhblack\uhblack\newarrowhead{littleblack}{\mkern-1mu%
\scriptaxis\blacktriangleright}{\scriptaxis\blacktriangleleft\mkern-2mu}%
\dhlblack\uhlblack

\def\rhla{\hbox{\setbox0=\lnchar55\dimen0=\wd0\kern-.6\dimen0\ht0\z@\raise
\axisheight\box0\kern.1\dimen0}}
\def\lhla{\hbox{\setbox0=\lnchar33\dimen0=\wd0\kern.05\dimen0\ht0\z@\raise
\axisheight\box0\kern-.5\dimen0}}
\def\dhla{\vboxtoz{\vss\rlap{\lnchar77}}}
\def\uhla{\vboxtoz{\setbox0=\lnchar66 \wd0\z@\kern-.15\ht0\box0\vss}}
\newarrowhead{LaTeX}\rhla\lhla\dhla\uhla

\def\lhlala{\lhla\kern.3em\lhla}
\def\rhlala{\rhla\kern.3em\rhla}
\def\uhlala{\hbox{\uhla\raise-.6ex\uhla}}
\def\dhlala{\hbox{\dhla\lower-.6ex\dhla}}
\newarrowhead{doubleLaTeX}\rhlala\lhlala\dhlala\uhlala

\def\hhO{\scriptaxis\bigcirc\mkern.4mu} \def\hho{{\circ}\mkern1.2mu}%
\newarrowhead{o}\hho\hho\circ\circ
\newarrowhead{O}\hhO\hhO{\scriptstyle\bigcirc}{\scriptstyle\bigcirc}

\def\rhtimes{\mkern-5mu{\times}\mkern-.8mu}\def\lhtimes{\mkern-.8mu{\times}%
\mkern-5mu}\def\uhtimes{\setbox0=\hbox{$\times$}\ht0\axisheight\dp0-\ht0%
\lower\ht0\box0 }\def\dhtimes{\setbox0=\hbox{$\times$}\ht0\axisheight\box0 }%
\newarrowhead{X}\rhtimes\lhtimes\dhtimes\uhtimes\newarrowhead+++++


\newarrowhead{Y}{\mkern-3mu\Yright}{\Yleft\mkern-3mu}\Yup\Ydown


\newarrowhead{->}\rightarrow\leftarrow\downarrow\uparrow

\newarrowhead{=>}\Rightarrow\Leftarrow{\@cmex7F}{\@cmex7E}

\newarrowhead{harpoon}\rightharpoonup\leftharpoonup\downharpoonleft
\upharpoonleft

\def\twoheaddownarrow{\rlap{$\downarrow$}\raise-.5ex\hbox{$\downarrow$}}
\def\twoheaduparrow{\rlap{$\uparrow$}\raise.5ex\hbox{$\uparrow$}}
\newarrowhead{->>}\twoheadrightarrow\twoheadleftarrow\twoheaddownarrow
\twoheaduparrow


\def\ltvee{\mkern-1mu{\lessthan}\mkern.4mu}
\newarrowtail{vee}\greaterthan\ltvee\vee\wedge

\newarrowtail{littlevee}{\scriptaxis\greaterthan}{\mkern-1mu\scriptaxis
\lessthan}{\scriptstyle\vee}{\scriptstyle\wedge}\ifx\boldmath\undefined
\newarrowtail{boldlittlevee}{\scriptaxis\greaterthan}{\mkern-1mu\scriptaxis
\lessthan}{\scriptstyle\vee}{\scriptstyle\wedge}\else\newarrowtail{%
boldlittlevee}{\boldscriptaxis\greaterthan}{\mkern-1mu\boldscriptaxis
\lessthan}{\boldscript\vee}{\boldscript\wedge}\fi

\newarrowtail{curlyvee}\succ{\mkern-1mu{\prec}\mkern.4mu}\curlyvee\curlywedge

\def\rttriangle{\mkern1.2mu\triangleright}
\newarrowtail{triangle}\rttriangle\lhtriangle\dhtriangle\uhtriangle
\newarrowtail{blacktriangle}\blacktriangleright{\mkern-1mu\blacktriangleleft
\mkern.4mu}\dhblack\uhblack\newarrowtail{littleblack}{\scriptaxis
\blacktriangleright\mkern-2mu}{\mkern-1mu\scriptaxis\blacktriangleleft}%
\dhlblack\uhlblack

\def\rtla{\hbox{\setbox0=\lnchar55\dimen0=\wd0\kern-.5\dimen0\ht0\z@\raise
\axisheight\box0\kern-.2\dimen0}}
\def\ltla{\hbox{\setbox0=\lnchar33\dimen0=\wd0\kern-.15\dimen0\ht0\z@\raise
\axisheight\box0\kern-.5\dimen0}}
\def\dtla{\vbox{\setbox0=\rlap{\lnchar77}\dimen0=\ht0\kern-.7\dimen0\box0%
\kern-.1\dimen0}}
\def\utla{\vbox{\setbox0=\rlap{\lnchar66}\dimen0=\ht0\kern-.1\dimen0\box0%
\kern-.6\dimen0}}
\newarrowtail{LaTeX}\rtla\ltla\dtla\utla

\def\rtvvee{\gg\mkern-3mu}
\def\ltvvee{\mkern-3mu\ll}
\def\dtvvee{\vbox{\hbox{$\vee$}\kern-.6ex \hbox{$\vee$}\vss}}
\def\utvvee{\vbox{\vss\hbox{$\wedge$}\kern-.6ex \hbox{$\wedge$}\kern\z@}}
\newarrowtail{doublevee}\rtvvee\ltvvee\dtvvee\utvvee

\def\ltlala{\ltla\kern.3em\ltla}
\def\rtlala{\rtla\kern.3em\rtla}
\def\utlala{\hbox{\utla\raise-.6ex\utla}}
\def\dtlala{\hbox{\dtla\lower-.6ex\dtla}}
\newarrowtail{doubleLaTeX}\rtlala\ltlala\dtlala\utlala

\def\utbar{\vrule height 0.093ex depth0pt width 0.4em}
\let\dtbar\utbar
\def\rtbar{\mkern1.5mu\vrule height 1.1ex depth.06ex width .04em\mkern1.5mu}%
\let\ltbar\rtbar
\newarrowtail{mapsto}\rtbar\ltbar\dtbar\utbar
\newarrowtail{|}\rtbar\ltbar\dtbar\utbar


\def\rthooka{\raisehook{}+\subset\mkern-1mu}
\def\lthooka{\mkern-1mu\raisehook{}+\supset}
\def\rthookb{\raisehook{}-\subset\mkern-2mu}
\def\lthookb{\mkern-1mu\raisehook{}-\supset}

\def\dthooka{\shifthook{}+\cap}
\def\dthookb{\shifthook{}-\cap}
\def\uthooka{\shifthook{}+\cup}
\def\uthookb{\shifthook{}-\cup}

\newarrowtail{hooka}\rthooka\lthooka\dthooka\uthooka\newarrowtail{hookb}%
\rthookb\lthookb\dthookb\uthookb

\ifx\boldmath\undefined\newarrowtail{boldhooka}\rthooka\lthooka\dthooka
\uthooka\newarrowtail{boldhookb}\rthookb\lthookb\dthookb\uthookb\newarrowtail
{boldhook}\rthooka\lthookb\dthooka\uthooka\else\def\rtbhooka{\raisehook
\boldmath+\subset\mkern-1mu}
\def\ltbhooka{\mkern-1mu\raisehook\boldmath+\supset}
\def\rtbhookb{\raisehook\boldmath-\subset\mkern-2mu}
\def\ltbhookb{\mkern-1mu\raisehook\boldmath-\supset}
\def\dtbhooka{\shifthook\boldmath+\cap}
\def\dtbhookb{\shifthook\boldmath-\cap}
\def\utbhooka{\shifthook\boldmath+\cup}
\def\utbhookb{\shifthook\boldmath-\cup}
\newarrowtail{boldhooka}\rtbhooka\ltbhooka\dtbhooka\utbhooka\newarrowtail{%
boldhookb}\rtbhookb\ltbhookb\dtbhookb\utbhookb\newarrowtail{boldhook}%
\rtbhooka\ltbhookb\dthbooka\utbhooka\fi

\newarrowtail{hook}\rthooka\lthookb\dthooka\uthooka\newarrowtail{C}\rthooka
\lthookb\dthooka\uthooka

\newarrowtail{o}\hho\hho\circ\circ
\newarrowtail{O}\hhO\hhO{\scriptstyle\bigcirc}{\scriptstyle\bigcirc}

\newarrowtail{X}\lhtimes\rhtimes\uhtimes\dhtimes\newarrowtail+++++


\newarrowtail{Y}\Yright\Yleft\Yup\Ydown




\newarrowfiller{=}=={\@cmex77}{\@cmex77}
\def\vfthree{\mid\!\!\!\mid\!\!\!\mid}
\newarrowfiller{3}\equiv\equiv\vfthree\vfthree

\def\vfdashstrut{\vrule width0pt height1.3ex depth0.7ex}
\def\vfthedash{\vrule width\Cd@qE height0.6ex depth 0pt}
\def\hfthedash{\Cd@MI\vrule\horizhtdp width 0.26em}
\def\hfdash{\mkern5.5mu\hfthedash\mkern5.5mu}
\def\vfdash{\vfdashstrut\vfthedash}
\newarrowfiller{dash}\hfdash\hfdash\vfdash\vfdash


\newarrowmiddle+++++




\iffalse
\newarrow{To}----{vee}
\newarrow{Arr}----{LaTeX}
\newarrow{Dotsto}....{vee}
\newarrow{Dotsarr}....{LaTeX}
\newarrow{Dashto}{}{dash}{}{dash}{vee}
\newarrow{Dasharr}{}{dash}{}{dash}{LaTeX}
\newarrow{Mapsto}{mapsto}---{vee}
\newarrow{Mapsarr}{mapsto}---{LaTeX}
\newarrow{IntoA}{hooka}---{vee}
\newarrow{IntoB}{hookb}---{vee}
\newarrow{Embed}{vee}---{vee}
\newarrow{Emarr}{LaTeX}---{LaTeX}
\newarrow{Onto}----{doublevee}
\newarrow{Dotsonarr}....{doubleLaTeX}
\newarrow{Dotsonto}....{doublevee}
\newarrow{Dotsonarr}....{doubleLaTeX}
\else
\newarrow{To}---->
\newarrow{Arr}---->
\newarrow{Dotsto}....>
\newarrow{Dotsarr}....>
\newarrow{Dashto}{}{dash}{}{dash}>
\newarrow{Dasharr}{}{dash}{}{dash}>
\newarrow{Mapsto}{mapsto}--->
\newarrow{Mapsarr}{mapsto}--->
\newarrow{IntoA}{hooka}--->
\newarrow{IntoB}{hookb}--->
\newarrow{Embed}>--->
\newarrow{Emarr}>--->
\newarrow{Onto}----{>>}
\newarrow{Dotsonarr}....{>>}
\newarrow{Dotsonto}....{>>}
\newarrow{Dotsonarr}....{>>}
\fi

\newarrow{Implies}===={=>}
\newarrow{Project}----{triangle}
\newarrow{Pto}----{harpoon}

\newarrow{Eq}=====
\newarrow{Line}-----
\newarrow{Dots}.....
\newarrow{Dashes}{}{dash}{}{dash}{}

\newarrowhead{cmexbra}{\@cmex7B}{\@cmex7C}{\@cmex3B}{\@cmex38}
\newarrowtail{cmexbra}{\@cmex7A}{\@cmex7D}{\@cmex39}{\@cmex3A}
\newarrowmiddle{cmexbra}{\braceru\bracelu}{\bracerd\braceld}{\@cmex3D}{\@cmex
3C}
\newarrow{@brace}{cmexbra}-{cmexbra}-{cmexbra}
\newarrow{@parenth}{cmexbra}---{cmexbra}
\def\rightBrace{\d@brace[cmex,thick,midvshaft]}
\def\leftBrace{\u@brace[cmex,thick,midvshaft]}
\def\upperBrace{\r@brace[cmex,thick,midhshaft]}
\def\lowerBrace{\l@brace[cmex,thick,midhshaft]}
\def\rightParenth{\d@parenth[cmex,thick]}
\def\leftParenth{\u@parenth[cmex,thick]}
\def\upperParenth{\r@parenth[cmex,thick]}
\def\lowerParenth{\l@parenth[cmex,thick]}



\let\lInto\lIntoA
\let\hEq\rEq
\let\vEq\uEq









\def\labelstyle{
\ifincommdiag
\textstyle
\else
\scriptstyle
\fi}
\let\objectstyle\displaystyle

\newdiagramgrid{pentagon}{0.618034,0.618034,1,1,1,1,0.618034,0.618034}{1.%
17557,1.17557,1.902113,1.902113}

\newdiagramgrid{perspective}{0.75,0.75,1.1,1.1,0.9,0.9,0.95,0.95,0.75,0.75}{0%
.75,0.75,1.1,1.1,0.9,0.9}

\diagramstyle[
dpi=300,
vmiddle,nobalance,
loose,
thin,
pilespacing=10pt,%
shortfall=4pt,
]

\ifx\ProcessOptions\undefined\else\Cd@QJ\ProcessOptions\relax\Cd@lE\Cd@b\fi
\fi

\cdrestoreat

\dimen0 200pt \dimen1 210pt \dimen2 220pt \dimen3 230pt \dimen4 240pt \dimen5
250pt \dimen6 260pt \dimen7 270pt \dimen8 280pt \dimen9 290pt
%
\NoBlackBoxes
%
\def\openC{\Bbb C}

\def\openP{\Bbb P}
\def\PP{\Bbb P}
\def\Pn#1{\openP^{#1}}
\def\openQ{\Bbb Q}
\def\openZ{\Bbb Z}
\def\ZZ{\Bbb Z}
\def\Hom{\operatorname{Hom}}

\def\rk{\operatorname{rk}}
\def\Im{\operatorname{Im}}
\def\dim{\operatorname{dim}}
\def\ker{\operatorname{ker}}
\def\id{\operatorname{id}}
\def\Hom{\operatorname{Hom}}
\def\Ext{\operatorname{Ext}}
\def\deg{\operatorname{deg}}
\def\Sing{\operatorname{Sing}}
\def\Trisec{\operatorname{Trisec}}
\def\Sec{\operatorname{Sec}}
\def\Grass{\operatorname{Grass}}
\def\NS{\operatorname{NS}}
\def\Num{\operatorname{Num}}
\def\Pic{\operatorname{Pic}}
\def\O#1{{\scr O}(#1)}
\def\rto{\raise.5ex\hbox{$\scriptscriptstyle ---\!\!\!>$}}
%
%
\topmatter
\title
Syzygies of abelian and bielliptic surfaces in ${\openP}^4$
\endtitle
\rightheadtext{A. Aure, W. Decker, K. Hulek, S. Popescu, K. Ranestad }
\leftheadtext{Syzygies of abelian and bielliptic surfaces in ${\openP}^4$}

\author
\hbox{Alf Aure\dag, Wolfram Decker, Klaus Hulek, Sorin Popescu,
Kristian Ranestad}
\endauthor
%
%
\address
Wolfram Decker\hfil\break
Fachbereich Mathematik,
Universit\"at des Saarlandes,
D 66041 Saarbr\"ucken,
Germany
\endaddress
\email
decker\@math.uni-sb.de
\endemail
\address
Klaus Hulek\hfil\break
Institut f\"ur Mathematik,
Universit\"at Hannover,
D 30167 Hannover,
Germany
\endaddress
\email
hulek\@math.uni-hannover.de
\endemail
\address
Sorin Popescu\hfil\break
Mathematics Department,
Brandeis University,
Waltham, MA 02254,
USA
\endaddress
\email
popescu\@oscar.math.brandeis.edu
\endemail
\address
Kristian Ranestad\hfil\break
Matematisk Institutt,
Universitetet i Oslo,
Box 1053 Blindern,
N-0316 Oslo 3,
Norway
\endaddress
\email
ranestad\@math.uio.no
\endemail
%
%
%
\endtopmatter
\toc
\widestnumber\head{9.}
\head 0. Introduction     \endhead
\head 1. The Heisenberg group of level 5    \endhead
\head 2. Surfaces and syzygies    \endhead
\head 3. The minimal abelian and bielliptic surfaces    \endhead
\head 4. Moore matrices     \endhead
\head 5. Modules obtained by concatenating three Moore matrices \endhead
\head 6. Syzygy construction of bielliptic surfaces     \endhead
\head 7. Non-minimal abelian surfaces obtained via linkage \endhead
\head 8. A new family of non-minimal abelian surfaces \endhead
\head 9. Degenerations of abelian and bielliptic surfaces \endhead
\widestnumber\head{Appendix}
\head Appendix \endhead
\endtoc
\endtopmatter
\document
%
%
\head{0. Introduction}
\endhead
So far six families of smooth irregular surfaces are known to exist in
$\openP^4$ (up to pullbacks of these families by suitable finite covers
$\openP^4\longrightarrow \openP^4$). These are the elliptic quintic scrolls,
the minimal abelian and bielliptic surfaces (of degree 10), two different
families of non-minimal abelian surfaces of degree 15 and one family of
non-minimal bielliptic surfaces of degree 15. The elliptic quintic scroll is
classical. The degree 10 abelian surfaces arise as sections of the
Horrocks-Mumford bundle \cite{HM}, and the degree 10 bielliptic surfaces were
found by Serrano \cite{Ser} using Reider's method \cite{Rei}. It was noticed by
Ellingsrud and Peskine (unpublished, compare also \cite{Au}) that the
Horrocks-Mumford surfaces can be
(5,5)-linked to non-minimal abelian surfaces of degree 15. The other type of
non-minimal abelian surfaces of degree 15 was found by the fourth author
in his thesis \cite{Po} using the syzygy approach of \cite{DES}.
Whereas the abelian surfaces obtained by liaison lie on three quintics,
the latter abelian surfaces are contained in one quintic hypersurface only.
Finally, the existence of the non-minimal bielliptic surfaces of degree 15
was established in \cite{ADHPR} while studying the geometry of the embedding
of the minimal bielliptic surfaces. The unifying element which governs the
properties of all these surfaces is their Heisenberg symmetry. In fact, it
is only this symmetry which makes it possible to construct the irregular
surfaces of degree 15 (where Reider's method does not apply). In all cases
every interesting property of these surfaces we know is related to or depends
on the presence of the Heisenberg symmetry. We do not know of any irregular
surfaces in $\openP^4$ which do not possess a Heisenberg symmetry (again up
to pullbacks under finite covers $\openP^4\longrightarrow \openP^4$), but we
also do not know of an a priori reason why such surfaces should not exist.
\medskip
The main purpose of this paper is to describe the structure of the
Hartshorne-Rao modules and the syzygies of the known smooth irregular surfaces
in $\openP^4$. At the same time we explain how to construct these surfaces
via syzygies, thus giving a unified construction method for the smooth
irregular surfaces known so far. In fact, this means that one can now
construct all known smooth non-general type surfaces in $\openP^4$ via the
syzygy approach (compare \cite{DES}, \cite{Po}, \cite{Sch}).
In particular, one can compute explicit equations for all these surfaces
with the help of Macaulay \cite{Mac}. We refer to \cite{DES} for remarks
concerning these computations in general and to
http://yosemite.math.uni-sb.de/software.html
for the corresponding programs. In fact, the main results of this paper
(6.3 and 6.6) were first predicted and checked by using  Macaulay.
The degree of the smooth non-general type surfaces in
$\openP^4$ is bounded \cite{EP}, and it is conjectured that the precise
bound is 15. The degree 15 surfaces considered in this paper are the only
smooth non-general type surfaces known of degree 15.
\medskip
The abelian and bielliptic surfaces have in common that either their first
Hartshorne-Rao module (in the degree 10 cases) or the dual of their second
Hartshorne-Rao module (in the degree 15 cases) has up to twists a minimal
free presentation of type
$$
0\leftarrow M\leftarrow 5R\overset{\alpha_1}\to\longleftarrow 15 R(-1)
\leftarrow\ldots,
$$
where $R$ is the homogeneous coordinate ring of $\openP^4$.
In fact, by Heisenberg invariance, $\alpha_1$ is of type
$$
\alpha_1=(M_{p_1}(y), M_{p_2}(y),M_{p_3}(y)),
$$
where each building block $M_{p_i}$ is what we call
a {\it{Moore matrix}} (compare Section
4 for a precise definition). The matrices $M_{p_i}$ are $5\times
5$-matrices with linear
entries and depend linearly on points $p_i\in\openP^4$. Thus the
module $M$ corresponds to a plane $\Pi=<p_1, p_2, p_3>\subset\openP^4$. The
plane
$\Pi$ is called the {\it{representing plane}} of $M$. In order
to obtain surfaces in this way the representing plane has to be chosen very
carefully, in particular the general plane does not give rise to
surfaces.\medskip
The abelian surfaces are in addition invariant under a certain involution
$\iota$ of $\openP^4$. There are basically only 3 possibilities to construct
a module as above which is Heisenberg- and $\iota$-invariant. This can
be stated precisely using the eigenspaces $\openP^2_+$ and $\openP^1_-$ of
$\iota$ (compare (5.1)). By \cite{De2} the plane $\openP^2_+$ gives the
first cohomology module of the Horrocks-Mumford bundle and leads thus to
the abelian surfaces of degree 10 (and hence also to the linked abelian
surfaces of degree 15). The second of these possibilities, in a special case,
leads to the union of two elliptic quintic scrolls (see 6.7).
Finally, the third possibility gives rise to the new abelian surfaces of degree
15 found in \cite{Po} (compare Section 8). The bielliptic surfaces
are only Heisenberg- but not $\iota$-invariant. On the other hand, we know
from \cite{ADHPR} that every minimal and non-minimal bielliptic surface
lies on a unique quintic hypersurface, namely either the trisecant scroll
of a quintic elliptic scroll or the secant variety of an elliptic
normal curve in $\openP^4$. In fact, every such quintic hypersurface
contains 8 bielliptic surfaces paired under $\iota$ and a pencil of abelian
surfaces of degree either 10 or 15 arising from the Horrocks-Mumford bundle.
Since the representing plane $\openP^2_+$ of the module of these abelian
surfaces can be seen as the plane spanned by the non-trivial 2-torsion points
of any given elliptic quintic normal curve, and since there are 8 non-trivial
3-torsion points on such a curve, a natural idea for the bielliptic surfaces is
to study planes which depend on the 3-torsion points of an elliptic normal
curve (compare 6.3 and 6.6).
These torsion points are given as the intersection of the curve
with a certain hypersurface (see 4.13). For explicit calculations
on a computer it is convenient to work over a finite prime
field $\openZ_p$ (compare again \cite{DES}). For this purpose one has to
choose the characteristic $p$ such that at least some of the torsion points are
rational over $\openZ_p.$ \medskip
In the case of abelian surfaces the module $M$ does not determine the surface.
It either gives rise to the Horrocks-Mumford bundle (in the degree 10 case) or
to certain reflexive sheaves (for the non-minimal abelian surfaces)
whose sections give the corresponding surfaces. On the
other hand both the minimal and the non-minimal bielliptic surfaces are
determined by their module $M$ (see again 6.3 and 6.6).\medskip
The building blocks of the syzygies of $M$ are again matrices associated to
homomorphisms of Heisenberg modules. The building blocks $L_{q}(y)$ of the
first order linear syzygies of $\alpha_1$ for example are $5\times 5$-matrices
of a type which we call {\it{syzygy matrices}}. They also depend linearly on
parameter points $q\in\openP^4$.
\medskip
The paper is organized as follows: \smallskip\noindent
In Section 1 we recall the representation theory which we will need later on.
Section 2 contains a quick review of the syzygy
approach of \cite{DES}, which is in fact an application of Beilinson's
theorem \cite{Bei}. In Section 3 we study the Horrocks-Mumford surfaces
from a Beilinson point of view, and we determine the shape of the syzygies
of the minimal bielliptic surfaces. We are, however, at this point not yet
able to determine the differentials of the minimal free resolution of the
module $M$ leading to the minimal bielliptic surfaces.
\smallskip\noindent
Sections 4, 5 and 6 contain the main new ideas and arguments of this paper.
In the sequel we will give a short outline of these sections (see the
beginning of Section 4 for more details).
\smallskip\noindent
Moore matrices and syzygy matrices make their first appearance in Section 4,
where they arise naturally as matrices representing certain Heisenberg
invariant homomorphisms coming from the Koszul complex on $\openP^4$. It is
also in this section that we first establish various
relations between Moore matrices and geometric objects. Indeed, it is well
known and goes back to Bianchi \cite{Bi}, Klein \cite{Kl}, and others, that
Moore matrices are the crucial ingredients in the minimal free resolution of
elliptic normal curves in $\openP^4$. Similarly, syzygy matrices, i.e.,
matrices of type $L_{q}(y)$, appear in the
minimal free resolution of elliptic quintic scrolls (see Lemma 4.4).
For certain parameter points the determinants of our matrices define the
secant variety of an elliptic normal curve and the trisecant variety of a
quintic elliptic scroll resp. (see 4.12). In order to describe the first
order linear syzygies of our modules we need to study products of Moore
matrices and syzygy matrices. It turns out that these products define
an incidence scroll when we restrict the parameter point $p$ of the Moore
matrix to an elliptic quintic normal curve. This scroll is closely related
to the degree 15
scroll which is the union of the vertices of the rank 3 quadrics through an
elliptic quintic normal curve. It is the geometry of these scrolls which
will allow us later on to get a hold on the syzygies of the modules $M$.
We will also determine the 3-torsion points and the 6-torsion points of
an elliptic quintic normal curve. For this we need to recall
some facts about Shioda's modular surface $S(5)$ and certain maps
related to the Horrocks-Mumford bundle.
\smallskip\noindent
In Section 5 we study modules obtained by concatenating three Moore matrices,
in
particular in the cases described in (5.1).
We investigate when these modules are artinian (5.5 and 5.6) and we study
their linear and quadratic syzygies. Using the results of Section 4 we can show
that these syzygies correspond to trisecants of certain scrolls.
This is another crucial point where we establish a correspondence between
algebraic properties of the syzygies and geometric properties of certain
scrolls.
\smallskip\noindent
Section 6 contains the main results of this paper, namely a complete
description of the syzygies of the bielliptic surfaces of degree 10
(see 6.3) and degree 15 (see 6.6).
\smallskip\noindent
In Section 7 we briefly review the abelian surfaces of degree 15 which
arise by liaison, whereas Section 8 contains a short discussion of the new
abelian surfaces of degree 15 which were originally found in \cite{Po}.
Finally,
in Section 9 we discuss some interesting degenerations of our surfaces.
Our considerations underline  in particular the fact that the
Hilbert scheme is connected \cite{Ha}.
\medskip
The first author, Alf Aure, died 17.10.94.  His death has been a great loss
to all of us. While this paper was finished after his death, Alf Aure
contributed to all the main ideas and supplied many of the crucial
arguments.
\medskip
\noindent
{\it{Acknowledgement}}. The authors would like to thank the DFG
(Grant Po 514/1-1 and Schwer\-punkt\-programm ''Komplexe Mannigfaltigkeiten'',
Grant Hu 337/4-3) and the Center for Advanced Studies, Oslo,
for partial support.
\bigskip

\head{1. The Heisenberg group of level 5}
\endhead
In this section we recall the representation theory of the
Heisenberg group $H_5$. The geometric significance of $H_5$
comes from the fact that its normalizer $N_5$ in $SL(5,\openC)$ is the
symmetry group of the Horrocks-Mumford bundle on $\openP^4$ \cite{HM},
\cite{De1}. In turn, the Horrocks-Mumford bundle provides a link between
some of the $H_5$-invariant varieties in $\openP^4$ which we study in
this paper. Some of these varieties are additionally invariant under
the involution $\iota\in N_5$ defined below. We also recall the
representation theory of the subgroup $G_5\subset N_5$ generated by $H_5$
and $\iota$.\medskip
Let $V$ be a $5$-dimensional $\openC$-vector space with
dual bases $e_0,\dots ,e_4$ of $V$ and $y_0,\dots,y_4$
of $V^*$. The Heisenberg group of level 5 in its Schr\"odinger representation
is the subgroup $H_5\subset SL(V)$ generated by the transformations
$$
\vbox{
\halign{&$#$\hfil\cr
\sigma : e_i\rightarrow e_{i-1}\cr
\tau : e_i\rightarrow \xi^i e_i\ ,&\quad \xi = e^{2\pi i/5}\ .\cr
}}
$$
$H_5$ has order 125 and is a central extension
$$
1\rightarrow \mu_5\rightarrow H_5\rightarrow \openZ_5 \times \openZ_5
\rightarrow 1\ ,
$$
where $\mu_5$ is the group of $5^{\text{th}}$ roots of unity. Let $\theta$ be
the generator of Gal $(\openQ(\xi):\openQ)$ given by
$\theta(\xi)=\xi^2$.\quad The Schr\"odinger representation of $H_5$ on $V$ is
an irreducible representation $V_0$ of $H_5$ and it gives rise to three
more: Let $V_i$ be the representation given by the composition
$$
H_5\overset{\theta^i}\to\rightarrow H_5 \rightarrow GL(V), \quad
i=\overline{1,3}.
$$
$V_0,\dots,V_3$ plus the 25 characters of $\openZ_5\times\openZ_5$
exhaust the irreducible representations of $H_5$.\par
\noindent The normalizer $N_5$ of $H_5$ in $SL(V)$ has order 15,000.
In fact,
$$
N_5/H_5\cong SL(2,\openZ_5),
$$
and
$$
N_5 = H_5 \rtimes SL(2,\openZ_5) \subset SL(V)
$$
is a semi-direct product. Let
$$
G_5 = H_5 \rtimes\openZ_2 \subset N_5
$$
be the subgroup generated by $H_5$ and the involution
$$
\iota:e_i\rightarrow e_{-i}\ .
$$
We recall the the character table of $G_5$ from \cite{Ma1}:
%
\par
%
$$
\vbox{\offinterlineskip
\halign{
\vrule height10.5pt depth 5.5pt#
&\hbox to 50pt{\hfil$#$\hfil}
&\vrule#
&\hbox to 100pt{\hfil$#$\hfil}
&\vrule#
&\hbox to 50pt{\hfil$#$\hfil}
&\vrule#
&\hbox to 50pt{\hfil$#$\hfil}
&\vrule#
\cr
\multispan 6\hrulefill\cr
&\{\alpha\}&&C_{m,n}&&C_\alpha&\cr
\noalign{\hrule}
&1&&1&&1&&I&\cr
\noalign{\hrule}
&5\theta^i(\alpha)&&0&&\theta^i(\alpha)&&V_i&\cr
\noalign{\hrule}
&1&&1&&-1&&S&\cr
\noalign{\hrule}
&5\theta^i(\alpha)&&0&&-\theta^i(\alpha)&&V_i^\#&\cr
\noalign{\hrule}
&2&&\xi^{sn+tm}+\xi^{-sn-tm}&&0&&Z_{s,t}&\cr
\noalign{\hrule}
}}
$$
%
\par\noindent
where $\{\alpha\}$ is the class containing only the central element
$\alpha\in\mu_5$,
$$
C_{m,n} = \{\alpha\, \xi^{2mn} \sigma^m \tau^n,\ \alpha\, \xi^{2mn}
\sigma^{-m} \tau^{-n} \mid\ \alpha\in\mu_5\}
$$
(there are 12 different $C_{m,n}$) and
$$
C_\alpha = \{\alpha\, \sigma^{2mn} \sigma^m \tau^n \iota \mid\ m,n \in
\openZ_5 \}
$$
(there are 5 classes $C_\alpha$). We denote by $Z$ the direct sum of the
12 different irreducible representations $Z_{s,t}$. If we restrict $Z$ to
$H_5$ we obtain the direct sum of the non-trivial
characters of $H_5$. We also need the following formulae from \cite{Ma1}:
$$V_i\otimes V_i\cong3V_{i+1}\oplus 2V^{\sharp}_{i+1},\qquad
V_i\otimes V_{i+1}\cong3V_{i+3}\oplus 2V^{\sharp}_{i+3},\qquad
V_i\otimes V_{i+2}\cong I \oplus Z$$
$${(V_i)}^{\ast}\cong V_{i+2},\qquad {(V^{\sharp}_{i})}^{\ast}\cong
V^{\sharp}_{i+2},\qquad V_i\otimes S\cong V^{\sharp}_i,\qquad V_i\otimes
Z\cong 12V_{i}\oplus 12V^{\sharp}_{i},$$
where $i\in\openZ_4$. For the exterior powers we have
$$\Lambda^{2}V_i\cong 2V^{\sharp}_{i+1},\qquad\Lambda^{3}V_i\cong
2V^{\sharp}_{i+3},\qquad \Lambda^{4}V_i\cong V_{i+2}.$$
\medskip
$\openP^4 = \openP^4(y) =\openP (V)$ will stand for the projective space of
lines in $V$ and
$$
R = \openC [y_0,\dots, y_4] = \underset {m\geq 0}\to{\textstyle\bigoplus}
S^m V^*
$$
for its homogeneous coordinate ring. If not otherwise mentioned, $R$ carries
the action of $G_5$ induced by the dual of the Schr\"odinger representation.
The well-known formula for the characters of the symmetric powers
$$
\chi_{\scriptscriptstyle S^m V}=\sum_{i=1}^{5}{(-1)}^{i+1}
\chi_{\scriptscriptstyle S^{m-i}V\otimes\Lambda^{i}V}
$$
yields
$$
H^0(\openP^4, \scr O_{\openP^4}(1)) \cong V_2,\qquad H^0(\openP^4,
\scr O_{\openP^4}(2)) \cong 3V_3,\qquad H^0(\openP^4, \scr O_{\openP^4}(3))
\cong 5V_1\oplus 2V^{\sharp}_1 $$
$$H^0(\openP^4, \scr O_{\openP^4}(4)) \cong 10V_0\oplus 4V^{\sharp}_0,\qquad
H^0(\openP^4, \scr O_{\openP^4}(5))\cong 6I\oplus 5Z.
$$
We write
$$
\openP^2_+ =\openP (W_+)\ ,\quad\openP^1_- =\openP (W_-)\subset
\openP (V)
$$
for the eigenspaces of $\iota$. Note that a subgroup of $N_5$ isomorphic
to the icosahedral group
$$
A_5\cong SL(2,\openZ_5)/<\iota>
$$
acts on $\openP^2_+$ and $\openP^1_-$ (compare \cite {BHM2, (1.1)}).
We will write $e_0=(1:0:\cdots :0), \dots, e_4=(0:\cdots :0:1)$
for the vertices of the standard simplex in $\openP^4$ with respect
to our fixed basis of $V$.\medskip
In Section 4 we consider different $\openP^4$'s with coordinates
denoted by $y$, $x$, $z$ and $w$. Certain natural rational maps
between these $\openP^4$'s are $G_5$-equivariant, if we consider
different $H_5$-actions as follows: $\openP^4(y)=\openP^4(V_0)$ as above,
$\openP^4(x)=\openP^4(V_1)$, $\openP^4(z)=\openP^4(V_0)$ and
$\openP^4(w)=\openP^4(V_2)$. Note, that a subscheme $\Sigma\subset\Pn 4$
which is $G_5$-invariant with respect to one of these actions is also
$G_5$-invariant with respect to any of the others.
In the context of different $\openP^4$'s, say e.g. $\openP^4(y)$ and
$\openP^4(x)$, we will use the following notation:
If $\Sigma=\Sigma(y)\subset\Pn 4(y)$ is a subscheme, then let
$\Sigma(x)\subset\Pn 4(x)$ denote the same subscheme $\Sigma$
with coordinates $y$ replaced by coordinates $x$.\par
\head
{2. Surfaces and syzygies}
\endhead
In this section we recall from \cite{DES} how to
construct surfaces in $\openP^4$ via
the syzygies of their Hartshorne-Rao modules. In fact, this is an
application of  Beilinson's spectral sequence.\medskip
$X$ will always denote a smooth surface in $\openP^4$ and\par
\noindent -\quad $d$ its degree\par
\noindent -\quad $\pi$ its sectional genus\par
\noindent -\quad $p_g$ its geometric genus\par
\noindent -\quad $q$ its irregularity\ .\medskip
In order to construct $\scr J_X$ and thus $X$ via syzygies one first
has to construct the {\it Hartshorne-Rao modules} of $X$:
$$
H^i_*\scr J_X = \underset{m\in\openZ}\to{\textstyle\bigoplus}
H^i (\openP^4,\scr J_X (m)),\quad i=1,2.
$$
Starting from these graded $R$-modules of finite length one is looking
for vector
bundles $\scr F$ and $\scr G$ on $\openP^4$ with rk $\scr G$ = rk
$\scr F +1$, and a morphism $\varphi\in\Hom (\scr F,\scr G)$, which
drops rank precisely along $X$. In fact, it is convenient to involve
suitable syzygy bundles of the
$H^i_*\scr J_X$ as direct summands of $\scr F$ and $\scr G$. Recall:
\proclaim
{Proposition 2.1}Let $M=\underset{k\in \openZ}\to\bigoplus M_k$ be a
graded $R$-module of finite length and let
$$
0\leftarrow M\leftarrow L_0\overset{\alpha_1}\to\leftarrow
L_1\leftarrow\dots \overset{\alpha_5}\to\leftarrow L_5\leftarrow 0
$$
be its minimal free resolution. Then for $1\leq i\leq 3$ the
sheafified syzygy module
$$
\scr F_i = \scr Syz_i(M) = (\ker \alpha_i)^\sim = (\Im \alpha_{i+1})^\sim
$$
is a vector bundle on $\openP^4$ with intermediate cohomology
$$
\underset {k\in\openZ}\to{\textstyle\bigoplus} H^j \left(\openP^4, \scr F_i
(k)\right) \cong \left\{
\aligned
M&\quad j=i\\
0&\quad j\neq i,\quad 1\leq j\leq 3\quad .
\endaligned
\right .
$$
Conversely, any vector bundle $\scr F$ on $\openP^4$ with this
intermediate cohomology is stably\par
\noindent equivalent with $\scr F_i$, i.e.,
$$
\scr F \cong \scr F_i\oplus \scr L\ ,\qquad \scr L \text{\quad a
direct sum of line bundles.}\qquad \qed
$$
\endproclaim
\proclaim
{Example 2.2} The Koszul complex
{\eightpoint
$$
0\leftarrow\openC (m)\leftarrow R(m)\leftarrow V^*\otimes R(m-1)\leftarrow
\Lambda^2 V^* \otimes R(m-2)\leftarrow \quad \dots \quad \leftarrow
\Lambda^5 V^*\otimes  R(m-5)\leftarrow 0\
$$}\noindent
resolves $\openC (m)$, where $\openC$ is considered as a graded
$R$-module sitting in degree $0$. It follows that
$$
\scr Syz_i(\openC (m)) = \Omega^i(m), 1\leq i\leq 3,
$$
where $\Omega^i= \Lambda^i T^*_{\openP (V)}$. By using the Koszul complex
one proves that
$$
\Hom(\Omega^i(i), \Omega^j(j))\cong\Lambda^{i-j}V,
$$
the isomorphism being defined by contraction.
\qed
\endproclaim
\noindent Which syzygy bundles to choose in order to construct a surface $X$
with given invariants comes out by analyzing  Beilinson's spectral sequence
for $\scr J_X (m)$, $m=3$ or $4$. Recall:
\proclaim{Theorem 2.3}\cite{Bei}\quad
For any coherent sheaf $\scr S$ on $\openP^4$ there is a spectral
sequence with $E_1$-terms
$$
E_1^{pq} = H^q\left(\openP^4, \scr S(p)\right)\otimes \Omega^{-p}(-p)
$$
converging to $\scr S$, i.e., $E_\infty^{pq} = 0$ for $p+q\neq 0$ and
$\bigoplus E_\infty^{-p,p}$ is the associated graded sheaf of a suitable
filtration of $\scr S$.\qquad\qed
\endproclaim
\noindent This fits nicely with the approach via the
Hartshorne-Rao modules since the $d_1$-differentials
$$
\split
d_1^{pq}
& \in \Hom \left(H^q\left(\openP^4, \scr S(p)\right)\otimes\Omega^{-p}(-p),
H^q\left(\openP^4, \scr S(p+1)\right)\otimes\Omega^{-p-1}(-p-1)\right)\\
&\cong\Hom \left(V^*\otimes H^q\left(\openP^4, \scr S(p)\right),
H^q\left(\openP^4, \scr
S(p+1)\right)\right)\
\endsplit
$$
coincide with the natural multiplication maps. In our case this typically
means that syzygy bundles of the $H_*^i\scr J_X$ will appear in the
$E_2$-diagram of the spectral sequence.
We will interpret one part of  Beilinson's spectral sequence for
$\scr J_X (m)$  as the spectral sequence of a vector
bundle $\scr F$, the other part as that of a vector
bundle $\scr G$ with rk $\scr G$ = rk $\scr F +1$. The differential between
the two parts will define a morphism $\varphi\in\Hom (\scr F,\scr G)$ whose
maximal minors vanish along the desired surface:
$$
0\rightarrow\scr F\overset\varphi\to\rightarrow\scr G\rightarrow \scr
J_X(m)\rightarrow 0\ .
$$
\noindent Then we can compute the explicit equations of $X$ from a
mapping cone between the minimal free resolutions of $\scr F$ and
$\scr G$.\par
\medskip
Some a priori information on the dimensions $h^i\scr J_X (m)$
can be obtained from Riemann-Roch:
\proclaim{Proposition 2.4 }
Let $X\subset\openP^4$ be a smooth surface. Then
$$
\chi(\scr J_X(m)) = \chi(\scr O_{\openP^4}(m))-\binom{m+1}2 d+m(\pi-1)-1+q-p_g\
.\qquad\qed
$$
\endproclaim
\noindent Moreover one knows:
\proclaim{Proposition 2.5} \cite{DES}
Let $X\subset \openP^4$ be a smooth, non-general type surface which is
not contained in any quartic hypersurface. Then we have the following
table for the $h^i\scr J_X (m)$:
%
$$
\vbox{\offinterlineskip
\halign{\hbox to 15pt{\hfil}
\vrule height10.5pt depth 5.5pt#
&\hbox to 50pt{\hfil$#$\hfil}
&\vrule#
&\hbox to 50pt{\hfil$#$\hfil}
&\vrule#
&\hbox to 50pt{\hfil$#$\hfil}
&\vrule#
&\hbox to 50pt{\hfil$#$\hfil}
&\vrule#
&\hbox to 50pt{\hfil$#$\hfil}
&\vrule#
&\hbox to 50pt{\hfil$#$\hfil}
&\vrule#
&\hbox to 15pt{\hfil$#$\hfil}
&#
\cr
\omit&&\omit\hbox{\hskip-3pt\hbox{$\bigg\uparrow$}\raise5pt\hbox{$i$}\hfil}\cr
\multispan{14}\hrulefill\cr
&0&&0&&0&&0&&0&&0&& \cr
\multispan{14}\hrulefill\cr
&N+1&&p_g&&0&&0&&0&&0&\cr
\multispan{14}\hrulefill\cr
&0&&q&&h^2\scr J_X(1)&&h^2\scr J_X(2)&&h^2\scr J_X(3)&&h^2\scr J_X(4)&\cr
\multispan{14}\hrulefill\cr
&0&&0&&0&&h^1\scr J_X(2)&&h^1\scr J_X(3)&&h^1\scr J_X(4)&\cr
\multispan{14}\hrulefill\cr
&0&&0&&0&&0&&0&&0&\cr
\multispan{14}\hrulefill
&\vbox to0pt{\vss\vskip5.75pt\hbox{$\!-\negthickspace\negmedspace@>> m
>$}\vss}\cr
}}
$$
\par\noindent
where
$$
N = \pi - q + p_g - 1\ .\qquad\qed
$$
\endproclaim\noindent
In the sequel we will represent a zero in a cohomology table
by an empty box.\par\medskip
We finally mention two classification results concerning the invariants
of those surfaces we are interested in here. Both results are an
immediate consequence of adjunction theory \cite{So}, \cite{SV}.
\proclaim
{Proposition 2.6} \cite{Ra, Proposition 9.1}\quad Let $X\subset\openP^4$
be a smooth surface with $d=10$ and $\pi=6$. Then $X$ is minimal and
either abelian, or bielliptic.\qquad\qed
\endproclaim
\noindent For the converse statement compare the next section.
\proclaim
{Proposition 2.7} \cite{Po, Lemma 7.1}\quad Let $X\subset \openP^4$ be a
smooth surface with $d=15,\ \pi=21$, and $\chi(\scr O_X)=0$. Then $X$ is
one of the following:\par
\noindent (i)\quad A ruled surface over an elliptic curve.\par
\noindent (ii)\quad A non-minimal abelian surface.\par
\noindent (iii)\quad A non-minimal bielliptic surface.\par
\noindent
In the cases (ii) and (iii) $X$ is embedded by
$$
H = H_{min} - \sum_{i=1}^{25} E_i\ .
$$
(with obvious notations). \qquad\qed
\endproclaim
\proclaim
{Remark 2.8}
It is not known, whether surfaces of type (i) exist.\qquad\qed
\endproclaim
\head
{3. The minimal abelian and bielliptic surfaces}
\endhead
In this section we review the minimal abelian surfaces in $\openP^4$ from
a syzygy point of view. These surfaces are classified as the smooth
zero-schemes of the Horrocks-Mumford bundle in a suitable twist, and they
may be constructed from the $H^1$-cohomology module of the bundle. It
follows, that this module is the first Hartshorne-Rao module for every
such surface. Its presentation matrix is known from \cite{De2}.
It is the type of this matrix which motivates our
considerations later on. We provide  further motivation by computing
the shape of the minimal free resolution of the first Hartshorne-Rao module
of a minimal bielliptic surface in $\openP^4$. From this computation we
may expect and will in fact prove in Section 6, that the presentation matrix
is of the same type as above. In this case, however, the surface is uniquely
determined by the module.
\medskip
We will abbreviate HM for Horrocks-Mumford.
Let $\scr E$ be the HM-bundle on $\openP^4 = \openP (V_0)$
normalized such that $c_1(\scr E)=-1$, $c_2(\scr E)=4$. By its construction
in \cite{HM} the bundle is acted on by the normalizer $N_5$ of
$H_5$ in $SL(V_0)$. $\scr E(3)$ has a four-dimensional space of sections.
Since $H^0 \scr E(2) = 0$ every $0\neq s\in H^0 \scr E(3)$ vanishes along
a surface, say $X_s$:
$$
0\rightarrow \scr O\overset s\to\rightarrow \scr E(3)\rightarrow\scr
J_{X_s}(5)\rightarrow 0\ .\tag3.1
$$
Every such $X_s$ is called a {\it HM-surface}. It is invariant under the
action of $G_5$.
$\scr E(3)$ and thus $\scr J_{X_s}(5)$ are globally generated outside
the 25 {\it HM-lines} $L_{ij}$ obtained from
$$
L_{00}=\openP_-^1
$$
by translating with the elements of
$$
H_5/\mu_5\cong \openZ_5\times\openZ_5\ .
$$
$L = \bigcup L_{ij}$ is cut out by the $H_5$-invariant quintics,
i.e., by the elements of
$$
H^0({\scr O}_{\openP^4}(5))^{H_5}\cong \Lambda^2 H^0({\scr E}(3)).
$$
Note that every $H_5$-invariant quintic is in fact also
invariant under $\iota$ (an explicit basis of
$H^0({\scr O}_{\openP^4}(5))^{H_5}$ will be given in Section 4).
Horrocks and Mumford show that for $s$ generic $X_s$ is smooth in the
points of $L$ too, and thus it is a minimal abelian surface with $d=10$ and
$\pi = 6$ (compare Proposition 2.6). Conversely, every minimal abelian
surface in $\openP^4$ arises in this way (up to coordinate transformations).
By (3.1) there are $G_5$-equivariant isomorphisms
$$
H_*^i\scr J_{X_s}(2) \cong H_*^i\scr E\,,\qquad i=1,2\ .
$$
In particular, the Hartshorne-Rao modules of $X_s$ are independent of $s$.
The $H^2$-module is just the vector space
$$
H^2\scr J_{X_s} \cong  H^2\scr E(-2) \cong  2S\ .
$$
\proclaim
{Proposition 3.2}\cite{HM}, \cite{De2}, \cite{Ma3}, \cite{Ma4}.
\quad $N=\bigoplus\limits_{k=0}^3 N_k = H_*^1\scr E$ is a graded
$R$-module with Hilbert function $(5,10,10,2)$, generated by
$N_0\cong V_3$. It has a minimal free resolution of type

$$
\vbox{%
\halign{&\hfil\,$#$\,\hfil\cr
0&\leftarrow&N&\leftarrow&5R&\overset\alpha_1\to\leftarrow
&15 R(-1)&&10R(-2)&&2R(-3)\cr
&&&&&&&\vbox to
10pt{\vskip-0pt\hbox{$\nwarrow$}\vss}&\oplus&\leftarrow&\oplus\cr
&&&&&&&&F_0&&F_1&\vbox to 10pt{\vskip-4pt\hbox{$\nwarrow$}\vss}&F_2&
\leftarrow&F_3&\leftarrow&0\ ,\cr
}}
$$
where the induced complex

%
$$
\vbox{%
\halign{&\hfil\,$#$\,\hfil\cr
&&&&&4R(-3)\cr
0&\leftarrow&&H_*^0\scr E\leftarrow&F_0 = &\oplus\cr
&&&&&15R(-4)&\vbox to 10pt{\vskip-4pt\hbox{$\nwarrow$}\vss}
&F_1 = 35 R(-5)&\leftarrow&F_2 = 20R(-6)
&\leftarrow&F_3 = 2R(-8)\leftarrow 0\cr
}}
$$

\noindent is the minimal free resolution of $H_*^0\scr E$.\qquad\qed
\endproclaim
\noindent Therefore via (3.1) we may compute from $N$ the explicit equations
of the minimal abelian surfaces in $\openP^4$ and their degenerations.
\proclaim
{Corollary 3.3} Let $X_s$ be a HM-surface. Then $\scr J_{X_s}$ has syzygies of
type
$$
\vbox{%
\halign{&\hfil\,$#$\,\hfil\cr
&&&&3\scr O(-5)\cr
0&\leftarrow&\scr J_{X_s}&\leftarrow&\oplus\cr
&&&&15\scr O(-6)&\vbox to
10pt{\vskip-4pt\hbox{$\nwarrow$}\vss}&35\scr O(-7)&\leftarrow&20\scr
O(-8)\cr
&&&&&&&&&\vbox to
10pt{\vskip-4pt\hbox{$\nwarrow$}\vss}&2\scr O(-10)&\leftarrow &0\ .\cr
}}
$$
In particular, $X_s$ is cut out by 3 quintics and 15 sextics. The quintics
alone cut out $X_s$ and the 25 HM-lines. These are the $6$-secant lines to
$X_s$.\qquad\qed\endproclaim
\noindent To compare later on with the bielliptic surfaces we also mention:
\proclaim
{Corollary 3.4} Let $X_s$ be a HM-surface. Then:\par
\noindent (i)\quad $\scr J_{X_s}$  has the following cohomology table:
$$
\vbox{\offinterlineskip
\halign{\hbox to 10pt{\hfil#}
&\vrule height10.5pt depth 5.5pt#
&\hbox to 25pt{\hfil$#$\hfil}
&\vrule#
&\hbox to 25pt{\hfil$#$\hfil}
&\vrule#
&\hbox to 25pt{\hfil$#$\hfil}
&\vrule#
&\hbox to 25pt{\hfil$#$\hfil}
&\vrule#
&\hbox to 25pt{\hfil$#$\hfil}
&\vrule#
&\hbox to 25pt{\hfil$#$\hfil}
&\vrule#
&\hbox to 25pt{\hfil$#$\hfil}
&\vrule#
&\hbox to 10pt{#\hss}
\cr
&\omit&&\omit\hbox{\hskip-3.015pt\hbox{$\bigg\uparrow$}
\raise5pt\hbox{$i$}\hfil}\cr
\multispan{17}\hrulefill\cr
&& && && && && && && &\cr
&\multispan{14}\hrulefill\cr
&&5&&1&& && && && && &\cr
&\multispan{14}\hrulefill\cr
&& &&2&& && && && && &\cr
&\multispan{14}\hrulefill\cr
&& && && &&5&&10&&10&&2&\cr
&\multispan{14}\hrulefill\cr
&& && && && && && &&3&\cr
\multispan{16}\hrulefill
&\vbox to 0pt{\vss\vskip 3.7pt\hbox to 40pt{\hskip-.6pt
\rightarrowfill}
\hbox to 33pt{\hfil$m$}\vss}\cr
}}
$$
\medskip
\noindent (ii)\quad Beilinson's spectral sequence for $\scr J_{X_s}(3)$
yields an exact sequence
$$
0\rightarrow \scr F = \scr O(-2)\oplus
2\Omega^3(3)\overset\varphi\to\rightarrow \scr G = \scr Syz_1
(N(1))\rightarrow \scr J_{X_s}(3)\rightarrow 0\ .\tag3.5
$$
\endproclaim
\demo
{Proof} (i)\quad We refer to \cite{HM} for the cohomology groups of $\scr E$.
\smallskip
\noindent (ii)\quad Comparing with the spectral sequence for $\scr E(1)$
we find that $5\Omega^1(1)\overset d_1^{-1,1}\to\rightarrow 10\scr O$
is surjective. Hence also $5\scr O(-1)\overset d_1^{-4,3}\to\rightarrow
\Omega^3(3)$
is surjective with kernel $\scr O(-2)$. By Beilinson's theorem we obtain a
resolution
$$
0\rightarrow \scr O(-2)\oplus 2\Omega^3(3)\rightarrow\scr
G\rightarrow\scr J_{X_s}(3)\rightarrow 0\ ,
$$
where $\scr G$ is the kernel
$$
0\rightarrow \scr G\rightarrow 5\Omega^1(1)\overset d_1^{-1,1}\to
\rightarrow 10\scr O\rightarrow 0\ .
$$
By Proposition 2.1
$$
\scr G = \scr Syz_1 (N(1))\
$$
since both bundles have the same rank and intermediate cohomology.
\qquad\qed
\enddemo
\proclaim
{Remark 3.6}
By Beilinson's spectral sequence applied to $\scr J_{X_s}(4)$ we reobtain
the exact sequence (3.1) tensored by $-1$.\qquad\qed
\endproclaim
We will next discuss the shape of the syzygies of the minimal
bielliptic surfaces in $\openP^4$. For this purpose we need to  recall
from \cite{ADHPR} some of the geometry of
the embedding of these surfaces. Let $W$ be a
$3$-dimensional $\openC$-vector space, let
$$
E = E_\lambda = \left\{x_0^3+x_1^3+x_2^3+\lambda x_0x_1x_2 =
0\right\}\subset \openP^2 = \openP(W)
$$
be a smooth elliptic curve in the Hesse pencil with origin $p_0 = (0:1:-1)$,
and let $z_0,\dots ,z_{14}$ be the canonical basis of $H^0\scr O_E (15
p_0)$. The Heisenberg group $H_{15}$ acts on $H^0\scr O_E(15 p_0)$ by
$$
\vbox{%
\halign{&\hfil\,$#$\,\hfil\cr
\sigma_{15} :z_i\rightarrow z_{i-1},&\cr
\tau_{15} : z_i\rightarrow \xi^{-i}_{15} z_i,&\qquad
\xi_{15} = e^{2\pi i/15}\ .\cr
}}
$$
We may identify the Heisenberg group $H_3 \subset SL(W)$ in its
Schr\"odinger representation
with the subgroup of $H_{15}$ generated by $\sigma_3 = \sigma_{15}^5$
and $\tau_3 = \tau_{15}^5$. Indeed,
$$
[\sigma_{15}^5,\tau_{15}^5] = \xi_3\cdot \id, \qquad\xi_{3} = e^{2\pi i/3} .
$$
So $H_3\times H_3$ acts naturally on
$W_E =\scr O_E(15 p_0)\otimes W^*$. The center of the diagonal $\Delta\subset
H_3\times H_3$ acts trivially on $W_E$. Hence $\scr E_E =
W_E^\vee /\Delta$ is a rank 3 vector bundle over
$E/\openZ_3\times\openZ_3 \cong E$. Consider the corresponding
projective quotient map
$$
E\times\openP^2 \rightarrow \openP^2_E = \openP (\scr
E_E)\ .
$$
The line bundle $\scr L' = \scr O_E (15 p_0)\boxtimes \scr
O_{\openP^2}(1)$ descends to a line bundle $\scr L$ on
$\openP^2_E$. Its sections
$$
\vbox{%
\halign{&\hfil\,$#$\,\hfil\cr
s_0&=&z_0\otimes x_0 &+& z_5\otimes x_1 &+& z_{10}\otimes x_2\cr
s_1&=&z_3\otimes x_0 &+& z_8\otimes x_1 &+& z_{13}\otimes x_2\cr
s_2&=&z_6\otimes x_0 &+& z_{11}\otimes x_1 &+& z_1\otimes x_2\cr
s_3&=&z_9\otimes x_0 &+& z_{14}\otimes x_1 &+& z_4\otimes x_2\cr
s_4&=&z_{12}\otimes x_0 &+& z_2\otimes x_1 &+& z_7\otimes x_2\cr
}}
$$
are invariant under $\Delta$ and define a basis of $H^0\scr L$.
We may identify $H_5\subset SL(V_0^*)$ in the dual Schr\"odinger
representation with the subgroup of $H_{15}$ generated by
$\sigma = \sigma_{15}^3$ and $\tau=\tau_{15}^6$. Indeed,
$$
[\sigma_{15}^3,\tau_{15}^6] = \xi^{-1}\cdot \id\ .
$$
The induced action of $H_5$ on $\scr L'$, with $H_5$ acting
trivially on $\scr O_{\openP^2}(1)$, commutes with the action of
$H_3$. On $H^0\scr L$ it defines the action
$$
\sigma(s_i) = s_{i-1}\,\qquad \tau(s_i) = \xi^{-i}s_i\ ,
$$
i.e., $H^0\scr L = V_0^*$ as $H_5$-modules.
The sections $s_0,\dots ,s_4$ map $\openP^2_E$ $H_5$-equivariantly to
a quintic hypersurface which is the trisecant variety
$$
\Trisec S^2 E\subset \openP^4 = \openP(V_0)\ .
$$
$\openP^2_E$ contains a pencil $|A_K|$ of abelian surfaces
isogeneous to a product and eight bielliptic surfaces
$B_{(i,j)}$, indexed by the non-trivial characters
$(i, j)\in \openZ_3 \times \openZ_3\setminus \{(0, 0)\}$. These surfaces are
embedded as smooth surfaces with
$d=10$ and $\pi=6$. Every minimal bielliptic surface in
$\openP^4$ arises in this way (up to coordinate
transformations). By construction, the surfaces $B_{(i,j)}$ are invariant
under the Schr\"odinger action of $H_5$ but not under the action
of $\iota$. Actually
$$
\iota\left(B_{(i,j)}\right) = B_{(-i,-j)}\ .
$$
\proclaim
{Proposition 3.7}
Let $B\subset\openP^4$ be a minimal bielliptic surface embedded
as above. Then:
\medskip
\noindent (i)\quad $\scr J_B$ has the cohomology table
$$
\vbox{\offinterlineskip
\halign{\hbox to 15pt{\hfil#}
&\vrule height10.5pt depth 5.5pt#
&\hbox to 25pt{\hfil$#$\hfil}
&\vrule#
&\hbox to 25pt{\hfil$#$\hfil}
&\vrule#
&\hbox to 25pt{\hfil$#$\hfil}
&\vrule#
&\hbox to 25pt{\hfil$#$\hfil}
&\vrule#
&\hbox to 25pt{\hfil$#$\hfil}
&\vrule#
&\hbox to 25pt{\hfil$#$\hfil}
&\vrule#
&\hbox to 25pt{\hfil$#$\hfil}
&\vrule#
&\hbox to 25pt{\hfil$#$\hfil}
&\vrule#
&\hbox to 10pt{#\hss}\cr
&\omit&\omit&\omit&\omit&\omit\hbox{\hskip-3.015pt\hbox{$\bigg\uparrow$}
\raise5pt\hbox{$i$}\hfil}\cr
\multispan{17}\hrulefill
&\vbox to 0pt{\vss\vskip 5pt\hbox to 15pt{\hskip-5pt\hrulefill}}\cr
&& && && && && && && && &&\cr
&\multispan{16}\hrulefill\cr
&&{20}&&5&& && && && && && &&\cr
&\multispan{16}\hrulefill\cr
&& && &&1&& && && && && &&\cr
&\multispan{16}\hrulefill\cr
&& && && && &&5&&10&&10&& &&\cr
&\multispan{16}\hrulefill\cr
&& && && && && && && &&1&&\cr
\multispan{17}\hrulefill
&\vbox to 0pt{\vss\vskip 3.75pt\hbox to 30pt{\hskip-.8pt\rightarrowfill}
\hbox to 20pt{\hfil$m$}\vss}\cr
}}
$$
\medskip
\noindent (ii)\quad Beilinson's spectral sequence for $\scr J_B(3)$
yields an exact sequence
$$
0\rightarrow \scr F = 5\scr O(-1)\oplus\Omega^3(3)\rightarrow
\scr G_B =\scr Syz_1 \left(N^B(1)\right)\rightarrow\scr
J_B(3)\rightarrow 0,\tag3.8
$$
where $N^B= \bigoplus\limits_{k=0}^2 N_k^B = H_*^1 \scr J_B(2)$
is a graded $R$-module with Hilbert function $(5,10,10)$,
generated by $N_0^B\cong  V_3$.
\medskip
\noindent (iii)\quad $N^B$ has a minimal free resolution of type
%
$$
\vbox{%
\halign{&\hfil\,$#$\,\hfil\cr
0\leftarrow N^B\leftarrow 5R\leftarrow 15R(-1)&&10R(-2)&&R(-3)\ \cr
&&\oplus\cr
&\vbox to 10pt{\vskip-4pt\hbox{$\nwarrow$}\vss}\ &R(-3)&\leftarrow &\oplus\cr
&&\oplus\cr
&&25 R(-4)&&55 R(-5)
&\vbox to 10pt{\vskip-12pt\hbox{$\nwarrow$}\vss}\ &40 R(-6)
\leftarrow 10 R(-7)\leftarrow 0\cr
}}.
$$
Equivalently, $\scr J_B$ has syzygies of type
%
$$
\vbox{%
\halign{&\hfil\,$#$\,\hfil\cr
&&&&\scr O(-5)\cr
0&\leftarrow&\scr J_B&\leftarrow&\oplus\cr
&&&&25\scr O(-6)&\vbox to 10pt{\vskip-4pt\hbox{$\nwarrow$}\vss}
&55\scr O(-7)&\leftarrow&40\scr O(-8)
&\leftarrow&10\scr O(-9)\leftarrow 0\cr
}}.
$$
In particular, $B$ is cut out by 1 quintic and 25 sextics, and it is
uniquely determined by $N^B$.
\endproclaim
\demo
{Proof}
(i)\quad The assertion follows from Riemann-Roch and Proposition 2.5.
First notice that $h^2\scr J_B(1) = \chi (\scr J_B(1)) = 0$. So by taking
cohomology in the exact sequences
$$
0\rightarrow\scr J_B(m-1)\rightarrow\scr J_B(m)\rightarrow\scr
J_H(m)\rightarrow 0
$$
associated to the generic hyperplane section $H$ we see that
$h^2\scr J_B(m) = 0$ for $m\geq 1$. This implies that
$$
h^1\scr J_B(5) = \chi (\scr J_B(5)) -1 = 0,
$$
since $B$ is contained in precisely one quintic hypersurface \cite{ADHPR}
.\smallskip
\noindent
(ii)\quad Beilinson's spectral sequence for $\scr J_B(3)$
yields a resolution
$$
0\rightarrow 5\scr O(-1)\oplus\Omega^3(3)\rightarrow\scr
G_B\rightarrow\scr J_B(3)\rightarrow 0\ ,
$$
where $\scr G_B$ is a kernel
$$
0\rightarrow \scr G_B\rightarrow 5\Omega^1(1)\rightarrow 10\scr
O\rightarrow 0\ .
$$
Then
$$
{\text{rank}}\,\scr G_B = 10 \leq {\text{rank}}\,\scr Syz_1
(N_B(1))\ .
$$
By Proposition 2.1 equality holds and
$$
\scr G_B = \scr Syz_1 (N^B(1))\ ,
$$
since both bundles have the same intermediate cohomology. Moreover,
Beilinson's spectral sequence represents $\scr J_B(2)$
as the cohomology of a monad
$$
0\rightarrow 20\scr O(-1)\rightarrow \scr B =
5\Omega^3(3)\oplus\Omega^2(2)\rightarrow N_0^B\otimes\scr
O\rightarrow 0\ ,
$$
and Bott's theorem gives $H^1_*\scr B = 0$. Hence by taking cohomology
in the display of the monad we see that $N_0^B$ generates $N^B$.
Part (i) implies that $N^B$ has Hilbert function $(5,10,10)$.
It remains to show that $N_0^B \cong  V_3$. Since $H^0\scr O_{\openP^4}(2)
\cong  3V_3$ this is equivalent to showing that $H^0\scr O_B(2) \cong  4V_3$.
With the notations introduced after Remark 3.6 it suffices to show
that $H^0\left(\scr O_E(30 p_0)\boxtimes \scr O_{\openP^2}(2)\right)$, or
equivalently that $H^0\scr O_E(30 p_0)$ consists of a direct sum of
copies of $V_3$ only. But since $E$, embedded by $\scr
O_E(15 p_0)$ in $\openP^{14}$, is projectively normal
\cite{Hu1}, $H^0\scr O_E (30 p_0)$ is a quotient of $H^0\scr
O_{\openP^{14}}(2) \cong  24 V_3$.
\smallskip
\noindent
(iii)\quad Finally we will compute the shape
of the minimal free resolution of $N^B$. Since
$$
{\text{rank}}\ \scr Syz_1 (N^B(1)) = 10
$$
the syzygies of $N^B$ are of type
%
$$
\vbox{%
\halign{&\hfil\,$#$\,\hfil\cr
0\leftarrow N^B\leftarrow 5R\leftarrow 15R(-1)&&10R(-2)&&aR(-3)&&bR(-4)\ \cr
&&\oplus&&\oplus&&\oplus\cr
&\vbox to 10pt{\vskip-4pt\hbox{$\nwarrow$}\vss}\ &aR(-3)&\leftarrow &cR(-4)&
\leftarrow&\dots&\leftarrow&\dots\cr
&&\oplus &&\oplus&&\oplus\cr
&&(25-b+c)R(-4)&&\dots&&\dots
\cr
}}.
$$
\noindent
{}From the mapping cone between the minimal free resolutions of $\scr F$
and $\scr G_B$ we see that necessarily $a = h^0 \scr J_B(5) =1$,
hence $b=0$ and  $c = h^0 \scr J_B(6)-25=\chi (\scr J_B(6))-25=0$ by
Riemann-Roch. The remaining graded Betti
numbers are uniquely determined by the Hilbert function.\qquad\qed
\enddemo
So far we have only discussed the shape of the syzygies of the modules
involved. On our way we will explicitly compute the presentation matrices
of these modules. For the abelian surfaces this was already done in
\cite{De2}. Namely, by using the monad of the HM-bundle it is shown that
$$
\alpha_1 =
\pmatrix
y_0& 0& 0& 0& 0& 0& y_4& 0& 0& y_1& 0& 0& y_3& y_2& 0\\
0& 0& 0& 0& y_2& y_1& 0& 0& y_3& 0& 0& y_0& y_4& 0& 0\\
0& 0& 0& y_4& 0& 0& 0& y_0& 0& y_3& y_2& y_1& 0& 0& 0\\
0& 0& y_1& 0& 0& 0& y_2& 0& y_0& 0& y_3& 0& 0& 0& y_4\\
0& y_3& 0& 0& 0& y_4& 0& y_2& 0& 0& 0& 0& 0& y_1& y_0
\endpmatrix
\
$$
is the presentation matrix of $N$. $\alpha_1$ is obtained by
concatenating three $5\times5$-matrices whose type we
will study next.
\head
{4. Moore matrices}
\endhead
In this section, which is the heart of our paper, we introduce
matrices which where first considered by R. Moore \cite{Moo}. These
matrices, which we call Moore matrices and syzygy matrices resp.,
arise naturally from the representation theory of $H_5$ and are the
crucial building blocks in the minimal free resolutions of
$G_5$-invariant elliptic normal curves and elliptic quintic
scrolls in $\openP^4$ resp. More general Moore matrices play also a
crucial role in \cite{GP1} where they provide equations for
$(1,d)$-polarized abelian surfaces. We recall in the sequel some of the known
geometry of elliptic curves and their associated varieties, but we also
prove some new results which are of a geometric nature. In particular,
Proposition 4.13, Lemma 4.18 and Proposition 4.19 are new. Furthermore, we
study the interplay between varieties associated to elliptic curves
and our matrices.
This is new and crucial for everything which comes later. Therefore it is
helpful to keep the following geometric picture in mind.  In this picture
there are two kinds of $G_5$-invariant elliptic curves and two kinds of
$G_5$-invariant elliptic scrolls. The curves are the elliptic normal curves
$E$ of degree 5 and certain elliptic curves $E^{\prime}$ of degree 15, while
the scrolls are the elliptic quintic scrolls $Q$ and certain singular elliptic
scrolls $S$ of degree 15. These curves and these scrolls all depend on a
parameter point in $\Pn 1$ which we omit in our notation at this point,
see Propositions 4.2, 4.10, 4.3, and 4.16 resp. for precise descriptions.
There are certain 3-3 correspondences between these curves and surfaces.
Each scroll $Q$ contains three curves $E$ as bisections, and each curve $E$
is contained in three scrolls $Q$.  Similarly, each scroll $S$ contains three
curves $E^{\prime}$ as bisections, and each curve $E^{\prime}$ lies on three
scrolls $S$. On the other hand, each curve $E^{\prime}$ is a section of a
unique scroll $Q$ and vice versa, and each curve $E$ is the 4-tuple singular
locus of a unique scroll $S$ and vice versa.  Finally, the union of the curves
$E$ is a surface $S_{15}$ of degree 15 with a $G_5$-equivariant birational map
to a
surface $S_{45}$ of degree 45 which is the union of the curves $E^{\prime}$
(Proposition 4.10). The first and the last correspondence is described in
Propositions 4.3 and 4.16 resp., the proofs of the others, which are
similar, are omitted since we do not need them explicitly.  The
significance of these curves and surfaces for the Moore matrices $M_p$ and
the syzygy matrices $L_q$ defined below lies in the fact that these
matrices have linear entries which depend linearly on a point in $\Pn 4$.
When the point $p$ is a general point on a curve $E$, then we may identify
the hypersurface defined by $\det M_p$ with the secant variety of the curve
$E$ and the hypersurface defined by $\det L_p$ with the trisecant variety of
a scroll $Q$ over $E$. On the other hand, when the point $p$ is a general
point on $E^{\prime}$, then we may identify the hypersurface defined by
$\det M_p$ with the trisecant variety of the unique scroll $Q$ containing
$E^{\prime}$ and the hypersurface defined by $\det L_p$ with the secant
variety of a curve $E$ isomorphic to $E^{\prime}$ via the map from $S_{15}$
to $S_{45}$. This is made precise in Proposition 4.12.
The scroll $S$ comes into play when we consider the incidence
$$
I=\{(p,q)\mid \quad\text{The quadratic entries of}\quad M_pL_q\quad
\text{all lie in the ideal of}\quad E\}.
$$
This incidence variety turns out to be an elliptic scroll over the curve
$E$, which is mapped onto a scroll $S$ in the $\Pn 4$ of quadrics in the
ideal of $E$ (Lemmas 4.14, 4.15).\par\noindent
All surfaces in this paper may be constructed via syzygies from modules
presented by a matrix $\alpha_1$ which is obtained by concatenating three
Moore matrices $M_{p_i}$, $i=\overline{1,3}$. And any five first order
linear syzygies of $\alpha_1$ will be given by a block of syzygy matrices
$L_{q_i}$, $i=\overline{1,3}$. To show that the module presented by such a
matrix $\alpha_1$ gives rise to a surface via a syzygy construction, we need
to show that the module is artinian and that it has the correct number of
linear and quadratic syzygies.  This is done for our choices of points
$p_1$, $p_2$ and $p_3$ in Section 5 using the properties of the Moore matrices
and syzygy matrices obtained in this section. For the bielliptic
surfaces we choose the points $p_i$ on an elliptic curve $E$ and we show in
Section 5 how double points and trisecant lines to the scroll $S$ via the
above incidence give rise to ten first order linear syzygies for $\alpha_1$
in blocks $L_{q_{ij}}$, $i=\overline{1,3}$, $j=\overline{1,2}$ (Lemma 5.7
and Proposition 5.8). In fact, given
the points $p_i$ on $E$ we can identify the corresponding points $q_{ij}$ on
a scroll $Q$. More precisely, these points $q_{ij}$ lie on curves of type $E$
or $E^{\prime}$ on $Q$. In Section 6 we use this identification to show that
the matrix $\alpha_1$ gives rise to a surface on the trisecant variety
of a scroll $Q$ or on the secant variety of a curve $E$ (Theorems 6.3 and 6.6).
Now the 3-3 correspondences above come into this picture. The module presented
by $\alpha_1$ depends only on the plane spanned by the points $p_i$, the
{\it{representing plane}} of this module. For the minimal bielliptic surfaces
this plane meets a unique scroll
$Q$ in a line and three points, and it is a trisecant plane to all three curves
$E$ on this scroll.  For the non-minimal bielliptic surfaces the representing
plane meets three scrolls $Q$ in a line and three points. These three scrolls
intersect along a curve $E$ and the plane is trisecant to this unique curve
$E$. For our purposes it is easier to work with a fixed curve $E$, so this
3-3 correspondence will not be explicit in our arguments.
Finally, each scroll $Q$ may be identified with the symmetric product
$S^2(E)$ of one of the curves $E$ (Remark 4.7, (iii)), so that one may
use the arithmetic on $E$ to make
computations on $Q$ (Proposition 4.16 ff., cf. [Hu1], [ADHPR]).
\medskip
Matrices of type
$$
M(x,y) := (x_{3i-3j} y_{3i+3j})_{i,j\in\openZ _5}\quad {\text {and}}\quad
L(z,y) := (z_{i-j} y_{2i-j})_{i,j\in\openZ _5}\quad {\text {resp.}},
$$
will be called {\it Moore matrices} and {\it syzygy matrices} resp. We will
also consider these matrices with the $x$'s, $y$'s or $z$'s being replaced
by the coordinates of a fixed point in $\Pn 4(x)$, $\Pn 4(y)$ or $\Pn 4(z)$,
where we will understand,
that these matrices are determined up to scalars only. For example, we will
write
$$
M_p(y) = M(p,y)\quad {\text {and}}\quad L_q(y) = L(q,y),
$$
if $p = (p_0:\dots :p_4)\in \Pn 4(x)$ and $q = (q_0:\dots :q_4)\in \Pn 4(z)$
are fixed parameter points.\par\noindent
Blocks of such matrices show up if we tensor the Koszul complex resolving
$\openC (m)$ by $V_3$ and decompose into irreducible $H_5$-modules:
$$
... \leftarrow V_3\otimes  R(m)\overset \tilde M\to\leftarrow 5V_1\otimes
R(m-1)\overset \tilde L\to\leftarrow 10V_0\otimes R(m-2)\leftarrow ...\quad .
$$
More precisely, let us use the notation of Section 1 and choose decompositions
of
$V_3\otimes V_0^*$ and $V_1\otimes V_0^*$ resp. into
irreducible $H_5$-modules $V_1^{(0)}, \dots ,V_1^{(4)}$ of type $V_1$ and
$V_0^{(0)},
\dots ,V_0^{(4)}$ of type $V_0$ resp. Namely, let
$$
V_1^{(k)} = \text{span}(f_0^{(k)}, \dots ,f_4^{(k)})\quad \text{with} \quad
f_i^{(k)} = e_{2k+i}\otimes y_{k+i}
$$
and
$$
V_0^{(k)} = \text{span}(g_0^{(k)}, \dots ,g_4^{(k)})\quad \text{with} \quad
g_i^{(k)} = e_{k+i}\otimes y_{2k+i}.
$$
With respect to the first decomposition every element in
$$
\openP(\Hom_{H_5}(V_1\otimes R(m-1),V_3\otimes R(m)))
\cong \openP(\Hom_{H_5}(V_1,V_3\otimes V_0^*))
$$
is given by the choice of a $H_5$-invariant subspace
$$
V^p\subset V_3\otimes V_0^*\cong V_1^{(0)}\oplus \dots
\oplus V_1^{(4)}, \quad p = (p_0:\dots :p_4)\in
\Pn 4:
$$
The $i^{\text{th}}$ basis vector of $V^p$ is the linear combination
of $f_i^{(0)}, \dots ,f_i^{(4)}$ with $p_0,\dots, p_4$ as coefficients.
Similarly, we may use the second decomposition to describe the elements in
$$
\openP(\Hom_{H_5}(V_0\otimes R(m-2),V_1\otimes R(m-1)))
\cong \openP(\Hom_{H_5}(V_0,V_1\otimes V_0^*)).
$$
With respect to these choices and with respect to our fixed basis $e_{i}$,
$i\in\openZ_5$, of $V_0$, $V_1$ and $V_3$ we get isomorphisms
$$
\Pn 4(x) \cong \openP(\Hom_{H_5}(V_1\otimes R(m-1),V_3\otimes R(m)))
, \quad p\mapsto M_p(y)
$$
and
$$
\Pn 4(z) \cong \openP(\Hom_{H_5}(V_0\otimes R(m-2),V_1\otimes R(m-1)))
, \quad q\mapsto L_q(y).
$$
\medskip
Let us describe the determinants of our matrices.
Following \cite{Moo} one can choose a basis $s_1, ..., s_4$ of
$H^0 {\scr E}(3)$ such that
$$
\matrix
\gamma_0=s_1\wedge s_2=5y_0y_1y_2y_3y_4 &
\gamma_1=s_1\wedge s_3=\sum_{i\in\openZ_5}y_iy_{i+2}^2y_{i+3}^2 \\
\gamma_2=s_1\wedge s_4=-\sum_{i\in\openZ_5}y_i^3y_{i+2}y_{i+3} &
\gamma_3=s_2\wedge s_3=\sum_{i\in\openZ_5}y_i^3y_{i+1}y_{i+4} \\
\gamma_4=s_2\wedge s_4=-\sum_{i\in\openZ_5}y_iy_{i+1}^2y_{i+4}^2 &
\gamma_5=s_3\wedge s_4={1\over 5}\sum_{i\in\openZ_5}y_i^5 -
\prod_{i\in\openZ_5} y_i\\
\endmatrix
$$
\noindent
form a basis of
$$
H^0({\scr O}_{\openP^4}(5))^{H_5}\cong \Lambda^2 H^0({\scr E}(3)).
$$
The rational map
$$
\Theta : \openP^4 \dashrightarrow {{\openP^{5}}^\ast}=
{\openP(H^0({\scr O}_{\openP^4}(5))^{H_5})}^{\ast},
\quad y\mapsto (\gamma_0(y):\gamma_1(y):\dots :\gamma_5(y))
$$
is well-defined outside the union $L$ of the HM-lines and has as
image a smooth hyperquadric $\Omega$ whose dual quadric
$\bar\Omega=\{h_2h_3-h_1h_4+h_0h_5=0\}\subset{\openP^{5}}(h)$
may be identified with the Pl\"ucker embedding of
$\Grass(2,H^0({\scr E}(3)))$ in
$\openP(\Lambda^2 H^0({\scr E}(3)))$.
By composing $\Theta$ with the Gauss map we obtain
$$
\bar\Theta : \Pn 4\setminus L \rightarrow
\bar\Omega\subset\openP(H^0({\scr O}_{\openP^4}(5))^{H_5})
\cong \openP(\Lambda^2 H^0({\scr E}(3))).
$$
The tangent hyperplane $T_{\bar\Omega,\bar\Theta(y)}$ parametrizes those
$H_5$-invariant quintics which contain $y\in\openP^4$, and the quintics
in $T_{\bar\Omega,\bar\Theta(y)}$ $\cap$ $\bar\Omega$ correspond
to decomposable tensors $s\wedge t \in\Lambda^2 H^0({\scr E}(3))$ with
$s(y)=0$. Hence $\bar\Theta$ can be seen as the map which
associates to a point
$y\in\openP^4\setminus L$ the pencil of sections of $\scr E(3)$ vanishing in
$y$. In particular, $\Theta$ is generically $100:1$ \cite{Moo}.
So via $p\mapsto\Theta^{-1}(T_{\Omega,{\Theta(p)}})$,
where $p\in\openP^4(y)$, we obtain all $H_5$-invariant quintics which
correspond to decomposable tensors, i.e., those
which correspond to points in the Grassmannian of lines in
$\openP(H^0({\scr E}(3)))$. Every such quintic will be called
a {\it HM-quintic}. The HM-quintics are exactly the determinantal
hypersurfaces associated to Moore matrices and syzygy matrices resp.:
\proclaim {Remark 4.1}
(i)\quad For Moore matrices we observe that\medskip
$\det M(x,y)=\det M(y,x)=$
$$
=\gamma_5(x)\gamma_0(y)- \gamma_4(x)\gamma_1(y)+
\gamma_3(x)\gamma_2(y)+ \gamma_2(x)\gamma_3(y)-
\gamma_1(x)\gamma_4(y)+ \gamma_0(x)\gamma_5(y).
$$
Therefore $\{\det M_p(y)=0\}=\Theta^{-1}(T_{\Omega,
(\gamma_0(p):\gamma_1(p):\dots :\gamma_5(p))})$ if $p$ is not on a HM-line.
\par\noindent
(ii)\quad For the syzygy matrices we get\medskip
$\det L(z,y)=$
$$
=\gamma_5(z)\gamma_0(y)+\gamma_2(z)\gamma_1(y)+
\gamma_4(z)\gamma_2(y)+ \gamma_1(z)\gamma_3(y)+
\gamma_3(z)\gamma_4(y)- \gamma_0(z)\gamma_5(y).
$$
Hence
$\{\det L_q(y)=0\}=\Theta^{-1}(T_{\Omega,
(-\gamma_0(q):-\gamma_3(q):\gamma_1(q):\gamma_4(q):-\gamma_2(q):\gamma_5(q))})$
if $q$ is not on a HM-line.\qed\endproclaim
\medskip
We will show next that the Moore matrices and the syzygy matrices are the
crucial building blocks in the minimal free resolutions of $G_5$-invariant
elliptic normal curves and elliptic quintic scrolls in $\Pn 4$ resp.
Let us first recall in the following two propositions a number of standard
facts concerning the geometry of these varieties. For proofs we refer to
\cite{Hu1}, \cite{BHM1}, \cite{BHM2} and \cite{Moo}.
\proclaim {Proposition 4.2}
(i)\quad
Elliptic curves with a level-5 structure correspond bijectively to
$G_5$-invariant elliptic normal curves in $\Pn 4$, i.e., to
elliptic normal curves in $\Pn 4$ on which $H_5$ operates by
translations with 5-torsion points and $\iota$ as the reflection in the
origin given by the intersection of the curve with $\openP^1_-$.
There exists a 1-dimensional family $E_{(\lambda:\mu)}$
of $G_5$-invariant elliptic normal curves in $\Pn 4$ (possibly degenerate).
$E_{(\lambda:\mu)}$ intersects $\Pn 1_-$ in $o_{E_{(\lambda:\mu)}}=
(\lambda:\mu)$.\par\noindent
(ii)\quad
Let $\Lambda = \{(0:1),(1:0),((1\pm\sqrt{5})\xi^k:2), k\in\openZ_5\}$.
Then for every $(\lambda:\mu)\in\Pn 1_-\setminus\Lambda$,
the curve $E_{(\lambda:\mu)}$ is smooth. For
$(\lambda:\mu)\in\Lambda$, however, $E_{(\lambda:\mu)}$ is a connected cycle of
5 lines, i.e., a pentagon. Notice that we may identify the elements of
$\Lambda$ with the vertices of  an icosahedron sitting inside
$S^2\cong \Pn 1$ (cf. \cite{Kl}).\par\noindent
(iii)\quad Two curves of type $E_{(\lambda:\mu)}$ intersect iff the
corresponding parameters belong to opposite vertices of the icosahedron.
In this case the two singular curves have common vertices and thus form
a complete pentagon. Under the action of $N_5$ these pentagons are
equivalent to the one with vertices $e_0,\dots ,e_4$, i.e., to
$E_{(0:1)}\cup E_{(1:0)}$.\par\noindent
(iv)\quad
The union $S_{15}=\bigcup_{(\lambda:\mu)\in\Pn 1_-} E_{(\lambda:\mu)}$ is an
irreducible determinantal surface $S_{15}$ of degree 15. It is smooth outside
the 30 vertices of the 6 complete pentagons formed by the singular curves
$E_{(\lambda:\mu)}$, ${(\lambda:\mu)}\in\Lambda$. The normalization of
$S_{15}$ is isomorphic to Shioda's modular surface $S(5)$. It thus has
a natural structure of an elliptic surface over the modular curve
$X(5)\cong\Pn 1$ of level 5 . The smooth fibres are isomorphic to the
smooth curves $E_{(\lambda:\mu)}$.\par\noindent
(v)\quad The $5$-torsion points in the smooth fibres give rise to
25 sections of $S(5)$ which are mapped 1:1 under $S(5) \rightarrow S_{15}$
onto the Horrocks-Mumford lines. In particular, $X(5)$ may be canonically
identified with $L_{00}=\Pn 1_-$. We will denote these sections also
by $L_{ij}$, $(i,j)\in\openZ_5\times\openZ_5$.\par\noindent
(vi)\quad The $3$-section of $S(5)$ which intersects
each smooth fibre in its non-trivial $2$-torsion points
is mapped under $S(5) \rightarrow S_{15}$ onto a plane, nodal, irreducible
sextic curve $B\subset\Pn 2_+$ (Brings curve). The equation of $B$ in
$\Pn 2_+$ is
$$
y_0^4y_1y_2-y_0^2y_1^2y_2^2-y_0(y_1^5+y_2^5)+2y_1^3y_2^3=0.
$$
$B$ has double points at the six points
$$
(1:0:0),\quad (1:\xi^i:\xi^{-i}),\quad i\in \openZ_5,
$$
which form the unique minimal $A_5$-orbit in $\Pn 2_+$.\par\noindent
(vii)\quad
Every smooth curve $E_{(\lambda:\mu)}$ meets $B$ in three points. These are
the vertices of a triangle tangent to the fixed conic
$C_+=\{y_0^2+4y_1y_2=y_1-y_4=y_2-y_3=0\}\subset\Pn 2_+$. We will identify
$X(5)\cong\Pn 1_-$ with $C_+$ via the $A_5$-equivariant map
$\psi_5:(\lambda:\mu)\mapsto (2\lambda\mu:\mu^2:-\lambda^2)$.\qed
\endproclaim
\noindent
\proclaim{Proposition 4.3}
(i)\quad  For $(\lambda:\mu)\in\Pn 1_-\setminus\Lambda$
and $\tau_i$ a non-trivial $2$-torsion point on  $E_{(\lambda:\mu)}$ the
translation scroll
$$
\Sigma_{(\lambda:\mu)}(\tau_i):=
\bigcup_{p\in E_{(\lambda:\mu)}}{\overline{p,p+\tau_i}}\subset\Pn 4
$$
is a smooth $G_5$-invariant elliptic quintic scroll in $\Pn 4$ over the
curve $E_{(\lambda:\mu)}/<\tau_i>$. $E_{(\lambda:\mu)}$ is
contained in $\Sigma_{(\lambda:\mu)}(\tau_i)$ as a bisection of the ruling
of this scroll. Conversely, any $G_5$-invariant elliptic quintic
scroll $\Sigma$ in $\Pn 4$ arises in this way.\par\noindent
(ii)\quad  Given $\Sigma$ there exist exactly three pairs $(E_i,\tau_i)$,
$i = \overline{0,2}$ of elliptic curves and $2$-torsion points as above,
such that $\Sigma$ is the $\tau_i$-translation scroll of $E_i$. In
particular, $\Sigma$ meets $\Pn 2_+$ along the disjoint union of
the $\tau_i$, $i = \overline{0,2}$ and the ruling of $\Sigma$ over the
origin of the base curve.\par\noindent
(iii)\quad  Analogously, given a smooth $G_5$-invariant elliptic normal curve
$E_{(\lambda:\mu)}\subset\Pn 4$ there exist exactly three
$G_5$-invariant elliptic quintic scrolls containing it,
each two of them meeting transversally along the elliptic curve.
Namely, these are the translation scrolls
$\Sigma_{(\lambda:\mu)}(\tau_i)$, $\tau_i$ a non-trivial $2$-torsion
point on $E_{(\lambda:\mu)}$.\smallskip\noindent
We may rephrase these results as follows:\smallskip\noindent
(iv)\quad  There exists a 1-dimensional family $Q_{(\lambda:\mu)}$,
$(\lambda:\mu)\in\Pn 1_-$, of $G_5$-invariant elliptic quintic scrolls
in $\Pn 4$ (possibly degenerate). The ruling of a smooth
$Q_{(\lambda:\mu)}$ over the origin of the base curve is the line
$l_{(\lambda:\mu)}=\{y_1-y_4=y_2-y_3=
\lambda \mu y_0 + {\mu}^2 y_2 - {\lambda}^2y_1=0\}
\subset \Pn 2_+$. This line is tangent to the conic $C_+$ at the point
$(2\lambda\mu:\mu^2:-\lambda^2)$, which corresponds via
$\psi_5$ to $(\lambda:\mu)\in\Pn 1_-\cong X(5)$.\par\noindent
(v)\quad  For every $(\lambda:\mu)\in\Pn 1_-\setminus\Lambda$ the
surface $Q_{(\lambda:\mu)}$ is smooth. Otherwise $Q_{(\lambda:\mu)}$
is a cycle of 5 planes.\par\noindent
(vi) $Q_{(\lambda:\mu)}$ contains precisely the curves
$E_{(\tilde\lambda:\tilde\mu)}$ such that $(\tilde\lambda:\tilde\mu)$
is a root of the polynomial
$$
P_{3,3} = -\lambda^2\mu\tilde\lambda^3+\mu^3\tilde\lambda^2\tilde\mu
+\lambda^3\tilde\lambda\tilde\mu^2+\lambda\mu^2\tilde\mu^3
$$
with $(\lambda:\mu)$ fixed. If $Q_{(\lambda:\mu)}$ is
singular then it contains 3 pentagons with 2 of them coming together.
Conversely, $P_{3,3}$ determines the surfaces $Q_{(\lambda:\mu)}$ containing
a fixed curve $E_{(\tilde\lambda:\tilde\mu)}$.
\endproclaim
\demo{Proof} To prove (v) we observe that  $Q_{(\lambda:\mu)}$
is singular iff it contains a pentagon $E_{(\tilde\lambda:\tilde\mu)}$.
The assertion follows from the corresponding smoothness statement for
the elliptic curves and (vi) above.\qed
\enddemo
\proclaim
{Lemma 4.4}
(i)\quad The Moore matrix $M(x,p)$ is skew-symmetric (up to column
permutations) for $p = (\lambda :\mu)\in\openP_-^1(y)$ and its
$4\times 4$-Pfaffians define the $G_5$-invariant elliptic normal curve
$E_p = E_p(x)$. In particular, $E_p$ has a minimal free resolution of type
$$
0\leftarrow{\scr O}_{E_p}\longleftarrow{\scr O}
\overset {I_{E_p}}\to\longleftarrow V_3\otimes\O {-2}
\overset {M(x,p)}\to\longleftarrow V_1\otimes \O {-3}
\overset {I_{E_p}^t}\to\longleftarrow \O {-5}\leftarrow 0.
$$
\smallskip\noindent
(ii)\quad For $p = (\lambda :\mu)\in\openP_-^1(y)$
the $4\times 4$-minors of the syzygy matrix $L(\iota (z),p)$ define the
$G_5$-invariant elliptic quintic scroll
$Q_p=Q_p(z)\subset\Pn 4(z)$. In particular, $Q_p$
has a minimal free resolution of type
$$
0\leftarrow{\scr O}_{Q_p}\longleftarrow{\scr O}
\overset {I_{Q_p}}\to\longleftarrow V_1\otimes\O {-3}
\overset {L(\iota (z),p)^t}\to\longleftarrow V_0\otimes\O {-4}
\overset {[z]^t}\to\longleftarrow \O {-5}\leftarrow 0
$$
with $[z]=(z_0,\dots ,z_4)$.
\endproclaim
\demo
{Proof}
(i)\quad is well-known, compare e.g. \cite{Bi}, \cite{Kl}.\par\noindent
(ii)\quad Let $p\in\openP_-^1(y)$ and $\tilde Q_p=
\{z\mid\rk L(\iota(z),p)\le 3\}$. A straightforward calculation shows that
the saturated ideal of $\tilde Q_p$ contains $5$ cubics which define a ruled
surface of degree $5$ (see \cite{Moo, 4.2}). On the other hand,
$(z_0,\dots,z_4)\cdot L(\iota(z),p)=0$. Hence $L(\iota(z),p)^t$ defines
a morphism from $\scr F=\Omega^3(3)$ to $\scr G=V_1\otimes{\scr O}$ which
drops rank along $\tilde Q_p$. Altogether $\tilde Q_p$ is a  $G_5$-invariant
elliptic quintic scroll with the claimed minimal free resolution (compare e.g.
\cite{Ok1}). Moreover, $L(\iota(z),p)$ has rank 3 along the line $l_p\subset
\openP_+^2$ which implies that $\tilde Q_p = Q_p$ by Proposition
4.3, (iv). \qed
\enddemo
\proclaim{Corollary 4.5} Let $p = (\lambda :\mu) \in \openP^1_-$.
A set of defining equations for
$E _{(\lambda:\mu)}(x)$ is given by
$$
q^{(\lambda:\mu)}_i (x) = -\lambda\mu x^2_i -
\mu^2 x_{i+1} x_{i+4}+\lambda^2 x_{i+2}x_{i+3},
$$
and a set of defining equations for $Q_{(\lambda:\mu)}(z)$ is given by
$$\align
c^{(\lambda:\mu)}_i (z) =&
\lambda^2\mu^2 z_i^3+\lambda^3\mu(z_{i+1}^2z_{i+3}+ z_{i+2}z_{i+4}^2)\\
&-\lambda\mu^3(z_{i+1}z_{i+2}^2+z_{i+3}^2z_{i+4})
-\lambda^4 z_iz_{i+1}z_{i+4}-\mu^4 z_iz_{i+2}z_{i+3},
\endalign
$$
where $i\in\openZ_5$.\qed
\endproclaim
\proclaim{Remark 4.6}
As in Lemma 4.4 we obtain: $M(p,y)$ is skew-symmetric
for $p=(\lambda:\mu)\in\openP_-^1(x)$
and its $4\times 4$-Pfaffians define $E_{(\lambda:\mu)}
=E_{(\lambda:\mu)}(y)$ in $\Pn 4(y)$. On the other side,
$L(q,y)(y_0,\dots , y_4)^t = 0$ for $q=(\lambda:\mu)\in\openP_-^1(z)$
and in this case the $4\times 4$-minors of the matrix $L(q,y)$ define
$Q_{(\mu:-\lambda)}\subset\Pn 4(y)$.\qed
\endproclaim
\proclaim{Remark 4.7} Let $p = (\lambda :\mu) \in \openP^1_-$.\par
\noindent (i)\quad By applying Beilinson's theorem to ${\scr I}_{Q_p}(2)$
we get a resolution of type
$$
0\longleftarrow{{\scr I}_{Q_p}(2)}
\longleftarrow \Omega^2(2)\overset {\varphi}\to
\longleftarrow 5\O {-1}\leftarrow 0.
$$
$\varphi$ is a $1\times5$-matrix with entries in $\Lambda^2 V$
(see Example 2.2). By comparing the mapping cone over $\varphi$
with the minimal free resolution in Lemma 4.4 we see that
$$
\varphi=(\mu f^+(e_{i})+\lambda f^-(e_{i}))_{i\in\openZ_5},
$$
where $e_0,\dots, e_4$ is the chosen basis for $V$ and $f^{\pm}
:V_1\to\Lambda^2 V_0$ are the Horrocks-Mumford maps \cite{HM} defined by
$$
f^{+}(\sum x_ie_i)=\sum x_i e_{i+2}\wedge e_{i+3} \quad and \quad
f^{-}(\sum x_ie_i)=\sum x_i e_{i+1}\wedge e_{i+4}.
$$
Notice that $f^\pm$ correspond to the pencil of trivial characters in
$\Hom_{G_5}(V_1, \Lambda^2 V_0^\sharp)\cong 2I\oplus 2Z.$
(ii)\quad The Pl\"ucker equations for $\mu f^+(x)+
\lambda f^-(x)\in\Lambda^2 V_0$ to be a decomposable tensor, i.e., to
correspond to a line $L_x\subset\Pn 4$,
are just the equations of the elliptic
curve $E_{(\lambda :\mu)}$.\par\noindent
(iii)\quad The lines $L_x$, $x\in E_{(\lambda :\mu)}$, are the rulings
of the scroll $Q_{(\lambda :\mu)}$.
$L_{o_{E_{(\lambda :\mu)}}}$ is exactly $l_{(\lambda :\mu)}$,
the ruling over the origin of the base curve. Any scroll $Q_{(\lambda
:\mu)}$ may be identified with the symmetric product of an elliptic curve
(cf. \cite{BHM2} p.744).
Therefore, in our case,
$Q_{(\lambda:\mu)}$ may be canonically identified with
$S^2(E_{(\lambda :\mu)}).$\qed
\endproclaim
\medskip Next we will discuss some facts concerning $S_{15}$ and
other projective models of the Shioda modular surface $S(5)$.
First we recall from \cite{BH} a  natural choice for a
$\openQ$-basis of the N\'eron-Severi group of $S(5)$.
As mentioned in Proposition 4.2, (v), $S(5)$ has 25 sections
$L_{ij}$, $(i,j)\in\openZ_5\times\openZ_5$,
which consist of the $5$-torsion points in each smooth
fibre. Furthermore  $S(5)$ has 12 singular fibres of type $I_5$.
We will order the components $F_{kl}$, $k=\overline{1,12}$ and
$l\in\openZ_5$ of the singular fibres $F_k$ in a coherent way, such that
$L_{00}\cdot F_{k0}=1$ and $F_{kl}\cdot F_{k,l+1}=1$. If $F$ denotes
the class of a smooth fibre,
then $L_{00}$, $F$ and $F_{kl}$, $k=\overline{1,12}$, $l=\overline{1,4}$,
are a basis for $\NS(S(5))\otimes\openQ$ (cf. \cite{BH, Prop.1}).
A result of Inoue and Livn\'e (cf. \cite{BH, Theorem 1}) says, that there
exists a unique divisor class $I\in\Pic(S(5))$ with $5I\equiv\sum_{i,j}L_{ij}$.
In terms of the above basis it is given by
$$
I=5L_{00}+24F-\sum_{k=1}^{12}(2F_{k1}+3F_{k2}+3F_{k3}+2F_{k4}).
$$
The complete linear system $|I+2F|$ defines an immersion $\varphi_1
:S(5)\to S_{15}\subset\Pn 4$ \cite{BHM1, Proposition 9},
which is nothing else than the normalization
morphism described in Proposition 4.2, (iv).
In the sequel we will be mainly concerned
with the linear system $|3I+5F-\sum P_{kl}|$, where $P_{kl}$
are the vertices of the singular fibres of $S(5)$,
$k=\overline{1,12}$ and $l\in\openZ_5$. In \cite{BH, Theorem 4}
the following result is proved:
\proclaim{Proposition 4.8} Let $\tilde S(5)$ denote the blow-up of $S(5)$
in the 60 points $P_{kl}$. The linear system $|3I+5F-\sum P_{kl}|$
defines a regular map $\varphi_4:{\tilde S(5)}\to\Pn 4$, which is
birational onto a surface $S_{45}$ of degree 45. The singular fibres of
${\tilde S(5)}$ are mapped to complete pentagons. Two such pentagons
coincide if the parameters of the fibres belong to opposite vertices
of the icosahedron, and are disjoint otherwise.\qed
\endproclaim
\noindent The proof in  \cite{BH} relies on some properties of the maps
$v^+$ and $v^-$ introduced in \cite{BHM1}:
$$
v^+: V_1 \to V_0,\qquad v^+(\sum_{i\in\openZ_5} x_ie_i)=
\sum_{i\in\openZ_5}(x_{i+2}x_{i+4}^2-x_{i+1}^2x_{i+3})e_{i},
$$
\noindent
and
$$
v^-:V_1 \to V_0,\qquad v^-(\sum_{i\in\openZ_5} x_ie_i)=
\sum_{i\in\openZ_5}(x_{i+1}x_{i+2}^2-x_{i+3}^2x_{i+4})e_{i}.
$$
$v^\pm$ correspond to the pencil of trivial characters in
$\Hom_{G_5}(V_1, V_0 \otimes V_0^\sharp)\cong 2I\oplus 3S\oplus 5Z.$
We recall from \cite{BHM1}:
\proclaim{Remark 4.9} Let $(\lambda :\mu) \in \openP^1_-$.
\par\noindent
(i)\quad $v^+(x)=v^-(x)=0$ exactly at the singular points of $S_{15}$.
\par\noindent
(ii)\quad If $x\in\Pn 4$ is none of these singular points then
$\mu v^+(x)-\lambda v^-(x)=0$ iff $\mu f^+(x)+\lambda f^-(x)$ is a
decomposable tensor in $\Lambda^2 V_0$. When this occurs the point
$\openP(v^{\pm}(x))$ lies on the corresponding line $L_x$.\par\noindent
(iii)\quad The five cubic equations arising from $\mu v^+(x)-\lambda v^-(x)=0$
cut out the union of $E_{(\lambda :\mu)}$ and the singular points of $S_{15}$
(in case $E_{(\lambda :\mu)}$ is a pentagon its vertices show up as embedded
points).
In particular, $v^+$ and $v^-$ are proportional exactly over $S_{15}$.\qed
\endproclaim
\noindent
In the sequel we will consider $v^+$ as the $G_5$-equivariant rational map
$$v^+:\Pn 4(x)=\Pn 4 (V_1)\rto\Pn 4(z)=\Pn 4 (V_0)$$ defined in coordinates by
$z_{i}\to x_{i+2}x_{i+4}^2-x_{i+1}^2x_{i+3}$. Then the base locus
of $v^+$ is the disjoint union of the pentagon
$E_{(0:1)}(x)=\{x_{i+1}x_{i+4}=0\}_{i\in\openZ_5}$ (scheme theoretically
the vertices show up as embedded points) and the 25 vertices of the other
5 complete pentagons. In particular, $v^+$ is dominant and a standard
computation shows that it is generically $16:1$.
Moreover, since $v^+$ and $v^-$ are given by cubic polynomials, and since
every linear combination of them vanishes on a fibre $E_{(\lambda:\mu)}$
of $S_{15}$ and on the vertices $P_{kl}$ we have a commutative diagram
\diagram[tight,width=4em,height=3em]
S(5) & \rTo^{\varphi_1} & S_{15}\subset\Pn 4\\
   & \rdDotsto[fixed]_{\varphi_4}&\dDashto>{v^+}\\
   && S_{45}\subset\Pn 4.
\enddiagram
\proclaim{Proposition 4.10}
(i)\quad $\varphi_4$ maps $L_{00}$ birationally onto the
conic $C_+=\{z_0^2+4z_1z_2=z_1-z_4=z_2-z_3=0\}\subset\Pn 2_+$.
\par\noindent
(ii)\quad $v^+$ maps $\Pn 2 _+$ onto $\Pn 1_-$.
The $3$-section $\scr D$ of $S(5)$ which intersects
each smooth fibre in its non-trivial $2$-torsion points
is mapped $3:1$ by $\varphi_4$ onto $\Pn 1_-$.\par\noindent
(iii)\quad $v^+$ embeds a curve $E_{(\lambda:\mu)}$,
$(\lambda:\mu)\in\openP^1_-\setminus\Lambda$, as a  $G_5$-invariant
elliptic curve $E^{\prime}_{(\lambda:\mu)}\subset\Pn 4$
of degree 15. The curve $E^{\prime}_{(\lambda:\mu)}$ is
the section of the scroll $Q_{(\lambda:\mu)}$ which corresponds
to the set $\{\{p,q\}\mid 3p+2q=0\}$ under the identification
$Q_{(\lambda:\mu)}=S^2(E_{(\lambda:\mu)})$.
\endproclaim
\demo{Proof}(i)\quad is clear since ${v^+}_{\mid_{\Pn 1_-}}$ is the map
$$
{v^+}_{\mid_{\Pn 1_-}}: (x_1:x_2)\to(2x_1^2x_2:x_1x_2^2:-x_1^3
:-x_1^3:x_1x_2^2).
$$
\smallskip\noindent
(ii)\quad follows readily since ${v^+}_{\mid_{\Pn 2_+}}$ is the map
$$
{v^+}_{\mid_{\Pn 2_+}}: (x_0:x_1:x_2)\to
(0:x_1x_2^2-x_0^2x_2:x_1^3-x_0x_2^2:x_0x_2^2-x_1^3:x_0^2x_2-x_1x_2^2),
$$
and since $v^+$ restricted to the Brings curve $B\subset
\Pn 2_+$ has as base locus a scheme of length 15.
Indeed, it follows from the above discussion that set-theoretically this
base locus coincides with
$$
\Pn 2_+\cap\{ v^+=v^-=0\}\cup \{(0:0:1:1:0)\},
$$
i.e., that it consists of the double points $e_0=(1:0:0)$,
$(1:\xi^i:\xi^{-i})$,
$i\in \openZ_5,$ of $B$ and the point $(0:0:1)$, and that these points
count with different multiplicities as follows. The node $e_0$ has
multiplicity 4 since there are exactly two pentagons, namely $E_{(0:1)}$
and $E_{(1:0)}$ having $e_0$ as a vertex where two non-trivial 2-torsion
points coalesce, while all the other five nodes have multiplicity 2 since
for each of them there are precisely two pentagons meeting $\Pn 2_+$
transversely in such a point (compare with Lemma 9.1, (ii) for more details).
Similarly, $(0:0:1)$ counts with multiplicity 1 in the base locus.
\smallskip\noindent
(iii)\quad It follows from Remark 4.7 and Remark 4.9 that the image
$E^{\prime}_{(\lambda:\mu)}$ of $E_{(\lambda:\mu)}$ under ${v^+}$ is
a section of $Q_{(\lambda:\mu)}$. This section meets the ruling over
the origin of this scroll in the point $v^+((\lambda:\mu))$, i.e., in its
point of tangency to the conic $C_+$. On the other hand,
$E^{\prime}_{(\lambda:\mu)}$ has degree either 15 or 5 by $G_5$-invariance.
Moreover, it is of numerical class $C_0+bf$ with $b=12$ or $b=2$ since it meets
the ruling in one point. Here $C_0$ is the minimal section and $f$ a ruling.
There is no $G_5$-invariant curve in the numerical class
of $C_0+2f$, thus $E^{\prime}_{(\lambda:\mu)}$ has degree 15 and is of
numerical class $C_0+12f$. A $G_5$-invariant curve in this class is of type
${2p+3q=\tau_i}$ on $S^2(E_{(\lambda:\mu)})$, $\tau_i$ a 2-torsion point on
$E_{(\lambda:\mu)}$. By looking again at the intersection point with the
ruling over the origin we see that necessarily $\tau_i=0$.\qed
\enddemo
\noindent
For later reference we recall from \cite{ADHPR} the following:
\proclaim{Remark 4.11} For general $(\lambda:\mu)\in\Pn 1_-$
the curve $E^{\prime}_{(\lambda:\mu)}$ coincides with the base locus
of the pencil $|A_K|$ of degree 10 abelian surfaces on
$\Trisec Q_{(\lambda:\mu)}$.
\endproclaim
\proclaim{Proposition 4.12}
(i)\quad If $p$ is a point on a smooth fibre $E_{(\lambda:\mu)}(x)$ of
$S_{15}$ which is not a 5-torsion point, i.e., which does not lie on
a HM-line, then the quintic hypersurface $\{\det M_p(y)=0\}$
is the secant variety  $\Sec E_{(\lambda:\mu)}$
and $\{\det L_p(y)=0\}$ is the trisecant variety
$\Trisec Q_{(\mu:-\lambda)}$.\par\noindent
(ii)\quad  If $q$ is a point on a smooth fibre $E^{\prime}_{(\lambda:\mu)}(x)$
of $S_{45}$ which does not lie on a HM-line,
then the quintic hypersurface $\{\det M_q(y)=0\}$
is the trisecant variety $\Trisec Q_{(\lambda:\mu)}$
and $\{\det L_q(y)=0\}$ is the secant variety
$\Sec E_{(\lambda:\mu)}$.\qed
\endproclaim
\demo{Proof} In view of the correspondence described in
\cite{ADHPR, Proposition 12} between secant varieties
of elliptic normal curves and trisecant varieties of elliptic
quintic scrolls, and of the identity
$$L(\iota(x),u)^t\cdot
\pmatrix
y_0\\
y_1\\
y_2\\
y_3\\
y_4\\
\endpmatrix=
M(x,y)\cdot
\pmatrix
u_0\\
u_1\\
u_2\\
u_3\\
u_4\\
\endpmatrix
$$
it is enough to prove only the statements made for matrices of type $M$.
\smallskip\noindent
(i)\quad There are 5 independent $H_5$-invariant
quintic hypersurfaces through a $G_5$-invariant elliptic normal
curve (cf. \cite{Moo, 2.4}), hence $\Theta$ contracts
$E_{(\lambda:\mu)}(x)$ to a point $P$. It is shown in \cite{Moo, \S 3}
that $\Theta^{-1}(T_{\Omega,P})=\Sec E_{(\lambda:\mu)}$. This
quintic in turn coincides with $\{\det M_p(y)=0\}$ by Remark 4.1, (i).
\smallskip\noindent
(ii)\quad Analogously, one shows that there are precisely 5 $H_5$-invariant
quintic hypersurfaces through an elliptic curve of type
$E^{\prime}_{(\lambda:\mu)}(x)$. Namely, $E^{\prime}_{(\lambda:\mu)}(x)$ lies
on the scroll $Q_{(\lambda:\mu)}$ by Proposition 4.10, (iii), and since there
are 3 $H_5$-invariant quintic hypersurfaces through an elliptic quintic scroll
\cite{Hu1, VIII.3.5}, since $|-2K_{Q_{(\lambda:\mu)}}|$
moves in a pencil \cite{BHM2, Proposition 5.4}, and since
$h^1({\scr I}_{Q_{(\lambda:\mu)}}(5))=0$ the claim follows by taking
the $H_5$-invariant part of the cohomology of the exact sequence
$$
o\to {\scr I}_{Q_{(\lambda:\mu)}}(5)
\rightarrow {\scr I}_{E^{\prime}_{(\lambda:\mu)}}(5)
\rightarrow{\scr O}_{Q_{(\lambda:\mu)}}(-2K_{Q_{(\lambda:\mu)}})
\to 0.
$$
Therefore to show (ii) it is  enough to check only one point $q$ on
$E^{\prime}_{(\lambda:\mu)}(x)\setminus L$. This is done in \cite{ADHPR,
Corollary 14}
for $q$ chosen as the intersection point of $E^{\prime}_{(\lambda:\mu)}(x)$
and $C_+$.\qed
\enddemo
\proclaim{Proposition 4.13} Let $\Gamma=\{x_0^2+(x_1+x_4)(x_2+x_3)=0\}$
be the quadric cone over the conic $C_+\subset\openP^2_+$ with
vertex $\openP^1_-$. Then $\Gamma$ meets each smooth fibre
$E_{(\lambda:\mu)}$ of the surface $S_{15}$ in twice its
origin and simply in the non-trivial 3-torsion points.
\endproclaim
\demo{Proof}
First note that $E_{(\lambda:\mu)}$ is not contained in
$\Gamma$. Since $\Gamma$ is singular along $\openP^1_-$ it
will be sufficient to show that the non-trivial
3-torsion points of $E_{(\lambda:\mu)}$ lie in $\Gamma$.
\smallskip\noindent
For this we use the 4:4 correspondence between
$G_5$-invariant elliptic normal curves and
trisecant varieties of $G_5$-invariant elliptic quintic scrolls.
Recall from \cite{ADHPR} that a trisecant scroll $\Trisec Q$
contains a pencil $|A_K|$ of abelian surfaces. Among these
there are precisely 4 translation scrolls coming from
the 4 triangles from the Hesse pencil with
respect to the twisted $H_3$-action in each plane of $\Trisec Q$.
Clearly, as translation scrolls they are given by a 3-torsion point on
the corresponding elliptic curve. Notice that these 4
elliptic  curves do not lie on the scroll $Q$.\smallskip\noindent
Conversely, given an elliptic curve, a 3-torsion point on it
defines a translation scroll of the above type. The eight non-trivial
3-torsion points on the curve give rise to only 4 translation scrolls,
hence the claimed correspondence.\smallskip\noindent
Now both the elliptic scroll $Q$ and the desingularization of the
trisecant scroll  $\Trisec Q$, are projective bundles over an elliptic
curve $E$, and the projections are $G_5$-equivariant.  We denote
the plane in $\Trisec Q$ over the origin $o_E$ of $E$
by $\pi_{o_E}$. It contains a plane cubic curve of the scroll $Q$
which we denote by $C_{o_E}$. Note that $\pi_{o_E}$ is also the
plane spanned by the origin of $E_{(\lambda:\mu)}$ and a pair
of points $(\rho,-\rho)$, where $\rho$ is a non-zero 3-torsion point of
$E_{(\lambda:\mu)}$ and that we obtain any such plane
in this way. \par\noindent
The proposition now follows from the

\proclaim{Claim} The plane $\pi_{o_E}$ is one of the planes
of the cone $\Gamma$.\endproclaim
\demo{Proof of claim} The plane cubic curve $C_{o_E}$ meets each of
the three elliptic normal curves on $Q$ in its origin,
i.e., in a point on $\Pn 1_-$. Thus $\Pn 1_-$ is contained in
$\pi_{o_E}$. On the other side, $C_{o_E}$ meets the ruling $l_{o_E}$
of $Q$ over $o_E$ in its point of tangency
to $C_+$, hence $\pi_{o_E}\cap\Pn 2_+\subset C_+$ and the claim
follows.\qed
\enddemo
\enddemo
\medskip
Later on we will consider $H_5$-invariant graded $R$-modules whose
presentation matrix is obtained by concatenating three Moore matrices
$M_{p_1}(y)$, $M_{p_2}(y)$ and $M_{p_3}(y)$. It will follow from the
representation theory described at the beginning of this section that
the first order linear syzygies of such modules come in blocks of
$L$-matrices. Therefore we are interested in the product matrices
$$
P_{x,z} (y) = M(x,y)\,L(z,y).
$$
Each entry of $P_{x,z} (y)$ is a quadric in the span of
$$
y^2_i,\ y_{i+1} y_{i+4},\ y_{i+2} y_{i+3}
$$
for some $i\in\openZ_5$ with coefficients that are bilinear forms in
$x$'s and $z$'s. We collect these forms in the matrix
$$
\pmatrix
x_0 z_0 & x_1 z_1 & x_2 z_2 & x_3 z_3 & x_4 z_4\\
x_1 z_3 + x_4 z_2 & x_2 z_4 +x_0 z_3 & x_3 z_0 + x_1 z_4
& x_4 z_1 + x_2 z_0 & x_0 z_2 + x_3 z_1 \\
x_2 z_1 + x_3 z_4 & x_3 z_2 + x_4 z_0 & x_4 z_3 + x_0 z_1
&x_0 z_4 + x_1 z_2 & x_1 z_0 + x_2 z_3
\endpmatrix ,
$$
which, after column permutations, can be rewritten in a compact way as
$$
T(x,z)=(t_{ji})=\pmatrix
x_{3i}z_{3i}\\
x_{3i+1}z_{3i+3}+x_{3i+4}z_{3i+2}\\
x_{3i+2}z_{3i+1}+x_{3i+3}z_{3i+4}\\
\endpmatrix_{i\in\openZ_5}.
$$
The entries of $T$ define a morphism
$$
\Phi_T : \openP^4(x) \times \openP^4(z) \rightarrow \openP^{14}(t).
$$
Let us fix a point $(\lambda:\mu)\in\Pn 1_-\setminus\Lambda$. The variety
$$
I_{(\lambda:\mu)} = \left\{ (x,z)\mid P_{x,z}(y) = 0
\text{\ for\ any\ } y\in
E_{(\lambda:\mu)} (y)\right\}\subset \openP^4(x) \times \openP^4(z)
$$
parametrizes the product matrices whose entries are in the linear
system of quadrics defining the smooth elliptic curve
$E_{(\lambda:\mu)}(y)$. $I_{(\lambda:\mu)}$ will be used in Section 5
to compute the first order linear syzygies of a module as above in
the case that all 3 points $p_1$, $p_2$ and $p_3$ are on the elliptic curve
$E_{(\lambda:\mu)}(x)$. For this purpose we will need the following results.
\proclaim
{Lemma 4.14} $I_{(\lambda:\mu)}$ is cut out by the 10 bilinear forms
$$
\align
\lambda x_{i} z_{i} + \mu \big(x_{i+2}z_{i+1}+x_{i+3}z_{i+4}\big) =
\mu x_{i} z_{i} - \lambda\big(x_{i+1}z_{i+3}+x_{i+4}z_{i+2}\big) = 0,
\quad& i\in\openZ_5\ .
\endalign
$$
It is an elliptic scroll, which maps by projection to the elliptic
curve $E_{(\lambda:\mu)} (x)$ in the first factor (this map defining
the ruling) and isomorphically to $Q_{(\lambda:\mu)}(z)$ in the
second factor.
\endproclaim
\demo
{Proof} By writing down $P_{x,z}(y)$ explicitly and by comparing with
the equations $q^{(\lambda:\mu)}_i(y)$ of $E_{(\lambda:\mu)}(y)$
given in Corollary 4.5 we see that
$I_{(\lambda:\mu)}$ is defined by the condition
$$
\pmatrix
\lambda& 0& \mu\\
\mu& -\lambda& 0
\endpmatrix
\cdot T (x,z) = 0\ ,
$$
i.e., by the 10 equations above. A straightforward calculation shows
that the saturation of these forms with respect to $z_0,\dots,z_4$\
$(x_0,\dots,x_4)$ contains the equations of $E_{(\lambda:\mu)}(x)$\
\,$\big(Q_{(\lambda:\mu)}(z)\big)$.\par\noindent
Conversely, by plugging in our fixed point $p = (\lambda:\mu)\in\openP^1_-$,
we see that $P_{p,z}(p) = 0$ iff $z$ is on the line $l_p$, the ruling of
$Q_{(\lambda:\mu)}(z)$ over the origin. By $H_5$-invariance
the image of $I_{(\lambda:\mu)}$ under the projection to $x$ (to $z$)
contains $25$ points (lines) of $E_{(\lambda:\mu)}(x)$\
\,$\big(Q_{(\lambda:\mu)}(z)\big)$. By Bezout this image contains the
whole elliptic curve (scroll).\qed
\enddemo
\noindent
We have shown that $I_{(\lambda:\mu)}$ is contained in $E_{(\lambda:\mu)}(x)
\times Q_{(\lambda:\mu)}(z)$ and that it is a smooth elliptic scroll. We now
restrict the morphism $\Phi_T$ to this scroll. Let $w_i=x_{3i}z_{3i}$,
$i=\overline{0,4}$, be the entries of the first row of $T$. Then $\Phi_T$ and
the rational map defined by the $w_i$ into $\openP^4(w)=\openP^4(V_2)$ coincide
on
$I_{(\lambda:\mu)}$, if $\openP^4(w)$  is suitably embedded in
$\openP^{14}(t)$. In fact, under the identification of $\Pn {14}(t)$ with the
space of quadrics in $y$'s, we embed $\Pn 4(w)$ as the space of quadrics in the
ideal of $E_{(\lambda:\mu)}(y)$. The image of $I_{(\lambda:\mu)}$ in
$\openP^4(w)$ proves to be very singular along the curve
$E_{(\lambda:\mu)}(w)$. More precisely, we have:
\proclaim
{Lemma 4.15}
Let $S_{(\lambda:\mu)}$ be the image of $I_{(\lambda:\mu)}$ under $\Phi_{T}$.
Then $S_{(\lambda:\mu)}$ is contained in the linear subspace
$\openP^4_{(\lambda:\mu)}(w)\subset \openP^{14}(t)$ which is cut out by the
linear forms
$$
\align
\lambda t_{0i} + \mu t_{2i} = \mu t_{0i} - \lambda t_{1i}  = 0,
\quad& i=\overline{0,4}\ .
\endalign
$$
In fact, $S_{(\lambda:\mu)}$ is the degree 15 scroll in
$\openP^4_{(\lambda:\mu)}(w)$ defined by the singular
lines of the rank 3 quadrics through $E_{(\lambda:\mu)}(w)$.
\endproclaim
\demo
{Proof}
Consider the rational map
$$
\Phi_{T_0} : \openP^4(x) \times \openP^4(z) \dashrightarrow
\openP^4_{(\lambda:\mu)}(w) \subset \openP^{14}(t),
$$
defined by the entries of the first row of $T$. Clearly, the
restrictions of $\Phi _T$ and $\Phi _{T_0}$ to $I_{(\lambda:\mu)}$
coincide. The first assertion of the lemma follows.  For the second
one we first remark that $\Phi_{T_0}$ maps the ruling
$\{p\}\times l_p(z)$, where $p=(\lambda:\mu)$, $1:1$
to $\Pn 1_-(w)$. Thus $S_{(\lambda:\mu)}$ is
two-dimensional. Recall now from \cite {Hul, VIII.3.5} that
$H^0(\scr J^2_{E_{(\lambda:\mu)}(w)}(5))^{H_5}$ is a net whose base
locus is exactly the scroll of the singular lines
of the rank 3 quadrics through $E_{(\lambda:\mu)}(w)$.
So it suffices to show that $S_{(\lambda:\mu)}$
is contained in every $H_5$-invariant quintic which is singular along
$E_{(\lambda:\mu)}(w)$. By \cite{Hu1, p.82} a basis for these quintics
is given by
$$
\align
Q_0(w) &= \sum\limits_{i=0}^4 w_i q_i^2\ ,\\
Q_1(w) &= \sum\limits_{i=0}^4 w_i q_{i+1} q_{i+4}\ ,\\
Q_2(w) &= \sum\limits_{i=0}^4 w_i q_{i+2} q_{i+3},
\endalign
$$
where $q_i = q_i^{(\lambda:\mu)}(w)$ are the equations
of $E_{(\lambda:\mu)}(w)$. A
straightforward calculation shows that for all $i$
$$
Q_i(x_0z_0, x_3z_3, x_1z_1, x_4z_4, x_2z_2) = 0
$$
modulo the equations of $I_{(\lambda:\mu)}$\ .\qed
\enddemo
\noindent
We recall the geometry of the degree 15 scroll
defined by the singular lines of the rank 3 quadrics through
an elliptic normal curve in $\Pn 4$.
\proclaim{ Proposition 4.16}
(i)\quad The scroll $S_{(\lambda:\mu)}$ is a birational image
of the symmetric product $S^2(E_{(\lambda:\mu)})$ via the
identification map which associates to a proper secant line $\overline{pq}$
of $E_{(\lambda:\mu)}$ the point of intersection of this secant with
the singular line of the rank 3 hyperquadric uniquely determined
by the fact that it contains the curve and that its singular locus
meets $\overline{pq}$ outside the curve.\par\noindent
(ii)\quad $S_{(\lambda:\mu)}$ is 4-tuple along
$E_{(\lambda:\mu)}$ and smooth elsewhere.\par\noindent
(iii)\quad The union $X_{p+q}$ of all the secants of
$E_{(\lambda:\mu)}$ spanned by points which add up to
$p+q$ in the group law is a smooth rational cubic scroll,
whose directrix $d_{p+q}$ is a ruling in $S_{(\lambda:\mu)}$.
As a ruling in the symmetric product $S^2(E_{(\lambda:\mu)})$
this is the ruling over the point $p+q$ on the base curve.
This ruling hits the curve $E_{(\lambda:\mu)}\subset\Pn 4$ in exactly
one point corresponding to $-2(p+q)$.\par\noindent
(iv)\quad The curve $\{\{p,q\} \mid 3p+2q=0\}$
on the symmetric product $S^2(E_{(\lambda:\mu)})$ is mapped 4:1
onto the curve $E_{(\lambda:\mu)}\subset S_{(\lambda:\mu)}$.
The image of $\{p,q\}$ is the intersection of the secant
$\overline{pq}$ with $d_{p+q}$.
\endproclaim
\demo{Proof} See \cite{Hu1}, \cite{BHM1}, \cite{BHM2}.\qed\enddemo
\proclaim{Remark 4.17}
(i)\quad
The scroll $S_{(\lambda:\mu)}$ contains three elliptic curves , which on
the symmetric product $S^2(E_{(\lambda:\mu)})$ are numerically equivalent
to the anticanonical divisor. Namely, these can be identified with the set
of unordered pairs of points $\{\{p,q\} \mid p=q+\tau_i\}$ on
$E_{(\lambda:\mu)}$,
$\tau_i$ a non-trivial $2$-torsion point on $E_{(\lambda:\mu)}$.  The group law
is given by $\{p,q\}+\{r,s\}=\{p+r,q+r\}=\{p+s,q+s\}$, which is well-defined
since $p-q=r-s=\tau_i$ is 2-torsion.

\par\noindent
(ii)\quad
Through any point of $E_{(\lambda:\mu)}$ there are exactly four
rulings of $S_{(\lambda:\mu)}$. In particular, the rulings through
the origin of $E_{(\lambda:\mu)}$ correspond to the directrices $d_{\tau_i}$,
$\tau_i$ a 2-torsion point on $E_{(\lambda:\mu)}$.\qed
\endproclaim
\proclaim{Lemma 4.18}
Let $\tau_i$, $i=\overline{1,3}$ be the non-trivial 2-torsion points
of $E_{(\lambda:\mu)}$ and let $\Pi_{\tau_1\tau_2}$ be the plane
spanned in $\Pn 4$ by the directrices $d_{\tau_1}$ and $d_{\tau_2}$.
Then $\Pi_{\tau_1\tau_2}$ meets $S_{(\lambda:\mu)}$ in $d_{\tau_1}$,
$d_{\tau_2}$ and the points corresponding on the symmetric
product to $\{\{p,q\}\mid 3(p+q)=p-q=\tau_3\}$. In particular:
\par\noindent
(i)\quad There are exactly 8 rulings not contained in
$\Pi_{\tau_1 \tau_2}$ which meet this plane. These are the rulings
over the points $x\in E_{(\lambda:\mu)}$ with
$3x=\tau_1+\tau_2$, $x\ne\tau_3$.\par\noindent
(ii)\quad These rulings meet the plane in non-trivial 6-torsion
points of the numerically anticanonical curve on $S_{(\lambda:\mu)}$
corresponding to $\tau_3$ as in Remark 4.17, (i).
\endproclaim
\demo{Proof} Assume that the plane $\Pi_{\tau_1 \tau_2}$ and the ruling
over a point $x\in E_{(\lambda:\mu)}$ intersect. Let $H$ be the hyperplane
spanned by $\Pi_{\tau_1 \tau_2}$ and this ruling.
The curve $E_{(\lambda:\mu)}$ hits this ruling
in the point $-2x$ and the other two rulings in the origin.
Therefore we may write
$$
 H\cap E_{(\lambda:\mu)}=o_{E_{(\lambda:\mu)}}+(-2x)+a+b+c.
$$
The chosen 3 rulings are directrices of rational cubic scrolls which
contain the elliptic curve, and $H$ intersects each of these scrolls
along their directrix and two rulings.  From Proposition 4.16, (iii)
we get the following relations:\par\noindent
On the scroll with directrix over $\tau_1$
$$a+b=c+(-2x)=\tau_1.$$
On the scroll with directrix over $\tau_2$
$$a+c=b+(-2x)=\tau_2.$$
Note that we have essentially defined the points $a,b,c$ via
these relations, hence for the scroll over $x$ we necessarily get
$$b+c=a+o_{E_{(\lambda:\mu)}}=x.$$
Adding the second relations in each case we get
$a+b+c+(-4x)=x+\tau_1+\tau_2$, and since $a+b+c+(-2x)=0$
we deduce $3x=\tau_3$.
\par\smallskip\noindent
Let us consider now the hyperplane $H'$ spanned by
$\Pi_{\tau_1\tau_2}$ and the ruling in the rational scroll
over the directrix of $x$. This ruling hits the curve in points $p$ and $q$
say, and thus we can write $H'\cap E_{(\lambda:\mu)}=
o_{E_{(\lambda:\mu)}}+p+q+a+b$. As above we get
that $a+p=b+q=\tau_1$ and $a+q=b+p=\tau_2$, which
implies $ p-q=\tau_3$. The rest of the lemma follows now from
Remark 4.17, (i).\qed
\enddemo
\medskip
We will finally determine the 6-torsion points of a $G_5$-invariant
elliptic normal curve in $\Pn 4$. For this purpose we will again use
the rational map  $v^+:\Pn 4(x) \rto\Pn 4(z)$, and the morphism
$\Phi_T$. Notice that the map
$$\Pn 4(x)=\Pn 4 (V_1)\rto\Pn 4(w)=\Pn 4 (V_2)$$ defined in coordinates by
$w_{2i}\to x_i(x_{i+2}x_{i+4}^2-x_{i+1}^2x_{i+3})$ is $G_5$-equivariant.
It corresponds to a trivial character in
$\Hom_{G_5}(V_1, V_0 \otimes V_3^\sharp)\cong 2I\oplus 3S\oplus 5Z.$
\proclaim{Proposition 4.19} Let $(\lambda:\mu)\in\Pn 1_-\setminus\Lambda$.
\par\noindent
(i)\quad The graph $\Gamma_{(\lambda:\mu)}\subset\Pn 4(x)\times\Pn 4(z)$
of the restriction of $v^+$ to the elliptic normal curve
$E_{(\lambda:\mu)}$ is contained in
the incidence scroll  $I_{(\lambda:\mu)}$ .\par\noindent
(ii)\quad $\Phi_T$ maps the graph $\Gamma_{(\lambda:\mu)}$ $4:1$ onto the
curve $E_{(\lambda:\mu)}(w)\subset\openP^4_{(\lambda:\mu)}(w)$.
Moreover, the composition
$$\Phi_{T}\circ(\id_{E_{(\lambda:\mu)}(x)}\times {v^+}_{\mid
_{E_{(\lambda:\mu)}(x)}}):E_{(\lambda:\mu)}(x)\to E_{(\lambda:\mu)}(w)$$
acts as multiplication by 2 in the group law of the elliptic curve
$E_{(\lambda:\mu)}$.\par\noindent
iii)\quad The octic hypersurface $\{w_0^2+(w_1+w_4)(w_2+w_3)=0\}$,
where $w_{2i}=x_{i}(x_{i+2}x_{i+4}^2-x_{i+1}^2x_{i+3})$,
$i\in\openZ_5$, meets each smooth fiber of $S_{15}\subset\Pn 4(x)$ in twice
the 2-torsion points and simply in the remaining 6-torsion points .
\endproclaim
\demo{Proof}
(i)\quad One can easily check that
$$
\pmatrix
\lambda& 0& \mu\\
\mu& -\lambda& 0
\endpmatrix
\cdot T (x,v^+(x)) = 0\
$$
modulo the equations $q_i^{(\lambda:\mu)}$, $i\in\openZ_5$, of
$E_{(\lambda:\mu)}(x)$. Hence we deduce from the proof of Lemma 4.14
that the graph $\Gamma_{(\lambda:\mu)}$
is contained in the incidence scroll $I_{(\lambda:\mu)}$.\smallskip\noindent
(ii)\quad By Lemma 4.15 the image of $I_{(\lambda:\mu)}$ under
$\Phi_{T}$ is the singular scroll $S_{(\lambda:\mu)}$. Therefore
all claims in (ii) follow from the fact that through any point
$p$ of $E_{(\lambda:\mu)}(w)$ there pass exactly four rulings of
$S_{(\lambda:\mu)}(w)$. Namely, these are the directrices $d_q$ of the
cubic scrolls $X_q$ with $2q=-p$ in the group law (cf. Proposition 4.16, (iii)
and (iv), and Remark 4.17, (ii)).\smallskip\noindent
(iii) is then clear from Proposition 4.13.\qed
\enddemo
For the next section we summarize some of the results of 4.14-4.19.  First
these results give rise to the diagrams

$$
\diagram[tight,width=5em,height=4em]
{I_{(\lambda:\mu)}} &\rTo^{\Phi_T} &{S_{(\lambda:\mu)}(w)}\\
\dTo^{\pi_1}&\rdTo^{\pi_2}&\uTo_{\rho}\\
{E_{(\lambda:\mu)}(x)} &\lTo^{\pi}&{Q_{(\lambda:\mu)}(z)}\\
\enddiagram
 $$
\noindent
and

$$
\diagram[tight,width=5em,height=4em]
{\Gamma_{(\lambda:\mu)}} &\rTo^{\Phi_T} &{E_{(\lambda:\mu)}(w)}\\
\uTo^{\psi }&\rdTo^{\pi_2}&\uTo_{\rho}\\
{E_{(\lambda:\mu)}(x)} &\rTo^{v^+}&{E^{\prime}_{(\lambda:\mu)}(z)}.\\
\enddiagram
 $$
The varieties of the second diagram are subvarieties of the corresponding
varieties of the first diagram.   In the first diagram $\rho =\Phi_T\circ
\pi^{-1}_2$, and in the second diagram $$\psi
=(\id_{E_{(\lambda:\mu)}(x)}\times {v^+}_{\mid
_{E_{(\lambda:\mu)}(x)}}). $$  The maps $\Phi_T$, $\pi_2$ and $\rho$ of the
second diagram are restrictions of those in the first diagram, while $\psi$
and $v^+$ are sections of the corresponding maps in the first diagram.
4.14, 4.15, 4.16 and 4.19 say that both diagrams commute.  In the next
section we will need to know which pairs of points $(p,q)\in
I_{(\lambda:\mu)}$, $p\in E_{(\lambda : \mu )}$, $q\in Q_{(\lambda : \mu )}$
are mapped by $\Phi_T$ to the points of $E_ {(\lambda : \mu )}(w)$ and to
the plane $\Pi_{\tau_1\tau_2}$ of Lemma 4.18.  Since $v^+$ is an
isomorphism we can identify points on $E^{\prime}_{(\lambda : \mu )}(z)$
with points on $E_{(\lambda : \mu )}(x)$.
\proclaim {Remark 4.20} Under the identification of $Q_{(\lambda :\mu)}$
with $S^2(E_{(\lambda:\mu)})$ each numerically
anticanonical curve  on $S^2(E_{(\lambda:\mu)})$ corresponds to one of the
three $G_5$-invariant elliptic quintic curves on $Q_{(\lambda :\mu)}$. The
curve corresponding to $\{\{p,q\}\mid p=q+\tau_i\}$ will be denoted by
$E_{(\lambda :\mu)}({\tau_i})$.\
\qed\endproclaim
\proclaim {Corollary 4.21} (i)\quad The map $\rho$ of the second diagram is
$4:1$ and acts as multiplication by 2 in the group law of the elliptic
curve $E^{\prime}_{(\lambda:\mu)}(z)$.\par\noindent
  (ii)\quad The HM-lines in $\openP^4(z)$ are proper 3-secants to
$Q_{(\lambda : \mu )}(z) $. In fact, they are all 3-secants to the curve
$E^{\prime}_{(\lambda : \mu )}(z)$ and intersect the curve in the
10-torsion points which are different from the 5-torsion points. \par\noindent
 (iii)\quad  The numerically anticanonical curve $E_{(\lambda
:\mu)}({\tau_3})\subset
Q_{(\lambda :\mu)}(z)$ is mapped by $\rho$ isomorphically to an elliptic curve
on $S_{(\lambda : \mu )}(w)$ which meets the plane $\Pi_{\tau_1 \tau_2}$
of Lemma 4.18 in 8 non-trivial 6-torsion points.
\endproclaim
\demo{Proof} (i) follows from 4.14 and 4.19.  (ii) follows from 4.10. (iii)
follows from and 4.18.\qed
\enddemo
\head
{5. Modules obtained by concatenating three Moore matrices.}
\endhead
In this section we use the results of Section 4 to analyse
graded $R$-modules whose presentation matrix is a block of three
Moore matrices. We are particularly interested in those modules
listed in (5.1). We check under which conditions the modules in (5.1)
are artinian and we compute their linear and quadratic syzygies.
This provides enough information for the following sections, where
we will construct surfaces from these modules and determine
their syzygies completely.\medskip
Let $\Pi = <p_1,p_2,p_3>\subset \openP^4(x)$ be a
plane, let
$$
\alpha_1 = (M_{p_1}(y), M_{p_2}(y), M_{p_3}(y)),
$$
and let $M = \bigoplus\limits_{k\geq 0} M_k$ be the $H_5$-invariant graded
$R$-module presented by $\alpha_1$:
$$
0\leftarrow M\leftarrow V_3 \otimes R \overset{\alpha_1}\to
\longleftarrow 3 V_1 \otimes R(-1)\ .
$$
We call $\Pi$ the {\it representing plane} of $M$.
The minimal free resolution of $M$ is of type
$$
0\leftarrow M\leftarrow F_{0,0}\otimes R\overset\alpha_1\to\leftarrow
F_{1,1}\otimes R(-1)\overset\alpha_2\to\leftarrow
{\textstyle\bigoplus\limits_{l\geq 2}} F_{2,l}\otimes
R(-l)\overset\alpha_3\to\leftarrow \cdots
\overset\alpha_5\to\leftarrow {\textstyle\bigoplus\limits_{l\geq 5}}
F_{5,l}\otimes R(-l)\leftarrow 0
$$
with $F_{0,0} = M_0 = V_3$ and $F_{1,1} = 3V_1$.
The syzygies $F_{k,k+l}$ are canonically isomorphic to the cohomology
groups of the Koszul complexes
$$
\cdots\leftarrow\Lambda^{k-1}V_0^*\otimes M_{l+1}\leftarrow
\Lambda^kV_0^*\otimes
M_l\leftarrow \Lambda^{k+1}V_0^*\otimes M_{l-1}\leftarrow\cdots\ ,
$$
defined by the multiplication maps $M_{l+1}\leftarrow V_0^*\otimes M_l$
(see \cite{Gr}). These Koszul complexes tensorized by $R(-k-l)$ fit
as vertical lines into a double complex with $H_5$-equivariant
differentials. The horizontal lines are obtained by tensoring $M_l$
with the Koszul complex resolving $\openC(-l)$.\par\noindent
Following \cite{Po} we first consider the case where $M$ is additionally
invariant under $\iota$. In terms of $G_5$-modules the double complex reads:

%
{\eightpoint
$$
\vbox{%
\halign{&\hfil\,$#$\,\hfil\cr
&&&\downarrow&&&&\downarrow&&&&\downarrow&&&&\downarrow\cr
0&\leftarrow&V_3&\otimes&R&\overset\alpha\to\leftarrow&(3V_1\oplus
2V_1^\#)&\otimes&R(-1)&\leftarrow&(4V_0\oplus 6V_0^\#)&\otimes&
R(-2)&\leftarrow&(2S\oplus 2Z)&\otimes&R(-3)&\leftarrow\cr
&&&\downarrow&&&&\downarrow&&&&\downarrow&&&&\downarrow\cr
&&&0&&\leftarrow &\hfill M_1&\otimes&R(-1)&\leftarrow&\hfill M_1\otimes
V_2&\otimes&R(-2)&\leftarrow&M_1\otimes 2V_3^\#&\otimes&R(-3)&\leftarrow\cr
&&&\downarrow&&&&\downarrow&&&&\downarrow&&&&\downarrow\cr
&&&0&&\leftarrow&&0&&\leftarrow&\qquad\qquad
M_2&\otimes&R(-2)&\leftarrow&M_2\otimes
V_2&\otimes&R(-3)&\leftarrow\cr
&&&\downarrow&&&&\downarrow&&&&\downarrow&&&&\downarrow\cr
}}
$$
}
%
\noindent
For every $k$ the induced maps $F_{k-1,l'}\otimes R(-l')\leftarrow
F_{k,l}\otimes R(-l)$ are just the components of $\alpha_k$.
In order to study $\alpha_1$ we use the notations of the beginning of
Section 4. We want to compute the restriction of the differential $\alpha$
of the double complex to a 5-dimensional $G_5$-invariant
subspace of $V_3\otimes V_0^*$. Therefore we choose a decomposition of
$V_3\otimes V_0^* \cong 3V_1\oplus 2V_1^\#$  into
irreducible $G_5$-modules (the subspaces selected in Section 4 are not
$\iota$-invariant). Namely, we consider the subspaces whose $i^{\text{th}}$
basis vector is given by $f_i^{(0)}$, $f_i^{(1)}+f_i^{(4)}$, $f_i^{(2)}+
f_i^{(3)}$ and $f_i^{(2)}-f_i^{(3)}$, $f_i^{(1)}-f_i^{(4)}$ resp.
Every $a=(a_0:\dots: a_4)\in \openP^4$ with either $a_3=a_4=0$ or
$a_0=a_1=a_2=0$
defines a $G_5$-invariant subspace
$$
V^a\subset V_3\otimes V_0^* \cong 3V_1\oplus 2V_1^\# :
$$
The $i^{\text{th}}$ basis vector of $V^a$ is the linear combination of the
$i^{\text{th}}$ basis vectors above with
$a_0,\dots, a_4$ as coefficients. The restriction to $V^a$ of the
differential $\alpha$ of the double complex is given by the Moore
matrix
$$
M_p(y)\quad {\text {with}}\quad
p = (a_0:a_1+a_4: a_2+a_3: a_2-a_3: a_1 -a_4)\in \openP^4.
$$
Then $p\in \openP^2_+$ iff $V^a\cong V_1$ and $p\in\openP^1_-$ iff
$V^a\cong V_1^\#$\ .\par\noindent
Suppose for the moment that $M_1$ is $10$-dimensional. Then there are
three possibilities for $M_1$, namely $M_1 \cong 2 V_1^\#$,
$M_1 \cong V_1\oplus V_1^\#$ or $M_1 \cong 2 V_1$.
Correspondingly, there are the following three
possibilities for the representing plane of $M$:\medskip\noindent
(5.1), (i)\quad $<p_1,p_2,p_3> = \openP^2_+$\
.\medskip\noindent
(5.1), (ii)\quad $<p_1,p_2>\subset\openP^2_+\
,\quad p_3\in\openP^1_-$\ .\medskip\noindent
(5.1), (iii)\quad $p_1=(a_0:a_1:a_2:a_2:a_1)\in\openP^2_+\ ,\quad
<p_2,p_3>=\openP^1_-\ .$
\medskip
\noindent (5.1), (i) gives the module of Section 3 corresponding to the
HM-surfaces. By Proposition 4.2 we can choose $p_1$, $p_2$ and $p_3$
to be the non-trivial 2-torsion points on any smooth elliptic curve
$E_{(\lambda:\mu)} (x)$. Due to the results on the incidence scroll
in Section 4 we will have some control over the module $M$ in the case where
all three points are on a smooth $E_{(\lambda:\mu)} (x)$. For example, the
following is a special case of (5.1), (ii):\medskip\noindent
(5.1), (ii)'\quad $<p_1,p_2>\subset\openP^2_+$ is the line spanned by two
non-trivial 2-torsion points of $E_{(\lambda:\mu)} (x)$ and
$p_3\in\openP^1_-$ is its origin.\medskip\noindent
Let us fix one smooth curve $E_{(\lambda:\mu)} (x)$ and denote by $\tau_i$,
$i=\overline{1,3}$ the non-trivial 2-torsion points and by $\rho$ one of
the non-trivial 3-torsion points of this curve. The following cases,
in which $M$ is not invariant under $\iota$, are of particular interest
because of Lemma 4.18 and Remark 4.17:\medskip\noindent
(5.1), (iv)\quad $p_i = \tau_i$, $i=1,2$ and $p_3 = \tau_3 + \rho$.
\medskip\noindent
(5.1), (v)\quad $p_i = \tau_i + \rho$, $i=\overline{1,3}.$
\bigskip
In the sequel we will analyse graded  $R$-modules $M$ with a presentation
matrix $\alpha_1$ defined by a representing plane $\Pi = <p_1,p_2,p_3>\subset
\openP^4(x)$ with particular interest in the cases of (5.1).\medskip
Let us first check under which conditions $M$ is artinian.
We need the following notation. Let $N(u,y)$ be the $5\times 5$-matrix with
bilinear entries in $u$'s and $y$'s defined by the relation
$$
\pmatrix
u_0, u_1, u_2, u_3, u_4\\
\endpmatrix\cdot M(x,y)=
\pmatrix
x_0, x_1, x_2, x_3, x_4\\
\endpmatrix\cdot N(u,y).\tag 5.2
$$
In other words,
$$
N(u,y)=(u_{2i+j}y_{i+j})_{i,j \in \openZ_5}.
$$
\proclaim{Lemma 5.3} If $y\in\openP^4(y)$ is in the support of the module $M$,
then the  matrix $M(x,y)$ has rank at most 4 on the plane $\Pi$ and
there exists a point $u\in\openP^4(u)$ such that $p_i \cdot N(u,y)=0$,
$i=\overline{1,3}$.
In particular, $M$ is artinian iff its representing plane
$\Pi$ is not in the kernel of the linear map
$$
N(u,y) : \openP^4(x) \rightarrow \openP^{4}(x),
\quad x \mapsto x \cdot N(u,y)
$$
for every $(u,y) \in \openP^4(u) \times \openP^4(y)$.
\endproclaim
\demo{Proof} If the module $M$ has support at $y$, then there exists a point
$u$
defining a common left side syzygy for all the three Moore blocks of $\alpha_1$
at that point $y$. The lemma now follows from the defining relation (5.2).\qed
\enddemo
\proclaim{Lemma 5.4} Let $y$ be in the support of the module $M$. Then we
have the following two possibilities:\par\noindent
(i)\quad If $y$ is not on a HM-line, then $\Pi$ is contained in the
HM-quintic $\{\det M(x,y)=0\}$. In particular, the image of $\Pi$ under the
rational map
$$
\Theta :\Pn 4(x)\dashrightarrow {\Pn 5}^\ast
\quad x\mapsto (\gamma_0(x):\gamma_1(x):\dots :\gamma_5(x))
$$
does not span all of ${\Pn 5}^\ast$.\par\noindent
(ii)\quad If $y$ is on one of the HM-lines, then the plane $\Pi$ contains two
lines of a $G_5$-invariant pentagon. In particular, $\Pi$ does not intersect
a smooth elliptic curve $E_{(\lambda :\mu)}(x)$ nor does it contain a HM-line.
\endproclaim
\demo{Proof} (i)\quad If $y$ is not on a HM-line then by Remark 4.1
$\det M(x,y)$ is not vanishing identically and  $\{\det M(x,y)=0\}$
is a HM-quintic which must contain $\Pi$ by Lemma 5.3.\smallskip\noindent
(ii)\quad By $H_5$-invariance we may assume that $y\in \openP^1_-$.
Then $M(x,y)$ is skew-symmetric (up to column permutations), it has the
left side syzygy
$$
I_{E_y}(x)\cdot M(x,y) = 0,
$$
and it has rank 4 outside $E_y$ by Lemma 4.4. Lemma 5.3 and the
relation (5.2) tell us that $I_{E_y}(x)$ is a multiple of a common
$u$ for every $x\in\Pi$. The entries of $I_{E_y}(x)$ are the
quadrics which cut out the elliptic curve $E_y$. These define a
quadro-cubic Cremona transformation from $\openP^4$ to $\openP^4$
which by the above must contract the plane $\Pi$. There
are no such planes if $E_y$ is smooth since the Cremona transformation
contracts only $\Sec E$ which does not contain any plane (see
\cite {Seg}, \cite {Sem} and also \cite {ADHPR} for the properties
of the quadro-cubic Cremona transformation). If $E_y$ is
a pentagon, then precisely the planes spanned by two lines
of the pentagon are contracted. The second assertion in (ii) can be checked
from the explicit equations of the elliptic curves and the HM-lines resp.
\qed
\enddemo
\proclaim{Proposition 5.5} The module $M$ defined by the plane
$\Pi =<p_1,p_2,p_3>$ is artinian in the cases (i), (ii)', (iv)
and (v) of (5.1) where all three points lie on a smooth $G_5$-invariant
elliptic normal curve $E\subset\openP^4$.
\endproclaim
\demo{Proof} Suppose that the module $M$ has support at some point
$y\in\openP^4(y)$. Then $\Pi$ is contained in the HM-quintic
$\{\det M(x,y)=0\}$ by Lemma 5.4, (ii) and (i). We will show that this quintic
coincides with $\Sec E$, a contradiction,
since $\Sec E$ contains no planes.\smallskip\noindent
Let more generally $U$ be any HM-quintic containing
all three lines $l_{ij}=\overline {p_i p_j}$, $i < j.$ The $l_{ij}$ are
secant lines to $E$, and we denote by $T_{ij}$ the translation scrolls
of $E$ generated by these lines (compare \cite{BHM2}).
$T_{ij}$ has degree 5 if $p_i-p_j$ is 2-torsion and degree 10 otherwise.
For every $i < j$ there exists precisely one section $s_{ij}$ of the
HM-bundle $\scr E(3)$ whose reduced zeroscheme is $T_{ij}$ {\cite{BHM2}}.
For our purpose it suffices to show that $U$ contains at least two different
$T_{ij}$'s since then $U$ is uniquely determined as the wedge product of the
corresponding sections of the HM-bundle. In particular,  $U$ must coincide
with $\Sec E$, because this HM-quintic contains all translation scrolls of
$E$ as above.\smallskip\noindent
The orbit of $l_{ij}$ under $G_5$ consists of 25 lines if
$l_{ij}$ is $\iota $-invariant, i.e., if $p_i+p_j=0$ or $p_i$ and $p_j$
are 2-torsion, and of 50 lines otherwise.
By Bezout $U$ contains $E$ and thus $T_{ij}$ unless
possibly if $l_{ij}$ is $\iota $-invariant and $T_{ij}$ has degree 10.
In the cases (i) and (ii)' $p_1$, $p_2$ and $p_3$ are 2-torsion points.
Therefore the $T_{ij}$ are the three distinct quintic elliptic scrolls which
are the 2-torsion translation scrolls of $E$ and they are all contained
in $U$. In case (iv) $T_{12}$ has degree 5, whereas $T_{13}$ and $T_{23}$
are two different scrolls of degree 10, and the lines $l_{13}$ and $l_{23}$
are not $\iota$-invariant.
In case (v) non of the lines is $\iota$-invariant and the three scrolls of
degree 10 are different. In any case $U$ contains three different $T_{ij}$'s.
\qed\enddemo
\proclaim{Proposition 5.6} Let $p_1=(a_0:a_1:a_2:a_2:a_1)\in\Pn 2_+$,
$p_2 = (0:1:0:0:-1)$ and $p_3 = (0:0:1:-1:0)$ as in case (5.1), (iii).
Then the  module $M$ defined by the plane $\Pi=<p_1,p_2,p_3>$ is artinian
iff $p_1\in\Pn 2_+\setminus\Delta$, where $\Delta$ is the union of the
6 lines
$$
x_0=0 \quad {\text and} \quad x_0+\xi^i x_1+\xi^{-i} x_2=0,\quad i\in\openZ_5,
$$
corresponding to the unique minimal $A_5$-orbit in $(\Pn 2_+)^*$.
\endproclaim
\demo{Proof} By Lemma 5.3 we have to show that $\Pi$ is contained
in the kernel of the linear map
$$
N(u,y) : \openP^4(x) \rightarrow \openP^{4}(x),
\quad x \mapsto x \cdot N(u,y)
$$
for some $(u,y) \in \openP^4(u) \times \openP^4(y)$ iff $p_1\in\Delta$.
The variety
$$
I_- = \{(u,y)\in \openP^4(u) \times \openP^4(y)\mid \Pn 1_-\subset
\ker N(u,y)\}
$$
is defined by the 10 bilinear forms
$$
u_{i+2} y_{i+1} - u_{i+3} y_{i+4} = u_{i+4} y_{i+2} -
u_{i+1} y_{i+3} = 0,
\quad i\in\openZ_5.
$$
The saturation of these forms with respect to $y_0, \dots, y_4$ contains
the 10 cubic equations
$$
v^+(u)=v^-(u)=0,
$$
which cut out the 30 vertices of the 6 complete pentagons
in $S_{15}\subset\openP^4(u)$ (compare Remark 4.9, (i)). Conversely, every
such vertex is contained
in the image of $I_-$ under the projection to $\openP^4(u)$.
Under the action of the normalizer $N_5$ of the
Heisenberg group these pentagons are equivalent to the one with
vertices $e_0, \dots, e_4$. For $(u=e_i,y)\in I_-$ the plane $\Pi$ is
contained in $\ker N(u,y)$ iff $p_1$ is on the line $x_0=0$. The lemma
follows by the action of $N_5$.\qed
\enddemo
\noindent
 From now on we will always suppose in case (5.1), (iii) that
$p_1\in\Pn 2_+\setminus\Delta$.\medskip
Next we will study the first order linear syzygies of $\alpha_1$.
It is clear from the representation theory described at the
beginning of Section 4 and from the double complex of the Koszul
cohomology that these come in
blocks of $L$-matrices. We recall from Section 4 that
the product matrices $M_{p}(y)L_{q}(y)$ have quadratic entries in
$y$'s with coefficients that are bilinear forms in $x$'s and $z$'s
defining a morphism
$$
\Phi_T : \openP^4(x) \times \openP^4(z) \rightarrow \openP^{14},
$$
which is a projection of the Segre embedding.
Five first order linear syzygies of $\alpha_1$ are given by a non-trivial
relation
$$
M_{p_1}(y)L_{q_1}(y)+M_{p_2}(y)L_{q_2}(y)+M_{p_3}(y)L_{q_3}(y) = 0
$$
with $q_i\in\openP^4(z)$ or $L_{q_i}(y) = 0$, $i=\overline{1,3}$.
We will express such a relation in terms of $\Phi_T$ and thus in terms of
the three linear subspaces defined as $R_{p_i} = \Phi_T(\{p_i\}\times
\openP^4(z))$, $i=\overline{1,3}$ (clearly, these subspaces are $\openP^4$'s).
First notice that the following case cannot occur,
since $\Phi_T$ is a morphism: \smallskip\noindent
(I)\quad $(p_i,q_i)$ is a base point of $\Phi_T$ for one $i$, while the
other two $L$-matrices are zero.\smallskip\noindent
Therefore every  five first order linear syzygies of $\alpha_1$ are a
linear combination of those corresponding to one of the following cases.
\smallskip\noindent
(II)\quad $\Phi_T(p_i,q_i) = \Phi_T(p_j,q_j)$ is a double point of $R_{p_i}\cup
R_{p_j}$ for one pair $i<j$, while the third $L$-matrix is zero.
\smallskip\noindent
(III)\quad The three different points $\Phi_T(p_i,q_i)$ are collinear, i.e.,
these points define a proper trisecant to $R_{p_1}\cup R_{p_2}\cup R_{p_3}$.
\smallskip\noindent
\proclaim
{Lemma 5.7}
Let $c\in\{1,2\}$ and suppose that $R_{p_1}$, $R_{p_2}$ and $R_{p_3}$
span a $\openP^{14-c}$. Then $\alpha_1$ has precisely $5c$ first order
linear syzygies.
\endproclaim
\demo
{Proof}
(i)\quad If $R_{p_1}$, $R_{p_2}$ and $R_{p_3}$ span a $\openP^{13}$,
then there are two cases:\par\noindent
a)\quad Two of the $\openP^4$'s meet in a point and their span is disjoint
from the third one.\par\noindent
b)\quad The $\openP^4$'s are pairwise disjoint and any of the three intersects
the span of the two others in a point.\par\noindent
In case a) there is no proper trisecant. In case b) let $\{i,j,k\} =
\{1,2,3\}$.
We spot one trisecant by projecting from the point of intersection
of $R_{p_k}$ with the span of $R_{p_i}$ and $R_{p_j}$ to a $\openP^{8}$,
where the images of $R_{p_i}$ and $R_{p_j}$ have to meet in a point. This
trisecant is the unique proper trisecant (independent of the choice
of $i,j,k$). So in both cases there are precisely 5 first order linear
syzygies.\par
\smallskip\noindent
(ii)\quad If $R_{p_1}$, $R_{p_2}$ and $R_{p_3}$ span a $\openP^{12}$,
then there are five cases:\par\noindent
a)\quad Two of the $\openP^4$'s meet in a line and their span is disjoint
from the third one.\par\noindent
b)\quad One of the $\openP^4$'s meets any of the others in a point, the
two others do not meet.\par\noindent
c)\quad $R_{p_1}$, $R_{p_2}$ and $R_{p_3}$ meet in a common point.\par\noindent
d)\quad Two of the $\openP^4$'s meet in a point and their span intersects
the third one in another point.\par\noindent
e)\quad The $\openP^4$'s are pairwise disjoint and any of the three intersects
the span of the two others in a line.\par\noindent
Arguing as above we see that there is no proper trisecant in the cases a), b)
and c), one in case d) and a pencil of them in case e). Altogether, there are
precisely 10 first order linear syzygies in each case.\qed
\enddemo
\proclaim
{Proposition 5.8}
For generic choices in (5.1), (ii) $\alpha_1$ has precisely 5 first order
linear syzygies. In the other cases of (5.1) there are precisely 10 such
syzygies $(L_{q_{ij}})_{{1\le i\le 3}\atop{1\le j\le 2}}$.
The points $q_{ij}$ can be chosen as
follows:\par\noindent
(i) $L_{q_{31}} = L_{q_{22}} = 0,\quad q_{ij}\in\openP^1_-$
otherwise, $q_{11} = q_{12}$.\par\noindent
(ii)' $L_{q_{31}} = L_{q_{22}} = 0$, $q_{11}$, $q_{12}$,
$q_{21}\in\openP^1_-$, $q_{11} = q_{12}$ and $q_{32} =
\psi_5(o_{E_{(\lambda:\mu)}})\in\openP^2_+.$\par\noindent
(iii) $q_{11}=(0:1:0:0:-1)$, $q_{21}=(0:-a_1:a_0:a_0:-a_1)$,
$q_{31}=(-2a_1:a_2:0:0:a_2)$, $q_{12}=(0:0:1:-1:0)$,
$q_{22}=(2a_2:0:-a_1:-a_1:0)$ and
$q_{32}=(0:a_0:-a_2:-a_2:a_0)$.\par\noindent
(iv) $L_{q_{31}} = 0$, $q_{i1}\in\openP^1_-$, $i=1,2$ and $q_{32}
\in E_{(\lambda:\mu)}(\tau_3)$ not on a HM-line.\par\noindent
(v) $L_{q_{31}} = L_{q_{22}} = 0$, and $q_{ij}\in E^{\prime}_
{(\lambda:\mu)}$ not on a HM-line otherwise, $q_{11} = q_{12}$.
\endproclaim
\demo
{Proof} We first study the situation, where all three points $p_1$, $p_2$
and $p_3$ are on the elliptic curve $E_{(\lambda:\mu)} (x)$. Let
$\openP^4_{(\lambda:\mu)} \subset \openP^{14}$ denote the $\openP^4$
of quadrics defining $E_{(\lambda:\mu)}$. Then $\openP^4_{(\lambda:\mu)}$
meets $W := \Phi_T(E_{(\lambda:\mu)} \times \openP^4(z))$
along the degree 15 scroll $S_{(\lambda:\mu)}$ as in Lemma 4.15.
Now define $W'$ as the $\openP^2$-scroll obtained by projecting $W$
from $\openP^4_{(\lambda:\mu)}$.\par
\proclaim{Lemma 5.9}
(i)\quad $W'$ is a linearly normal, smooth
$\openP^2$-scroll in $\Pn 9$ of degree 10.\par\noindent
(ii)\quad $W'$ has no proper trisecants.
\endproclaim
\demo{Proof}
Recall from Proposition 4.16 that $S_{(\lambda:\mu)}$ is a birational image of
the symmetric product $S^2(E_{(\lambda:\mu)})$. More precisely, if
${\scr E}_{(\lambda:\mu)}$ denotes the unique non-split extension
$$
0\rightarrow{\scr O}_{E_{(\lambda:\mu)}}(-o_{E_{(\lambda:\mu)}})
\rightarrow{\scr E}_{(\lambda:\mu)}\rightarrow
{\scr O}_{E_{(\lambda:\mu)}}
\rightarrow 0,
$$
where $o_{E_{(\lambda:\mu)}}$ is the origin of $E_{(\lambda:\mu)}$,
then the correspondence in Proposition 4.16 maps
$\openP({\scr E}_{(\lambda:\mu)}(-7o_{E_{(\lambda:\mu)}}))$
to $S_{(\lambda:\mu)}$ via the morphism induced by
${\scr O}_{\openP({\scr E}_{(\lambda:\mu)}(-7o_{E_{(\lambda:\mu)}}))}(1)$.
Then the identification in Lemma 4.15 of $S_{(\lambda:\mu)}$ with the image
of $I_{(\lambda:\mu)}$ under $\Phi_T$ accounts for the
exact sequence of vector bundles on $E_{(\lambda:\mu)}$:
$$
0\rightarrow{\scr E}_{(\lambda:\mu)}(-7o_{E_{(\lambda:\mu)}})
\rightarrow H^0({\scr I}_{E_{(\lambda:\mu)}}(2))\otimes
{\scr O}_{E_{(\lambda:\mu)}}(-1)\rightarrow
{\scr Q}
\rightarrow 0,
$$
where ${\scr Q}$ is a rank 3 vector bundle such that $W'$ coincides
with the image of $\openP({\scr Q})$ under the map induced by
${\scr O}_{\openP({\scr Q})}(1)$. Twisting this sequence by ${\scr
O}_{E_{(\lambda:\mu)}}(1)$ we get
$$
0\rightarrow{\scr E}_{(\lambda:\mu)}(-2o_{E_{(\lambda:\mu)}})
\rightarrow H^0({\scr I}_{E_{(\lambda:\mu)}}(2))\otimes
{\scr O}_{E_{(\lambda:\mu)}}\rightarrow
{\scr Q}(1)
\rightarrow 0.
$$
Together with its dual this sequence presents naturally $\openP({\scr
Q^{\vee}}(-1))\subset ({\Pn 4})^{\ast}\times E_{(\lambda:\mu)}$ as the dual
of the scroll $\openP({\scr E}_{(\lambda:\mu)}(-2o_{E_{(\lambda:\mu)}}))
\subset \Pn 4\times E_{(\lambda:\mu)}$.  The projections to $\Pn 4$ and its
dual resp. map the $\Pn 1$-bundle to an elliptic quintic scroll and the
$\Pn 2$ -bundle to the dual quintic  trisecant  scroll resp.  Thus we know
from [ADHPR] that ${\scr Q}(3o_{E_{(\lambda:\mu)}})$ is indecomposable of
degree -1. Furthermore ${\scr Q^{\vee}}(-3o_{E_{(\lambda:\mu)}})$ has one
section.
In fact, $h^0({\scr Q^{\vee}}\otimes \scr L^{\vee})=1$
for all degree 3 line bundles $\scr L$ on $E_{(\lambda:\mu)}$ and
$h^0({\scr Q^{\vee}})=10$.  Thus $W'$ is linearly normal of degree 10 in
$\Pn 9$,
and every triple of planes in $W'$ spans a hyperplane.  Therefore any two
planes must span a $\Pn 5$ and no three planes can have a trisecant.   We
conclude that $W'$ is smooth with no trisecant.\qed
\enddemo
\noindent
As a consequence of the lemma we see that $W$ is smooth outside
$\openP^4_{(\lambda:\mu)}$ and that any proper trisecant to
$W$ lies in $\openP^4_{(\lambda:\mu)}$, and thus must be
a trisecant to $S_{(\lambda:\mu)}$. In particular, up to linear
combinations, 5 first order linear syzygies of $\alpha_1$ correspond either
to a point of intersection of two of the three rulings
$\tilde R_{p_i}=R_{p_i}\cap\openP^4_{(\lambda:\mu)}$, or
to a proper trisecant of these rulings.\smallskip\noindent
In the cases (i) and (ii)'
$\tilde R_{p_1}\cap\tilde R_{p_2}\cap\tilde R_{p_3}=\{o_{E_{(\lambda:\mu)}}\}$,
while in case (v) the rulings intersect in the chosen 3-torsion point $\rho$
(compare Proposition 4.16, (iii) and Remark 4.17, (ii)). Here we identify
as always points on $E_{(\lambda :\mu)}(x)$ with points on
$E_{(\lambda :\mu)}(w)$.
In case (iv), however, $\tilde R_{p_3}$ hits the plane spanned by
$\tilde R_{p_1}$ and $\tilde R_{p_2}$ in the point
$\Phi_T(p_3, q_{32})$, where $q_{32}\in E_{(\lambda :\mu)}({\tau_3})$
by Lemma 4.18 and Corollary 4.21. Notice more precisely that
$\Phi_T(p_3, q_{32})$ lies on a ruling of $S_{(\lambda :\mu)}$ which meets
$E_{(\lambda :\mu)}$ in a non-trivial 6-torsion point. In both cases (iv)
and (v) the points $q$ are not on a HM-line, since any point $(p,q)$  with
$q$ in $S^2(E_{(\lambda :\mu)})$ on a HM-line would map to a ruling of
$S_{(\lambda :\mu)}$  which meets $E_{(\lambda :\mu)}$ in
a non-trivial 5-torsion point again by Corollary 4.21. Therefore there are
precisely 10 first order linear syzygies in all cases where $p_1$, $p_2$
and $p_3$ lie on the elliptic curve $E_{(\lambda:\mu)}(x)$, and the
points $q$ can be chosen as claimed.
\smallskip\noindent
Of course (i) is known from \cite{De2}, \cite{Ma3} and \cite{Ma4} and
can also be seen by just computing the equations
of $R_{p_1}$, $R_{p_2}$ and $R_{p_3}$, where, say, $p_1 = (1:0:0:0:0)$,
$p_2 = (0:1:0:0:1)$ and $p_3 = (0:0:1:1:0)$.
A similar calculation, depending on $p_1 = (a_0:a_1:a_2:a_2:a_1)
\in\Pn 2_+\setminus\Delta$, and with $p_2 = (0,1,0,0,-1) $ and
$p_3 = (0,0,1,-1,0)$, shows the existence of precisely ten first order linear
syzygies in case (iii). It is easy to compute that
$\sum_{i=1}^{3} M_{p_i}\cdot L_{q_{ij}}=0$, $j=1,2$, when the $q_{ij}$
are chosen as claimed (see also \cite{Po}).\smallskip\noindent
We finally need to show that for generic choices in (5.1), (ii)
$\alpha_1$ has precisely 5 first order linear syzygies.
For this, we consider first $S_-=\Phi_T (\Pn 1_-(x) \times \Pn 2_+(z))$ and
$S_+=\Phi_T (\Pn 2_+(x)\times \Pn 1_-(z))$.
Then $S_-$ is contained in the $\Pn 5\subset \Pn {14}(t)$ defined by
$$
t_{00}=t_{10}=t_{20}=t_{02}+t_{03}=t_{12}+t_{13}=t_{22}+t_{23}=t_{01}+t_{04}
=t_{11}+t_{14}=t_{21}+t_{24}=0,
$$
while $S_+$ is contained in the $\Pn 3$ inside this $\Pn 5$ cut out by the
additional
$$
t_{11}+t_{03}=t_{01}+t_{22}=0.
$$
Thus for $p_1,p_2\in \Pn 2_+$ and $p_3\in \Pn 1_-$ the subspaces $R_{p_1}$
and $R_{p_2}$ will meet $S_+$ in a line while $R_{p_3}$ will meet $S_-$ in
a plane.  Since this latter plane will meet the $\Pn 3$ of $S_+$ in at least
a point there must be at least one not necessarily proper trisecant line for
the  $R_{p_i}$, $i=\overline{1,3}$. Thus the corresponding $\alpha_1$ has at
least 5 linear syzygies. It now suffices to find three points such that the
corresponding $\openP^{4}$'s span a $\openP^{13}$. E.g., we may take
$p_1=(1:0:0:0:0)$, $p_2=(0:1:0:0:1)$ and $p_3=(0:1:1:-1:-1)$ in which case
$R_{p_1}$, $R_{p_2}$ and $R_{p_3}$ span the hyperplane defined by
$$
2t_{20}+t_{11}-t_{02}-t_{03}+t_{14}=0
$$
in $\Pn {14}(t)$.
\qed
\enddemo
\noindent
Note that in the $G_5$-invariant cases the first order linear syzygies are
decomposed into $G_5$-representations. This will be the case also for the rest
of the syzygies below, although we do not make it explicit since
$H_5$-invariance is enough to make the computations.
\proclaim
{Proposition 5.10}
In the above cases of (5.1) with 10 first order linear syzygies
$(L_{q_{ij}})$ the linear part of the minimal free resolution of $M$,
written in terms of $H_5$-modules, is as follows:
$$
V_3\otimes  R\overset {(M_{p_i})}\to\leftarrow 3V_1\otimes R(-1)
\overset {(L_{q_{ij}})}\to\leftarrow 2V_0\otimes R(-2)\overset K
\to\leftarrow aI\otimes  R(-3)\leftarrow 0,
$$
where $K$ is a direct sum of Koszul maps $[y]^t=(y_0,\dots , y_4)^t$.
Moreover $a = 2$ in the case (i), $a = 1$ in the cases (ii)' and (iv)
and $a = 0$ in the cases (iii) and (v).
\endproclaim
\demo{Proof}
We will take into account the results of Proposition
5.8 concerning the points $q_{ij}$.\par\noindent
Notice that $L_q(y)[y]^t = 0$ iff $q\in\openP^1_-.$ In this case
$[y]^t$ defines the only linear syzygy of $L_q(y)$
(compare Remark 4.6). If however $q$ is a point not on one
of the HM-lines, then $\det L_q(y)$ does not vanish identically
by Remark 4.1, and $L_q(y)$ has no syzygies at all.
This settles the cases (i), (ii)', (iv) and (v). In case (iii)
we compute that the linear syzygies of $(L_{q_{11}},L_{q_{12}})$ are
precisely a direct sum of two copies of $[y]^t$. These do not give rise
to a linear syzygy of $(L_{q_{ij}})_{{2\le i\le 3}\atop{1\le j\le 2}}$.\qed
\enddemo
\proclaim
{Proposition 5.11}
In the above cases of (5.1) with 10 first order linear syzygies
$(L_{q_{ij}})$ there are 4 first order quadratic syzygies of $\alpha_1$
in the case (i), 2 in the case (iii) and 1
in the cases (ii)', (iv) and (v).
\endproclaim
\demo{Proof}
The decomposition of the space of quadrics in $\Pn 4$ into irreducible
$H_5$-modules
is $S^2V_0^* \cong 3V_3$. Explicitely such a decomposition is given by
distinguishing the following subspaces:
$$
<y_i^2>_{i\in\openZ_5},\quad <y_{i+1}y_{i+4}>_{i\in\openZ_5},
\quad <y_{i+2}y_{i+3}>_{i\in\openZ_5}.
$$
Moreover, $S^2V_0^*\otimes V_1$ decomposes into characters of $H_5$.
 From the double complex of the Koszul cohomology we see that
quadratic syzygies of $\alpha_1$ are elements in
$$
\Hom_{H_5}(\chi \otimes \O{-3}, 3V_1 \otimes \O{-1})
\cong \Hom_{H_5}(\chi, S^2V_0^* \otimes 3V_1),
$$
where $\chi$ is some character of $H_5$.\medskip\noindent
We first study those first order quadratic syzygies which involve the
trivial character. These come in blocks of column vectors of type
$$
L^{\prime}_b(y) =\pmatrix b_0 y^2_{i} + b_1 y_{i+2}y_{i+3}+ b_2 y_{i+1}
y_{i+4}\\
\endpmatrix_{i\in\openZ_5},
$$
where $b = (b_0:b_1:b_2)$ is a parameter point in $\openP^{2}$.
The $i^{\text{th}}$ entry of the  $5\times 1$-matrix
$M_x(y) L^{\prime}_b(y)$ is a cubic in the span of
$$
y_{i}^3,\ y_{i}y_{i+1}y_{i+4},\ y_{i}y_{i+2}y_{i+3},\ y_{i+1}^2y_{i+3},
\ y_{i+1}y_{i+2}^2,\ y_{i+2}y_{i+4}^2,\ y_{i+3}^2y_{i+4},
$$
$i = \overline{0,4}$, whereas the coefficients are bilinear forms in $x$'s
and $b$'s. We collect these forms in the column vector
$$
T^{\prime}(x,b) =
\pmatrix
x_0 b_0\\
x_0 b_2 + (x_1 + x_4) b_1\\
x_0 b_1 + (x_2 + x_3) b_2\\
x_2 b_0 + x_4 b_2\\
x_3 b_1 + x_4 b_0\\
x_1 b_2 + x_3 b_0\\
x_1 b_0 + x_2 b_1
\endpmatrix.
$$
$T^{\prime}$ defines a rational map
$$
\Phi_{T^{\prime}} : \openP^4(x) \times \openP^2(b) \rto \openP^{6}(s),
$$
well-defined outside the codimension 5 set
$$
B_{\Phi_{T^\prime}}=(\Pn 1_-(x) \times \openP^2(b)) \cap
\{b_0^2+b_1b_2=x_2b_0-x_1b_2=x_2b_1+x_1b_0=0\},
$$
which is a rational curve projecting bijectively onto $\Pn 1_-$ in the first
factor, and onto the conic $\{b_0^2+b_1b_2=0\}\subset\openP^2(b)$ in the
second factor.\smallskip\noindent
So a first order quadratic syzygy for $\alpha_1$ involving the trivial
character is given by a non-trivial relation
$$
M_{p_1}(y)L^{\prime}_{r_1}(y)+M_{p_2}(y)L^{\prime}_{r_2}(y)+M_{p_3}(y)
L^{\prime}_{r_3}(y)=0
$$
with $r_i\in\openP^2(b)$ or $L^{\prime}_{q_i}(y) = 0$, $i=\overline{1,3}$.
We may express such a relation in terms of $\Phi_{T^{\prime}}$ and thus
in terms of the three linear subspaces
$S_{p_i}:=\Phi_{T^{\prime}}(\{p_i\}\times \openP^2(b))\subset
\openP^{6}(s)$, $i=\overline{1,3}$. Then we see that every quadratic syzygy
is a linear combination of those corresponding to one of the
following cases.\smallskip\noindent
(I)\quad A basepoint, i.e., $(p_i,r_i)\in B_{\Phi_{T^\prime}}$ for some $i$,
while the other two $L^\prime$-matrices are zero.
\smallskip\noindent
(II)\quad A double point of $S_{p_i}\cup S_{p_j}$, i.e.,
$\Phi_{T^\prime}(p_i,r_i)=\Phi_{T^\prime}(p_j,r_j)$ for one pair $i<j$,
while the third $L^\prime$-matrix is zero.\smallskip\noindent
(III)\quad A proper trisecant line of $S_{p_1}\cup S_{p_2}\cup S_{p_3}$,
i.e., the three different points $\Phi_{T^\prime}(p_1,r_1)$,
$\Phi_{T^\prime}(p_2,r_2)$ and
$\Phi_{T^\prime}(p_3,r_3)$ are  collinear.\par\smallskip\noindent
The genuine quadratic syzygies which appear in the minimal resolution
of the module $M$ are determined only modulo those quadratic syzygies
which depend on the linear ones. These in turn appear as non-trivial
relations
$$
(M_{p_1}(y)L_{q_1}(y)+M_{p_2}(y)L_{q_2}(y)+M_{p_3}(y)L_{q_3}(y))
(y_0,\dots,y_4)^t=0.
$$
Each summand is a product matrix $M_p(y)L_q(y)(y_0,\dots,y_4)^t$
with cubic entries as above, but  this time the coefficient vector is
$$
T^{\prime\prime}(x,z) =
\pmatrix
x_0z_0 \\
x_0(z_1+z_4)+(x_1+x_4)(z_2+z_3)\\
x_0(z_2+z_3)+(x_2+x_3)(z_1+z_4)\\
x_2z_0 + x_4(z_1+z_4)\\
x_3(z_2+z_3) + x_4z_0\\
x_1(z_1+z_4)+x_3z_0\\
x_1z_0+x_2(z_2+z_3)
\endpmatrix.
$$
$T^{\prime\prime}$ defines a rational map
$$
\Phi_{T^{\prime\prime}}: \openP^4(x) \times \openP^4(z) \rto \openP^{6}(s),
$$
which is a projection of $\Phi_T$. Therefore the linear map $\rho$ defined by
$$
b_0\mapsto z_0,\qquad b_1\mapsto(z_2+z_3),\qquad
b_2\mapsto(z_1+z_4)
$$
gives rise to a commutative diagram
$$
\diagram[tight,width=4em,height=3em]
{\Pn 4(x)\times\Pn 2(b)} &\rDashto^{\Phi_{T^{\prime}}} &{\Pn 6(s)}\\
\uDashto^{\rho}&\ruDotsto^{\Phi_{T^{\prime\prime}}}&\uDashto_{\pi_T}\\
{\Pn 4(x)\times\Pn 4(z)} &\rTo^{\Phi_{T}}&{\Pn {14}(t)}.\\
\enddiagram
 $$
\noindent
The projection $\pi_T$ restricted to $\Phi_T(\Pn 4(x)\times\Pn 4(z))$ is
well-defined outside the image via $\Phi_T$ of the base locus
$B_{\Phi_{T^{\prime\prime}}}$ of $\Phi_{T^{\prime\prime}}$.
The first component of $B_{\Phi_{T^{\prime\prime}}}$ is precisely the base
locus ${\Pn 4(x)\times\Pn 1_-(z)}$ of the map
$\rho:\openP^4(x)\times \openP^4(z)\rto\openP^4(x)\times \openP^2(b),$
while the second component
$$
{((\Pn 1_-(x)\times\Pn 4(z))}\cap
\{z_0^2+(z_1+z_4)(z_2+z_3)=x_2z_0-x_1(z_1+z_4)=x_2(z_2+z_3)+x_1z_0=0\}
$$
is a $\openP^2$-bundle over $\openP^1_-(x)$ which maps onto the quadric cone
$\Gamma=\{z^2_0+(z_1+z_4)(z_2+z_3)=0\}\subset\openP^4(z)$.
Under this projection the plane over a point $o\in \openP^1_-(x)$ is
mapped to the plane on $\Gamma$ over the point $\psi_5(o)$ on
$C_+\subset \openP^2_+(z)$.\par\smallskip\noindent
We will need one further fact, namely, that $S_p=\Phi_{T^{\prime}}(\{p\}
\times \openP^2(b))\subset
\openP^{6}(s)$ coincides with the fixed line
$$S_-=\{s_0=s_1=s_2=s_3+s_5=s_4+s_6=0\}\subset\Pn 6(s)$$
for all $p\in \openP^1_-(x)$, whereas $S_p$ is a plane
inside the 4-dimensional linear subspace
$$S_+=\{s_3-s_5=s_4-s_6=0\}\subset\openP^6(s)$$
for all $p\in \openP^2_+$. \par\medskip\noindent
Now we may proceed to study the cases of (5.1), in each case determining all
quadratic syzygies and comparing with Proposition 5.8 to see which of them are
dependent on the linear ones.\par\smallskip\noindent
In case (i) we get three planes $S_p$ inside $S_+$. Each pair of planes
intersect
in a point inducing thus three quadratic syzygies of type (II), while there is
exactly one proper trisecant to all three planes
giving rise to a fourth syzygy of type (III).
In this case the points giving linear syzygies are
base points for the map $\rho$, so $\alpha_1$ has 4 genuine quadratic syzygies.
\smallskip\noindent
In case (ii)' the planes $S_{p_1}$ and $S_{p_2}$ span $S_+$ while
$S_{p_3}=S_-$. Since $S_+$ and $S_-$ span $\openP^6(s)$, there are 2 candidates
for quadratic
syzygies: one of type (II) from the double point of the first two planes, and
one of
type (I) from the basepoint of $\Phi_{T^{\prime}}$ over the point $p_3$.
Of the linear syzygies in this case, the points
$(p_1,q_{11}),(p_1,q_{12}),(p_2,q_{21})$
lie in the base locus $(\Pn 4(x)\times\Pn 1_-(z))$ of $\rho$, while
$(p_3,q_{32})$ is mapped
to the basepoint of $\Phi_{T^{\prime}}$ over $p_3$.  Thus only the
quadratic syzygy of type (I) is dependent on the linear ones, and
$\alpha_1$ has exactly one genuine quadratic syzygy, which we may
choose to be of type (II). \smallskip\noindent
In case (iii) $S_{p_1}$ is a plane in $S_+$, while $S_{p_2}=S_{p_3}=S_-$.
Thus there are 2 quadratic syzygies of type (I) from the basepoints over
$p_2$ and $p_3$ and 2 quadratic syzygies of type (II)
from the doublepoints on $S_-$.  Of the linear syzygies,
the points $(p_1,q_{11})$ and $(p_1,q_{21})$ are basepoints of $\rho$, while
the other points are not. Therefore the linear syzygies give rise to
2 quadratic syzygies of type (I) or (II), depending on the choice of
$p_1 = (a_0:a_1:a_2:a_2:a_1)\in\Pn 2_+\setminus\Delta$. In any case,
$\alpha_1$ has precisely 2 genuine quadratic syzygies.
\par\smallskip\noindent
For the last two cases we will take a closer look at the  condition
that all three points $p_1$, $p_2$ and $p_3$ lie on a smooth elliptic curve
$E_{(\lambda :\mu)}(x)$. First we look at the projection via $\pi_T$ of the
$\openP^4$- bundle $W = \Phi_T(E_{(\lambda :\mu)}(x) \times
\openP^4(z))\subset\openP^{14}(t)$ into $\openP^6(s)$.
Inside $W$ lies the image $\tilde U$ of $E_{(\lambda :\mu)}(x) \times
\openP^2_+(z)$. Since $\Phi_T(\openP^4(x) \times \openP^1_-(z))$
is in the base locus of the projection,
the images of $W$ and $\tilde U$ in $\openP^6(s)$ coincide.
Furthermore $\tilde U$ is projected from a point on the
plane over the origin, hence its image
$U$ has degree one less.\par\noindent
We have seen that the linear syzygies all come from a double point
or from a proper trisecant to the scroll $S_{(\lambda :\mu)}$
inside $\openP^4_{(\lambda :\mu)} \subset \openP^{14}(t)$. This
$\openP^4_{(\lambda :\mu)}$ meets the projection center in the
origin $o_{(\lambda :\mu)}$ of the elliptic curve $E_{(\lambda :\mu)}$
in $S_{(\lambda :\mu)}$.  $S_{(\lambda :\mu)}$ is 4-tuple in this point so
its image $T_{(\lambda :\mu)}$ under $\pi_T$ must be a scroll of degree 11 in a
$\openP^3_{(\lambda :\mu)}\subset \openP^6(s)$.
This scroll is 4-tuple along an elliptic quartic curve $F_{(\lambda :\mu)}$,
which is the image of $E_{(\lambda :\mu)}$ under the projection from the
origin.\par\noindent
The scroll $T_{(\lambda :\mu)}$ sits inside the $\openP^2$-scroll $U$,
which in turn spans $\openP^6(s)$ since the line
$S_{o_{(\lambda :\mu)}}=S_-$ over the origin and  the three planes
over the nontrivial 2-torsion points span $\openP^6(s)$. \par\noindent
The degree of $\tilde U$ is 15, whence the degree of $U$ is 14, and any
hyperplane
which contains $T_{(\lambda :\mu)}$ must intersect $U$ in three more planes.
Three planes in $U$ over the points $p_i$ do have a double point or a proper
trisecant outside the span of $T_{(\lambda :\mu)}$ iff they span
at most a hyperplane in $\openP^6(s)$. This happens exactly when the three
planes together with $T_{(\lambda :\mu)}$ span a hyperplane.
Now $\tilde U$ is the $\openP^2$-bundle over $E = E_{(\lambda:\mu)}$  of a
rank 3 vector bundle with determinant ${\scr O}_E(15o)$,
$o =o_{E_{(\lambda:\mu)}}$, and since
one gets $U$ by projecting from a point over the origin its rank
3 bundle has determinant  ${\scr O}_E(14o)$.  The rank 2
vector bundle associated to $S_{(\lambda :\mu)}$ has also determinant
${\scr O}_E(15o)$, so the bundle associated to $T_{(\lambda :\mu)}$
has determinant ${\scr O}_E(11o)$. Therefore three planes in $U$
over the points $p_i$ have a double point or a trisecant outside the span of
$T_{(\lambda :\mu)}$ precisely when the three points $p_i$
add up to  the origin in the group law of the elliptic curve.
\par\smallskip\noindent
We can now settle the last two cases.\smallskip\noindent
In case (iv) the three $S_{p_i}$'s are planes which meet
$T_{(\lambda:\mu)}$ along lines. The three points do not add
up to the origin, thus there are no trisecant lines
outside the $\openP^3_{(\lambda :\mu)}$ spanned by $T_{(\lambda :\mu)}$.
On the other hand, the two lines over 2-torsion points do
intersect on $T_{(\lambda :\mu)}$, while all three lines
have a pencil of trisecant lines on $T_{(\lambda :\mu)}$. This accounts
for one quadratic syzygy of type (II) and one of type (III).
Of the linear syzygies, the points $(p_i, q_{i1})$ are
basepoints of $\rho$, while the points $(p_i, q_{i2})$ are mapped to a
proper trisecant line. Therefore $\alpha_1$ has exactly one genuine
quadratic syzygy, which we may choose to be of type (II).
\smallskip\noindent
In case (v) the three $S_{p_i}$'s meet in $T_{(\lambda :\mu)}$ as they did
in $S_{(\lambda :\mu)}$, and there is one trisecant line to the above planes
which lies outside $T_{(\lambda :\mu)}$ (the three points in this case do
add up to the origin). So there is only one genuine quadratic syzygy which is
is of type (III).
\par\smallskip\noindent
This concludes the analysis of quadratic syzygies corresponding to the
trivial character. \par\medskip\noindent
For the non-trivial characters the calculations are similar.
Again we multiply the matrix $M_x(y)$ with a column vector with quadratic
entries in the $y$'s, namely, with the column vector for the trivial character
acted on by $\sigma^j\tau^k$.
\smallskip\noindent
We look at the case $(j,k)=(1,2)$, the other cases are similar.
\smallskip\noindent
Now the $i^{\text{th}}$ entry of the product matrix is a cubic in the span of
$$
y_{i+1}^3,\ y_{i}y_{i+1}y_{i+2},\ y_{i+1}y_{i+3}y_{i+4},
\ y_{i+2}^2y_{i+4},
\ y_{i+2}y_{i+3}^2,\ y_{i}^2y_{i+3},\ y_{i}y_{i+4}^2,
$$
$i = \overline{0,4}$, and its coefficient column vector is
$$
T^{\prime}_{(1,2)}(x,b) =
\pmatrix
x_4 b_0\\
x_4 b_2 + (\xi^3 x_0 + \xi^2 x_3) b_1\\
x_4 b_1 + (\xi x_1 + \xi^4 x_2) b_2\\
\xi x_1 b_0 + \xi^2 x_3 b_2\\
\xi^4 x_2 b_1 + \xi^2 x_3 b_0\\
\xi^3 x_0 b_2 + \xi^4 x_2 b_0\\
\xi^3 x_0 b_0 + \xi x_1 b_1
\endpmatrix
$$
multiplied by $\xi ^{i+2}.$
$T^{\prime}_{(1,2)}$ defines a rational map
$$
\Phi_{T^{\prime}_{(1,2)}} : \openP^4(x) \times \openP^2(b) \rto \openP^{6}(s),
$$
well-defined outside the codimension 5 set
$$
B_{\Phi_{T^{\prime}_{(1,2)}}}=(L_{12}(x) \times \openP^2(b)) \cap
\{b_0^2+b_1b_2=\xi x_1b_0-\xi^3 x_0b_2=\xi x_1b_1+\xi^3 x_0b_0=0\},
$$
which is a rational curve projecting bijectively onto the
HM-line
$$
L_{12}=\sigma\tau^2\Pn 1_-
=\{x_4=\xi^3 x_0+\xi^2x_3=\xi x_1+\xi^4 x_2=0\}\subset\Pn 4(x)
$$
in the first factor, and onto the conic $\{b_0^2+b_1b_2=0\}
\subset\openP^2(b)$ in the second factor.\smallskip\noindent
The quadratic syzygies which are induced by the linear ones
appear with summands of type $M_p(y)L_q(y)(y_3,\xi y_4,
\xi^2 y_0,\xi^3 y_1,\xi^4 y_2)^t$, with cubic entries as above.
This time the coefficient vector of the $i^{\text{th}}$ entry is
$$
T^{\prime}_{(1,2)}(x,z) =
\pmatrix
x_4 z_4\\
x_4 (\xi^4 z_0 +\xi z_3) + (\xi^3 x_0 + \xi^2 x_3)(\xi^3 z_1 +\xi^2 z_2)\\
x_4 (\xi^3 z_1 +\xi^2 z_2) + (\xi x_1 + \xi^4 x_2)(\xi^4 z_0 +\xi z_3) \\
\xi x_1 z_4 + \xi^2 x_3 (\xi^4 z_0 +\xi z_3)\\
\xi^4 x_2 (\xi^3 z_1 +\xi^2 z_2) + \xi^2 x_3 z_4\\
\xi^3 x_0 (\xi^4 z_0 +\xi z_3) + \xi^4 x_2 z_4\\
\xi^3 x_0 z_4 + \xi x_1 (\xi^3 z_1 +\xi^2 z_2)
\endpmatrix
$$
multiplied by $\xi ^{i+3}.$ The linear map
$$
b_0\mapsto z_4,\qquad b_1\mapsto (\xi^3 z_1 +\xi^2 z_2) ,\qquad
b_2\mapsto (\xi^4 z_0 +\xi z_3)
$$
gives rise to a commutative diagram similar
to the one for the trivial character. Only this time the map from
$\openP^{14}(t)$ to $\openP^6(s)$ factors through the
projection from the HM-line $L_{12}$, which corresponds to our character.
\par\noindent
Since none of the planes spanned by the $p_i$'s of the cases listed
in (5.1) meets any HM-line except $\openP^1_-$,
the only possible quadratic syzygies of this character come either
from double points, or from trisecants of the images of the three planes
$\{p_i\}\times \openP^2(b)$ in $\Pn 6(s)$,
which are not inherited from $\Pn {14}(t)$.  The calculations run
as above, except that now we have no new trisecants or double points.
For instance, if all the $p_i$'s  lie on an elliptic curve,
then this is because addition on the elliptic curve should now sum up to a
non-trivial 5-torsion point instead of the origin.
In the $\iota$-invariant cases the three planes will span all of $\Pn 6(s)$,
without intersecting each other, and therefore they will have a pencil of
trisecant lines, all of which are inherited from $\Pn {14}(t)$.
\qed\enddemo

\proclaim
{Proposition 5.12}
In the above cases of (5.1) with 10 first order linear syzygies
$(L_{q_{ij}})$ there are no second order quadratic syzygies of $\alpha_1$.
\endproclaim
\demo{Proof} The number of such syzygies, if they
exist at all, is a multiple of 5. More precisely,
from Proposition 5.10, Proposition 5.11 and the double complex
of the Koszul cohomology we may deduce that the linear and quadratic part of
the minimal free resolution of $M$ contains a subcomplex
$$
\vbox{%
\halign{&\hfil\,$#$\,\hfil\cr
V_3\otimes R&\overset {(M_{p_i})}\to\leftarrow&3V_1\otimes R(-1)
&\quad \overset {(L_{q_{ij}},L^{\prime}_i)}\to\quad &2V_0\otimes R(-2)
&&aI\otimes  R(-3)\cr
&&&\vbox to 10pt{\vskip-3pt\hbox{$\nwarrow$}\vss}&\oplus&\leftarrow&\oplus\cr
&&&&(a+b)I\otimes  R(-3)&&cV_2\otimes R(-4)&\cr
}},
$$
\noindent
where $a,b\le 2$ and $L^{\prime}_i = (L^{\prime}_{r_{ik}})$,
$i=\overline{1,3}$, is as in the proof of Proposition 5.11.
Hence 5 second order quadratic syzygies are given
by a relation
$$
L_{q_{i1}}(y)P_1(y)+L_{q_{i2}}(y)P_2(y)+
L^{\prime}_i(y)P_3(y)=0,\quad i=\overline{1,3},\tag 5.13
$$
where $P_1$ and $P_2$ are $5\times 5$-matrices with quadratic entries
and $P_3$ is an $(a+b)\times 5$-matrix with linear entries.
We want to show that c=0 or equivalently that such a relation is trivial
or depends on the second order linear syzygies.
\smallskip\noindent
We consider first the cases (ii)', (iv) and (v).
In these cases $L_{q_{31}}=0$ by Proposition 5.8 and $(a+b) = 1$
by Proposition 5.10. Thus the third of the above matrix equations is just
$$
L_{q_{32}}(y)P_2(y)+L^{\prime}_3(y)P_3(y)=0,\tag 5.14
$$
and $L^{\prime}_3(y)P_3(y)$ has generic rank $\le 1$. Hence
$L_{q_{32}}(y)P_2(y)$ must also have generic rank $\le 1$.
But in all the above cases  $\det L_{q_{32}}$ is a genuine quintic
polynomial, thus $L_{q_{32}}$ has generic rank 5 which means that $P_2(y)$
must have rank $\le 1$. If $P_2(y)$ is non-zero it represents an
element in
$$
\openP(\Hom_{H_5}(V_2\otimes R(-4),V_0\otimes R(-2)))
\cong \openP(\Hom_{H_5}(V_2,V_0\otimes S^2 V_0^*)),
$$
and is thus of type
$$
P_2(y)=b_{i-j}y_{3i+3j}^2+c_{i-j}y_{3i+3j+2}y_{3i+3j+3}+
d_{i-j}y_{3i+3j+1}y_{3i+3j+4},
\quad i,j\in\openZ_5,
$$
where $(b_0: \dots: d_4)$ is a parameter point in $\openP^{14}$. In this
case a check on the
$2\times 2$-minors of $P_2(y)$ shows that $P_2(y)$ factors through
$\chi \otimes R(-3)$ for some character $\chi$ of $H_5$, i.e.,
$P_2(y)$ is a product $P_2(y)=P^{\prime}_2(y)P^{\prime\prime}_2(y)$,
where $P^{\prime}_2(y)=(\xi^{ji}y_{4k+i})_{i\in\openZ_5}^t$ and
$P^{\prime\prime}_2(y)=(\xi^{4ji}y_{k+i})_{i\in\openZ_5}$
for some $(j,k)\in\openZ_5\times\openZ_5$ (of course this follows also
from the factoriality of the polynomial ring). By comparing coefficients
in (5.14) we see that in fact $\chi$ must be the trivial character, i.e.,
$P^{\prime\prime}_2(y)=P_3(y)=[y]=(y_0,\dots,y_4)$ and
$P^{\prime}_2(y)=[y]^t$. \smallskip\noindent
In case (v) $L^{\prime}_3$ is non-zero by the proof of Proposition 5.11.
It follows that $P_2(y)$ and $P_3(y)$ are zero, since otherwise
the first order quadratic syzygies would depend on the linear ones. Indeed,
this can be seen by rewriting (5.14) as
$$
L^{\prime}_{q_{32}}(y)P^{\prime\prime}_2(y)+L^{\prime}_3(y)P_3(y)=0
$$
with $L^{\prime}_{q_{32}}=L_{q_{32}}P^{\prime}_2(y)$. Then also
$P_1(y)=0$, since
$\det L_{q_{11}}$ (and also $\det L_{q_{21}}$) are not vanishing
identically. Notice that in this case there are no second order
linear syzygies.\smallskip\noindent
In the cases (ii)' and (iv) there is precisely one
second order linear syzygy. We may suppose by the proof of Proposition
5.11 that $L^{\prime}_3(y)$ and thus also $P_2(y)$ are zero. Hence (5.13)
yields
$$(\lambda L_{q_{11}}(y)+\mu L_{q_{21}}(y))P_1(y)+
(\lambda L^{\prime}_1(y)+\mu L^{\prime}_2(y))P_3(y)=0$$
for all $(\lambda:\mu)\in\Pn 1$.
$\lambda L_{q_{11}}(y)+\mu L_{q_{21}}(y)$ has generic rank 4
(being skew-symmetric) but nevertheless moving kernel with
respect to $(\lambda:\mu)\in\Pn 1$. Therefore a similar argument
as above allows us to conclude that $P_1(y)=0$ or $P_1(y)=[y]^t[y]$. In
the latter case $(\lambda L_{q_{11}}(y)+\mu L_{q_{21}}(y))[y]^t=0$
implies $P_3(y)=0$, so $(P_1,P_2,P_3)^t=([y]^t[y],0,0)$
depends on the second order linear syzygies. In the first case
$P_1=P_2=P_3=0$.
\par\medskip
\noindent
Case (i) is well known (see \cite{De2}, \cite{Ma3} and \cite{Ma4}) and can
be worked out directly. In case (iii) we may suppose that $L_1^{\prime}=0$
by the proof of Proposition 5.11 (compare the proof of Proposition 8.2
for $L_2^{\prime}$ and $L_3^{\prime}$). Since everything is very
explicit, this case can also be worked out directly.
\qed\enddemo
\head
{6. Syzygy construction of bielliptic surfaces}
\endhead
 In this section we prove the main results of this paper. Namely, we
describe the structure of the Hartshorne-Rao modules and the syzygies
of the minimal bielliptic surfaces in $\Pn 4$ and of the bielliptic surfaces
of degree 15 studied in \cite{ADHPR}. At the same time we explain
how to construct these surfaces via syzygies. We are not able to show directly
that our constructions do lead to smooth surfaces. To overcome this difficulty
we make use of a detailed classification of the $H_5$-invariant surfaces
of small degree lying on either the trisecant variety of an elliptic quintic
scroll in $\Pn 4$, or on the secant variety of an elliptic normal curve in
$\Pn 4$.\medskip
\proclaim{Proposition 6.1} Let $V_{(\lambda:\mu)}=\Trisec Q_{(\lambda:\mu)}
\subset\Pn 4$, $(\lambda:\mu)\in\Pn 1_-\setminus\Lambda$, be the trisecant
variety  of the elliptic scroll $Q_{(\lambda:\mu)}$ and let
$X\subset V_{(\lambda:\mu)}$ be a reduced $H_5$-invariant surface
of degree $d\le 10$ lying on this quintic hypersurface. Then
either\hfill\break
(i) $d=5$ and $X$ is the elliptic scroll $Q_{(\lambda:\mu)}$, or\hfill\break
(ii) $d=10$ and $X$ is either a minimal abelian surface in the
pencil $|A_K|$, or one of the four singular translation scrolls in
this pencil (see \cite{ADHPR} and Remark 4.11), or a minimal
bielliptic surface $B_{(i,j)}$.
\endproclaim
\demo{Proof} Let $\tilde V$ denote the normalization of $V_{(\lambda:\mu)}$.
Recall from \cite{ADHPR, Proposition 24 and Proposition 26} that $\tilde V$
can be identified with a $\Pn 2$-bundle $\Pn 2_E$ over $E$, the elliptic curve
which is the base of the scroll  $Q_{(\lambda:\mu)}$.
The group $G_5$ acts on $V_{(\lambda:\mu)}$ and thus also on $\tilde V$,
and furthermore the projection  $p:\tilde V\to E$ is equivariant under this
action.
The generators of $\Num \tilde V$ are $F$, the class of a fiber of $p$, and
$H$,
the class of ${\scr O}_{\Pn 2_E}(1)$. The normalization map onto
$V_{(\lambda:\mu)}\subset\Pn 4$ is induced by  ${\scr O}_{\Pn 2_E}(1)$,
so we compute that $H^3=5$, $H^2\cdot F=1$ and $H\cdot F^2=F^3=0$.
On the other side, the canonical class $K$ on $\tilde V$ is numerically
equivalent
to $-3H+5F$, and in fact $|-K|$ corresponds
to the pencil $|A_K|$ on $V_{(\lambda:\mu)}$ (see \cite{ADHPR, Theorem 30}).
\smallskip\noindent
Assume first that $X$ is not contained inside the non-normal locus
of $V_{(\lambda:\mu)}$, and let $X_0$ denote the strict transform
of $X$ on $\tilde V$. By Heisenberg invariance the degree of $X$ is divisible
by 5,
thus we need to find the integral solutions of
$X_0\equiv aH+bF$, with $H^2\cdot X_0=5m$, $m=1,2$. The coefficients a and
b satisfy the inequalities
$$
a=H\cdot X_0\cdot F\ge 0, \quad{\text {and}}
\quad 10a+3b= H\cdot X_0\cdot(3H-5F)\ge 0.
$$
Moreover, $X_0$ must intersect the pullback on $\tilde V$ of any curve
in the antibicanonical pencil on $Q_{(\lambda :\mu)}$ in a scheme whose
length is a multiple of 25, since these are all $H_5$-invariant.  If $C_K$
is the base curve $K^2$ of the pencil of anticanonical surfaces on $\tilde V$,
then
the numerical class of the curves in the antibicanonical pencil on
$Q_{(\lambda :\mu)}$ is $5H-C_K$. So their pullback to $\tilde V$ is equivalent
to $(3H-5F)\cdot (5H-(3H-5F))=(3H-5F)\cdot (2H+5F)$.  Thus we get that
$35a+6b$ is a multiple of 25.  Combining this with the above
$H^2\cdot X_0=5a+b=5m$, for $m=1,2$, we get that
$$
10a+b \quad {\text{is a multiple of}}\quad 25.
$$
For $m=1$ there are no numerical solutions. In case $m=2$
we are left with $(a,b)=(3,-5)$. This is the numerical class of
the anticanonical bundle on $\tilde V$ and from \cite{ADHPR, Lemma 3, Lemma 4
and
Theorem 30} it follows that the only reduced surfaces of degree 10 in this
class are
those listed at (ii) in the statement of the lemma. \smallskip\noindent
At last (i) follows from the fact that the scroll $Q_{(\lambda:\mu)}$
coincides with the non-normal locus of $V_{(\lambda:\mu)}$.\qed
\enddemo
\medskip
\proclaim{Proposition 6.2} Let $W_{(\lambda:\mu)}=\Sec E_{(\lambda:\mu)}
\subset\Pn 4$, $(\lambda:\mu)\in\Pn 1_-\setminus\Lambda$,
be the secant variety of the elliptic curve
$E_{(\lambda:\mu)}\subset\Pn 4$, and let $X\subset W_{(\lambda:\mu)}$ be a
reduced $H_5$-invariant surface of degree $d\le 15$ lying on this quintic
hypersurface. Then either\hfill\break
(i) $d=5$ and $X$ is one of the 3 translation scrolls
$\Sigma_{(\lambda:\mu)}(\tau_i)$,
where $\tau_i\in E_{(\lambda:\mu)}$ is a non-trivial 2-torsion point,
or \hfill\break
(ii) $d=10$ and $X$ is either a singular translation scroll (see \cite{BHM2}
and \cite{ADHPR, proof of Proposition 36}), or the union
$\Sigma_{(\lambda:\mu)}(\tau_1)\cup\Sigma_{(\lambda:\mu)}(\tau_2)$,
where $\tau_i\in E_{(\lambda:\mu)}$, $i=1,2$,  are distinct
non-trivial 2-torsion points, or\hfill\break
(iii) $d=15$ and $X$ is either a non-minimal abelian surface
$\tilde A_K$ (see \cite{ADHPR, Theorem 32}), or a non-minimal bielliptic
surface $\tilde B_{(i,j)}$ (see again \cite{ADHPR}), or the union of a
degenerate abelian surface of degree 10 and an elliptic quintic scroll
$\Sigma_{(\lambda:\mu)}(\tau_i)$, or the singular scroll $S_{(\lambda:\mu)}$
(see Proposition 4.16), or the union
$\Sigma_{(\lambda:\mu)}(\tau_1)\cup\Sigma_{(\lambda:\mu)}(\tau_2)\cup\Sigma_
{(\lambda:\mu)}(\tau_3)$, where $\tau_i\in E_{(\lambda:\mu)}$,
$i=\overline{1,3}$, are the non-trivial 2-torsion points. Among these surfaces
the bielliptic are the only ones that are not also $\iota$-invariant.
\endproclaim
\demo{Proof}
Let $\tilde W$ be the blow-up of $W_{(\lambda:\mu)}$ along the curve
$E_{(\lambda:\mu)}$.  Then $\tilde W$ has the structure of a $\Pn 1$-bundle
over the symmetric product $\Sigma$ of $E_{(\lambda:\mu)}$ (see
\cite{ADHPR, proof of Proposition 36} for more details).  Let $X_0$ be the
strict transform of $X$ on $\tilde W$. Now  $G_5$ acts on $\tilde W$ and both
projections $p:\tilde W\to\Sigma$ and $q:\tilde W\to E_{(\lambda:\mu)}$
are equivariant under this action. Let $F$ denote the class of a fiber of $q$.
Such a fiber is mapped to a rational cubic scroll in $\openP^4$. Let $B$ denote
the scroll in $\tilde W$ over a section of $\Sigma$ with self-intersection 1.
The image of $B$ on $W_{(\lambda:\mu)}$ is the cone of
secants to $E_{(\lambda:\mu)}$ through a
given point on it. In particular, it is a cone over an elliptic curve of
degree 4.  Therefore $B$ has a unique section with self-intersection $-4$
which is contracted to a point on $E_{(\lambda:\mu)}$.  Let
$S$ be the section in $\tilde W$ over $\Sigma$ which consists of the
directrices
in each of the $F$'s. The image of $S$ in $\Sec E_{(\lambda:\mu)}$
is the scroll
$S_{(\lambda:\mu)}$.\par\smallskip
\noindent
If $H$ on $\tilde W$ is the pullback of the class of a hyperplane section of
$W_{(\lambda:\mu)}$, then $H\equiv S+2B-6F$ in $\Num \tilde W$. Furthermore, if
$K$ is the class of the
dualizing sheaf of $\tilde W$ and $K_\Sigma$ the pullback of the
canonical class on $\Sigma$, then $K\equiv -2S-3B+14F$ and
$K_\Sigma\equiv -2B+F$.  The exceptional surface $Y$ on $\tilde W$ is
anticanonical and has the structure of a product $E_{(\lambda:\mu)}\times
E_{(\lambda:\mu)}$ (see \cite{ADHPR, p. 886}). In particular,
every curve on it has non-negative selfintersection.\smallskip\noindent
The intersection numbers on $\tilde W$ are:
$H^2\cdot F=3,\quad H^2\cdot S=15,\quad H^2\cdot B=4$. Moreover,
$B^2\cdot F=B\cdot F^2=H\cdot F^2=F^3=B^3=0$,
while $H\cdot B\cdot F=H\cdot B^2=1$,
and $B\cdot S^2=12,\quad B\cdot F\cdot S
=-F\cdot S^2=B^2\cdot S=1,\quad S^3=-25$.
We are looking for integral solutions of
$X_0\equiv aF+bB+cS$, $H^2\cdot X_0=5m$, $m=\overline{1,3}$,
with coefficients satisfying the inequalities
$$c=X_0\cdot B^2\geq 0,\qquad b-c=
X_0\cdot S\cdot F\geq 0\qquad {\text {or}}\quad
a=b=0,\ c=1,$$
and
$$a+b+4c= X_0\cdot B\cdot (S+2B-10F)\geq 0,$$
\noindent
where the last inequality comes from
the fact that $B\cdot (H-4F)= B\cdot (S+2B-10F)$ is the section with
selfintersection $-4$ on any surface in the class $B$. Furthermore,
$$
2ab+ac+2b+13bc=X_0^2\cdot (2S+3B-14F)\geq 0
$$
since every curve on the exceptional surface has
non-negative selfintersection.  Finally,
$$
2a+b\quad {\text{is a multiple of}}\quad 25.
$$
Indeed, consider the 3 quintic
elliptic scrolls which are the 2-torsion translation scrolls of
$E_{(\lambda:\mu)}$. The pullback of any antibicanonical curve
on one of these scrolls belong to the class
$(2B-F)\cdot (4S+2B-26F)$, and it intersects $X_0$ in a
multiple of 25 points, unless $X_0$ is the union of the three scrolls.
\smallskip\noindent
For $m=1$ only the triple $(a,b,c)=(-1, 2, 0)$ satisfies our conditions.
The surfaces $X_0\equiv -F+2B$ are  the pullbacks of numerical
anticanonical curves on $\Sigma$. It follows from \cite{CC} that
$X$ is necessarily one of the three elliptic quintic
scrolls $\Sigma_{(\lambda:\mu)}(\tau_i)\subset W_{(\lambda:\mu)}$, where
$\tau_i\in E_{(\lambda:\mu)}$ is a non-trivial 2-torsion point.
\smallskip\noindent
For $m=2$ we are left with $(a,b,c)=(-2,4,0)$. Invariant divisors in the
corresponding numerical class are either singular translation scrolls of
degree 10 (see \cite{HM}, or \cite{Hu1, Theorem VII.3.1}), or the union of
two elliptic quintic scrolls $\Sigma_{(\lambda:\mu)}(\tau_1)\cup
\Sigma_{(\lambda:\mu)}(\tau_2)$.\smallskip\noindent
Finally, for $m=3$ the numerical solutions are
$$
(a,b,c)\in\{(-3,6,0), (0,0,1), (-14,3,3)\}.
$$
A surface in the first numerical class
is a  pullback of an antitricanonical curve on $\Sigma$,
and so it is an $\iota$-invariant scheme (``translation scroll'') of
degree 15. Such a surface is either
the union of a degenerate abelian surface of degree 10
(see \cite{BHM2}) and a quintic elliptic scroll
$\Sigma_{(\lambda:\mu)}(\tau_i)$, or
the union of all the three quintic elliptic scrolls
$\Sigma_{(\lambda:\mu)}(\tau_i)$.  The second class corresponds to
the surface $S_{(\lambda:\mu)}$ of degree 15,  while the last
solution corresponds to the pullback of numerically anticanonical surfaces
on the variety $\tilde V$ of the proof of Proposition 6.1 via the blow-up map
$\tilde W\to \tilde V$ with exceptional divisor $S$ (see also \cite{ADHPR,
Proposition
28 and 36}. Here the bielliptic surfaces are the only ones
that are not also $\iota$-invariant.\qed
\enddemo
\proclaim{Theorem 6.3} Let $M=\underset{k\ge 0}\to\bigoplus M_k$
be the $H_5$-invariant graded $R$-module presented by
$$\alpha_1 = (M_{\tau_1}(y), M_{\tau_2}(y), M_{\tau_3+\rho}(y)),$$
where $\tau_i$, $i=\overline{1,3}$, are the non-trivial 2-torsion
points and $\rho$ is a non-trivial 3-torsion point on the elliptic
normal curve $E_{(\lambda:\mu)}\subset\Pn 4$,
$(\lambda:\mu)\in\Pn 1_-\setminus\Lambda$. Then $Syz_1(M(1))$ is a
vector bundle of rank 10 and there exists a sheaf monomorphism
$$
\varphi: 5{\scr O}(-1)\oplus\Omega^3(3)\to\scr Syz_1(M(1)).
$$
Such a morphism drops rank along a smooth minimal bielliptic surface $B$
of degree 10, and it is uniquely determined up to isomorphisms, hence $B$
is uniquely determined by $M$. Moreover, $B$ is contained in precisely
one quintic hypersurface. This is the trisecant variety
$\Trisec Q=\{\det L_q(y)=0\}$ defined by a point $q$ on the numerically
anticanonical
curve $E_{(\lambda:\mu)}(\tau_3)\subset Q_{(\lambda:\mu)}$ not on a HM-line
(compare 4.3, 4.12, 4.20 and 4.21).
Conversely, any smooth minimal bielliptic surface $B\subset\Pn 4$
can be obtained in this way (up to coordinate transformations).
In particular, there are 8 different bielliptic surfaces on
$\Trisec Q$, corresponding to the eight different
choices of a non-trivial 3-torsion point on $E_{(\lambda:\mu)}$.
\endproclaim
\demo{Proof}
First of all $M$ is artinian by Proposition 5.5, whence
$\scr Syz_1(M(1))$ is a vector bundle. Now $\alpha_1$ has exactly
ten first order linear syzygies by Proposition 5.8, one first order
quadratic syzygy by Proposition 5.11, one second order
linear syzygy by Proposition 5.10 and no second order quadratic syzygies
by Proposition 5.12. In particular, $M$ has Hilbert function $(5,10,10)$
and a minimal free resolution of type
$$
\vbox{%
\halign{&\hfil\,$#$\,\hfil\cr
0\leftarrow M\leftarrow 5R\leftarrow 15R(-1)&&10R(-2)&&R(-3)\ \cr
&&\oplus\cr
&\vbox to 10pt{\vskip-4pt\hbox{$\nwarrow$}\vss}\ &R(-3)&\leftarrow &\oplus\cr
&&\oplus\cr
&&25 R(-4)&&55 R(-5)
&\vbox to 10pt{\vskip-12pt\hbox{$\nwarrow$}\vss}\ &40 R(-6)
\leftarrow 10 R(-7)\leftarrow 0\cr
}}.
$$
It follows readily that there exists, up to isomorphisms, a uniquely
determined sheaf mono-\par\noindent morphism $\varphi\in
\Hom(5{\scr O}(-1)\oplus\Omega^3(3),\scr Syz_1(M(1)))$, and that
this morphism drops rank in the expected codimension 2. Thus the
maximal minors of $\varphi$ cut out a locally Cohen-Macaulay surface
$B$ whose ideal sheaf $\scr J_B$ has a resolution
$$
0\rightarrow 5\scr O(-1)\oplus\Omega^3(3)\overset\varphi\to\rightarrow
\scr Syz_1 \left(M(1)\right)\rightarrow\scr
J_B(3)\rightarrow 0.
$$
By Porteus' formula $\deg B=c_2(Syz_1\left(M(1)\right))-
c_2(5{\scr O}(-1)\oplus\Omega^3(3))=10$, while, in case $B$ is reduced
and irreducible, Riemann-Roch on $\Pn 3$ yields $\pi(B)=6$ for the
arithmetic sectional genus. On the other side, $\chi({\scr O}_B)=0$ and
$h^1({\scr O}_B)=h^3(\Omega^3)=1$, thus $B$ has the same numerical
invariants as a minimal bielliptic surface in $\Pn 4$. In order to see
that $B$ is indeed a smooth bielliptic surface we will use Proposition 6.1.
\smallskip\noindent
We will first show that $B$ is contained in $\Trisec Q$.
{}From the mapping cone
$$
\diagram[width=4em,height=3em,PS]
{10{\scr O}\oplus 5{\scr O}(-1)} & \lInto &
{{\Omega^3(3)\oplus 5{\scr O}(-1)}} &\lOnto
& {5{\scr O}(-1)\oplus 5{\scr O}(-1)} &\lTo& {{\scr O}(-2)}\\
\dDashto^{\tilde{\varphi}} &  &\dTo^{\varphi} &\ldTo(4,2)_{\tilde
\varphi_0} \ldTo&\dTo^{\varphi_0}&&\dTo^{\varphi_1}\\
15{\scr O} & \lInto & {\scr Syz_1(M(1))} & \lOnto &
10{\scr O}(-1)\oplus{\scr O}(-2)\oplus 25{\scr O}(-3)&
\lTo&{\scr O}(-2)\oplus\dots\\
\enddiagram
$$
between the minimal free resolutions of $5\scr O(-1)\oplus
\Omega^3(3)$ and $\scr Syz_1\left(M(1)\right)$ we see that $B$ is
contained in precisely one quintic hypersurface in $\Pn 4$. In fact
we can say more. The composition $\tilde\varphi_0:5{\scr O}(-1)\oplus
5{\scr O}(-1)\to15{\scr O}$ in the above diagram
is given by the first order linear syzygies of $\alpha_1$. Since
$h^1(\Omega^2(3))=0$ we may lift $\tilde\varphi_0$ to a morphism
$\tilde\varphi\in\Hom(10{\scr O}\oplus 5{\scr O}(-1), 15{\scr O})$ whose
constant part is of maximal rank. After moding out constants in both sides
we are left with a $5\times 5$-matrix $\tilde L$ with linear entries whose
determinant is exactly the equation of the unique quintic containing $B$.
We may suppose by Proposition 5.8, (iv) that the first order linear syzygies
of $\alpha_1$ are given by a block of syzygy matrices
$$
\left( \matrix
L_{q_{11}}(y) & L_{q_{12}}(y)\\
L_{q_{21}}(y) & L_{q_{22}}(y)\\
0 & L_{q_{32}}(y)
\endmatrix \right)
$$
with $q_{i1}\in\Pn 1_-$, $i=1,2$, and $q_{32}\in E_{(\lambda:\mu)}(\tau_3)$
not on a HM-line. Moreover, by Proposition 5.10, $\alpha_1$ has only one second
order linear syzygy and this is of type $([y]^t, 0)$.  We therefore may suppose
that our matrix $\tilde L$ coincides with $L_{q_{32}}(y)$, whence $B$ lies on
the quintic hypersurface $\Trisec Q=\{\det L_{q_{32}}(y)=0\}$
(compare Proposition 4.12, (i)).
\smallskip\noindent
We will show next that $B$ is reduced. Assume the contrary and let
$B_{\text{red}}$ denote the reduced support of $B$.
Then $B_{\text{red}}\subset\Trisec Q$ is a $H_5$-invariant
surface of degree 5, hence $B_{\text{red}}$ coincides with $Q$
by Proposition 6.1, and $B$ is a locally Cohen-Macaulay
double structure on $Q$. It follows that $B$ is necessarily a local complete
intersection (see e.g. \cite{Ma2, proof of Theorem 1}), and that it
comes through a Ferrand doubling \cite{Fe}. In other words, there
is an exact sequence
$$
0\to\scr J_B/\scr J_Q^2\rightarrow\scr J_Q/\scr J_Q^2 \rightarrow \scr L
\to 0,\tag 6.4
$$
where $\scr L$ is a line bundle on $Q$.
Since $\scr L=\scr J_Q/\scr J_B$ we also have an exact sequence
$$
0\to\scr L \rightarrow\scr O_B \rightarrow\scr O_Q\to 0.      \tag 6.5
$$
Let $C_0$ and $f$ be the generators of $\Num Q$, and write
$\scr L\equiv a C_0 +b f$. If follows from (6.4) that the twisted normal bundle
$\scr N_Q\otimes \scr L$ of $Q$ in $\openP^4$ has a nowhere vanishing
section. In particular,
$$
0=c_2(\scr N_Q\otimes \scr L)=25 +14a +3b +a^2 +2ab.
$$
On the other hand, $\chi(\scr L)= \chi(\scr O_B) -\chi(\scr O_Q)=0$ by (6.5)
since $B$ has the numerical type of a bielliptic surface. Thus
Riemann-Roch implies that
$$
\scr L^2=K_{Q}\cdot \scr L, \quad {\text{i.e.}},\quad a^2+2ab=-a-2b,
$$
where $K_{Q}$ denotes the canonical class on $Q$. The last equation
gives that either $a=-2b$, or $a=-1$. In conclusion,
there are exactly 2 solutions for these equations, namely
$$
\scr L\equiv -2C_0+f \quad {\text{and}}\quad  \scr L\equiv -C_0-12f.
$$
Since we already know all the line bundles of numerical type
$\scr L\equiv -2C_0+f$ it is easy to deduce (cf. the proof of
\cite{HV, Proposition 4}) that in the first case
the only such doublings are the HM-surfaces in \cite{HV}.
This is a contradiction, since the HM-surfaces are of the numerical type of
an abelian surface and thus in particular have different cohomology.
The second case can be ruled out by taking global sections in the exact
sequence
$$
0\to \scr J_B/\scr J_Q^2(5)\rightarrow\scr J_Q/\scr J_Q^2 (5)\rightarrow
\scr L(5)\to 0.
$$
$\scr J_B/\scr J_Q^2(5)$ has no sections since both $B$ and the first
infinitesimal neighborhood of the scroll $Q$ are contained in the
same unique quintic, while $\scr L(5)$ has at most 2 sections
since $\scr L(5)\equiv -2K$. But this is a contradiction since
$H^0\scr J_Q/\scr J_Q^2 (5)=50.$
\smallskip\noindent
In conclusion, $B$ is a reduced $H_5$-invariant surface of degree 10 on
$\Trisec Q$ which is not $\iota$-invariant. Thus it must
be a smooth minimal bielliptic surface by Proposition 6.1.\smallskip\noindent
The converse statement is clear from \cite{ADHPR} and Proposition 3.7.
\qed\enddemo
\proclaim{Theorem 6.6} Let $M^{\prime}=\underset{k\ge 0}\to\bigoplus
M^{\prime}_k$ be the $H_5$-invariant graded $R$-module presented by
$$\alpha_1 = (M_{\tau_1+\rho}(y), M_{\tau_2+\rho}(y), M_{\tau_3+\rho}(y)),$$
where $\tau_i$, $i=\overline{1,3}$, are the non-trivial 2-torsion
points and $\rho$ is a non-trivial 3-torsion point on the elliptic
normal curve $E_{(\lambda:\mu)}\subset\Pn 4$,
$(\lambda:\mu)\in\Pn 1_-\setminus\Lambda$. Then:\par\noindent
(i)\quad $M^{\prime}$ has Hilbert function $(5,10,10,1)$
and a minimal free resolution of type
$$
\vbox{%
\halign{&\hfil\,$#$\,\hfil\cr
0\leftarrow M^{\prime}\leftarrow 5R\leftarrow 15R(-1)&&10R(-2)&&\ \cr
&&\oplus\cr
&\vbox to 10pt{\vskip-4pt\hbox{$\nwarrow$}\vss}\ &R(-3)&\cr
&&\oplus\cr
&&20 R(-4)&\vbox to 10pt{\vskip-12pt\hbox{$\nwarrow$}\vss}\ &45 R(-5)
&\leftarrow&30 R(-6)& &5 R(-7)\cr
&&&&&&&\vbox to 10pt{\vskip-5pt\hbox{$\nwarrow$}\vss}\ &\oplus\cr
&&&&&&&&R(-8)&\vbox to 10pt{\vskip-4pt\hbox{$\nwarrow$}\vss}\ 0\cr
}}.$$
(ii)\quad
$Syz_2\left((M^{\prime})^{\vee}(3)\right)$ is a
vector bundle of rank 21 and there exists a sheaf monomorphism
$$\varphi: 20{\scr O}(-1)\to\scr Syz_2\left((M^{\prime})^{\vee}(3)\right).$$
Such a morphism drops rank along a smooth non-minimal
bielliptic surface $B^{\prime}$ of degree 15, and it is uniquely determined
up to isomorphisms, hence $B^{\prime}$ is uniquely determined by $M^{\prime}$.
Moreover, $B^{\prime}$ is contained in precisely one quintic hypersurface,
namely $\Sec E_{(\lambda:\mu)}$, and the ideal sheaf ${\scr I}_{B'}$ has a
minimal free resolution of type
$$
\vbox{%
\halign{&\hfil\,$#$\,\hfil\cr
&&&&\scr O(-5)\cr
0&\leftarrow&\scr J_{B'}&\leftarrow&\oplus\cr
&&&&10\scr O(-6)&\vbox to 10pt{\vskip-4pt\hbox{$\nwarrow$}\vss}
&15\scr O(-7)&\leftarrow&5\scr O(-8)&\leftarrow 0.\cr
}}
$$\endproclaim
\demo{Proof}
As in the proof of Theorem 6.3 we deduce that $M^{\prime}$ is artinian,
that it has Hilbert function $(5,10,10,1)$ (use Gotzmann's persistence
theorem \cite{Go}), and a minimal free resolution as stated above.
In particular, $\scr Syz_2((M^{\prime})^{\vee}(3))$ is a rank 21 vector
bundle, and it follows that there exists, up to isomorphisms,
a uniquely determined sheaf monomorphism $\varphi\in\Hom(20{\scr O}(-1),
{\scr Syz_2((M^{\prime})^{\vee}(3))})$, whose degeneracy locus is a locally
Cohen-Macaulay scheme $B^{\prime}$ of codimension 2. From the exact sequence
$$
0\rightarrow 20\scr O(-1)\overset\varphi\to\rightarrow
{\scr Syz_2((M^{\prime})^{\vee}(3))}\rightarrow
\scr J_{B^{\prime}}(3)\rightarrow 0
$$
we deduce as in the proof of Theorem 6.3 that ${\scr I}_{B'}$ has the
syzygies as claimed, that $\deg B^{\prime}=15$, that, in case $B^{\prime}$
is reduced and irreducible, $\pi(B^{\prime})=21$, and moreover that
$\chi({\scr O}_{B^{\prime}})=0$ and $h^1({\scr O}_{B^{\prime}})=1$.
In order to see that $B^{\prime}$ is indeed a smooth non-minimal bielliptic
surface we will use Proposition 6.2. \smallskip\noindent
We will first show that $B^{\prime}$ is contained in $\Sec
E_{(\lambda:\mu)}$. We may suppose by Proposition 5.8, (v) that the first order
linear syzygies of $\alpha_1$ are given by a block of syzygy matrices
$$
\left( \matrix
L_{q_{1}}(y) & L_{q_{1}}(y)\\
L_{q_{2}}(y) & 0\\
0 & L_{q_{3}}(y)
\endmatrix \right),
$$
where the common determinant $\gamma$ of all 4 non-zero matrices $L_{q_i}(y)$
defines $\Sec E_{(\lambda :\mu)}$. The first order linear syzygies together
with the first order quadratic syzygy of $\alpha_1$ form a $15\times11$-matrix
whose transpose is the presentation matrix of the ideal $I_{B^{\prime}}$ of
$B^{\prime}$. It follows that $I_{B^{\prime}}$ contains all the
$10\times 10$-minors of this $15\times11$-matrix. In particular,
$I_{B^{\prime}}$
contains ${\text {det}}L_{q_2}\cdot {\text {det}} L_{q_3}$, which
is the  square of $\gamma$.  But we have more: $I_{B^{\prime}}$
contains the product of $\gamma$ by any $5\times 5$-minor
of any $10\times 5$-matrix formed by two of the three blocks $L_{q_i}$.
That the unique quintic containing  $B^{\prime}$ is $\Sec E_{(\lambda :\mu)}
=\{\det\gamma=0\}$ follows therefore immediately from the\smallskip\noindent
\proclaim{Claim} The $5\times 5$-minors of all three $10\times 5$-matrices
$$
\left( \matrix
L_{q_{i}}(y) \\
L_{q_{j}}(y)
\endmatrix \right)
$$
vanish at most on a curve.\endproclaim
\demo{Proof of Claim}
Let $y$ be a point where all these minors vanish.  We may
assume that $y$ lies on $\Sec E_{(\lambda:\mu)}=\{\gamma =0\}$ outside the
curve $E_{(\lambda:\mu)}$, thus we may assume that $L_{q_1}(y)$ has rank 4.
We argue as in the proofs of Lemmas 5.3, 5.4 and Proposition 5.5.
In this case let
$$
L(z,y)(\lambda_0,\dots,\lambda_4)^t=0
$$
be the unique relation between the columns of $L(z,y)$ for
$z=q_1$.  The rank condition on the $10\times 5$-matrices now says that
$$
L(z,y)(\lambda_0,\dots,\lambda_4)^t=0
$$
for $z=q_i$, $i=\overline{1,3}$.  Thus the linear span of the three
points $q_i$ is contained in $\det L(z,y)=0$ for the given point $y$.   Now
this linear span is a plane:
The points  $q_i=\rho+\tau_i$ on $E^{\prime}_{(\lambda:\mu)}$ may as points
on $Q_{(\lambda:\mu)}$ be identified with the pair of points
$(\rho,\tau_i)$ when we identify $Q_{(\lambda:\mu)}$  with the symmetric
product of $E_{(\lambda:\mu)}$. Therefore they lie on the same plane cubic
curve on the quintic scroll, namely the curve of points with one coordinate
equal to $\rho$. It follows from the argument at the end of Section 2 in
[ADHPR] that the three points on this curve are collinear iff the sum
$\tau_1+\tau_2+\tau_3=-\rho$.   Therefore the $q_i$ span the plane of the
cubic curve.  This is a plane in the trisecant scroll of
$Q_{(\lambda:\mu)}$, so in fact  this plane must lie in the quintic
$\det L(z,y)=0$.  Thus the 50 $G_5$-translates of this plane must lie in this
quintic also, so  by Bezout this quintic has to be the trisecant quintic.
But the map $\bar\Theta$ from Section 4 contracts only a curve to the point
corresponding to this quintic and the claim follows.
\qed\enddemo\smallskip\noindent
We will show next that $B^{\prime}$ is reduced. Assume the contrary, and let
$B^{\prime}_{\text{red}}\subset \Sec E_{(\lambda:\mu)}$ be the reduced
support of $B$. To derive a contradiction we take
cohomology in the exact sequence
$$
0\rightarrow \scr J\rightarrow
{\scr O}_{B^{\prime}}\rightarrow
{\scr O}_{B^{\prime}_{\text{red}}}\rightarrow 0.
$$
Since $h^2({\scr O}_{B^{\prime}})=0$ it follows that
$h^2({\scr O}_{B^{\prime}_{\text{red}}})=0$.
In fact, the same argument
shows that $h^2({\scr O}_{Y})=0$ for any subscheme $Y\subset B^{\prime}$.
We now consider two cases. In any case we can apply Proposition 6.2 since
$B^{\prime}_{\text{red}}$ is $H_5$-invariant.
\smallskip\noindent
If $\deg B^{\prime}_{\text {red}}=10$, then $B^{\prime}_{\text {red}}$ is
of type $B^{\prime}_{\text {red}}=\Sigma_{(\lambda:\mu)}(\tau_1)\cup
\Sigma_{(\lambda:\mu)}(\tau_2)$.
So $B^{\prime}$ is a locally Cohen-Macaulay surface which is the union of a
double structure $X$ on one scroll, say $\Sigma_{(\lambda:\mu)}(\tau_1)$,
and the other scroll $\Sigma_{(\lambda:\mu)}(\tau_2)$. On the other side,
$\Sec E_{(\lambda:\mu)}$ contains the unique HM-double structure on
$\Sigma_{(\lambda:\mu)}(\tau_1)$. As a double structure contained
in $\Sec E_{(\lambda:\mu)}$ this is uniquely determined outside the curve
$E_{(\lambda :\mu)}$ since the quintic is smooth outside this curve.
Since the HM-double structure is locally Cohen-Macaulay,
$X$ and therefore $B^{\prime}$  must contain it as a subscheme.
This contradicts the fact that $h^2({\scr O}_{Y})=0$ for
any subscheme $Y$ of $B^{\prime}$.\smallskip\noindent
If $\deg B^{\prime}_{\text {red}}=5$, then $B^{\prime}$ is a locally
Cohen-Macaulay triple structure on a scroll $\Sigma_{(\lambda:\mu)}(\tau_i)$.
Again, such a triple structure is unique outside the curve $E_{(\lambda :\mu)}$
since the quintic is smooth outside this curve. Therefore $B^{\prime}$ must
contain the HM-double structure on the scroll, a contradiction as
above.\smallskip\noindent
In conclusion, $B^{\prime}$ is a $H_5$-invariant, reduced surface of degree
15 on
$\Sec E_{(\lambda:\mu)}$ which is not $\iota$-invariant. Thus it is a
smooth non-minimal bielliptic surface by Proposition 6.2.
\qed\enddemo
At this point we want to mention one further type of $G_5$-invariant
surfaces in $\Pn 4$ with $d=10$ and $\pi=6$. Namely, as we will see below,
there exists a one-dimensional family of surfaces $B^{\prime\prime}$ which
are a union of two $G_5$-invariant elliptic quintic scrolls
meeting exactly along a $G_5$-invariant elliptic normal curve.
The topological realization of the dual graph of a type III degeneration
of an elliptic scroll is a M\"obius band, hence the topological realization
of the dual graph of a type III degeneration of two elliptic scrolls
meeting exactly along an elliptic curve is a Klein bottle, i.e., two
M\"obius bands glued along their boundaries. This already indicates
that the surfaces  $B^{\prime\prime}$ should be similar, at least in terms
of numerical invariants, to a minimal bielliptic surface in $\Pn 4$.\par
\noindent
On the other hand, the canonical bundle of the union of two scrolls is a
2-torsion bundle, whereas the canonical bundle of the bielliptic surfaces in
$\openP^4$ is 3-torsion. This already indicates that the union of two scrolls
is
not a degeneration of the smooth bielliptic surfaces. In fact, we shall see in
the appendix that the union of two scrolls corresponds to a smooth point in
the Hilbert scheme and that every (small) deformation is again a union of two
scrolls. On the other
side, D. Morrison gives in \cite{Mor1}, \cite{Mor2} a complete description of
all weakly projective semistable degenerations of bielliptic surfaces with
3-torsion canonical class. Unfortunately, none of the examples of
$H_5$-invariant degenerations to be discussed in Section 9 falls in his
list.
\proclaim{Theorem 6.7}  Let $M^{\prime\prime}=\underset{k\ge 0}\to\bigoplus
M^{\prime\prime}_k$ be the $H_5$-invariant graded $R$-module presented by
$$\alpha_1 = (M_{\tau_1}(y), M_{\tau_2}(y), M_{o_{E_{(\lambda:\mu)}}}(y)),$$
where $\tau_i$, $i=\overline{1,2}$, are two non-trivial 2-torsion
points and $o_{E_{(\lambda:\mu)}}$ is the origin on the elliptic
normal curve $E_{(\lambda:\mu)}\subset\Pn 4$,
$(\lambda:\mu)\in\Pn 1_-\setminus\Lambda$.
Then:\par\noindent
(i)\quad $M^{\prime\prime}$ has Hilbert function $(5,10,10)$
and a minimal free resolution of type
$$
\vbox{%
\halign{&\hfil\,$#$\,\hfil\cr
0\leftarrow M^{\prime\prime}\leftarrow 5R\leftarrow 15R(-1)&&10R(-2)&&R(-3)\
\cr
&&\oplus\cr
&\vbox to 10pt{\vskip-4pt\hbox{$\nwarrow$}\vss}\ &R(-3)&\leftarrow &\oplus\cr
&&\oplus\cr
&&25 R(-4)&&55 R(-5)
&\vbox to 10pt{\vskip-12pt\hbox{$\nwarrow$}\vss}\ &40 R(-6)
\leftarrow 10 R(-7)\leftarrow 0\cr
}}.
$$
(ii)\quad $Syz_1(M^{\prime\prime}(1))$ is a vector bundle of rank 10 and there
exists a sheaf monomorphism
$$
\varphi: 5{\scr O}(-1)\oplus\Omega^3(3)\to\scr Syz_1(M^{\prime\prime}(1)).
$$
Such a morphism drops rank along the surface $B^{\prime\prime}=\Sigma_
{(\lambda:\mu)}(\tau_1)\cup\Sigma_{(\lambda:\mu)}(\tau_2)$ of
degree 10, and it is
uniquely determined up to isomorphisms, hence $B^{\prime\prime}$ is
uniquely determined by $M^{\prime\prime}$. Moreover, $B^{\prime\prime}$ is
contained in precisely one quintic hypersurface, namely
$\Sec E_{(\lambda:\mu)}$, and the ideal sheaf ${\scr I}_
{B^{\prime\prime}}$ has a minimal free resolution of type
$$
\vbox{%
\halign{&\hfil\,$#$\,\hfil\cr
&&&&\scr O(-5)\cr
0&\leftarrow&\scr J_{B^{\prime\prime}}&\leftarrow&\oplus\cr
&&&&25\scr O(-6)&\vbox to 10pt{\vskip-4pt\hbox{$\nwarrow$}\vss}
&55\scr O(-7)&\leftarrow&40\scr O(-8)
&\leftarrow&10\scr O(-9)\leftarrow 0.\cr
}}
$$
Conversely, any union $B^{\prime\prime}\subset\Pn 4$ of two elliptic quintic
scrolls
meeting exactly along an elliptic normal curve is obtained in this way (up to
coordinate transformations).
\endproclaim
\demo{Proof} The proof of (i) and (ii) is similar to the one of Theorem
6.3. Details are left to the reader. For the converse statement let
$B^{\prime\prime}=\Sigma_{(\lambda:\mu)}(\tau_1)\cup\Sigma_{(\lambda:\mu)}
(\tau_2)$ be the union of the $\tau_i$-translation scrolls of
$E_{(\lambda:\mu)}$, $i=1,2$. By comparing the minimal free resolutions in
Lemma 4.4, (i) and (ii) we deduce easily the existence of a commutative
diagram:
\smallskip\noindent
$$
\diagram[width=4em,height=3em,PS]
{\scr O}&\lTo^{I_{E_{(\lambda:\mu)}}}&
V_3\otimes\O {-2}&\lTo^{M(y,o_{E_{(\lambda:\mu)}})}& V_1\otimes \O {-3}\\
&\luTo_{I_{\Sigma_{(\lambda:\mu)}(\tau_i)}}&\uTo_{M(y,\tau_i)}
&&\uTo^{L(\iota (y),\nu_i)^t}&\luTo^{I_{E_{(\lambda:\mu)}}^t}\\
&&V_1\otimes\O {-3}&\lTo^{L(\iota (y),\zeta_i)^t}&V_0\otimes\O
{-4}&\lTo^{[y]^t}
& \O {-5}\\
\enddiagram
$$
with suitably chosen $\nu_i$ and $\zeta_i\in\Pn 1_-$, $i=1,2$.
In particular, this provides an explicit description of the
first and second order linear syzygies  and of the  first
order quadratic syzygies of the first Hartshorne-Rao module of
$B^{\prime\prime}$, from which our claim follows as in
Proposition 3.7.\qed
\enddemo
\head {7. Non-minimal abelian surfaces obtained via linkage}
\endhead
In this section we review the non-minimal abelian surfaces
obtained from the minimal ones via linkage from a syzygy point of view.
This also serves as a further motivation for Section 8 where we will
present a new family of non-minimal abelian surfaces in $\Pn 4$.
\medskip
Let $A=X_s$ be an abelian surface of degree $d=10$ and
sectional genus $\pi=6$ as in Section 3. Then we can link $A$ in the complete
intersection of two quintic hypersurfaces to a smooth surface $A'$
with $d=15$ and $\pi=21$. It follows from a Bertini type argument on the
blow-up
of $\openP^4$ in $L$ that for generic choices $A'$ is indeed also
smooth in the points of the union $L$ of the HM-lines $L_{ij}$
(Ellingsrud and Peskine, unpublished). By
Proposition 2.7 (compare also \cite{Au}) $A'$ is a non-minimal abelian
surface with the 25 HM-lines $L_{ij}$ as exceptional lines and the canonical
divisor $K= \sum L_{ij}$. Moreover, $A'$ carries a polarization of type
(2,10) \cite{Au, Corollary 4.11}, \cite{ADHPR, Corollary 40} (compare
Proposition 8.4 below). From the exact sequence
$$
0\rightarrow \scr K\rightarrow 3\scr O(-5)\oplus 15\scr O(-6)\rightarrow
\scr J_A\rightarrow 0
$$
given by the syzygies of $\scr J_A$ we deduce via linkage \cite{PS} a
resolution
$$
0\rightarrow \scr F'=\scr O(-5)\oplus 15\scr O(-4)\rightarrow \scr G' =
\scr K^{\vee} (-10)\rightarrow \scr J_{A'} \rightarrow 0.\tag7.1
$$
Hence $H^2_*\scr J_{A'}$ is the $\openC$-dual module
$$
N' = H^2_*\scr J_{A'} = H^2_* \scr K^{\vee} (-10) = (H^1_*\scr
J_A(5))^* = (N(3))^*,
$$
and $H^1_*\scr J_{A'}$ is just the vector space
$$
H^1_*\scr J_{A'} = H^1\scr J_{A'}(5) = (H^2\scr J_A)^* = 2 S\
$$
(compare Section 3). In terms of $N'$
$$
\scr G' = \ker (2\scr O(-5)\overset\beta\to\leftarrow \scr
Syz_2(N'))\ ,
$$
where $\beta$ is induced by the minimal free resolution of $N'$. So we
may compute the equations of $A'$ directly from $N'$ and thus from $N$.
Of course this can also be seen from the exact sequence
$$
0\rightarrow \scr E^{\vee}(-8)\rightarrow 3\scr O(-5) \rightarrow
\scr J_{A'}\rightarrow 0\ ,
$$
obtained from (3.1) via linkage.
\proclaim
{Proposition 7.2} $\scr J_{A'}$ has syzygies of type
%
$$
\vbox{%
\halign{&\hfil\,$#$\,\hfil\cr
&&&&3\scr O(-5)\cr
0&\leftarrow&\scr J_{A'}&\leftarrow&\oplus\cr
&&&&5\scr O(-7)&\vbox to 10pt{\vskip-4pt\hbox{$\nwarrow$}\vss}&15\scr
O(-8)&\leftarrow&10\scr
O(-9)&\leftarrow&2\scr O(-10)\leftarrow 0\ .\cr
}}
$$
The quintics alone already cut out $A'$\ .\qquad\qed
\endproclaim
\proclaim
{Remark 7.3} We have the following cohomology table for $\scr
J_{A'}$:
$$
\vbox{\offinterlineskip
\halign{\hbox to 10pt{\hfil#}
&\vrule height10.5pt depth 5.5pt#
&\hbox to 25pt{\hfil$#$\hfil}
&\vrule#
&\hbox to 25pt{\hfil$#$\hfil}
&\vrule#
&\hbox to 25pt{\hfil$#$\hfil}
&\vrule#
&\hbox to 25pt{\hfil$#$\hfil}
&\vrule#
&\hbox to 25pt{\hfil$#$\hfil}
&\vrule#
&\hbox to 25pt{\hfil$#$\hfil}
&\vrule#
&\hbox to 25pt{\hfil$#$\hfil}
&\vrule#
&\hbox to 10pt{#\hss}
\cr
&\omit&&\omit\hbox{\hskip-3.015pt\hbox{$\bigg\uparrow$}
\raise5pt\hbox{$i$}\hfil}\cr
\multispan{17}\hrulefill\cr
&& && && && && && && &\cr
&\multispan{14}\hrulefill\cr
&&20&&1&& && && && && &&\cr
&\multispan{14}\hrulefill\cr
&& &&2&&10&&10&&5&& && &\cr
&\multispan{14}\hrulefill\cr
&& && && && && && &&2&\cr
&\multispan{14}\hrulefill\cr
&& && && && && && &&3&\cr
\noalign{\vskip-3pt}
\noalign{\hbox to 9.5cm{\rightarrowfill}}
\noalign{\hbox to 9cm{\hfil$m$}}
}}
$$
By applying Beilinson's spectral sequence to $\scr J_{A'}(3)$ we
reobtain the exact sequence (7.1).\qquad\qed
\endproclaim
\proclaim
{Remark 7.4} Generalized Serre correspondence \cite{Ok2, Theorem 2.2}
yields an exact sequence
$$
0\rightarrow \scr O\rightarrow \scr E'(3)\rightarrow \scr
J_{A'}(5)\rightarrow 0\ ,
$$
where $\scr E'$ is a stable rank 2 reflexive sheaf with Chern classes
$c_1=-1$, $c_2=9$, $c_3=25$ and $c_4=50$. Since
$$
\dim \Ext^1 (\scr J_{A'} (5), \scr O) = p_g (A') = 1\ ,
$$
this extension is uniquely determined. So it coincides with the one
given by Beilinson's spectral sequence applied to $\scr J_{A'}(4)$.
 From (7.1) we can derive an alternative description: $\scr E'(-2)$
is the cokernel
$$
0\rightarrow 15\scr O(-4)\overset\gamma\to\rightarrow\scr
G'\rightarrow \scr E'(-2)\rightarrow 0\ ,
$$
where $\gamma$ is given by the syzygies of $\scr G'$. So $\scr E'$
depends on $N'$ only and plays the role of the HM-bundle for the
linked surfaces.\qquad\qed
\endproclaim
\head {8. A new family of non-minimal abelian surfaces}
\endhead
In this section we present the construction of a new family of
non-minimal
abelian surfaces of degree 15 and sectional genus 21 in $\Pn 4$. This
construction and Proposition 8.2 are taken from \cite{Po, Section 7}. For
some of our results we need to rely on a computation using \cite{Mac}.
\par\medskip
Let $\Pi=<p_1,p_2,p_3>$ be a plane as in (5.1), (iii) and Proposition 5.6
with $p_1 =a=(a_0:a_1:a_2:a_2:a_1 )\in\openP^2_+\setminus\Delta$,
$p_2 = (0:1:0:0:-1)$ and $p_3 = (0:0:1:-1:0)$, and let
$M$ be the $G_5$-invariant graded $R$-module with representing plane $\Pi$.
Then $M$ is artinian by Proposition 5.6. By construction
$M_0\cong V_3$ and $M_1\cong 2V_1$ as $G_5$-modules. We will use the double
complex of the Koszul cohomology to compute $M_k$ for $k\geq 2$.
The kernel of the vertical differential $4V_0\oplus 6 V_0^\#\rightarrow
M_1\otimes V_2 \cong  6 V_0\oplus 4 V_0^\#$ contains $2V_0^\#$, and in fact
equality holds by Proposition 5.8.  Therefore $M_2\cong 2V_0^\#$.
$M$ has no second order linear syzygies by  Proposition 5.10 and
it has exactly 2 genuine first order quadratic syzygies  by Proposition
5.11. Thus a similar argument as above shows that $M_3\cong 2I$.
Moreover, $M$ has no second order quadratic syzygies by Proposition
5.12. \medskip
In the sequel we will claim correctness of our results only for a general
point $p_1=a\in\openP^2_+\setminus\Delta$. In fact, the module $M$ with
representing plane $\Pi$ will lead to a smooth surface
only if we choose $p_1$ general enough (compare Section 9). So it will be
sufficient for our purposes to check in an example that some Zariski open
conditions are non empty.
E.g., we can check in an example, say $p_1=(1:0:0:0:0)$, that
$M_k=0$ for $k\geq 4$. Therefore, for generic choices,  $M$ has Hilbert
function $(5,10,10,2)$ as in the case of the first Hartshorne-Rao module
of the HM-surfaces (compare Proposition 3.2). But this time the minimal
free resolution of $N_a = (M)^{\vee}(3)$ is of the form
%
$$
\vbox{%
\halign{&\hfil\,$#$\,\hfil\cr
0\leftarrow N_a\leftarrow 2R\ \cr
&\vbox to 10pt{\vskip-4pt\hbox{$\nwarrow$}\vss}\ 20R(-2)\leftarrow 35
R(-3)&&15 R(-4)\cr
&&&\oplus\cr
&&\vbox to 10pt{\vskip-4pt\hbox{$\nwarrow$}\vss}&2 R(-5)\cr
&&&\oplus&\vbox to 10pt{\vskip-2pt\hbox{$\nwarrow$}\vss}\cr
&&&10 R(-6)&&15R(-7)\leftarrow 5 R(-8)\leftarrow 0\ .\cr
}}
$$
Comparing with (7.1) we define
$$
\scr F =\scr O(-5)\oplus 15\scr O(-4)\ ,\quad
\scr G_a =\scr Syz_2(N_a).
$$
As one can check with \cite{Mac} in an example the generic
$\varphi\in\Hom(\scr F,\scr G_a)$ gives rise to a
smooth surface $A_a\subset\openP^4$ with $d=15$, $\pi=21$ and
$\chi=0$. By construction $\scr J_{A_a}$ has a minimal
free resolution of type
%
$$
\vbox{%
\halign{&\hfil\,$#$\,\hfil\cr
&&&&\scr O(-5)\cr
0&\leftarrow&\scr J_{A_a}&\leftarrow&\oplus\cr
&&&&10\scr O(-6)&\vbox to 10pt{\vskip-4pt\hbox{$\nwarrow$}\vss}
&15\scr O(-7)&\leftarrow&5\scr O(-8)&\leftarrow 0\cr
}}
$$
and the cohomology table
$$
\vbox{\offinterlineskip
\halign{\hbox to 10pt{\hfil#}
&\vrule height10.5pt depth 5.5pt#
&\hbox to 25pt{\hfil$#$\hfil}
&\vrule#
&\hbox to 25pt{\hfil$#$\hfil}
&\vrule#
&\hbox to 25pt{\hfil$#$\hfil}
&\vrule#
&\hbox to 25pt{\hfil$#$\hfil}
&\vrule#
&\hbox to 25pt{\hfil$#$\hfil}
&\vrule#
&\hbox to 25pt{\hfil$#$\hfil}
&\vrule#
&\hbox to 25pt{\hfil$#$\hfil}
&\vrule#
&\hbox to 10pt{#\hss}
\cr
&\omit&&\omit\hbox{\hskip-3.015pt\hbox{$\bigg\uparrow$}
\raise5pt\hbox{$i$}\hfil}\cr
\multispan{17}\hrulefill\cr
&& && && && && && && &\cr
&\multispan{14}\hrulefill\cr
&&20&&1&& && && && && &&\cr
&\multispan{14}\hrulefill\cr
&& &&2&&10&&10&&5&& && &\cr
&\multispan{14}\hrulefill\cr
&& && && && && && && &\cr
&\multispan{14}\hrulefill\cr
&& && && && && && &&1&\cr
\noalign{\vskip-3pt}
\noalign{\hbox to 9.5cm{\rightarrowfill}}
\noalign{\hbox to 9cm{\hfil$m$}}
}}
$$
Hence by Proposition 2.7 $A_a$ is a non-minimal abelian surface
with 25 exceptional lines. Generalized Serre correspondence yields an exact
sequence
$$
0\rightarrow\scr O\rightarrow\scr E_a(3)\rightarrow
\scr J_{A_a}(5)\rightarrow 0,\tag 8.1
$$
where $\scr E_a$ is a stable rank 2 reflexive sheaf on $\Pn 4$ with
Chern classes
$c_1=-1$, $c_2=9$, $c_3=25$ and $c_4=50$. As in Remark 7.4 one
sees that (8.1) can also be obtained by applying Beilinson's spectral
sequence to $\scr J_{A_a}(4)$, and hence deduce that $\scr E_a$
is the cokernel
$$
0\rightarrow 15\scr O(-4)\overset\gamma\to\rightarrow\scr
G_a\rightarrow \scr E_a(-2)\rightarrow 0\ ,
$$
where $\gamma$ is induced by the syzygies. In any case $H^0(\scr
E_a(3))$ has dimension 2.
\proclaim
{Proposition 8.2}
Fix a general point $p_1=a$ as above and a basis $s_1,s_2$ of
$H^0(\scr E_a(3))$. Then there exists a pencil
$$
A_{(\lambda:\mu)} = \{\lambda s_1 +\mu s_2 = 0\}
$$
of $G_5$-invariant, smooth, non-minimal abelian surfaces with $d=15$, $\pi=21$,
all lying on the unique quintic hypersurface
$$
U = \{s_1\wedge s_2=0\}\subset\openP^4\ .
$$
Moreover
$$
\Sing(\scr E_a) = \bigcup E_{ij}
$$
consists of the 25 lines $E_{ij}$ obtained as the $H_5$-translates of
the line
$$
E_{00} = \{y_1-y_4 = y_2-y_3 = a_0y_0+2a_1y_1+2a_2y_2 = 0\}
\subset\openP^2_+\ .
$$
These lines are in the base locus of the pencil and hence  are
the exceptional lines of the smooth members in the pencil.
\endproclaim
\demo
{Proof}
$\scr E_a$ has the syzygies
%
$$
\vbox{%
\halign{&\hfil\,$#$\,\hfil\cr
&&&&2\scr O(-3)\cr%
0&\leftarrow&\scr E_a&\leftarrow&\oplus&\overset{\beta =
(\alpha_{22}^t,\alpha_{21}^t)}\to{\hbox to 2cm{\leftarrowfill}}&15\scr
O(-5)\ \overset{\alpha_1^t}\to{\hbox to 1truecm{\leftarrowfill}}\ 5\scr
O(-6)\leftarrow 0,\cr
&&&&10\scr O(-4)\ .\cr
}}
$$
where the maps are given by the minimal free resolution of $N_a$.
Recall that
$$
\alpha_{21} = \left( \matrix
L_{q_{11}}(y) & L_{q_{12}}(y)\\
L_{q_{21}}(y) & L_{q_{22}}(y)\\
L_{q_{31}}(y) & L_{q_{32}}(y)
\endmatrix \right)
$$
with $q_{ij}$, $1\le i\le 3$, $1\le j\le 2$, as in Proposition 5.8.
$\alpha_{22}$ is given by the 2 genuine first order quadratic syzygies
of $\alpha_1$ which we have to choose among all first order quadratic
syzygies
corresponding to the trivial character of $H_5$ as in the proof of
Proposition 5.11. There are 4 such syzygies which we collect in the matrix
$$
S = \left( \matrix
0 &0&0&0\\
(y_{i+2}y_{i+3})_{i\in\openZ_5} & 0&(y_{i+1}y_{i+4})_{i\in\openZ_5}
&(y_{i}^2)_{i\in\openZ_5}\\
0 & (y_{i+1}y_{i+4})_{i\in\openZ_5}&(y_{i}^2)_{i\in\openZ_5}&
-(y_{i+2}y_{i+3})_{i\in\openZ_5}
\endmatrix \right).
$$
Modulo the first order linear syzygies we have the relations
$$
S\cdot \left( \matrix
a_1&a_0\\
-a_0&a_2\\
0&-a_1\\
-a_2&0\endmatrix \right)=0,
$$
which tell us how to choose $\alpha_{22}$ depending on $a$.
For the general point $p_1=a$ one can check
that $\beta$ drops rank on $E_{00}$ and thus also on all the
$H_5$-translates $E_{ij}$. Hence $\bigcup E_{ij}\subset \Sing (\scr E_a)$
and equality follows from
$$
\deg\ \Sing (\scr E_a) = c_3 (\scr E) = 25.
$$
That $\bigcup E_{ij}$ is contained in the base locus is clear
from generalized Serre correspondence.\qed
\enddemo
\proclaim
{Remark 8.3}
(i)\quad One can show that in fact $\scr E_a(3)$ is globally generated
outside $\bigcup E_{ij}$ and that $\bigcup E_{ij}$ is precisely the
base locus of the pencil $A_{(\lambda:\mu)}$.\par
\noindent (ii)\quad The quintic hypersurface $U$ is interesting in
its own and has
a rich geometry. For details compare \cite{Po, Section 7} and \cite{GP2}.\qed
\endproclaim
\proclaim
{Proposition 8.4} (i)\quad If $A_a$ is a new abelian surface, then $H$ is of
type (1,20).\par\noindent
(ii)\quad If $A'$ is a linked abelian surface as in Section 7,
then $H$ is of type (2,10).
\endproclaim
\demo
{Proof} (i)\quad One can check using \cite{Mac} that $A_a$ intersects
$\openP^1_-$ in 3 points and $\openP^2_+$ in the line $E_{00}$ and 12
further points.
Choose any of these 12 points as the origin of $A_a$ (strictly speaking of
the minimal model $A_a^{\text {min}}$ of $A_a$, of course). The involution
$\iota$ leaves the union of the 25 lines $E_{ij}$ invariant and descends to
an involution $\iota$ on $A_a^{\text {min}}$. $\iota$ has exactly 16 fixed
points on $A_a^{\text {min}}$, namely the images of the points where $A_a$
intersects $\openP^1_-$ and $\openP^2_+$ resp., and the point which arises
from blowing down the line which is contained in the intersection of $A_a$
and $\openP^2_+$. Moreover, $\iota$ induces an involution on $\scr L =
\scr O_{A_a^{\text {min}}}(H)$, i.e., $\iota$ is an involution of the
polarized abelian surface $(A_a^{\text {min}},\scr L)$ and $\scr L$ is a
symmetric line bundle, i.e.,
$\iota^* \scr L = \scr L$. Since $\iota$ has exactly 16 fixed points it
must be the involution $x\mapsto - x$ on $A_a^{\text {min}}$. This follows
from the Lefschetz fixed point formula (see e.g. \cite{BL}).
The involution $\iota$ has 12 even 2-torsion points (they come from
the 12 isolated points of $(\openP^2)^+\cap A_a)$ and 4 odd
2-torsion points (these come from $(\openP^1)^-\cap A_a$ and the
invariant line that was blown down). Hence by \cite {LB, Proposition 4.7.5}
the polarization given by $\scr L$ must be of type (1,20).\smallskip\noindent
(ii)\quad In this case we get by the same argument that there are 16
even 2-torsion points and no odd 2-torsion points. Again by
\cite {LB, Proposition 4.7.5} the polarization must be the
square of a primitive polarization, i.e., it must be of type (2,10).\qed
\enddemo
\head
{9. Degenerations of abelian and bielliptic surfaces.}
\endhead
In this section we describe several types of $H_5$-invariant
degenerations of abelian and bielliptic surfaces of degree 10 and 15.
The degenerations of the smooth HM-surfaces are already completely
classified and very well understood
(cf. \cite{HM}, \cite{BHM2}, \cite{HV}, \cite{BM}, \cite{Hu2}).
Therefore we will focus in the sequel only on degenerations of the
other abelian and bielliptic surfaces described in this paper. Our
results are in the spirit of those in \cite{BHM2}, but we do not
aim at complete lists here. Nevertheless our results underline
the fact that the Hilbert scheme is connected \cite{Ha}.\medskip
In order to construct degenerations we specialize the representing
planes of the modules which give rise to the smooth surfaces. In (5.1),
(ii'), (iv) and (v) these planes are spanned by torsion points
on an elliptic normal curve. Therefore we first need to
understand the limits of such points on a singular fibre of the
Shioda modular surface $S(5)$. We recall that all such fibres are mapped
to pentagons in $\Pn 4$ via the normalization map described in Proposition
4.2. For the degenerations we look at the limit points on the standard
pentagons $E_{(0:1)}$ and  $E_{(1:0)}$ with vertices $e_0,\dots, e_4$ only.
On $S(5)$ we have the following:
\proclaim{Lemma 9.1} Let $\scr D$, $\scr X$, and $\scr T$ resp. be the
3-, 8-, and 24-section resp. of the Shioda modular surface $S(5)$
defined by the torsion points of order 2, 3, and 6 resp.
on the smooth fibres. Furthermore, let $O$ denote the zero section,
and let $\Gamma$ be a singular fibre, say with vertices $e_0,\dots, e_4$,
such that the involution $\iota$ operates by $e_i\to e_{-i}$.
If the edges of $\Gamma$ are $\overline{e_{i+2}e_{i+3}}$, $i\in\openZ_5$,
then the following is true:\par\noindent
(i)\quad $O$ intersects $\Gamma$ in the fixed point of $\iota$ on
$\overline{e_2e_3}$ which corresponds to (1:-1).\par\noindent
(ii)\quad $\scr D$ passes transversely through the $\iota$-invariant
point on $\overline{e_2e_3}$ which corresponds to (1:1),
and it passes with multiplicity 2 through $e_0$.\par\noindent
(iii)\quad $\scr X$ passes with multiplicity 3 through $e_1$ and $e_4$ and
intersects transversely $\overline{e_2e_3}$ in 2
points corresponding to the 2 primitive roots of -1 of order 3.\par\noindent
(iv)\quad  $\scr T$ passes with multiplicity 4 through $e_0$,
with multiplicity 3 through $e_1$ and $e_4$, with multiplicity 6 through
$e_2$ and $e_3$, and intersects transversely $\overline{e_2e_3}$ in 2
points corresponding to the 2 primitive roots of unity of order 3.
\par\noindent
The analogous results hold for the pentagon with edges
$\overline{e_{i+1}e_{i+4}}$, $i\in\openZ_5$.
\endproclaim
\demo{Proof}
(i) and (ii) are well-known, see \cite{Sh}, \cite{Na} or
\cite{HKW, Proposition 2.33} for a proof and comments. For (iii) and (iv)
one can proceed as in the proof of \cite{HKW, Proposition 2.33},
along with a careful book-keeping of the various glueings in the
toroidal compactification. Another method is to make these computations
directly on a projective model of $S(5)$, namely in our case on
$S_{15}\subset\Pn 4$ with $\Gamma=E_{(0:1)}$.  In this setting (iii)
follows directly from Proposition 4.13.  In fact, on $\overline{e_2e_3}$
the limits of the proper 3-torsion points are defined by the
equation $y_2^2-y_2y_3+y_3^2=0$. For (iv) one uses the fact that the octic
hypersurface $w_0^2+(w_1+w_4)(w_2+w_3)$, where $w_{2i}=y_{i}(y_{i+2}
y_{i+4}^2-y_{i+1}^2y_{i+3})$, $i\in\openZ_5$, meets each smooth fiber of
$S_{15}$ in 6-torsion points only by Proposition 4.19.
Thus, for instance, on $\overline{e_2e_3}$
the limits of the proper 6-torsion points are defined by the
equation $y_2^6y_3^6(y_2^2+y_2y_3+y_3^2)=0$. The proof of the other
assertions in (iv) is similar.\qed
\enddemo
We are now ready to study degenerations of our surfaces.
Again a first question will be whether our limit planes represent artinian
modules with the expected number of syzygies, i.e., with the same graded Betti
numbers as the general members of our families. The planes contained in
coordinate hyperplanes are not of interest for us:
\proclaim{Remark 9.2} Any plane $\Pi$ inside a coordinate hyperplane
$\{x_i=0\}$
gives rise to a non-artinian module by Lemma 5.3, since $\Pi$ is contained in
the
kernel of $N(e_i, e_i)$.\qed\endproclaim\medskip
For the minimal bielliptic surfaces of degree 10 the first
Hartshorne-Rao module of the surface has a representing plane of type
$\Pi=<\tau_1,\tau_2,\tau_3+\rho>$, where $\tau_i$, $i=\overline{1,3}$ are
the non-trivial 2-torsion points and $\rho$ is a non-trivial 3-torsion point
on a $G_5$-invariant elliptic normal curve in $\Pn 4$. Remark 9.2 rules out
most of the possible limit planes. We may, however, consider the case where
the points $\tau_1$ and $\tau_2$  come together at $e_0$, while $\tau_3+\rho$
goes to one of the two possible limit points on $\overline{e_2e_3}$ (in the
case of $E_{(0:1)}$) or $\overline{e_1e_4}$ (in the case of $E_{(1:0)}$).
There are two branches of the Brings curve $B$ through $e_0$ with
tangents $\{y_2=y_3=y_1-y_4=0\}$ or $\{y_2-y_3=y_1=y_4=0\}$. This gives rise
to four limit planes which are spanned by $\{y_2=y_3=y_1-y_4=0\}$ and one of
the two points $(0:0:\theta_k:1:0)\in\overline{e_2e_3}$ or by
$\{y_2-y_3=y_1=y_4=0\}$ and one of the two points $(0:\theta_k:0:0:1)\in
\overline{e_1e_4}$, where $\theta_1$ and $\theta_2$ are two different primitive
third roots of unity.
One checks through a direct computation that these limit
planes define indeed four artinian graded modules with the Hilbert
function $(5,10,10)$ and with the expected syzygies. Moreover,
a determinantal construction as in Proposition 3.7 and
Theorem 6.3 leads to four locally Cohen-Macaulay surfaces
$Y_{(\theta_k:1)}$ or $Z_{(\theta_k:1)}$ in $\Pn 4$ with the
same numerical invariants as a smooth minimal bielliptic surface of degree 10.
\smallskip\noindent
We will describe the equations of these degenerations. In fact,
$Y_{(\theta_k:1)}$ and $Z_{(\theta_k:1)}$ can be seen as special points of
other components of the Hilbert scheme. Namely, consider the surfaces
$$
Y_{(\lambda:\mu)}=\bigcup_{i\in\openZ_5}Y^{(\lambda:\mu)}_i \quad {\text{and}}
\quad Z_{(\lambda:\mu)}=\bigcup_{i\in\openZ_5}Z^{(\lambda:\mu)}_i,
$$
where each $Y^{(\lambda:\mu)}_i$ and $Z^{(\lambda:\mu)}_i$ is a double
structure
(only generically Cohen-Macaulay) on a plane:
$$
Y^{(\lambda:\mu)}_i=\lbrace y_{i}^2=y_{i+2}^2=y_{i}y_{i+2}=
\lambda y_{i+2}y_{i+3}^2+\mu y_{i}y_{i+4}^2=0\rbrace
$$
and
$$
Z^{(\lambda:\mu)}_i=\lbrace y_{i}^2=y_{i+1}^2=y_{i}y_{i+1}=
\lambda y_{i+1}y_{i+4}^2+\mu y_{i}y_{i+2}^2=0\rbrace.
$$
The reduced support of $Y_{(\lambda:\mu)}$ or $Z_{(\lambda:\mu)}$ is
the union of 5 planes, which in turn is a degeneration of an elliptic
quintic scroll in $\Pn 4$. Therefore, in a way, these surfaces are double
structures on a degenerate elliptic quintic scroll.
$Y_{(1:1)}$ and $Z_{(1:1)}$ are HM-surfaces and lie on three
quintics, but for the general point $(\lambda:\mu)$ these surfaces lie
on the same unique degenerate quintic
$\lbrace\gamma_0=y_0y_1y_2y_3y_4=0\rbrace$. Correspondingly, the
plane which is spanned by the tangent $\{y_2=y_3=y_1-y_4=0\}$ and a point
$(0:0:\lambda:\mu:0)\in\overline{e_2e_3}$ or the tangent
$\{y_2-y_3=y_1=y_4=0\}$ and a point $(0:\lambda:0:0:\mu)\in
\overline{e_1e_4}$ is $\openP^2_+$ for $(\lambda:\mu)=(1:1)$, whereas for
the general point $(\lambda:\mu)$ this plane represents an artinian module
with the same number of syzygies as for the minimal bielliptics.
\smallskip\noindent
Now recall that the first Hartshorne-Rao module of the union of two
$G_5$-invariant elliptic quintic scrolls as in Theorem 6.7 has a
representing plane of type $\Pi=<\tau_1,\tau_2,o>$, where $\tau_i$,
$i=\overline{1,2}$, are two non-trivial 2-torsion points and $o$ is the
origin of a $G_5$-invariant elliptic normal curve in $\Pn 4$.
Therefore we can specialize as above such that $\tau_1$ and $\tau_2$ come
together to the vertex $e_0$, while $o$ goes to the point
$(0:0:1:-1:0)\in\overline{e_2e_3}$ or $(0:1:0:0-1)\in\overline{e_1e_4}$.
Thus $Y_{(1:-1)}$ and $Z_{(1:-1)}$ are degenerations of the union of two
elliptic quintic scrolls as in Proposition 6.7.\medskip
Next we study the  non-minimal bielliptic surfaces of degree 15 as in
Theorem 6.6. The dual of the second Hartshorne-Rao module of such a
surface has a representing plane of type
$\Pi=<\tau_1+\rho,\tau_2+\rho,\tau_3+\rho>$, where $\tau_i$,
$i=\overline{1,3}$ are the three non-trivial 2-torsion points
and $\rho$ is a non-trivial 3-torsion point of a
$G_5$-invariant elliptic normal curve in $\Pn 4$.
As in the case of the minimal bielliptic surfaces Remark 9.2 rules out most
of the possible limit planes. Some limit planes however give rise to
artinian modules. Let $\tau_1+\rho$ and $\tau_2 +\rho$ come together at
the vertex $e_0$ while the third point $\tau_3+\rho$ goes to one of the two
possible limit points on $\overline{e_2e_3}$ (in the case of $E_{(0:1)}$) or
$\overline{e_1e_4}$ (in the case of $E_{(1:0)}$). The image in $\Pn 4$ of the
curve $\scr T$ of points of order 6 has four branches through $e_0$ with
tangents spanned by $e_0$ and one of the four points corresponding to primitive
roots of unity on the lines $\overline{e_2e_3}$ and $\overline{e_1e_4}$. This
gives rise to four limit planes containing $e_0$ and meeting
$\overline{e_2e_3}$
and $\overline{e_1e_4}$ transversally at points of order 6. The corresponding
modules are artinian and give rise to four reduced but reducible surfaces.
These are the union of five cubic surfaces, namely
$$
\{y_i= \theta_k y_{i+1}^2y_{i+3}+y_{i+2}y_{i+4}^2=0\}, \quad i\in\openZ_5,
$$
or
$$
\{y_i= \theta_k y_{i+3}^2y_{i+4}+y_{i+1}y_{i+2}^2=0\}, \quad i\in\openZ_5,
$$
for $k=1,2$, where again $\theta_1$ and $\theta_2$ are two different primitive
third roots of unity. In each case the first five cubic surfaces glue along
the pentagon $E_{(0:1)}$, while the last five cubic surfaces glue along the
pentagon $E_{(1:0)}$. As in the minimal bielliptic case these degenerations
belong to pencils of surfaces, where this time the surface corresponding to
$(\lambda:\mu)\in\openP^1$ has the five components
$$
\{y_i= \lambda y_{i+1}^2y_{i+3}+\mu y_{i+2}y_{i+4}^2=0\}, \quad i\in\openZ_5,
$$
or
$$
\{y_i= \lambda y_{i+3}^2y_{i+4}+\mu y_{i+1}y_{i+2}^2=0\}, \quad i\in\openZ_5.
$$
For $(\lambda:\mu)=(1:1)$ these surfaces are (5,5)-linked to HM-surfaces
and lie on on three quintics, but for the general point $(\lambda: \mu)$
these surfaces lie on the same unique degenerate quintic
$\{\gamma_0=y_0y_1y_2y_3y_4=0\}$.  Correspondingly, the plane spanned by
$e_0$, $\lambda e_4 +\mu e_1$ and $\lambda e_2 +\mu e_3$ or $e_0$, $\lambda
e_1 +\mu e_4$ and $\lambda e_2 +\mu e_3$ is $\openP_+^2$ for
$(\lambda:\mu)=(1:1)$ while for the general point $(\lambda:\mu)$ this
plane represents a module with the same number of syzygies as for the
non-minimal bielliptic surfaces of degree 15. Recall from
\cite{Po} that there is a pencil of degenerations of the linked abelian
surfaces, such that the
surface corresponding to $(\lambda:\mu)\in\openP^1$ has the five components
$$
\lbrace y_i=\lambda(y_{i+1}^2y_{i+3}+y_{i+2}y_{i+4}^2)+
\mu(y_{i+1}y_{i+2}^2+y_{i+3}^2y_{i+4})=0\rbrace, \quad i\in\openZ_5.
$$
\smallskip\noindent
Now we will discuss briefly some degenerations of the new non-minimal
abelian surfaces of degree 15. The following results are
taken from \cite{Po, Section 7}.\smallskip\noindent
We recall that for these surfaces the dual of the second Hartshorne-Rao
module has a representing plane spanned by $\Pn 1_-$ and a general point
$a\in\Pn 2_+$. We have seen in Proposition 8.2 that such a module gives
rise to a pencil $A_{(\lambda:\mu)}$ of non-minimal
abelian surfaces, all lying on the same unique quintic hypersurface.
It is shown in \cite{Po, Lemma 7.6} that this quintic is
uniquely determined as the  $H_5$-invariant quintic containing the
25 exceptional lines of the smooth members of the pencil, i.e., the
25 lines $E_{ij}$ obtained as the $H_5$-translates of the line
$$
E_{00} = \{y_1-y_4 = y_2-y_3 = a_0y_0+2a_1y_1+2a_2y_2 = 0\}
\subset\openP^2_+
$$
(compare again Proposition 8.2). To obtain degenerations we can either
specialize the parameter point $a$ or look at special members of the
pencil $A_{(\lambda:\mu)}$.\smallskip\noindent
A first type of degenerations corresponds to the choice of
the parameter point $a\in\Pn 2_+$ as one of the nodes of
the Brings curve, say $a=(1:0:0:0:0)$. These degenerations are singular
surfaces $X_{(\lambda:\mu)}\subset\PP^4$,
$(\lambda:\mu)\in\PP^1$, which all lie on the
degenerate quintic hypersurface
$\lbrace\gamma_0=y_0y_1y_2y_3y_4=0\rbrace$.
A straightforward but tedious computation shows that
$$X_{(\lambda:\mu)}=\bigcup_{i\in\ZZ_5}X^{(\lambda:\mu)}_i,$$
where the $X^{(\lambda:\mu)}_i$ are the cubic surfaces defined by
$$X^{(\lambda:\mu)}_i=\lbrace y_i=\lambda(y_{i+1}^2y_{i+3}-y_{i+2}y_{i+4}^2)+
\mu(y_{i+1}y_{i+2}^2-y_{i+3}^2y_{i+4})=0\rbrace.$$
\noindent
We fix now a general $(\lambda:\mu)\in\PP^1$ and write
$X_i=X^{(\lambda:\mu)}_i$.
\proclaim{Proposition 9.3}(i)\quad $X_i$ is a smooth Del Pezzo cubic surface
in the hyperplane $H_i=\{y_i=0\}$. Moreover, it is invariant under
the action of $\tau$, whereas $\sigma(X_i)=X_{i-1}$, $i\in\ZZ_5$.\par
\noindent (ii)\quad  Two surfaces $X_i$ and $X_j$, $i\ne j$, meet along
a smooth conic and a point outside the conic. There are altogether five
such points, namely the vertices $e_i$, $i\in\ZZ_5$, of the complete
pentagon $E_{(0:1)}\cup E_{(1:0)}$. Through each point $e_i$ there pass exactly
four such Del Pezzo surfaces.\par
\noindent (iii)\quad Write $E_{(0:1)}=\bigcup_{i\in\ZZ_5}L_i$ with
$L_i=\lbrace y_{i+2}=y_{i+3}=y_{i+4}=0\rbrace$, $i\in\openZ_5$, and
$E_{(1:0)}=\cup_{i\in\ZZ_5}{L'}_i$ with
${L'}_i=\lbrace y_{i+1}=y_{i+3}=y_{i+4}=0\rbrace$, $i\in\openZ_5$. Then
$$X\cap H_i=X_i\cup{L'}_{i+1}\cup{L}_{i+2}\cup{L}_{i+3}\cup{L'}_{i+4}.$$
Furthermore, these four lines are exceptional on all $X_j$ with
$j\in\ZZ_5\setminus\lbrace i\rbrace$.
\endproclaim
\demo{Proof} The claims follow by a straightforward computation using the
explicit equations of $X_i$.\qed\enddemo
Observe now that for this choice of the parameter point the exceptional lines
of the degree 15 non-minimal abelian surfaces degenerate to
$\tilde E_{ij}=\sigma^i\tau^j \tilde E_{00}$, $i,j\in\ZZ_5$, where
$\tilde E_{00}=\lbrace y_0=y_1-y_4=y_2-y_3=0\rbrace$.
It is easily checked that $\tilde E_{i0}$, $\tilde E_{i1}$,
$\tilde E_{i2}$, $\tilde E_{i3}$
and $\tilde E_{i4}$ are $(-1)$-curves on the Del Pezzo surface $X_i$.
As a consequence we obtain the following geometric characterization
of our configuration:
\proclaim{Proposition 9.4} For each $i\in\ZZ_5$ the lines
$\tilde E_{i0}$, $\tilde E_{i1}$, $\tilde E_{i2}$, $\tilde E_{i3}$,
$\tilde E_{i4}$ and $L_{i+3}$, ${L'}_{i+1}$ can be completed
to a Schl\"afli double six configuration of lines
in the hyperplane $H_i$, which in turn determines
the Del Pezzo surface $X_i$ as the unique cubic surface
in $H_i$ containing the given double six.
\endproclaim
\demo{Proof} The intersection patterns are clear from the explicit
description of the configuration, and thus
the claim follows from \cite{HCV, \S 25}.
\qed\enddemo
\medskip
We will finally mention some degenerations, which we have checked in an
explicit example using \cite{Mac}. For some of them see \cite{Po} and
\cite{GP2} for more details.\smallskip\noindent
First of all each general pencil $A_{(\lambda:\mu)}$ of non-minimal
abelian surfaces as above contains degenerations which are ruled surfaces
${\openP}({\scr O}_E+{\scr L})$, where $E$ is an elliptic curve and
$\scr L\in \Pic^0(E)$, mapped to $\Pn 4$
either by $H\equiv C_0+20f_o-\sum_{i,j\in\openZ^5}E_{ij}$ or
by $H\equiv 4C_0+5f_o-\sum_{i,j\in\openZ^5}E_{ij}$, where in both cases the 25
exceptional lines form an $H_5$-orbit. In the first case the image in $\Pn 4$
is a scroll of degree 15 with singularities along the image of $C_0$, which is
a curve of degree 20, arithmetic genus 26 and with 25 nodes (the images
of the points where the curves $f_o-E_{ij}$ meet $C_0$). In the second case
the image is a ruled surface of degree 15, with singularities along the image
of
$C_0$, which is an elliptic curve of degree 5.\par
\smallskip\noindent
Another class of degenerations corresponds to the choice of the
parameter point $a$ as a general point on the Brings curve. These are
conic bundles ${\openP}({\scr O}_E+{\scr L})$, with $E$ and $\scr L$ as
above, mapped to $\Pn 4$  by $H\equiv 2C_0+10f_o-\sum_{i,j\in\openZ^5}E_{ij}$.
\smallskip\noindent
Finally, for a choice of the parameter point $a$ as a general point
on the modular conic $C_+\subset\Pn 2_+$, we obtain a pencil of surfaces
$X_{(\lambda:\mu)}$, which has an elliptic quintic scroll as base component.
Furthermore, $X_{(\lambda:\mu)}$ contains residual to the elliptic scroll
a surface $T_{(\lambda:\mu)}$ of degree $10$ and sectional genus $6$.
For  general $(\lambda:\mu)\in\Pn 1$ the surface $T_{(\lambda:\mu)}$ is a
smooth
minimal abelian surface isogenous to a product, and meeting the base
elliptic scroll along a section of degree 15, which is of the type described in
Proposition 4.10, (iii) and Remark 4.11. Moreover, a special member in the
pencil
$X_{(\lambda:\mu)}$ is the first infinitesimal neighborhood of the base
elliptic
quintic scroll.\par
%
\bigskip
\heading  Appendix\\
\\
\rm by C. Ciliberto and K. Hulek
\endheading
%
In this appendix we show that the union $X$ of two quintic elliptic
scrolls in $\openP^4$ intersecting along an elliptic normal curve defines
a smooth point in the Hilbert scheme $\scr H$ of surfaces of degree $d=10$,
sectional genus $\pi=6$ and $\chi ({\scr O}_X)=0$. The dimension of
$\scr H$ at $[X]$ is 25. As a corollary we obtain that every (small)
deformation of $X$ is again a union of two quintic elliptic scrolls.\par

Let $E$ be an elliptic normal curve of degree $5$ in $\openP^4$ with
origin
$0$ and non-zero 2-torsion points $\{Q_i; i=1, 2, 3\}$. The union of all
secants joining points $P$ and $P+Q_i$ where $P$ varies in $E$ is a smooth
elliptic ruled surface $S_i$ of degree 5. For $i\neq j$ the union
$$
X=S_i \cup S_j\quad (i\neq j)
$$
is a surface which is singular exactly along $E$ where the two surfaces
$S_i$ and $S_j$ meet transversally. In particular $X$ is a locally
complete intersection surface, and it is not difficult to compute that
$2K_X={\scr O}_X$.  As an abstract surface $X$ is a chain of two elliptic
ruled surfaces glued along a 2-section. Such surfaces appear in the
classification of degenerations of bielliptic surfaces \cite {FM}. Clearly
the degree of $X$ is 10. The hyperplane section consists of two elliptic
curves intersecting in 5 points, i.e.
$\pi (X)=6$. Moreover $q(X)=1,p_g(X)=0$, i.e. $\chi({\scr  O}_X)=0$.

\proclaim{Theorem 1}\par

\noindent(i)\quad
The Hilbert scheme $\scr H$ of surfaces of degree 10, sectional genus 6 and
$\chi ({\scr O}_X)=0$ is smooth at $[X]$.\par

\noindent(ii)\quad
The dimension of $\scr H$ at $[X]$ is 25.

\endproclaim

\proclaim{Corollary 2}
Every (small) deformation of $X$ is again a union of two quintic elliptic
scrolls.
\endproclaim

Before we can give a proof of this, we must have a brief look at the
geometry of the surfaces $S_i$. For this we again fix an elliptic curve
$E$ with origin $0$. (This is not the elliptic curve from above, but the
curve $E/\langle Q_i\rangle$.) We now denote the non-zero 2-torsion points
of
$E$ by
$\{P_i; i= 1, 2, 3\}$. The second symmetric product $S^2E$ of $E$ is, via
the map

$$
\vbox{
\halign{&$#$\hfil\cr
\pi:      S^2 E\rightarrow E\cr
     \quad   \{x, y\} \mapsto x + y\cr
}}
$$
a $\openP^1$-bundle over $E$. It was already observed in \cite{At, p. 451}
that this is the unique indecomposable $\openP^1$-bundle with $e=-1$. Let
$p:E\times E\rightarrow S^2 E$ be the natural quotient map. The curve
$$
E_0= p(E\times \{0\})
$$
is a section of $S^2 E$. We shall choose the point $p(0,0)$ as its origin
and, by abuse of notation, we shall denote it again by $0$. The curves
$$
\Delta_i =\{(x, x+P_i); x \in E\} (\cong E)\quad (i=1, 2, 3)
$$
are mapped 2:1 under $p$ onto 2-sections $E_i\subset S^2 E$. As abstract
curves $E_i=E/\langle P_i\rangle $. We shall choose the point $0_i=p(0,
P_i)=p(P_i, 0)$ as the origin of $E_i$. The group of 2-torsion points of
$\Delta_i$ is mapped to two points $\{ 0_i, Q_i\}$. Every fibre $f$ of
$S^2 E$ intersects   $E_i$ in two points which differ by $Q_i$. We shall
denote the fibre of $S^2E$ over $P\in E$ by $f_P$, and put $S=S^2E$. The
following formulae follow immediately from the above description:

$$
\vbox{
\halign{&$#$\hfil\cr
(1)\quad K_S = \scr O_S(-2E_0+f_0) \cr
(2)\quad \scr O_S(E_0)|_{E_0} = \scr O_{E_0}(0) \cr
(3)\quad \scr O_S(E_i) = \scr O_S(2E_0-f_{P_i}) \cr
(4)\quad \scr O_S(E_i)|_{E_i} = \scr O_{E_i}(0_i-Q_i)\cr
}}
$$

The line bundle $\scr O_S(H)=\scr O_S(E_0+2f_0)$ is very ample and
embeds $S$ as a scroll of degree 5 into $\openP^4$. Under this embedding
the 2-sections $E_i$ are mapped to quintic elliptic curves and $S$ can be
reconstructed from each of these quintic normal curves as translation
scrolls defined by the point $Q_i$.  Conversely fixing $E_i$, the three
non-zero 2-torsion points of $E_i$ give rise to three elliptic quintic
scrolls $S_i$. We choose two of them and consider the singular surface $X$
which is the object of our considerations. In order to prove Theorem 1 we
have to compute the cohomology groups $H^i(N_X)$ of the normal bundle of
$X$. For this we first have to investigate the normal bundle of $S$. By
\cite {HV, Proposition 4} there is an extension
$$
(5)\quad 0\rightarrow K^{-1}_S\rightarrow N_S\rightarrow
\scr O_S(5H)\otimes K^2_S\rightarrow 0.
$$
\proclaim{Lemma 3} (i) $\quad h^0(N_S)=25$\par
\noindent(ii)$\quad h^j(N_S)=0$ for $j\ge 1$.
\endproclaim
\demo{Proof}
By (1) $K^{-1}_S=\scr O_{S_i}(2E_0-f_0)$. It follows, e.g. since
$K^{-1}_S|_{E_i}\neq \scr O_{E_i}$ that $h^0(K_{S_i}^{-1})=0$.
Clearly $h^2(K_S^{-1})=h^0(K_S^2)=0$ and hence by Riemann-Roch
$h^1(K_S^{-1})=0$. Let  $ \scr M=\scr O_S(5H)\otimes K_S^2$. From the
definition of $H$ and (1) one obtains that $\scr M=\scr O_S(E_0+12f_0)$.
It is easy to see by the Nakai-Moishezon criterion that $\scr M\otimes
K_S^{-1}=\scr O_S(3E_0+11f_0)$ is ample and hence $h^j(\scr M)=0$ for
$j\ge 1$. But then Riemann-Roch gives $h^0(\scr M)=25$ and the
assertions follow from sequence (5) .
\qed
\enddemo
\proclaim{Lemma 4}
$h^j(N_S(-E_i))=0$ for $j\ge 0$.
\endproclaim
\demo{Proof}
Again setting $\scr M=\scr O_S(5H)\otimes K_S^2$ and twisting (5) by
$\scr O_S(-E_i)$ we obtain
$$
(6)\quad 0\rightarrow K_S^{-1}(-E_i)\rightarrow N_S(-E_i)\rightarrow \scr
M(-E_i)\rightarrow 0.
$$
Since $K_S^{-1}(-E_i)=\scr O_S(f_{P_i}-f_0)$ this line bundle has no
cohomology. By definition of $\scr M$ and (3) we obtain $\scr
M(-E_i)=\scr O_S(-E_0+12f_0+f_{P_i})$. Clearly $h^0(\scr
M(-E_i))=h^2(\scr M(-E_i))=0$ and once again by Riemann-Roch $h^1(\scr
M(-E_i))=0.$ The assertion now follows from the above sequence (6).
\qed
\enddemo

We can now compute the cohomology of $N_X$. Theorem 1 is an immediate
consequence of
\proclaim{Proposition 5}
(i)$\quad h^0(N_X)=25$\par
\noindent(ii)$\quad h^j(N_X)=0$ for $j\ge 1$.
\endproclaim

\demo{Proof}
Let $C=E_i$ and consider $T=N_{C/S_i}\otimes N_{C/S_j}=\scr
O_C(20_i-Q_i-Q_j)$. In particular $T$ is a non-trivial element of
order 2. As in \cite {CLM, section 2.2} we have exact sequences \par

$$
(7)\quad 0\rightarrow N_{S_i}\rightarrow N_X|_{S_i}\rightarrow T
\rightarrow 0
$$
$$
(8)\quad 0\rightarrow N_X|_{S_j}\otimes \scr O_{S_j}(-C)\rightarrow
N_X\rightarrow N_X|_{S_i}\rightarrow 0.
$$

It follows from Lemma 3 and sequence (7) that $h^0(N_X|_{S_i})=25$ and
$h^j(N_X|_{S_i})=0$ for $j\ge 1$. Twisting (7) by $\scr O_{S_i}(-C)$ we
get
$$
(9)\quad 0\rightarrow N_{S_i}(-C)\rightarrow N_X|_{S_i}\otimes \scr
O_{S_i}(-C)\rightarrow T\otimes\scr O_{S_i}(-C)\rightarrow 0.
$$
Using (3) we find that $T\otimes \scr O_{S_i}(-C)=\scr O_C(0_i-Q_j)$
which again has no cohomology. It therefore follows from Lemma 4 that
$h^j(N_X|_{S_i}\otimes \scr O_{S_i}(-C))=0$ for $j\ge 0$. The claim then
follows from sequence (8).
\qed
\enddemo

\demo{Proof of Corollary 2}
Let $X(2,5)$ be the modular curve parametrizing elliptic curves with a
level 5 structure and a non-zero 2-torsion point. We claim that every
point of $X(2,5)$ which is not a cusp gives rise to some $X$ as before.
Indeed the elliptic curves with level 5 structure are in 1:1
correspondence with Heisenberg invariant elliptic normal curves in
$\openP^4$. Given a non-zero 2-torsion point $P$ we have exactly 2 other
non-zero 2-torsion points. Hence we can use these 2 points to construct
$X$. Now consider some $X$ as before. Up to a change of coordinates we
can assume that $X$ is constructed as the union of two translation
scrolls of a Heisenberg invariant elliptic normal curve. Let $\scr H_0$ be
the unique component of $\scr H$ containing $X$. We have a map $\Phi:
X^0(2,5) \times $ PGL $(5,\openC)\rightarrow \scr H_0$ where $X^0(2,5)$ is
$X(2,5)$ without the cusps. Since every elliptic normal curve has only
finitely many automorphisms $\Phi$ is locally 1:1 and hence a local
(analytic) isomorphism. This proves the claim.
\qed
\enddemo

Clearly this situation can be generalized to higher dimensions in several
ways. We hope to return to this in the near future.\par

\enddocument